\documentclass[11pt]{article}
\usepackage{amsfonts}
\usepackage{amsmath}
\usepackage{amssymb}
\usepackage{amsthm,amscd}
\usepackage{tikz}
\usepackage{stmaryrd}
\newtheorem{corollary}{Corollary}[section]
\newtheorem{theorem}{Theorem}[section]
\newtheorem{lemma}{Lemma}[section]
\newtheorem{proposition}{Proposition}[section]
\newtheorem{definition}{Definition}[section]
\newtheorem{remark}{Remark}[section]
\newtheorem{convention}{Convention}[section]
\newtheorem{claim}{Claim}[section]

\newcommand{\C}{\mathbb{C}}
\newcommand{\Z}{\mathbb{Z}}
\newcommand{\N}{\mathbb{N}}
\newcommand{\D}{\mathbb{D}}
\renewcommand{\H}{\mathbb{H}}
\newcommand{\R}{\mathbb{R}}
\newcommand{\T}{\mathcal{T}}
\newcommand{\M}{\mathcal{M}}

\title{Goldman's Ergodicity Conjecture and Branched Developing Maps}

\author{Yiran Lin\footnote{Qiuzhen College, Tsinghua University. Email: lin-yr21@mails.tsinghua.edu.cn}}

\begin{document}

\maketitle

\begin{abstract}
    We prove that every representation $\rho:\pi_1(\Sigma_g)\rightarrow\mathrm{PSL}_2(\R)$ which sends every simple closed curve to a hyperbolic element admits a branched developing map. As a corollary, the mapping class group action on the Euler class zero component of the $\mathrm{PSL}_2(\R)$-character variety of a genus $g$ surface is ergodic if $g\geq3$. (The action is not ergodic if $g=2$ \cite{MW16}.) Moreover, almost every representation of nonzero Euler class admits a branched developing map.
\end{abstract}

\section{Introduction}

Let $\Sigma_g$ denote a closed oriented surface of genus $g\geq2$. The $\mathrm{PSL}_2(\R)$-character variety
\[\M(\Sigma_g)=\mathrm{Hom}(\pi_1(\Sigma_g),\mathrm{PSL}_2(\R))\sslash \mathrm{PSL}_2(\R)\]
has $4g-3$ connected components
\[\M(\Sigma_g)=\bigsqcup_{k=2-2g}^{2g-2}\M^k(\Sigma_g),\]
where $k\in\{2-2g,\cdots,2g-2\}$ is the Euler class. The components with extremal Euler class (i.e., $|k|=2g-2$) consist of Fuchsian representations, and are identified with the Teichm\"uller space.

Let $\mathrm{Mod}(\Sigma_g)=\mathrm{Homeo}^+(\Sigma_g)/\mathrm{Homeo}^+_0(\Sigma_g)$ be the mapping class group. Then $\mathrm{Mod}(\Sigma_g)$ naturally acts on $\M^k(\Sigma_g)$. It is well known that the action is discrete if $|k|=2g-2$. However, Goldman \cite{Goldman06} conjectured that the action is ergodic if $|k|<2g-2$.

\begin{remark}
    The genus $2$ case was studied in March\'e-Wolff \cite{MW16}, where the authors proved that the action on $\M^{\pm1}(\Sigma_2)$ is ergodic, while $\M^0(\Sigma_2)$ has two invariant open subsets, with the action ergodic on one of them.
\end{remark}

Closely related to Goldman's ergodicity conjecture is the question of which representations admit branched developing maps.

\begin{definition}
    Let $\Sigma$ be a connected, oriented surface, and let $\rho:\pi_1(\Sigma)\rightarrow\mathrm{PSL}_2(\R)$ be a representation. A developing map for $\rho$ is a continuous, $\rho$-equivariant map $f:\widetilde{\Sigma}\rightarrow\H$, where $\widetilde{\Sigma}$ is the universal covering of $\Sigma$.
\end{definition}

\begin{definition}
    A developing map $f$ is called branched if it is a local homeomorphism away from a discrete set of branch points, at each of which $f$ is locally given by $z\mapsto z^k$ for some integer $k\geq2$.
\end{definition}

If Goldman's ergodicity conjecture were true, then almost every representation of nonzero Euler class admits a branched developing map.

Another closely related question was posed by Bowditch \cite{Bow98} (Question C), which ask whether every non-Fuchsian representation sends a simple closed curve to a non-hyperbolic element of $\mathrm{PSL}_2(\R)$. In this paper we explore these questions and prove the following theorem

\begin{theorem}\label{thm1}
    Suppose that $\rho:\pi_1(\Sigma_g)\rightarrow\mathrm{PSL}_2(\R)$ is a representation that sends every simple closed curve to a hyperbolic element of $\mathrm{PSL}_2(\R)$. Then there exists a branched developing map for $\rho$.
\end{theorem}

There are several immediate consequences.

\begin{theorem}
    Every representation $\rho:\pi_1(\Sigma_g)\rightarrow\mathrm{PSL}_2(\R)$ of Euler class zero sends some simple closed curve to a non-hyperbolic element.
\end{theorem}

This answers Bowditch's question for all representations of Euler class zero.

\begin{proof}
    A non-trivial representation of Euler class zero does not admit a branched developing map.
\end{proof}

\begin{theorem}
    The action of $\mathrm{Mod}(\Sigma_g)$ on $\M^0(\Sigma_g)$ is ergodic if $g\geq3$.
\end{theorem}

\begin{proof}
    Combine Theorem \ref{thm1} with Theorem 1.6 in \cite{MW16}.
\end{proof}

\begin{remark}
    The reader should note that Theorem \ref{thm1} holds for all $g\geq2$. The assumption $g\geq3$ in the above Theorem comes from Theorem 1.6 in \cite{MW16}.
\end{remark}

\begin{theorem}
    Almost every representation $\rho\in\M^k(\Sigma_g)$ of Euler class $k\neq0$ admits a branched developing map.
\end{theorem}

\begin{proof}
    By Theorem 1.6 in \cite{MW16}, if $k\neq0$, $\mathrm{Mod}(\Sigma_g)$ acts ergodically on the subset $\mathcal{NH}^k(\Sigma_g)$ of $\M^k(\Sigma_g)$ consisting of representations which send some simple closed curve to a non-hyperbolic element. It follows that almost every $\rho\in\mathcal{NH}^k(\Sigma_g)$ admits a branched developing map. (See also \cite{Nguyen21}, Theorem 1.1.3.) On the other hand, by Theorem \ref{thm1}, every $\rho\in\M^k(\Sigma_g)\backslash\mathcal{NH}^k(\Sigma_g)$ admits a branched developing map. It follows that almost every $\rho\in\M^k(\Sigma_g)\ (k\neq0)$ admits a branched developing map.
\end{proof}

\begin{remark}
    Very recently, the paper of Deroin-Tholozan \cite{DT26} appeared, in which the authors proved that every representation with almost extremal Euler class $(|k|=2g-3)$ admits a branched developing map. The reader can find extensive background material in this inspiring paper.
\end{remark}

We now describe the proof of Theorem \ref{thm1} in detail. Recall that our goal is to find a branched developing map. As an approximation, we introduce the class of \textit{excellent developing maps}, which roughly speaking are branched away from a finite collection of short, disjoint simple closed curves, on which the developing maps fold. An excellent developing map without folding curve is a branched developing map. We first find an excellent developing map, which is much easier than finding a branched developing map. We then perform a sequence of surgeries to reduce the number of folding curves, so that eventually there is no folding curve. We start with the definition of excellent developing maps.

\subsection{Excellent developing maps}

\begin{definition}
    Let $f$ be a developing map on a closed surface $\Sigma$. A simple closed curve $\gamma\subset\Sigma$ is called a folding curve of $f$ if there is a neighborhood $N$ of $\gamma$ such that $f$ is a local homeomorphism on $N\backslash\gamma$ but changes orientation when crossing $\gamma$.
\end{definition}

\begin{definition}
    An equidistant line in the Poincar\'e plane $\H$ is a horocycle or a component of the set of points equidistant to a geodesic. If $A\in\mathrm{PSL}_2(\R)$ is a hyperbolic or parabolic element, an $A$-equidistant line is an $A$-invariant equidistant line. An $A$-equidistant curve is the quotient of an $A$-equidistant line by $A$.
\end{definition}

We define the area of a developing map $f$ by
\[A(f)=\int_{\Sigma}|f^*\omega_\H|\]
where $\omega_\H$ is the canonical volume form on $\H$.

Let $\rho:\pi_1(\Sigma)\rightarrow\mathrm{PSL}_2(\R)$ be a representation, where $\Sigma$ is a closed surface. define
\[a_0(\rho)=\inf\{A(f):f\text{ is a developing map for }\rho\}.\]

\begin{definition}\label{def-excellent}
    Let $\Sigma$ be a closed surface of genus at least $2$ and let $\rho:\pi_1(\Sigma)\rightarrow\mathrm{PSL}_2(\R)$ be a representation. Suppose that $\rho$ sends every simple closed curve to a hyperbolic element of $\mathrm{PSL}_2(\R)$. A developing map $f$ for $\rho$ is called $\epsilon$-excellent if it satisfies the following properties:
\begin{itemize}
    \item[(1)] $f$ is branched away from a finite collection of folding curves;
    \item[(2)] $f$ maps every folding curve homeomorphically to an equidistant curve of length $\leq\epsilon$;
    \item[(3)] $A(f)\leq a_0(\rho)+\epsilon.$
\end{itemize}
\end{definition}

The first non-trivial fact is that excellent developing maps do exist.

\begin{proposition}\label{prop-ede}
    Let $\Sigma$ be a closed surface of genus at least $2$. Suppose that $\rho:\pi_1(\Sigma)\rightarrow\mathrm{PSL}_2(\R)$ sends every simple closed curve to a hyperbolic element. Then for any $\epsilon>0$, there exists an $\epsilon$-excellent developing map for $\rho$ with at most $-3\chi(\Sigma)$ folding curves.
\end{proposition}

The next step is to reduce the number of folding curves. Let $\epsilon_0(\Sigma)=10^{100\chi(\Sigma)}$, where $\chi(\Sigma)$ denotes the Euler characteristic of $\Sigma$.

\begin{proposition}\label{prop1}
    Let $\rho:\pi_1(\Sigma)\rightarrow\mathrm{PSL}_2(\R)$ be a representation which sends every simple closed curve to a hyperbolic element. Suppose that there exists an $\epsilon$-excellent developing map with $n$ folding curves $(\epsilon\leq\epsilon_0,n\geq1)$. Then there exists a $100\epsilon$-excellent developing map with at most $n-1$ folding curves.
\end{proposition}

Assuming these two propositions, we now finish the proof of Theorem \ref{thm1}.

\paragraph{Proof of Theorem \ref{thm1}}

Let $\Sigma$ be a closed surface of genus at least $2$ and let $\rho:\pi_1(\Sigma)\rightarrow\mathrm{PSL}_2(\R)$ be a representation which sends every simple closed curve to a hyperbolic element. Let $\epsilon=100^{3\chi(\Sigma)}\epsilon_0$. By Proposition \ref{prop-ede}, there exists an $\epsilon$-excellent developing map $f$ for $\rho$ with at most $-3\chi(\Sigma)$ folding curves. If $f$ has at least one folding curve, we apply Proposition \ref{prop1} to get a $100\epsilon$-developing map $f_1$ with less folding curves. If $f_1$ still has at least one folding curve, we continue this process until eventually we get a $100^k\epsilon$-excellent developing map $f_k$ with no folding curves. Note that since $f$ has at most $-3\chi(\Sigma)$ folding curves, we have $k\leq-3\chi(\Sigma)$. This implies $100^k\epsilon\leq\epsilon_0$, hence Proposition \ref{prop1} is always applicable, justifying this process. $f_k$ is the desired branched developing map.

\subsection{Existence of excellent developing maps}

The proof of Proposition \ref{prop-ede} uses analytic methods, in particular the theory of harmonic maps. A natural approach to finding branched developing maps is by minimizing the energy functional. Recall that every representation $\rho:\pi_1(\Sigma)\rightarrow\mathrm{PSL}_2(\R)$ gives rise to an energy functional on the Teichm\"uller space $\T(\Sigma)$, which associates to every marked Riemann surface $X$ the energy of the unique $\rho$-equivariant harmonic map on $X$. If $f$ is a branched developing map, then the energy of $f$ with respect to the pullback conformal structure is equal to the absolute value of the Euler class, which is the global minimum of the energy functional. Conversely, if $X$ is a critical point of the energy functional, then the $\rho$-equivariant harmonic map on $X$ is branched.

Therefore, it is natural to consider a minimizing sequence $X_n$ of the energy functional. But it is well possible that the sequence diverges in the Teichm\"uller space. If the projection of the sequence in the moduli space $\M(\Sigma)$ stays in a compact subset then we can perturb the harmonic map on $X_n$ for large enough $n$ into a branched developing map. The real trouble comes when $X_n$ has short curves, cutting it into components. The harmonic map can be perturbed to be branched on each component (more precisely, those components on which the harmonic map doesn't converge to a constant), but the orientation may be different for different components. This suggests the notion of excellent developing maps.

The proof of Proposition \ref{prop-ede} is essentially a detailed version of the rough idea above. But more care is needed for those components on which the harmonic maps converges to a constant. We cut these components further into pairs of pants and construct branched developing maps on each of these pairs of pants. This is the content of Subsection \ref{subsec-ede-Z}.

\begin{remark}
	The analytic approach is also explored in the recent paper \cite{DT26}, which has some overlap with Subsection \ref{subsec-ede-find decomp} of the present paper.
\end{remark}

\subsection{Reducing folding curves}

The proof of Proposition \ref{prop1} roughly consists of two steps. Given an excellent developing map with at least one folding curve, we first deform it so that there are two folding curves close to each other. Then we perform certain surgeries either to cancel them or to merge them into one folding curve. Here, to make sense of the term `two folding curves close to each other', we need a metric on $\Sigma$, which is defined as follows.

Given an excellent developing map $f$ on $\Sigma$, the pullback metric $g_f:=f^*g_{\H}$ is a branched hyperbolic metric on each component of $\Sigma\backslash\mathrm{Fold}(f)$ (we denote by $\mathrm{Fold}(f)$ the union of all folding curves of $f$). Let $d_f$ denote the integrated distance function of $g_f$.

Let $\delta_0(\Sigma)=10^{10\chi(\Sigma)}$.

\begin{proposition}\label{prop2}
    Let $\rho:\pi_1(\Sigma)\rightarrow\mathrm{PSL}_2(\R)$ be a representation which sends every simple closed curve to a hyperbolic element. Suppose that there exists an $\epsilon$-excellent developing map with $n$ folding curves $(\epsilon\leq\epsilon_0,n\geq1)$. Then there exists a $2\epsilon$-excellent developing map $f$ with $n$ folding curves, such that there are two folding curves $\gamma_1,\gamma_2$ satisfying $d_{f}(\gamma_1,\gamma_2)\leq\delta_0$.
\end{proposition}

(In particular, an excellent developing map cannot have exactly one folding curve.)

\begin{proposition}\label{prop3}
    Let $\rho:\pi_1(\Sigma)\rightarrow\mathrm{PSL}_2(\R)$ be a representation which sends every simple closed curve to a hyperbolic element. Suppose that there exists an $\epsilon$-excellent developing map $f$ with $n$ folding curves $(\epsilon\leq2\epsilon_0,n\geq1)$, such that two folding curves $\gamma_1,\gamma_2$ satisfy $d_f(\gamma_1,\gamma_2)\leq\delta_0$. Then there exists a $50\epsilon$-excellent developing map with at most $n-1$ folding curves.
\end{proposition}

It is evident that these two propositions imply Proposition \ref{prop1}.

\subsubsection{Moving a folding curve close to another}

\begin{figure}
\centering

\tikzset{every picture/.style={line width=0.75pt}} %set default line width to 0.75pt        

\begin{tikzpicture}[x=0.75pt,y=0.75pt,yscale=-1,xscale=1]
%uncomment if require: \path (0,339); %set diagram left start at 0, and has height of 339

%Curve Lines [id:da007972574231985163] 
\draw    (102,113) .. controls (131,124) and (182,127) .. (218,111) ;
%Curve Lines [id:da032849686058611605] 
\draw    (102,166) .. controls (139,155) and (176,158) .. (220,168) ;
%Curve Lines [id:da036911774948993226] 
\draw    (162,122) .. controls (153,122) and (153,157) .. (161,159) ;
%Curve Lines [id:da44614823045437835] 
\draw  [dash pattern={on 4.5pt off 4.5pt}]  (162,122) .. controls (170,125) and (168,159) .. (161,159) ;
%Curve Lines [id:da7711679075133349] 
\draw    (205,115) .. controls (198,116) and (195,167) .. (206,165) ;
%Curve Lines [id:da4064013358461872] 
\draw  [dash pattern={on 4.5pt off 4.5pt}]  (205,115) .. controls (214,116) and (213,162) .. (206,165) ;
%Curve Lines [id:da06283770489221274] 
\draw    (114,117) .. controls (106,118) and (105,161) .. (114,162) ;
%Curve Lines [id:da394802082760081] 
\draw  [dash pattern={on 4.5pt off 4.5pt}]  (114,117) .. controls (122,119) and (121,161) .. (114,162) ;
%Curve Lines [id:da9981891218842199] 
\draw    (300,113) .. controls (329,120) and (371,125) .. (400,113) ;
%Curve Lines [id:da9613158051832599] 
\draw    (300,163) .. controls (328,153) and (378,158) .. (398,164) ;
%Curve Lines [id:da3987787526885106] 
\draw    (363,120) .. controls (356,121) and (356,159) .. (363,159) ;
%Curve Lines [id:da9587221580786878] 
\draw  [dash pattern={on 4.5pt off 4.5pt}]  (363,120) .. controls (371,120) and (370,159) .. (363,159) ;
%Curve Lines [id:da026618436028021164] 
\draw  [dash pattern={on 4.5pt off 4.5pt}]  (323,117) .. controls (315,120) and (314,158) .. (322,159) ;
%Curve Lines [id:da9617013621573849] 
\draw  [dash pattern={on 4.5pt off 4.5pt}]  (323,117) .. controls (332,118) and (329,157) .. (322,159) ;
%Curve Lines [id:da5586873701387585] 
\draw    (468,115) .. controls (500,123) and (528,126) .. (557,115) ;
%Curve Lines [id:da5890962384861609] 
\draw    (467,164) .. controls (495,154) and (532,160) .. (557,167) ;
%Curve Lines [id:da046558977740303376] 
\draw  [dash pattern={on 4.5pt off 4.5pt}]  (523,122) .. controls (516,124) and (515,158) .. (523,160) ;
%Curve Lines [id:da7930097701039448] 
\draw  [dash pattern={on 4.5pt off 4.5pt}]  (523,122) .. controls (530,122) and (531,161) .. (523,160) ;
%Curve Lines [id:da1836480007386445] 
\draw    (482,118) .. controls (473,119) and (476,162) .. (483,160) ;
%Curve Lines [id:da6836301037277971] 
\draw  [dash pattern={on 4.5pt off 4.5pt}]  (482,118) .. controls (491,120) and (490,159) .. (483,160) ;
%Straight Lines [id:da04952230243721478] 
\draw    (419,140) -- (449,140) ;
\draw [shift={(451,140)}, rotate = 180] [color={rgb, 255:red, 0; green, 0; blue, 0 }  ][line width=0.75]    (10.93,-3.29) .. controls (6.95,-1.4) and (3.31,-0.3) .. (0,0) .. controls (3.31,0.3) and (6.95,1.4) .. (10.93,3.29)   ;

% Text Node
\draw (151,203.4) node [anchor=north west][inner sep=0.75pt]    {$\Sigma $};
% Text Node
\draw (154,162.4) node [anchor=north west][inner sep=0.75pt]    {$\gamma $};
% Text Node
\draw (126,131.4) node [anchor=north west][inner sep=0.75pt]    {$C_{+}$};
% Text Node
\draw (173,132.4) node [anchor=north west][inner sep=0.75pt]    {$C_{-}$};
% Text Node
\draw (347,165.4) node [anchor=north west][inner sep=0.75pt]    {$f( \gamma )$};
% Text Node
\draw (460,170.4) node [anchor=north west][inner sep=0.75pt]    {$f_{1}( \gamma )$};
% Text Node
\draw (319,206.4) node [anchor=north west][inner sep=0.75pt]    {$\mathbb{H} /\rho (\gamma )$};
% Text Node
\draw (478,208.4) node [anchor=north west][inner sep=0.75pt]    {$\mathbb{H} /\rho (\gamma )$};

\end{tikzpicture}
\caption{Cancel the collar neighborhoods $C_\pm$. One can also think of this as sliding the folding curve $\gamma$ in $\H/\rho(\gamma)$.}
\label{fig: intro-cancel collar}
\end{figure}
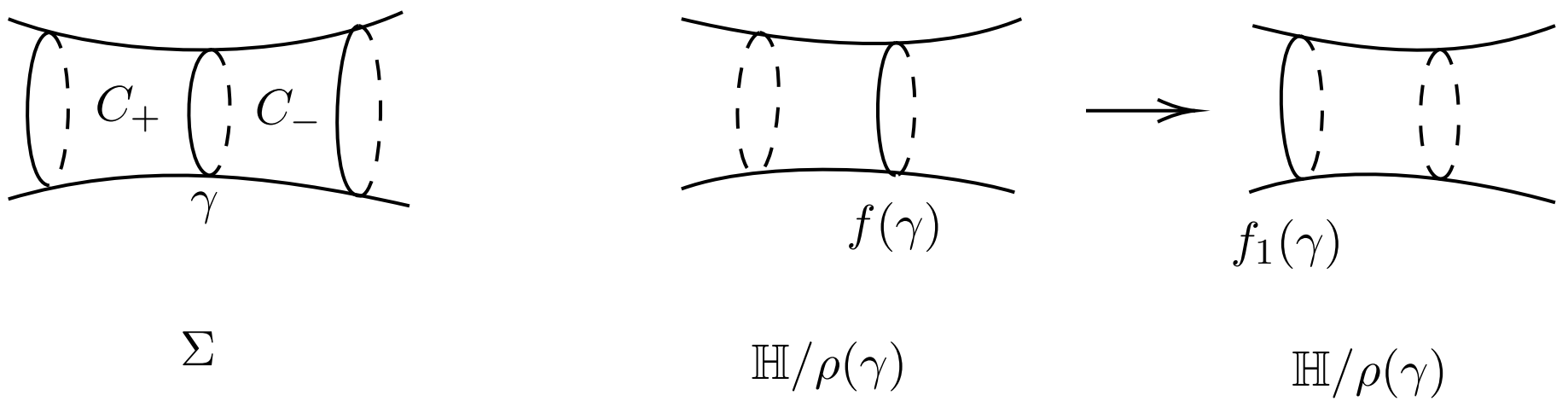

The proof of Proposition \ref{prop2} can be summarized as `moving branch points' and `cancelling collar neighborhoods'. We first describe `cancelling collar neighborhoods'. Let $\gamma$ be a folding curve of $f$, and let $\Sigma_+,\Sigma_-$ be the two components of $\Sigma\backslash\mathrm{Fold}(f)$ adjacent to $\gamma$ (throughout this paper, when we say `a component of $\Sigma\backslash\mathrm{Fold}(f)$', we mean a compact surface with boundary, with all boundary curves added). Denote by $S_\pm=(\Sigma_\pm,g_f)$ the surfaces equipped with the pullback branched hyperbolic metric. Let $l>0$. Let $C_\pm=N^{S_\pm}_l(\gamma)$ be the $l$-neighborhoods of $\gamma$ in $S_\pm$. If (the closure of) $C_\pm$ does not meet any other folding curve, does not contain any branch point, and is topologically an annulus, we call it the \textit{collar neighborhood} of $\gamma$ of length $l$. Equivalently, $C_\pm$ is the diffeomorphic image of $\gamma\times[0,l)$ in $S_\pm$ under the exponential map, for $l<\mathrm{inj}_\gamma$ (see Subsection \ref{subsec-c nbhd}). Suppose that $\gamma$ has a pair of collar neighborhoods $C_\pm$ of length $l$, then $f|_{C_\pm}$ are embeddings, with the same image in $\H/\rho(\gamma)$. Then we can remove the annulus $C_+\cup C_-$ from $\Sigma$ and identify its two boundary curves. The resulting surface $\Sigma_1$ is naturally homeomorphic to $\Sigma$, and $f$ restricts to a developing map $f_1$ on $\Sigma_1$, with the same monodromy as $f$ if we identify $\Sigma_1$ with $\Sigma$ via a natural homeomorphism. This is what we call \textit{removing a pair of collars} (see Subsection \ref{subsec-c sur}). Another way to think about this surgery is to fix the domain $\Sigma$, but deform the developing map $f$ by sliding the folding curve $\gamma$ in $\H/\rho(\gamma)$ (see Figure \ref{fig: intro-cancel collar}).

\begin{figure}
\centering

\tikzset{every picture/.style={line width=0.75pt}} %set default line width to 0.75pt        

\begin{tikzpicture}[x=0.75pt,y=0.75pt,yscale=-1,xscale=1]
%uncomment if require: \path (0,339); %set diagram left start at 0, and has height of 339

%Curve Lines [id:da960642303138076] 
\draw    (84,99) .. controls (110,109) and (164,108) .. (191,100) ;
%Curve Lines [id:da7973259400876378] 
\draw    (87,165) .. controls (124,157) and (158,157) .. (194,164) ;
%Curve Lines [id:da6208695332356338] 
\draw    (140,106) .. controls (150,106) and (148,161) .. (142,158) ;
%Curve Lines [id:da42204607630182456] 
\draw  [dash pattern={on 4.5pt off 4.5pt}]  (140,106) .. controls (133,105) and (134,159) .. (142,158) ;
%Curve Lines [id:da27610831052939644] 
\draw    (84,99) .. controls (98,98) and (97,165) .. (87,165) ;
%Curve Lines [id:da4123341594316564] 
\draw    (84,99) .. controls (74,111) and (75,162) .. (87,165) ;
%Curve Lines [id:da36576010152182536] 
\draw    (301,118) .. controls (328,109) and (360,96) .. (382,71) ;
%Curve Lines [id:da7974158938227603] 
\draw    (191,100) .. controls (199,102) and (201,123) .. (194,128) ;
%Curve Lines [id:da6289184012994868] 
\draw  [dash pattern={on 4.5pt off 4.5pt}]  (191,100) .. controls (186,102) and (187,124) .. (194,128) ;
%Curve Lines [id:da17360509775071298] 
\draw    (194,164) .. controls (199,165) and (204,138) .. (198,137) ;
%Curve Lines [id:da9577295123442662] 
\draw  [dash pattern={on 4.5pt off 4.5pt}]  (198,137) .. controls (194,137) and (188,162) .. (194,164) ;
%Curve Lines [id:da9980058791303438] 
\draw    (303,142) .. controls (333,144) and (372,165) .. (385,179) ;
%Curve Lines [id:da47335327592768073] 
\draw    (399,85) .. controls (376,108) and (376,144) .. (398,161) ;
%Curve Lines [id:da6612248562747625] 
\draw    (365,86) .. controls (392,113) and (378,158) .. (362,162) ;
%Curve Lines [id:da2315385000301483] 
\draw  [dash pattern={on 4.5pt off 4.5pt}]  (365,86) .. controls (372,85) and (385,89) .. (381,126) ;
%Curve Lines [id:da5562770659261326] 
\draw  [dash pattern={on 4.5pt off 4.5pt}]  (362,162) .. controls (372,166) and (383,156) .. (381,126) ;
%Curve Lines [id:da0032765169602089195] 
\draw    (301,118) .. controls (309,118) and (305,142) .. (303,142) ;
%Curve Lines [id:da74078011107087] 
\draw    (301,118) .. controls (298,119) and (296,142) .. (303,142) ;
%Curve Lines [id:da31387754835188275] 
\draw    (147,134) .. controls (163,135) and (182,133) .. (198,137) ;
%Curve Lines [id:da8616896468385119] 
\draw    (147,134) .. controls (162,133) and (180,133) .. (194,128) ;
%Curve Lines [id:da487484912107748] 
\draw    (382,71) .. controls (388,69) and (400,79) .. (399,85) ;
%Curve Lines [id:da7976445783466877] 
\draw  [dash pattern={on 4.5pt off 4.5pt}]  (382,71) .. controls (381,77) and (393,85) .. (399,85) ;
%Curve Lines [id:da8640329063807918] 
\draw    (398,161) .. controls (401,166) and (392,179) .. (385,179) ;
%Curve Lines [id:da6730470314628421] 
\draw  [dash pattern={on 4.5pt off 4.5pt}]  (398,161) .. controls (396,159) and (381,174) .. (385,179) ;

% Text Node
\draw (79,174.4) node [anchor=north west][inner sep=0.75pt]    {$\gamma $};
% Text Node
\draw (134,171.4) node [anchor=north west][inner sep=0.75pt]    {$\gamma _{1}$};
% Text Node
\draw (295,151.4) node [anchor=north west][inner sep=0.75pt]    {$\gamma $};
% Text Node
\draw (355,173.4) node [anchor=north west][inner sep=0.75pt]    {$\gamma _{1}$};
% Text Node
\draw (106,124.4) node [anchor=north west][inner sep=0.75pt]    {$C_{+}$};
% Text Node
\draw (337,119.4) node [anchor=north west][inner sep=0.75pt]    {$C_{+}$};

\end{tikzpicture}

\caption{There are two obstructions to moving $\gamma_1$ further. One is that $\gamma_1$ hits a branch point, the other is that $\gamma_1$ runs into a self-intersection. The latter case is impossible, since then the collar $C_+$ would have area $\geq\frac{1}{2}$, contradicting the area condition $A(f)\leq a_0+\epsilon$.}
\label{fig: intro-two obs}
\end{figure}
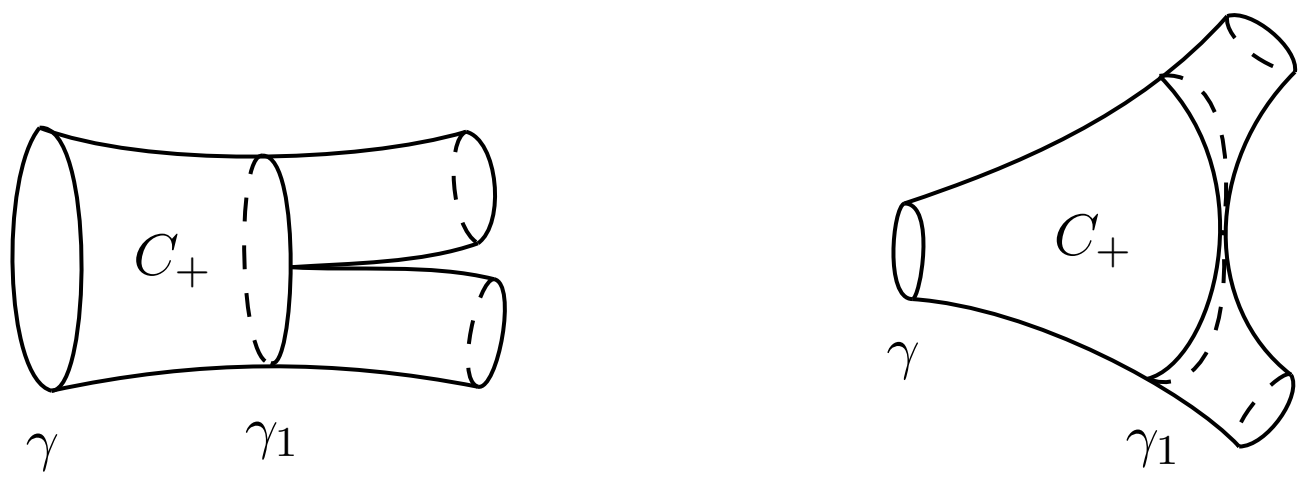

The above surgery can be think of as `moving the folding curve $\gamma$'. The proof of Proposition \ref{prop2} can be described as `moving $\gamma$ close to another folding curve'. More precisely, if there is a pair of collar neighborhoods $C_\pm$ such that the other boundary curve $\gamma_1$ of $C_+$ satisfies $d_{S_+}(\gamma_1,\gamma_2)\leq\delta_0$ for some other boundary curve $\gamma_2$ of $S_+$, then by removing $C_+\cup C_-$ we finish the proof of Proposition \ref{prop2}. So we look for maximal collar neighborhoods. There are two obstructions for extending a collar neighborhood, one is that $\gamma_1$ hits a branch point, the other is that $\gamma_1$ has a self-intersection (see Figure \ref{fig: intro-two obs}).

The latter case is in fact impossible. Indeed, if $\gamma_1$ has a self-intersection (which must be a tangential self-intersection), the assumption that $\rho$ send every simple closed curve to a non-elliptic element would imply that $\gamma_1$ has length $\geq1$. Then by elementary hyperbolic geometry one deduces that the collar $C_+$ has area $\geq\frac{1}{2}$. However, recall that we have the area condition $A(f)\leq a_0+\epsilon$. By removing the pair of collar neighborhoods $C_\pm$ we get a developing map $f_1$ with area $A(f_1)=A(f)-2\mathrm{Area}(C_+)$. This yields a contradiction.

Therefore branch points are the only obstacles. These obstacles are movable: we can move branch points away from a folding curve. The (local) surgery of moving branch points is well studied in the literature (see, for example, \cite{GBS14}). Here, we need to show that branch points can be moved as far as we want. This makes use of the assumption that the monodromy around every simple closed curve is non-elliptic. Indeed, without this assumption, it is possible that the surface runs into a singularity when moving a branch point (a curve around which the monodromy is elliptic might shrink to a point). Thanks to the non-elliptic assumption, we can give a lower bound of $\mathrm{inj}_P$ (Lemma \ref{lem-m2}), which shows that the surface does not degenerate in finite time. This enables us to move branch points as far away from $\gamma$ as possible, until they hit another folding curve; afterwards, we move $\gamma$ close to another folding curve.

\subsubsection{Merging or cancelling two folding curves}

The proof of Proposition \ref{prop3} is sketched as follows. Suppose that we are given two folding curves $\gamma_1,\gamma_2$ that are close to each other. We may assume that $\gamma_1,\gamma_2$ are boundary curves of a common component $\Sigma_0$ of $\Sigma\backslash\mathrm{Fold}(f)$. If $\Sigma_0$ is an annulus, then we can cancel $\gamma_1,\gamma_2$. This is the content of Subsection \ref{subsec-the case of annulus}. In this case the number of folding curves drops by $2$.

If $\Sigma_0$ is not an annulus, the dichotomy in Subsection \ref{subsec-d} shows that we are reduced to the so-called \textit{standard developing maps}. Roughly speaking, these are developing maps such that there is only one branch point $P$ in between $\gamma_1$ and $\gamma_2$. All other branch points are far away so that they do not hinder the surgery that we are going to perform.

\begin{figure}
    \centering

\tikzset{every picture/.style={line width=0.75pt}} %set default line width to 0.75pt        

\begin{tikzpicture}[x=0.75pt,y=0.75pt,yscale=-1,xscale=1]
%uncomment if require: \path (0,339); %set diagram left start at 0, and has height of 339

%Curve Lines [id:da666622023766586] 
\draw    (109,78) .. controls (145,85) and (179,99) .. (224,98) ;
%Curve Lines [id:da6181252373194502] 
\draw    (102,117) .. controls (118,117) and (135,125) .. (149,126) ;
%Curve Lines [id:da10178230336950478] 
\draw    (102,162) .. controls (141,154) and (194,152) .. (226,153) ;
%Curve Lines [id:da17044184277633434] 
\draw    (100,133) .. controls (116,132) and (133,129) .. (149,126) ;
%Curve Lines [id:da22024209633382819] 
\draw    (199,152) .. controls (189,148) and (192,100) .. (203,97) ;
%Curve Lines [id:da06898580522427633] 
\draw  [dash pattern={on 4.5pt off 4.5pt}]  (203,97) .. controls (214,104) and (209,150) .. (199,152) ;
%Curve Lines [id:da9815631865088188] 
\draw    (119,159) .. controls (111,155) and (110,135) .. (117,131) ;
%Curve Lines [id:da04562398337057261] 
\draw  [dash pattern={on 4.5pt off 4.5pt}]  (117,131) .. controls (123,134) and (125,154) .. (119,159) ;

%Curve Lines [id:da6814886235615397] 
\draw    (332,90) .. controls (362,103) and (389,112) .. (421,113) ;
%Curve Lines [id:da8480156606696256] 
\draw    (332,164) .. controls (358,154) and (390,161) .. (414,164) ;
%Curve Lines [id:da08499291951581067] 
\draw    (331,139) .. controls (346,138) and (357,137) .. (367,133) ;
%Curve Lines [id:da12452379739706632] 
\draw    (325,124) .. controls (343,126) and (356,128) .. (367,133) ;
%Curve Lines [id:da2600295721004914] 
\draw    (351,127) .. controls (345,124) and (352,100) .. (358,101) ;
%Curve Lines [id:da9472243029772387] 
\draw  [dash pattern={on 4.5pt off 4.5pt}]  (358,101) .. controls (365,104) and (356,129) .. (351,127) ;

%Straight Lines [id:da7156521388679785] 
\draw    (263,132) .. controls (264.67,130.33) and (266.33,130.33) .. (268,132) .. controls (269.67,133.67) and (271.33,133.67) .. (273,132) .. controls (274.67,130.33) and (276.33,130.33) .. (278,132) -- (282,132) -- (290,132) ;
\draw [shift={(292,132)}, rotate = 180] [color={rgb, 255:red, 0; green, 0; blue, 0 }  ][line width=0.75]    (10.93,-3.29) .. controls (6.95,-1.4) and (3.31,-0.3) .. (0,0) .. controls (3.31,0.3) and (6.95,1.4) .. (10.93,3.29)   ;
%Curve Lines [id:da6637273711185834] 
\draw    (112,223) .. controls (140,235) and (178,247) .. (213,249) ;
%Curve Lines [id:da6495986661654883] 
\draw    (98,298) .. controls (131,290) and (179,291) .. (210,293) ;
%Curve Lines [id:da3938750382672619] 
\draw    (100,268) .. controls (117,263) and (148,266) .. (161,266) ;
%Curve Lines [id:da4241257224429157] 
\draw    (104,253) .. controls (121,254) and (146,261) .. (161,266) ;
%Curve Lines [id:da9919140504123398] 
\draw    (123,294) .. controls (114,292) and (116,268) .. (121,265) ;
%Curve Lines [id:da5282679195423794] 
\draw  [dash pattern={on 4.5pt off 4.5pt}]  (121,265) .. controls (127,266) and (129,289) .. (123,294) ;
%Curve Lines [id:da5391572626186971] 
\draw    (123,255) .. controls (119,250) and (125,229) .. (132,230) ;
%Curve Lines [id:da6525715521496342] 
\draw  [dash pattern={on 4.5pt off 4.5pt}]  (132,230) .. controls (134,234) and (129,255) .. (123,255) ;
%Straight Lines [id:da4246850805764263] 
\draw    (265,273) .. controls (266.67,271.33) and (268.33,271.33) .. (270,273) .. controls (271.67,274.67) and (273.33,274.67) .. (275,273) -- (279,273) -- (287,273) ;
\draw [shift={(289,273)}, rotate = 180] [color={rgb, 255:red, 0; green, 0; blue, 0 }  ][line width=0.75]    (10.93,-3.29) .. controls (6.95,-1.4) and (3.31,-0.3) .. (0,0) .. controls (3.31,0.3) and (6.95,1.4) .. (10.93,3.29)   ;
%Curve Lines [id:da6392385010077283] 
\draw    (336,237) .. controls (357,246) and (384,256) .. (413,250) ;
%Curve Lines [id:da3066227848933243] 
\draw    (333,296) .. controls (359,288) and (386,292) .. (410,292) ;
%Curve Lines [id:da3088811975010527] 
\draw    (334,267) .. controls (343,264) and (352,269) .. (358,269) ;
%Curve Lines [id:da7527677603546244] 
\draw    (332,274) .. controls (342,273) and (351,272) .. (358,269) ;
%Curve Lines [id:da3256521587686505] 
\draw    (391,251) .. controls (383,256) and (382,286) .. (389,292) ;
%Curve Lines [id:da1667774779947736] 
\draw  [dash pattern={on 4.5pt off 4.5pt}]  (391,251) .. controls (397,255) and (398,288) .. (389,292) ;

% Text Node
\draw (111,163.4) node [anchor=north west][inner sep=0.75pt]    {$\gamma _{1}$};
% Text Node
\draw (191,157.4) node [anchor=north west][inner sep=0.75pt]    {$\gamma _{2}$};
% Text Node
\draw (359,74.4) node [anchor=north west][inner sep=0.75pt]    {$\beta $};
% Text Node
\draw (115,297.4) node [anchor=north west][inner sep=0.75pt]    {$\gamma _{1}$};
% Text Node
\draw (129,209.4) node [anchor=north west][inner sep=0.75pt]    {$\gamma _{2}$};
% Text Node
\draw (381,298.4) node [anchor=north west][inner sep=0.75pt]    {$\beta $};
% Text Node
\draw (159,111.4) node [anchor=north west][inner sep=0.75pt]    {$+$};
% Text Node
\draw (403,260.4) node [anchor=north west][inner sep=0.75pt]    {$+$};
% Text Node
\draw (190,260.4) node [anchor=north west][inner sep=0.75pt]    {$+$};
% Text Node
\draw (327,100.4) node [anchor=north west][inner sep=0.75pt]    {$+$};
% Text Node
\draw (93,137.4) node [anchor=north west][inner sep=0.75pt]    {$-$};
% Text Node
\draw (364,260.4) node [anchor=north west][inner sep=0.75pt]    {$-$};
% Text Node
\draw (102,229.4) node [anchor=north west][inner sep=0.75pt]    {$-$};
% Text Node
\draw (96,273.4) node [anchor=north west][inner sep=0.75pt]    {$-$};
% Text Node
\draw (218,116.4) node [anchor=north west][inner sep=0.75pt]    {$-$};
% Text Node
\draw (382,126.4) node [anchor=north west][inner sep=0.75pt]    {$-$};

\end{tikzpicture}

    \caption{$\gamma_1,\gamma_2$ merge into one folding curve $\beta$. `$+$' and `$-$' stand for orientation.}
    \label{fig: merge}
\end{figure}
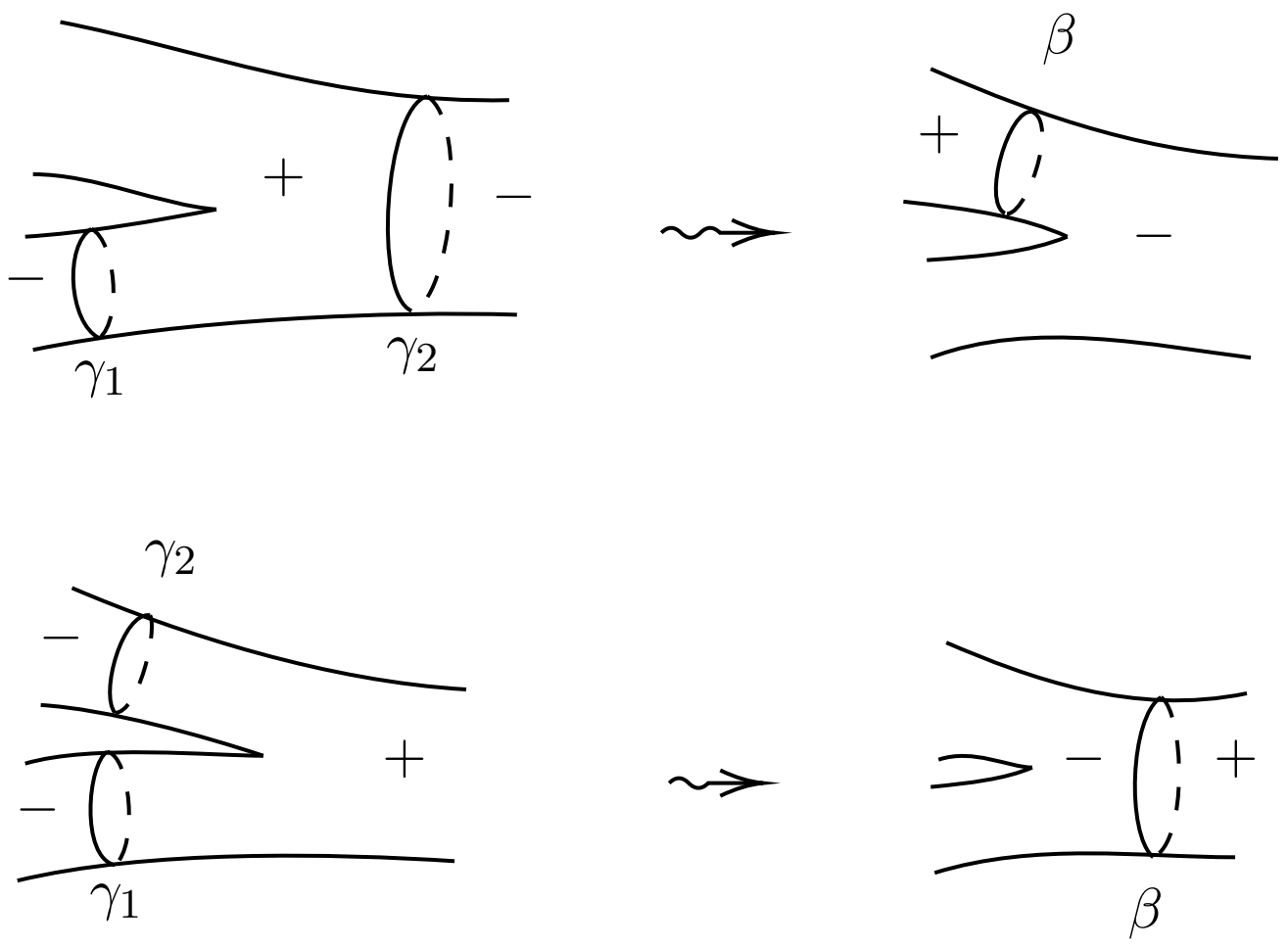

It remains to handle standard developing maps. Proposition \ref{prop-rm} further reduce standard developing maps to three models described in Section \ref{sec-three models}. We then perform specific surgeries for these three models to merge the folding curves $\gamma_1,\gamma_2$ into one. The surgeries for different models are simply different geometric realizations of the topological pictures suggested by Figure \ref{fig: merge}.

\begin{remark}
In fact all the arguments in this paper go through with `hyperbolic' replaced by `hyperbolic or parabolic', except the material of Section 8.
\end{remark}

\subsection{Conventions and some preliminaries}

Throughout this paper $\Sigma$ denotes a compact, connected, oriented surface (possibly with boundary), $\chi(\Sigma)$ denotes its Euler characteristic. $\epsilon_0(\Sigma)=10^{100\chi(\Sigma)},\delta_0(\Sigma)=10^{10\chi(\Sigma)},\delta_1(\Sigma)=10^{5\chi(\Sigma)}$ (so $\epsilon_0\ll\delta_0\ll\delta_1$).

We let $g_\H$ denote the hyperbolic metric on $\H$, $d_\H$ denote the distance function. If $(X,d)$ is a metric space and $A\subset X$ is a subset, $r>0$, define $N_r(A)=\{x\in X:d(x,A)<r\}$. If $A_1,A_2\subset X$ are two subsets, define $d(A_1,A_2)=\inf\{d(x_1,x_2):x_1\in A_1,x_2\in A_2\}$. If $B$ is a finite set we denote by $\#B$ the cardinality of $B$.

All maps are piecewise $C^1$, in particular all curves are piecewise $C^1$. When we say a \textit{homeomorphism}, we implicitly mean a piecewise $C^1$ homeomorphism.

\paragraph{Curves}

For a curve $\alpha:[0,l]\rightarrow\Sigma$ (or $\alpha:[0,l]\rightarrow\H$), we let $\overrightarrow\alpha$ denote the tangent vector $\alpha'(0)$ at $\alpha(0)$, and let $\overleftarrow\alpha$ denote the tangent vector $-\alpha'(l)$ at $\alpha(l)$. A curve is called \textit{simple} if it is injective. We usually do not distinguish a simple curve from its image. If $\Sigma$ is given a metric we denote by $l_\alpha$ the length of $\alpha$. Given two curves $\alpha:[0,l_1]\rightarrow\Sigma,\beta:[0,l_2]\rightarrow\Sigma$ such that $\alpha(l_1)=\beta(0)$, we define $\alpha*\beta:[0,l_1+l_2]\rightarrow\Sigma$ by
\begin{equation}
(\alpha*\beta)(t)=\begin{cases}&\alpha(t),\ t\leq l_1;\\ &\beta(t-l_1),\ t>l_1.\end{cases}
\end{equation}
If $\alpha,\beta$ are simple closed curves on a surface, $i(\alpha,\beta)$ denotes their geometric intersection number.

\paragraph{Developing maps and their restrictions}

An equivalent definition of a developing map is a continuous section of the fiber bundle $E^\rho\rightarrow\Sigma$, with fiber $\H$, defined by $E^\rho=\widetilde{\Sigma}\times\H/\pi_1(\Sigma)$, where $\pi_1(\Sigma)$ acts on $\widetilde{\Sigma}\times\H$ by $g\cdot(x,y)=(g\cdot x,\rho(g)\cdot y)$. Therefore for any subsurface $\Sigma_0\subset\Sigma$, we can restrict a developing map $f$ to $\Sigma_0$, denoted $f|_{\Sigma_0}$. In particular if $\Sigma_0$ is simply connected then $f$ well-defines a map $\Sigma_0\rightarrow\H$ (up to an automorphism of $\H$). If $\Sigma_0$ is an annulus with core curve $\gamma$ such that $\rho(\gamma)$ is hyperbolic or parabolic, then $f$ well-defines a map $\Sigma_0\rightarrow\H/\rho(\gamma)$.

More generally, given any continuous map $\Sigma_1\rightarrow\Sigma$, we can pull back $f$ to a developing map on $\Sigma_1$. A typical example is that $\Sigma_1$ is obtained by cutting $\Sigma$ along a simple closed curve. In this case we still call the pullback of $f$ to $\Sigma_1$ the restriction of $f$ to $\Sigma_1$.

Conversely, if we have a decomposition $\Sigma=\Sigma_1\cup\cdots\cup\Sigma_k$ of $\Sigma$ into subsurfaces with disjoint interior, and a section $f_i$ of $E^\rho|_{\Sigma_i}$ (i.e., a developing map for $\rho|_{\pi_1(\Sigma_i)}$ on $\Sigma_i$) is given for each $i$, such that $f_i$ and $f_j$ agree on $\partial\Sigma_i\cap\partial\Sigma_j$ for every pair $i,j$, then $\{f_i\}_{i=1}^k$ fit together into a developing map on $\Sigma$.

\paragraph{Branched hyperbolic metric}

A branched hyperbolic metric $g$ on $\Sigma$ is a hyperbolic metric away from finitely many branch points, at which $g$ is locally given by the pullback of $g_\H$ via $z\mapsto z^k$ for some integer $k\geq2$. We usually denote by $S=(\Sigma,g)$ a surface equipped with a branched hyperbolic metric, and denote by $d$ the integrated distance function of $g$.

For two tangent vectors $\overrightarrow{\alpha},\overrightarrow{\beta}$ at a smooth point $P$, $\measuredangle(\overrightarrow{\alpha},\overrightarrow{\beta})$ takes value in $[0,2\pi)$ unless otherwise stated. If $P$ is a branch point of multiplicity $m-1$, $\measuredangle(\overrightarrow{\alpha},\overrightarrow{\beta})$ takes value in $[0,2m\pi)$.

The set of branch points is denoted $\mathrm{Br}(g),\mathrm{Br}(S)$ or $\mathrm{Br}(f)$.

\paragraph{Metric v.s. developing map}

There is a one-to-one correspondence between branched hyperbolic metrics on $\Sigma$ and branched developing maps on $\Sigma$ up to automorphisms of $\H$, given by $f\mapsto g_f:=f^*g_\H$ and $g\mapsto\mathrm{dev}_g$. Therefore we usually think of them as the same structure on a surface $\Sigma$. It is useful to keep in mind two points of view.

\paragraph{Geodesics}

A geodesic on a branched hyperbolic surface is a curve without interior branch point, which is a geodesic in the hyperbolic metric. Geodesics are parametrized by arc length unless otherwise stated. By a geodesic loop, we mean a geodesic $\alpha:[0,l]\rightarrow S$ such that $\alpha(0)=\alpha(l)$. We call $\alpha(0)$ its base point, and write $\mathring\alpha=\alpha\backslash\{\alpha(0)\}$.

\paragraph{Surgery} By a surgery, we mean the following. Start with a developing map $f$ on a surface $\Sigma$, use some cut and paste procedure to a produce a surface $\Sigma_1$ that is naturally homeomorphic to $\Sigma$, and a developing map $f_1$ on $\Sigma_1$ with the same monodromy as $f$ if we identify $\Sigma_1$ with $\Sigma$ via a natural homeomorphism. Here naturally homeomorphic means that $\Sigma_1$ and $\Sigma$ are homeomorphic, with a well-defined isotopy class of homeomorphism.

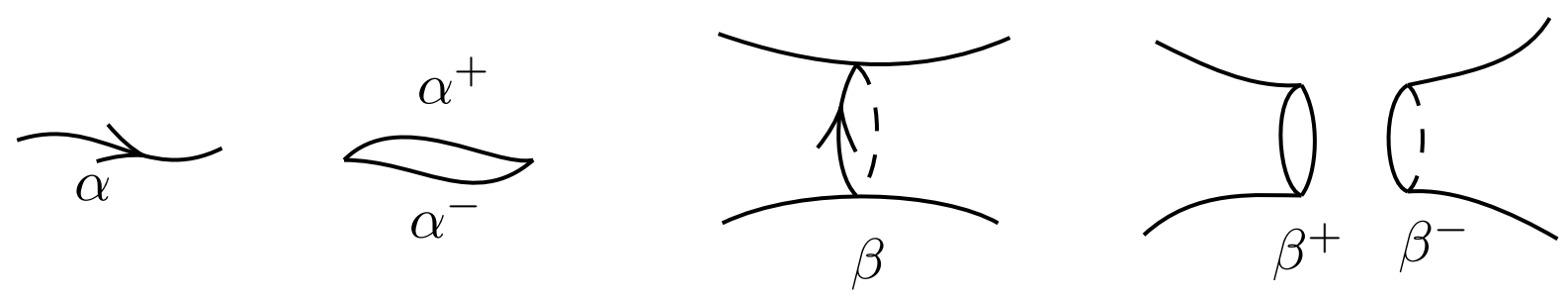
\begin{figure}
    \centering
    \tikzset{every picture/.style={line width=0.75pt}} %set default line width to 0.75pt        

\begin{tikzpicture}[x=0.75pt,y=0.75pt,yscale=-1,xscale=1]
%uncomment if require: \path (0,339); %set diagram left start at 0, and has height of 339

%Curve Lines [id:da30262398501133647] 
\draw    (101,132) .. controls (122,125) and (133,144) .. (153,134) ;
\draw [shift={(133.06,135.86)}, rotate = 197] [color={rgb, 255:red, 0; green, 0; blue, 0 }  ][line width=0.75]    (10.93,-4.9) .. controls (6.95,-2.3) and (3.31,-0.67) .. (0,0) .. controls (3.31,0.67) and (6.95,2.3) .. (10.93,4.9)   ;
%Curve Lines [id:da3934500994504231] 
\draw    (184,137) .. controls (198,123) and (222,139) .. (232,137) ;
%Curve Lines [id:da0678612247681617] 
\draw    (184,137) .. controls (203,137) and (217,150) .. (232,137) ;

%Curve Lines [id:da7482049809158364] 
\draw    (279,105) .. controls (302,113) and (328,117) .. (353,106) ;
%Curve Lines [id:da24884812008701762] 
\draw    (280,153) .. controls (301,143) and (336,145) .. (350,153) ;
%Curve Lines [id:da8588447675875823] 
\draw    (314,146) .. controls (308,140) and (308,122) .. (314,113) ;
\draw [shift={(310.12,123.53)}, rotate = 95.76] [color={rgb, 255:red, 0; green, 0; blue, 0 }  ][line width=0.75]    (10.93,-4.9) .. controls (6.95,-2.3) and (3.31,-0.67) .. (0,0) .. controls (3.31,0.67) and (6.95,2.3) .. (10.93,4.9)   ;
%Curve Lines [id:da493833693343228] 
\draw  [dash pattern={on 4.5pt off 4.5pt}]  (314,113) .. controls (321,118) and (321,141) .. (314,146) ;
%Curve Lines [id:da23729716898463227] 
\draw    (390,107) .. controls (402,113) and (414,119) .. (427,118) ;
%Curve Lines [id:da11386986481162686] 
\draw    (387,156) .. controls (400,144) and (416,146) .. (427,146) ;
%Curve Lines [id:da510751646102548] 
\draw    (427,118) .. controls (419,120) and (421,144) .. (427,146) ;
%Curve Lines [id:da9281379636924915] 
\draw    (427,118) .. controls (432,126) and (431,141) .. (427,146) ;
%Curve Lines [id:da806179148738447] 
\draw    (454,118) .. controls (468,115) and (483,113) .. (490,101) ;
%Curve Lines [id:da3671433294276293] 
\draw    (454,145) .. controls (466,144) and (480,150) .. (492,157) ;
%Curve Lines [id:da11206543783192291] 
\draw    (454,118) .. controls (447,122) and (448,143) .. (454,145) ;
%Curve Lines [id:da172362725365507] 
\draw  [dash pattern={on 4.5pt off 4.5pt}]  (454,118) .. controls (459,122) and (459,141) .. (454,145) ;

% Text Node
\draw (311,155.4) node [anchor=north west][inner sep=0.75pt]    {$\beta $};
% Text Node
\draw (418,151.4) node [anchor=north west][inner sep=0.75pt]    {$\beta ^{+}$};
% Text Node
\draw (450,149.4) node [anchor=north west][inner sep=0.75pt]    {$\beta ^{-}$};
% Text Node
%\draw (186,199) node [anchor=north west][inner sep=0.75pt]   [align=left] {Figure. Cut along a curve};
% Text Node
\draw (114,139.4) node [anchor=north west][inner sep=0.75pt]    {$\alpha $};
% Text Node
\draw (201,109.4) node [anchor=north west][inner sep=0.75pt]    {$\alpha ^{+}$};
% Text Node
\draw (199,143.4) node [anchor=north west][inner sep=0.75pt]    {$\alpha ^{-}$};

\end{tikzpicture}
    \caption{Cut along a curve}
    \label{fig: cut along a curve}
\end{figure}

\paragraph{Cut along a curve} If we cut a surface $\Sigma$ along a simple curve $\alpha$, the resulting surface $\Sigma_1$ has a new boundary curve composed of two copies of $\alpha$. We denote them by $\alpha^+$ and $\alpha^-$, where $\Sigma_1$ lies on the left of $\alpha^+$. Similarly, if we cut along a simple loop $\beta$, the resulting surface $\Sigma_1$ has two new boundary curves, both are copies of $\beta$. We denote them by $\beta^+$ and $\beta^-$, where $\Sigma_1$ lies on the left of $\beta^+$ (see Figure \ref{fig: cut along a curve}).

\subsection*{Organization}

Section 2 collects some elementary results on hyperbolic geometry. In section 3 we introduce the notion of collar, which plays a central role in this paper. We prove Proposition \ref{prop-ede} in Section 4, which is largely independent of other sections. The surgery of moving branch points is introduced in Section 5, and Section 6 contains some applications of moving branch points for our purpose. Proposition \ref{prop2} is proved in section 7. The rest of the paper is devoted to the proof of Proposition \ref{prop3}.

\subsection*{Acknowledgement}

I would like to sincerely thank Vladimir Markovic. I learned about Goldman’s conjecture from his lectures. We discussed the problem together, and he helped me a lot in the writing of this paper.

I would like to thank Jinyi Wang for teaching me how to draw figures.

The work was done in Qiuzhen College, Tsinghua University. I would like to thank Qiuzhen College for providing a stimulating research environment.

zbMATH and arXiv were used to search for references, and DeepSeek was used to correct grammatical errors.

\section{Some hyperbolic geometries}

In this section we collect some elementary lemmas on hyperbolic geometry. The spirit behind all these lemmas is that, in small scale, the hyperbolic space looks like the Euclidean space, and a hyperbolic element of $\mathrm{PSL}_2(\R)$ with small translation length looks like a translation.

Throughout this section, $d=d_\H$ denotes the hyperbolic distance on $\H$. $\epsilon$ denotes a positive real number less than $\frac{1}{10}$. 

\begin{lemma}\label{lem1}
    The curvature of an equidistant line has absolute value $\leq1$.
\end{lemma}

\begin{proof}
    The curvature is maximized by horocycles, in which case it is equal to $1$.
\end{proof}

\begin{lemma}\label{lem2}
    Let $\gamma\subset\H$ be an equidistant line. Let $P,Q$ be points on $\gamma$ of distance $\leq1$. Let $\nu_P,\nu_Q$ be unit normal vectors of $\gamma$ at $P,Q$ (pointing toward the same side of $\gamma$). Suppose $A\in\mathrm{PSL}_2(\R)$ maps $P$ to $Q$ and $A_*$ maps $\nu_P$ to $-\nu_Q$. Then $A$ is elliptic.
\end{lemma}

\begin{proof}
    The tangent geodesic lines of $\gamma$ at $P,Q$ intersect at a point, which is fixed by $A$.
\end{proof}

\begin{lemma}\label{lem-h11}
    Let $P,Q\in\H$. Let $\gamma$ be the geodesic through $P,Q$ and let $\alpha$ be an equidistant line through $P,Q$. Suppose that $d(P,Q)=\epsilon$. Then $\gamma$ and $\alpha$ form an angle less than $\frac{\epsilon}{2}$ at $P$.
\end{lemma}

\begin{proof}
    For fixed $P,Q$, this angle is maximized when $\alpha$ is a horocycle. The lemma then follows from straight forward computation.
\end{proof}

\begin{lemma}\label{lem-h1}
    Let $A\in\mathrm{PSL}_2(\R)$ be a hyperbolic or parabolic element. Suppose that $P,Q\in\H$ satisfies $A(P)=Q$ and $d(P,Q)=\epsilon$. Let $\mathcal{P}:T_P\H\rightarrow T_Q\H$ be the parallel transport along the geodesic connecting $P,Q$. Then $\mathcal{P}^{-1}\circ A_*$ acts on $T_P\H$ by a rotation of an angle less than $\epsilon$.
\end{lemma}

\begin{proof}
    Let $\gamma$ be the geodesic connecting $P,Q$ and let $\alpha$ be the $A$-equidistant line passing through $P,Q$. Then $\mathcal{P}^{-1}\circ A_*$ acts on $T_P\H$ by a rotation of twice the angle between $\gamma$ and $\alpha$ at $P$. The lemma then follows from the previous lemma.
\end{proof}

\begin{lemma}\label{lem-h2}
    Let $\gamma_1,\gamma_2:(-\epsilon,\epsilon)\rightarrow\H$ be $C^1$ curves such that $\gamma_1(0)\neq\gamma_2(0)$. Define $s(t)=d(\gamma_1(t),\gamma_2(t))$. Let $\alpha:[0,s(0)]\rightarrow\H$ be the unit speed geodesic from $\gamma_1(0)$ to $\gamma_2(0)$. Then
    \[s'(0)=\alpha'(s(0))\cdot\gamma_2'(0)-\alpha'(0)\cdot\gamma_1'(0)\]
    (where the dot product between tangent vectors is the Riemannian inner product). 
\end{lemma}

\begin{proof}
    Transform $\H$ to the strip $S=\{z\in\C:-\frac{\pi}{2}<\mathrm{Re}(z)<\frac{\pi}{2}\}$ so that $\gamma_1(0),\gamma_2(0)$ both lie on the vertical line $\{x=0\}$. $S$ is equipped with the hyperbolic metric $\frac{dx^2+dy^2}{\cos^2x}$, which dominates the Euclidean metric $dx^2+dy^2$. Denote the Euclidean distant between $\gamma_1(t)$ and $\gamma_2(t)$ by $s_E(t)=d_E(\gamma_1(t),\gamma_2(t))$. We have $s(t)\geq s_E(t),s(0)=s_E(0)$, while both functions are $C^1$. Therefore $s'(0)=s_E'(0)$. It is easy to compute $s'_E(0)$, equating it with $\alpha'(s(0))\cdot\gamma_2'(0)-\alpha'(0)\cdot\gamma_1'(0)$, thus proving the lemma.
\end{proof}

\begin{lemma}\label{lem-h3}
    Let $A\in\mathrm{PSL}_2(\R)$ be a non-elliptic element and let $\gamma:(-\epsilon,\epsilon)\rightarrow\H$ be a unit speed curve. Define $s(t)=d(\gamma(t),A(\gamma(t)))$. Then
    \[|s'(t)|\leq s(t)\]
\end{lemma}

\begin{proof}
    Let $\alpha:[0,s(0)]\rightarrow\H$ be the unit speed geodesic from $\gamma(0)$ to $A(\gamma(0))$. Then by Lemma \ref{lem-h2} we have
    \begin{align}\label{eq3}
        s'(0)=\alpha'(s(0))\cdot A_*(\gamma'(0))-\alpha'(0)\cdot\gamma'(0)=\alpha'(0)\cdot[\mathcal{P}^{-1}A_*(\gamma'(0))-\gamma'(0)],
    \end{align}
    where $\mathcal{P}$ is the parallel transport along $\alpha$. Note that $\gamma'(0)$ has unit length by assumption. Thus, by Lemma \ref{lem-h1}, we have $\|\mathcal{P}^{-1}A_*(\gamma'(0))-\gamma'(0)\|\leq s(0)$. Substituting this into (\ref{eq3}) and using $\|\alpha'(0)\|=1$ completes the proof.
\end{proof}

\begin{lemma}\label{lem-h4}
    Let $A_1,A_2\in\mathrm{PSL}_2(\R)$ be hyperbolic or parabolic elements such that $A_2\circ A_1^{-1}$ is also hyperbolic or parabolic. Let $P\in\H$ and let $\alpha_1,\alpha_2,\beta_1,\beta_2$ be the geodesic segments from $P$ to $A_1(P),A_2(P),A_1^{-1}(P),A_2^{-1}(P)$ respectively. Suppose that $\measuredangle(\overrightarrow{\alpha_2},\overrightarrow{\alpha_1}),\measuredangle(\overrightarrow{\beta_1},\overrightarrow{\beta_2})\in[0,\frac{1}{2}\pi)$, and that $l_{\alpha_2}\leq l_{\alpha_1}\leq\epsilon$. Let $T$ be the geodesic triangle with sides $\alpha_1,\alpha_2$ (assume that $\alpha_1,\alpha_2$ are not colinear). Then the interior angles of $T$ at the two vertices of $\alpha_1$ are both less than $2\epsilon$ (see Figure \ref{fig: lem-h4}).
\end{lemma}

\begin{proof}
    Let $\alpha_3$ be the geodesic segment from $A_1(P)$ to $A_2(P)$. Then $T$ has sides $\alpha_1,\alpha_2,\alpha_3$. Let $\theta_i$ be the interior angle of $T$ opposite to $\alpha_i$. Then we need to prove $\theta_2,\theta_3<2\epsilon$.
    
    Let $A_3=A_2\circ A_1^{-1}$. Then $A_3$ maps $A_1(P)$ to $A_2(P)$. Let $\phi=\measuredangle(A_{3*}(\overrightarrow{\alpha_3}),\overleftarrow{\alpha_3})$. By Lemma \ref{lem-h1}, we have $\phi\in(\pi-l_{\alpha_3},\pi+l_{\alpha_3})$. We show that
    \begin{align}\label{eq-h41}
        \measuredangle(\overrightarrow{\beta_1},\overrightarrow{\beta_2})\equiv\theta_1+\theta_2+\phi\mod2\pi.
    \end{align}
    Indeed, we have
    \begin{align}\label{eq-h42}
        \theta_1=\measuredangle(\overleftarrow{\alpha_3},\overleftarrow{\alpha_2}),\ \theta_2=\measuredangle(\overleftarrow{\alpha_1},\overrightarrow{\alpha_3})=\measuredangle(A_{3*}(\overleftarrow{\alpha_1}),A_{3*}(\overrightarrow{\alpha_3})).
    \end{align}
    Note that $\overrightarrow{\beta_1}=(A_{1}^{-1})_*(\overleftarrow{\alpha_1}),\ \overrightarrow{\beta_2}=(A_{2}^{-1})_*(\overleftarrow{\alpha_2})$, while $\measuredangle((A_{1}^{-1})_*(\overleftarrow{\alpha_1}),(A_{2}^{-1})_*(\overleftarrow{\alpha_2}))=\measuredangle((A_{2}\circ A_{1}^{-1})_*(\overleftarrow{\alpha_1}),\overleftarrow{\alpha_2})=\measuredangle(A_{3*}(\overleftarrow{\alpha_1}),\overleftarrow{\alpha_2})$, we get
    \begin{align}\label{eq-h43}
        \measuredangle(\overrightarrow{\beta_1},\overrightarrow{\beta_2})=\measuredangle(A_{3*}(\overleftarrow{\alpha_1}),\overleftarrow{\alpha_2}).
    \end{align}
    Combining (\ref{eq-h42}),(\ref{eq-h43}) and the definition of $\phi$ proves (\ref{eq-h41}).
    
    By assumption $\measuredangle(\overrightarrow{\beta_1},\overrightarrow{\beta_2})\in[0,\frac{1}{2}\pi)$. So (\ref{eq-h41}) implies $\theta_1+\theta_2+\phi\in[2\pi,\frac{5}{2}\pi)$, in particular $\theta_1+\theta_2+\phi\geq2\pi$. Combining this with $\phi\in(\pi-l_{\alpha_3},\pi+l_{\alpha_3})$, we get $\theta_1+\theta_2>\pi-l_{\alpha_3}$. Since $\theta_1+\theta_2+\theta_3\leq\pi$, we get
    \begin{align}\label{eq-h44}
        \theta_3<l_{\alpha_3}.
    \end{align}
    Since $l_{\alpha_3}\leq 2\epsilon$, (\ref{eq-h44}) implies $\theta_3<2\epsilon$. By the hyperbolic sine formula $\frac{\sin\theta_2}{\sinh l_{\alpha_2}}=\frac{\sin\theta_3}{\sinh l_{\alpha_3}}$, we deduce that $\frac{\sin\theta_2}{\sinh l_{\alpha_2}}=\frac{\sin\theta_3}{\sinh l_{\alpha_3}}\leq1$. Thus $\sin\theta_2\leq\sinh\epsilon$. The condition $l_{\alpha_2}\leq l_{\alpha_1}$ implies $\theta_2\leq\frac{1}{2}\pi$. We conclude that $\theta_2<2\epsilon$.
\end{proof}

\begin{figure}
    \centering

\tikzset{every picture/.style={line width=0.75pt}} %set default line width to 0.75pt        

\begin{tikzpicture}[x=0.75pt,y=0.75pt,yscale=-1,xscale=1]
%uncomment if require: \path (0,339); %set diagram left start at 0, and has height of 339

%Curve Lines [id:da8832926113251929] 
\draw    (103,152) .. controls (130,147) and (175,137) .. (204,134) ;
\draw [shift={(146.74,143.41)}, rotate = 349.05] [color={rgb, 255:red, 0; green, 0; blue, 0 }  ][line width=0.75]    (10.93,-4.9) .. controls (6.95,-2.3) and (3.31,-0.67) .. (0,0) .. controls (3.31,0.67) and (6.95,2.3) .. (10.93,4.9)   ;
%Curve Lines [id:da9020899014295596] 
\draw    (204,134) .. controls (242,139) and (270,145) .. (306,156) ;
\draw [shift={(261.59,144.04)}, rotate = 192.48] [color={rgb, 255:red, 0; green, 0; blue, 0 }  ][line width=0.75]    (10.93,-4.9) .. controls (6.95,-2.3) and (3.31,-0.67) .. (0,0) .. controls (3.31,0.67) and (6.95,2.3) .. (10.93,4.9)   ;
%Curve Lines [id:da9477037679263722] 
\draw    (131,119) .. controls (150,121) and (184,129) .. (204,134) ;
\draw [shift={(161.13,124.11)}, rotate = 11.89] [color={rgb, 255:red, 0; green, 0; blue, 0 }  ][line width=0.75]    (10.93,-4.9) .. controls (6.95,-2.3) and (3.31,-0.67) .. (0,0) .. controls (3.31,0.67) and (6.95,2.3) .. (10.93,4.9)   ;
%Curve Lines [id:da46413501767106047] 
\draw    (204,134) .. controls (226,128) and (257,122) .. (277,118) ;
\draw [shift={(246.3,124.19)}, rotate = 168.14] [color={rgb, 255:red, 0; green, 0; blue, 0 }  ][line width=0.75]    (10.93,-4.9) .. controls (6.95,-2.3) and (3.31,-0.67) .. (0,0) .. controls (3.31,0.67) and (6.95,2.3) .. (10.93,4.9)   ;
%Curve Lines [id:da948938074734543] 
\draw    (103,152) .. controls (114,139) and (123,128) .. (131,119) ;
\draw [shift={(120.86,130.76)}, rotate = 130.1] [color={rgb, 255:red, 0; green, 0; blue, 0 }  ][line width=0.75]    (10.93,-4.9) .. controls (6.95,-2.3) and (3.31,-0.67) .. (0,0) .. controls (3.31,0.67) and (6.95,2.3) .. (10.93,4.9)   ;
%Curve Lines [id:da09597361426463347] 
\draw    (161.04,191.42) .. controls (187.79,184.21) and (222.58,188.33) .. (245,201) ;
\draw [shift={(159,192)}, rotate = 343.5] [color={rgb, 255:red, 0; green, 0; blue, 0 }  ][line width=0.75]    (10.93,-4.9) .. controls (6.95,-2.3) and (3.31,-0.67) .. (0,0) .. controls (3.31,0.67) and (6.95,2.3) .. (10.93,4.9)   ;
%Curve Lines [id:da6838247230316516] 
\draw    (173.99,78.05) .. controls (195.56,89.11) and (220.66,90.7) .. (242,81) ;
\draw [shift={(172,77)}, rotate = 28.61] [color={rgb, 255:red, 0; green, 0; blue, 0 }  ][line width=0.75]    (10.93,-4.9) .. controls (6.95,-2.3) and (3.31,-0.67) .. (0,0) .. controls (3.31,0.67) and (6.95,2.3) .. (10.93,4.9)   ;

% Text Node
\draw (151,147.4) node [anchor=north west][inner sep=0.75pt]    {$\alpha _{1}$};
% Text Node
\draw (160,107.4) node [anchor=north west][inner sep=0.75pt]    {$\alpha _{2}$};
% Text Node
\draw (245,149.4) node [anchor=north west][inner sep=0.75pt]    {$\beta _{1}$};
% Text Node
\draw (242,105.4) node [anchor=north west][inner sep=0.75pt]    {$\beta _{2}$};
% Text Node
\draw (130,125.4) node [anchor=north west][inner sep=0.75pt]    {$T$};
% Text Node
\draw (96,119.4) node [anchor=north west][inner sep=0.75pt]    {$\alpha _{3}$};
% Text Node
\draw (204,54.4) node [anchor=north west][inner sep=0.75pt]    {$A_{2}$};
% Text Node
\draw (194,199.4) node [anchor=north west][inner sep=0.75pt]    {$A_{1}$};
% Text Node
\draw (198,141.4) node [anchor=north west][inner sep=0.75pt]    {$P$};

\end{tikzpicture}

    \caption{Lemma \ref{lem-h4}}
    \label{fig: lem-h4}
\end{figure}

\begin{lemma}\label{lem-h6}
    Let $A\in\mathrm{PSL}_2(\R)$ be a hyperbolic or parabolic element, and let $\gamma:(-\infty,+\infty)\rightarrow\H$ be a unit speed geodesic. Let $\alpha_t$ be the $A$-equidistant line passing through $\gamma(t)$. Let $\theta(t)\in[0,\frac{\pi}{2}]$ be the angle formed by $\gamma$ and $\alpha_t$ at $\gamma(t)$. Then
    \[|\theta'(t)|\leq2.\]
\end{lemma}

\begin{proof}
    We assume $A$ is hyperbolic. The parabolic case follows from taking limit. We may assume that $A$ is of the form $A(z)=az\ (a>0)$, so that the $A$-equidistant lines are rays $\R_+\cdot z_0\ (z_0\in\H)$. We may also assume that $\gamma$ does not intersect with $\partial\H$ at $\infty$, so that $\gamma$ is a semi-circle with center $x_0$ on the real axis. Then
    \begin{align}\label{eq-h1}
        |\theta'(t)|\leq\left|\frac{d}{dt}\arg\gamma(t)\right|+\left|\frac{d}{dt}\arg(\gamma(t)-x_0)\right|.
    \end{align}
    In general for any curve $\beta:(-\delta,\delta)\rightarrow\C\backslash\{0\}$ we have
    \[\left|\frac{d}{dt}\arg\beta(t)\right|\leq\frac{|\beta'(t)|}{|\beta(t)|}.\]
    Note that $\frac{|\beta'(t)|}{|\beta(t)|}\leq\|\beta'(t)\|_{g_\H}$. So if $\beta$ is unit speed, we have
    \[\left|\frac{d}{dt}\arg\beta(t)\right|\leq1.\]
    Combining this with inequality (\ref{eq-h1}) proves the lemma.
\end{proof}

The following is an immediate corollary of the above lemma.

\begin{lemma}\label{cor-h1}
    Let $\gamma:[0,l]\rightarrow\H$ be a geodesic of length $l$. Let $A\in\mathrm{PSL}_2(\R)$ be a hyperbolic or parabolic element. Let $\alpha_t$ be the $A$-equidistant line passing through $\gamma(t)$. Let $\theta\in[0,\frac{\pi}{2}]$ be the angle formed by $\alpha_0$ and $\gamma$ at $\gamma(0)$. Suppose that $\theta>2l$. Then for every $t\in[0,l]$, $\gamma$ intersects with $\alpha_t$ only at $\gamma(t)$ and the angle between them is $\in[\theta-2l,\theta+2l]$. In particular $\gamma$ projects onto a simple geodesic in $\H/A$.
\end{lemma}

The proofs of the following two lemmas are left to the reader.

\begin{lemma}\label{lem-h5}
    Let $T\subset\H$ be a geodesic triangle with vertices $P,Q,R$. Suppose that $\angle P\leq\frac{\pi}{2}$ and $d(P,R)\leq d(P,Q)\leq\frac{1}{100}$. Let $A\in\mathrm{PSL}_2(\R)$ be a hyperbolic or parabolic element such that $A(P)=Q$. Then the map $\overline{T}\backslash\{P\}\rightarrow\H/A$ is injective.
\end{lemma}

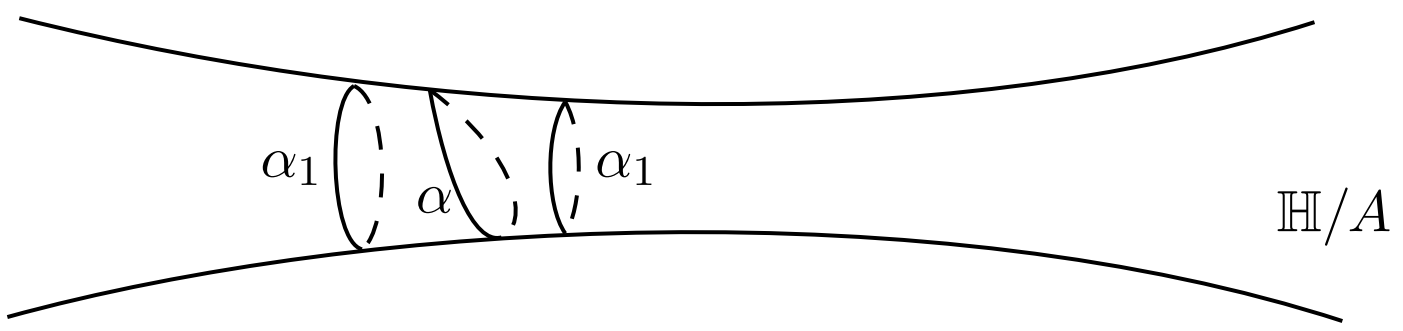
\begin{figure}
    \centering
    \tikzset{every picture/.style={line width=0.75pt}} %set default line width to 0.75pt        

\begin{tikzpicture}[x=0.75pt,y=0.75pt,yscale=-1,xscale=1]

%Curve Lines [id:da30016437422309605] 
\draw    (86,86) .. controls (189,112) and (317,117) .. (411,87) ;
%Curve Lines [id:da022088886394549423] 
\draw    (83,161) .. controls (174,136) and (316,129) .. (418,162) ;
%Curve Lines [id:da4863251931174214] 
\draw    (189,104) .. controls (193,126) and (200,143) .. (207,141) ;
%Curve Lines [id:da16755367956633227] 
\draw  [dash pattern={on 4.5pt off 4.5pt}]  (189,104) .. controls (212,121) and (214,138) .. (207,141) ;
%Curve Lines [id:da8352773416363617] 
\draw    (170,103) .. controls (163,107) and (164,142) .. (172,144) ;
%Curve Lines [id:da6454798083122758] 
\draw  [dash pattern={on 4.5pt off 4.5pt}]  (170,103) .. controls (179,107) and (179,139) .. (172,144) ;
%Curve Lines [id:da45899327054895367] 
\draw    (223,107) .. controls (218,114) and (218,133) .. (223,140) ;
%Curve Lines [id:da11055017895627073] 
\draw  [dash pattern={on 4.5pt off 4.5pt}]  (223,107) .. controls (228,116) and (227,134) .. (223,140) ;

% Text Node
\draw (184,127) node [anchor=north west][inner sep=0.75pt]   [align=left] {$ \alpha $};
% Text Node
\draw (400,127) node [anchor=north west][inner sep=0.75pt]   [align=left] {$ \H/A $};
% Text Node
\draw (145,118) node [anchor=north west][inner sep=0.75pt]   [align=left] {$ \alpha_1$};
% Text Node
\draw (229,118) node [anchor=north west][inner sep=0.75pt]   [align=left] {$\alpha_1$};
% Text Node
%\draw (80,178) node [anchor=north west][inner sep=0.75pt]   [align=left] {Figure. Equidistant curves close to $\displaystyle \overline{\alpha} $ on either sides.};

\end{tikzpicture}
    \caption{Equidistant curves close to $\alpha$ on either sides.}
    \label{fig: lem-h7}
\end{figure}

\begin{lemma}\label{lem-h7}
    Let $A\in\mathrm{PSL}_2(\R)$ be a hyperbolic or parabolic element. Let $P\in\H$ and let $\alpha$ be the geodesic segment from $P$ to $A(P)$. Then $\alpha$ projects onto a simple geodesic loop in $\H/A$. Moreover, if $l_{\alpha}\leq\epsilon$, then there exists an $A$-equidistant curve $\alpha_1$ of length $\leq\frac{5}{4}\epsilon$ on either side of $\alpha$ (disjoint from $\alpha$), such that the annulus bounded by $\alpha$ and $\alpha_1$ has area $\leq\frac{1}{10}\epsilon$ (see Figure \ref{fig: lem-h7}).
\end{lemma}

\section{Collars and surgery}

\subsection{Collars}

\begin{definition}\label{def-collar}
    A collar is a (non-complete) hyperbolic annulus such that both of its boundary curves are equidistant curves.
\end{definition}

\begin{definition}
    The length of a collar is the distance between its two boundary curves.
\end{definition}

\begin{definition}
    A collar is called $\epsilon$-thin if both of its boundary curves have length $\leq\epsilon$.
\end{definition}

\begin{lemma}\label{lem-c1}
    The area of an $\epsilon$-thin collar is $\leq2\epsilon$.
\end{lemma}

\begin{proof}
    By Lemma \ref{lem1} and the Gauss-Bonnet theorem.
\end{proof}

\begin{lemma}\label{lem-c2}
    Let $C$ be a collar with boundary curves $\gamma_1,\gamma_2$. Suppose that $l_{\gamma_1}\geq 2l_{\gamma_2}$. Then $\mathrm{Area}(C)>\frac{1}{3}l_{\gamma_1}$.
\end{lemma}

\begin{proof}
    Let $\kappa_1,\kappa_2\geq0$ be the curvature of $\gamma_1,\gamma_2$ respectively. Then the Gauss-Bonnet theorem implies $\mathrm{Area}(C)=\kappa_1l_{\gamma_1}\pm\kappa_2l_{\gamma_2}$. Note that $l_{\gamma_1}\geq l_{\gamma_2}$ implies $\kappa_1\geq\kappa_2$. We get $\mathrm{Area}(C)\geq\kappa_1(l_{\gamma_1}-l_{\gamma_2})\geq\frac{1}{2}\kappa_1l_{\gamma_1}$. By simple computation, $l_{\gamma_1}\geq 2l_{\gamma_2}$ implies $\kappa_1\geq\frac{\sqrt{3}}{2}>\frac{2}{3}$, proving the lemma.
\end{proof}

\begin{lemma}\label{lem-c3}
    Let $C$ be a collar of length $\leq \frac{1}{3}$. Then its boundary curves $\gamma_1,\gamma_2$ satisfy $l_{\gamma_1}<2l_{\gamma_2}$.
\end{lemma}

\begin{proof}
    Follows from Lemma \ref{lem-c2}.
\end{proof}

\begin{lemma}\label{lem-c31}
    let $C$ be a collar of length $\leq 1$. Then its boundary curves $\gamma_1,\gamma_2$ satisfy $l_{\gamma_1}<8l_{\gamma_2}$.
\end{lemma}

\begin{proof}
    Apply Lemma \ref{lem-c3} three times.
\end{proof}

\subsection{Collar neighborhoods}\label{subsec-c nbhd}

Let $S=(\Sigma,g)$ be branched hyperbolic surface, with distance function $d$. Let $\gamma$ be a boundary curve of $S$. Define
\[r_\gamma=d\left(\gamma,\mathrm{Br}(S)\cup(\partial S\backslash \gamma)\right).\]

Suppose that the holonomy of $g$ around $\gamma$ is hyperbolic or parabolic, and $\gamma$ is an equidistant curve (which means that $\mathrm{dev}_g$ maps $\widetilde{\gamma}$ homeomorphically to an equidistant line). Let $N\gamma$ denote the normal bundle of $\gamma$ in $\Sigma$. We trivialize $N\gamma$ by the inward pointing unit normal vector field, giving rise to an isomorphism $N\gamma\cong \gamma\times\R$. The exponential map
\[\exp_\gamma:\gamma\times[0,r_\gamma]\rightarrow S,\]
is a local diffeomorphism. Define
\[\mathrm{inj}_\gamma=\sup\left\{r:\exp_\gamma|_{\gamma\times[0,r]}\text{ is an embedding}\right\}.\]
For $r<\mathrm{inj}_\gamma$, we call $\exp_\gamma(\gamma\times[0,r))\subset S$ the \textit{collar neighborhood} of $\gamma$ of length $r$. These are collars in the sense of Definition \ref{def-collar}.

We equip $\gamma\times[0,\mathrm{inj}_\gamma]$ with the metric $\exp_\gamma^*(g)$, and denote it by $C_\gamma$ (the `maximal collar neighborhood'). Denote by $\gamma_1$ its boundary curve $\gamma\times\{\mathrm{inj}_\gamma\}$.

\begin{lemma}\label{lem-c4}
    Suppose that the holonomy of $g$ around every simple closed curve is non-elliptic. If $\mathrm{inj}_\gamma<r_\gamma$, then $l_{\gamma_1}\geq1$.
\end{lemma}

\begin{proof}
    If $\mathrm{inj}_\gamma<r_\gamma$, then $\exp_\gamma$ is a local diffeomorphism on a neighborhood of $\gamma_1$. Therefore $\exp_{\gamma}$ is not injective on $\gamma_1$. Otherwise, there would exist $r>\mathrm{inj}_{\gamma}$ such that $\exp_\gamma$ is injective on $\gamma\times[0,r)$, a contradiction. So there are two points $P\neq Q\in \gamma_1$ such that $\exp_\gamma(P)=\exp_\gamma(Q)$. Since $\exp_\gamma$ is injective on $\gamma\times[0,\mathrm{inj}_\gamma)$, we deduce that $d\exp_\gamma(\nu_P)=-d\exp_\gamma(\nu_Q)$, where $\nu_P,\nu_Q$ are the outward pointing unit normal vectors of $C_\gamma$ at $P,Q$ respectively. Let $\alpha\subset \gamma_1$ be an arc connecting $P,Q$. Then $\exp_\gamma(\alpha)$ represents a simple closed curve. If $l_{\gamma_1}<1$, then $l_{\alpha}<1$. By Lemma \ref{lem2}, the holonomy of $g$ around $\exp_\gamma(\alpha)$ is an elliptic element, a contradiction. Therefore $l_{\gamma_1}\geq1$.
\end{proof}

\begin{lemma}\label{cor-c2}
    Suppose that the holonomy of $g$ around every simple closed curve is non-elliptic. If $l_\gamma\leq\frac{1}{10}$, then $\mathrm{inj}_\gamma\geq\min\{1,r_\gamma\}$.
\end{lemma}

\begin{proof}
    If $\mathrm{inj}_\gamma<\min\{1,r_\gamma\}$, then Lemma \ref{lem-c4} implies $l_{\gamma_1}\geq1$. By Lemma \ref{lem-c31}, we get $l_\gamma>\frac{1}{8}$, a contradiction.
\end{proof}

\begin{lemma}\label{cor-c1}
    Suppose that the holonomy of $g$ around every simple closed curve is non-elliptic. If $l_\gamma\leq\frac{1}{10}$ and $r_\gamma\leq1$, then $\mathrm{inj}_\gamma=r_\gamma$.
\end{lemma}

\begin{proof}
    Follows from Lemma \ref{cor-c2}.
\end{proof}

\subsection{Attaching and removing collars}\label{subsec-c sur}

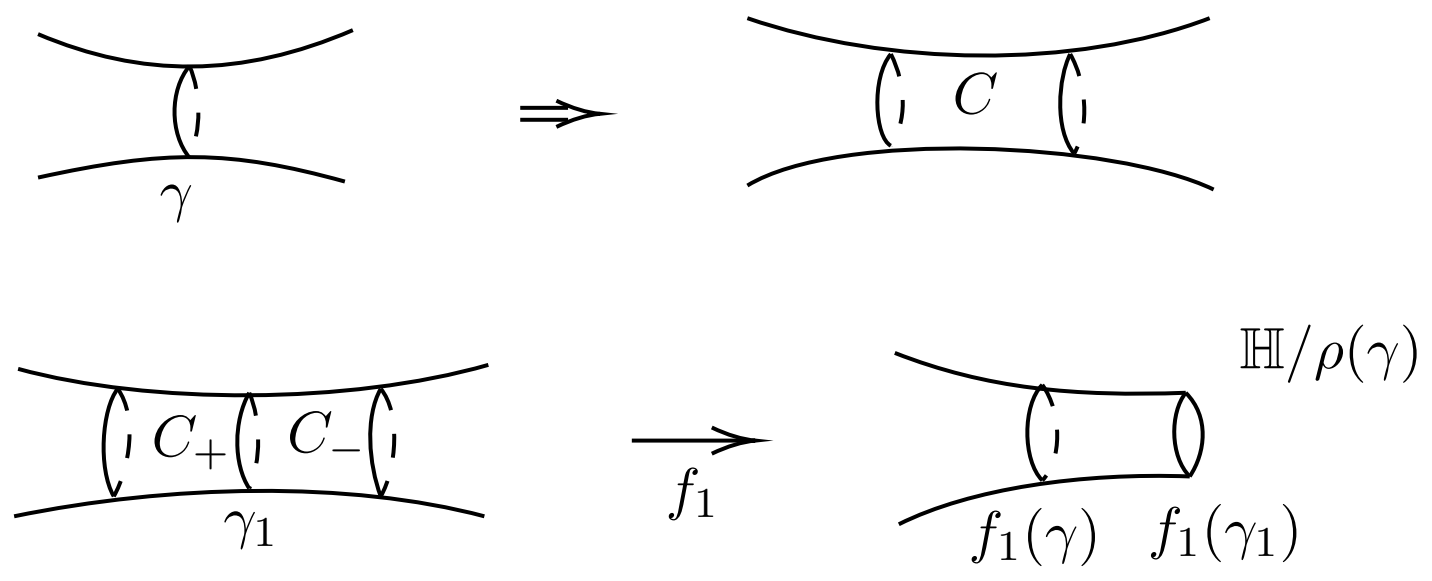
\begin{figure}
    \centering
    \tikzset{every picture/.style={line width=0.75pt}} %set default line width to 0.75pt        

\begin{tikzpicture}[x=0.75pt,y=0.75pt,yscale=-1,xscale=1]

%Curve Lines [id:da373184715545794] 
\draw    (95,63) .. controls (121,74) and (146,74) .. (174,62) ;
%Curve Lines [id:da5545936937804016] 
\draw    (95,99) .. controls (127,92) and (143,92) .. (172,100) ;
%Curve Lines [id:da853955770727319] 
\draw    (133,71) .. controls (128,77) and (128,88) .. (133,94) ;
%Curve Lines [id:da3918009507893817] 
\draw  [dash pattern={on 4.5pt off 4.5pt}]  (133,71) .. controls (136,78) and (136,87) .. (133,94) ;
%Straight Lines [id:da1231714741712372] 
\draw    (216,81.5) -- (228,81.5)(216,84.5) -- (228,84.5) ;
\draw [shift={(236,83)}, rotate = 180] [color={rgb, 255:red, 0; green, 0; blue, 0 }  ][line width=0.75]    (10.93,-3.29) .. controls (6.95,-1.4) and (3.31,-0.3) .. (0,0) .. controls (3.31,0.3) and (6.95,1.4) .. (10.93,3.29)   ;
%Curve Lines [id:da534843954130566] 
\draw    (273,59) .. controls (308,71) and (355,72) .. (389,59) ;
%Curve Lines [id:da7155232385573721] 
\draw    (273,101) .. controls (296,87) and (365,90) .. (390,102) ;
%Curve Lines [id:da41543328903092347] 
\draw    (309,68) .. controls (304,74) and (305,89) .. (309,91) ;
%Curve Lines [id:da7295209275773421] 
\draw  [dash pattern={on 4.5pt off 4.5pt}]  (309,68) .. controls (313,76) and (313,85) .. (309,91) ;
%Curve Lines [id:da2655257644605834] 
\draw    (354,68) .. controls (351,74) and (350,87) .. (355,93) ;
%Curve Lines [id:da8712046785600831] 
\draw  [dash pattern={on 4.5pt off 4.5pt}]  (354,68) .. controls (358,75) and (359,86) .. (355,93) ;

%Curve Lines [id:da31453066144899056] 
\draw    (90,147) .. controls (123,156) and (173,156) .. (208,146) ;
%Curve Lines [id:da03216370813730318] 
\draw    (89,184) .. controls (128,176) and (173,175) .. (207,184) ;
%Curve Lines [id:da2950805727080126] 
\draw    (115,152) .. controls (111,157) and (110,172) .. (114,179) ;
%Curve Lines [id:da00020784271227247686] 
\draw  [dash pattern={on 4.5pt off 4.5pt}]  (115,152) .. controls (120,159) and (118,172) .. (114,179) ;
%Curve Lines [id:da45587270536812663] 
\draw    (181,152) .. controls (177,159) and (178,171) .. (181,179) ;
%Curve Lines [id:da978573432276569] 
\draw  [dash pattern={on 4.5pt off 4.5pt}]  (181,152) .. controls (186,158) and (185,172) .. (181,179) ;
%Curve Lines [id:da028378571391129426] 
\draw    (148,153) .. controls (144,160) and (144,171) .. (148,177) ;
%Curve Lines [id:da9494612243277243] 
\draw  [dash pattern={on 4.5pt off 4.5pt}]  (148,153) .. controls (151,160) and (151,170) .. (148,177) ;

%Straight Lines [id:da8777063222227718] 
\draw    (244,165) -- (273,165) ;
\draw [shift={(275,165)}, rotate = 180] [color={rgb, 255:red, 0; green, 0; blue, 0 }  ][line width=0.75]    (10.93,-3.29) .. controls (6.95,-1.4) and (3.31,-0.3) .. (0,0) .. controls (3.31,0.3) and (6.95,1.4) .. (10.93,3.29)   ;
%Curve Lines [id:da9490333808083112] 
\draw    (310,143) .. controls (333,152) and (355,154) .. (383,153) ;
%Curve Lines [id:da8481581052746817] 
\draw    (311,186) .. controls (329,177) and (357,173) .. (384,174) ;
%Curve Lines [id:da848767622592576] 
\draw    (383,153) .. controls (379,158) and (379,169) .. (384,174) ;
%Curve Lines [id:da8717448311213095] 
\draw    (383,153) .. controls (389,159) and (388,168) .. (384,174) ;
%Curve Lines [id:da8891650106036576] 
\draw    (347,151) .. controls (342,156) and (342,170) .. (347,175) ;
%Curve Lines [id:da05196678475981065] 
\draw  [dash pattern={on 4.5pt off 4.5pt}]  (347,151) .. controls (352,158) and (352,170) .. (347,175) ;

% Text Node
\draw (124,99.4) node [anchor=north west][inner sep=0.75pt]    {$\gamma$};
% Text Node
\draw (323,71.4) node [anchor=north west][inner sep=0.75pt]    {$C $};
% Text Node
\draw (122,157.4) node [anchor=north west][inner sep=0.75pt]    {$C _{+}$};
% Text Node
\draw (156,156.4) node [anchor=north west][inner sep=0.75pt]    {$C _{-}$};
% Text Node
\draw (140,181.4) node [anchor=north west][inner sep=0.75pt]    {$\gamma_1$};
% Text Node
\draw (251,170.4) node [anchor=north west][inner sep=0.75pt]    {$f_1$};
% Text Node
\draw (395,134.4) node [anchor=north west][inner sep=0.75pt]    {$\mathbb{H} /\rho ( \gamma)$};
% Text Node
\draw (372,179.4) node [anchor=north west][inner sep=0.75pt]    {$f_1( \gamma_1)$};
% Text Node
\draw (327,180.4) node [anchor=north west][inner sep=0.75pt]    {$f_1( \gamma)$};
% Text Node
%\draw (140,221) node [anchor=north west][inner sep=0.75pt]   [align=left] {Figure. Insert $\displaystyle C $ and extend $\displaystyle f$.};

\end{tikzpicture}
    \caption{Insert $C$ and extend $f$.}
    \label{fig: sur-add collar}
\end{figure}

Let $f$ be an $\epsilon$-excellent developing map and let $\gamma\subset\Sigma$ be a folding curve. We introduce a surgery, called \textit{attaching a pair of collars of length $l$ along $\gamma$} (see Figure \ref{fig: sur-add collar}). Cut $\Sigma$ along $\gamma$ and insert an annulus $C$. The resulting surface $\Sigma_1$ is naturally homeomorphic to $\Sigma$. Let $\gamma_1\subset C$ be a core curve, which divides $C$ into two annuli $C_+,C_-$. We extend $f$ to $\Sigma_1$, mapping both $C_+$ and $C_-$ homeomorphically onto a collar $C_1\subset\H/\rho(\gamma)$ of length $l$, so that it is a local homeomorphism along the two boundary curves of $C$, and $\gamma_1$ becomes a new folding curve. This produces an $\epsilon_1$-excellent developing map $f_1$ on $\Sigma_1$, with the same monodromy as $f$ if we identify $\Sigma_1$ with $\Sigma$ via a natural homeomorphism, where $\epsilon_1=\max\{l_{f_1(\gamma_1)},\epsilon+2\mathrm{Area}(C_1)\}$. We always identify $\Sigma$ with $\Sigma_1$ via a homeomorphism $\Sigma\rightarrow\Sigma_1$ which is the identity away from a small neighborhood of $\gamma$, and maps $\gamma$ onto $\gamma_1$. Therefore, $f_1$ is an excellent developing map on $\Sigma$, with $\mathrm{Fold}(f_1)=\mathrm{Fold}(f)$.

The proof of the following lemma is left to the reader.

\begin{lemma}\label{lem-c6}
    Let $f$ be an $\epsilon$-excellent developing map. Then the developing map obtained by attaching a pair of collars of length $\leq\frac{1}{4n}$ along $n$ folding curves is $2\epsilon$-excellent.
\end{lemma}

The inverse surgery, \textit{removing a pair of collars of length $l$ at $\gamma$}, is as follows. Let $S_+,S_-$ be the two components of $\Sigma\backslash\mathrm{Fold}(f)$ adjacent to $\gamma$, equipped with the metric $g_f$. Let $0<l<\mathrm{inj}_\gamma$, where $\mathrm{inj}_\gamma=\min\{\mathrm{inj}_\gamma(S_+),\mathrm{inj}_\gamma(S_-)\}$. Let $C_+\subset S_+,\ C_-\subset S_-$ be the collar neighborhoods of $\gamma$ of length $l$. Denote by $\gamma_\pm$ the boundary curves of $C_\pm$ other than $\gamma$. Then $f$ maps $C_{\pm}$ homeomorphically to the same collar in $\H/\rho(\gamma)$, in particular $f$ maps $\gamma_\pm$ homeomorphically to the same equidistant curve, giving rise to an identification $\gamma_+\cong\gamma_-$. So we can remove $C_+\cup C_-$ from $\Sigma$ and identify $\gamma_+$ with $\gamma_-$ via the above identification. The resulting surface $\Sigma_1$ is naturally homeomorphic to $\Sigma$ and $f$ restricts to an $\epsilon_1$-excellent developing map on $\Sigma_1$, where $\epsilon_1=\min\{\epsilon,l_{\gamma_+}\}$. We always identify $\Sigma$ with $\Sigma_1$ via a homeomorphism $\Sigma\rightarrow\Sigma_1$ which is the identity way from a small neighborhood of $\overline{C_+\cup C_-}$, and maps $\gamma$ onto $\gamma_1$, where $\gamma_1$ is the curve on $\Sigma_1$ obtained by gluing $\gamma_{\pm}$. Therefore, $f_1$ is an excellent developing map on $\Sigma$, with $\mathrm{Fold}(f_1)=\mathrm{Fold}(f)$.

Define $r_\gamma(f)=\min\{r_\gamma(S_+),r_\gamma(S_-)\}$.

\begin{lemma}\label{lem-c7}
    Let $f$ be an $\epsilon$-excellent developing map and let $\gamma\subset\Sigma$ be a folding curve. Then $\mathrm{inj}_\gamma=r_\gamma$. Moreover, for every $0<l<\mathrm{inj}_\gamma$, the collar neighborhoods of $\gamma$ of length $l$ are $2\epsilon$-thin.
\end{lemma}

\begin{proof}
    Let $0<l<\mathrm{inj}_\gamma$ be arbitrary. Let $C_\pm\subset S_\pm$ be the collar neighborhoods of length $l$. We have
    \begin{align}\label{eq1}
        \mathrm{Area}(C_\pm)\leq\frac{1}{2}\epsilon.
    \end{align}
    Indeed, $\mathrm{Area}(C_+)=\mathrm{Area}(C_-)$ and by removing this pair of collar neighborhoods we get a developing map $f_1$ with area
    \[A(f_1)=A(f)-2\mathrm{Area}(C_+).\]
    Since by assumption $A(f)\leq a_0+\epsilon\leq A(f_1)+\epsilon$, this proves (\ref{eq1}).

    By (\ref{eq1}) and Lemma \ref{lem-c2}, the boundary curves of $C_\pm$ have length $<2\epsilon$. So $C_\pm$ are $2\epsilon$-thin. This proves the moreover part of the lemma. The equality $\mathrm{inj}_\gamma=r_\gamma$ then follows from Lemma \ref{lem-c4}.
\end{proof}

One can also attach collars along equidistant boundary curves of a branched hyperbolic surface. Let $S$ be a branched hyperbolic surface, with holonomy representation $\rho$, and let $\gamma$ be a boundary curve. Suppose that $\rho(\gamma)$ is hyperbolic or parabolic, and $\mathrm{dev}_S$ maps $\gamma$ homeomorphically onto an equidistant curve $\gamma_0$ in $\H/\rho(\gamma)$. Let $C\subset\H/\rho(\gamma)$ be a collar of length $l$ bounded by equidistant curves $\gamma_0,\gamma_1$, such that outward pointing vectors of $S$ along $\gamma$ are mapped (by $\mathrm{dev}_S$) to inward pointing vectors of $C$ along $\gamma_0$. Then we can attach $C$ to $S$, identifying $\gamma$ with $\gamma_0$ via $\mathrm{dev}_S$. We call this \textit{attaching a collar of length $l$ along $\gamma$}.

\section{Excellent developing map exists}

In this section we prove Proposition \ref{prop-ede}. The idea is to decompose $\Sigma$ into subsurfaces and construct branched developing maps on each of them.

\subsection{Build up excellent developing maps from subsurfaces}

\begin{definition}
    A developing map on a closed surface $\Sigma$ is called pseudo $\epsilon$-excellent if it satisfies all the conditions for being $\epsilon$-excellent except the area condition (condition (3) in Definition \ref{def-excellent}).
\end{definition}

\begin{definition}
    A developing map on a compact surface with boundary is called $\epsilon$-good if it is branched and maps every boundary curve homeomorphically to an equidistance curve of length $\leq\epsilon$.
\end{definition}

The next lemma tells us how to build a pseudo excellent developing map from good developing maps on subsurfaces.

\begin{lemma}\label{lem-ede1}
    Let $\Sigma$ be a closed surface. Let $\gamma_1,\gamma_2,\cdots,\gamma_n\subset\Sigma$ be a collection of disjoint simple closed curves, cutting $\Sigma$ into connected components $\Sigma_1,\cdots,\Sigma_k$. Suppose that $\rho|_{\pi_1(\Sigma_i)}$ has an $\epsilon$-good developing map $f_i$ for every $i=1,2,\cdots,k$, then there exists a pseudo $\epsilon$-excellent developing map $f$ for $\rho$ with at most $2n$ folding curves, such that its area satisfies
    \[A(f)\leq\sum_{i=1}^kA(f_i)+2n\epsilon.\]
\end{lemma}

\begin{proof}
    For every $1\leq i\leq n$, cut $\Sigma$ along $\gamma_i$ and insert an annulus $C_i$ (see Figure \ref{fig: insert}). Define $f$ to be $f_j$ on each $\Sigma_j$ and extend it across each $C_i$ by mapping $C_i$ homeomorphically onto the collar in $\H/\rho(\gamma_i)$ bounded by $f_{j_1}(\gamma_i)$ and $f_{j_2}(\gamma_i)$, where $\Sigma_{j_1}$ and $\Sigma_{j_2}$ are the two components adjacent to $\gamma_i$. (If it happens that $f_{j_1}(\gamma_i)=f_{j_2}(\gamma_i)$, we do not insert $C_i$.)

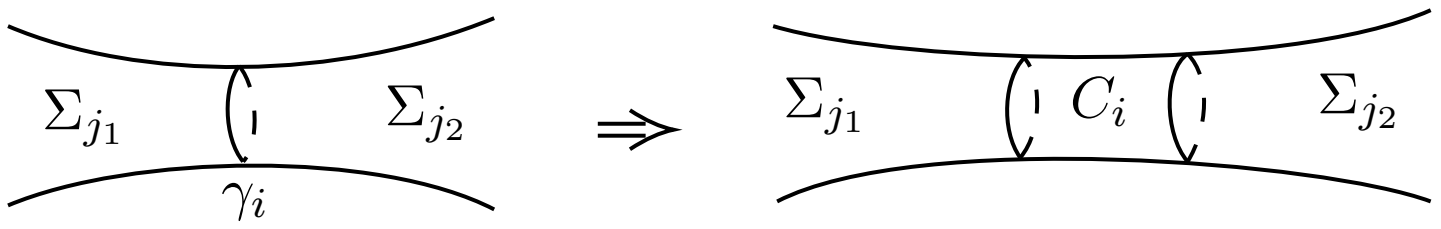
\begin{figure}
    \centering
    \tikzset{every picture/.style={line width=0.75pt}} %set default line width to 0.75pt        

\begin{tikzpicture}[x=0.75pt,y=0.75pt,yscale=-1,xscale=1]

%Curve Lines [id:da9320930467928172] 
\draw    (88,96) .. controls (121,110) and (171,110) .. (210,94) ;
%Curve Lines [id:da7230944093103555] 
\draw    (88,141) .. controls (119,128) and (181,127) .. (210,142) ;
%Curve Lines [id:da216173346203565] 
\draw    (147,130) .. controls (142,124) and (142,111) .. (146,106) ;
%Curve Lines [id:da7351698343240993] 
\draw  [dash pattern={on 4.5pt off 4.5pt}]  (146,106) .. controls (151,112) and (151,126) .. (147,130) ;
%Curve Lines [id:da30943549990363806] 
\draw    (280,96) .. controls (314,106) and (409,108) .. (445,92) ;
%Curve Lines [id:da30642527183616153] 
\draw    (281,140) .. controls (316,122) and (418,130) .. (445,140) ;
%Curve Lines [id:da9097150265702489] 
\draw    (342,129) .. controls (337,122) and (338,109) .. (343,104) ;
%Curve Lines [id:da06128241138183921] 
\draw  [dash pattern={on 4.5pt off 4.5pt}]  (342,129) .. controls (348,122) and (347,108) .. (343,104) ;
%Curve Lines [id:da9007164368946239] 
\draw    (384,130) .. controls (378,123) and (378,107) .. (384,103) ;
%Curve Lines [id:da7178139714500462] 
\draw  [dash pattern={on 4.5pt off 4.5pt}]  (384,103) .. controls (390,109) and (389,122) .. (384,130) ;

%Straight Lines [id:da9407875271750226] 
\draw    (236,120.5) -- (248,120.5)(236,123.5) -- (248,123.5) ;
\draw [shift={(255,122)}, rotate = 180] [color={rgb, 255:red, 0; green, 0; blue, 0 }  ][line width=0.75]    (10.93,-4.9) .. controls (6.95,-2.3) and (3.31,-0.67) .. (0,0) .. controls (3.31,0.67) and (6.95,2.3) .. (10.93,4.9)   ;

% Text Node
\draw (140,134) node [anchor=north west][inner sep=0.75pt]   [align=left] {$\gamma_{i}$};
% Text Node
\draw (95,110.4) node [anchor=north west][inner sep=0.75pt]    {$\Sigma _{j_{1}}$};
% Text Node
\draw (181,109.4) node [anchor=north west][inner sep=0.75pt]    {$\Sigma _{j_{2}}$};
% Text Node
%\draw (192,163) node [anchor=north west][inner sep=0.75pt]   [align=left] {Figure. Insert $\displaystyle \Pi _{i} $.};
% Text Node
\draw (415,106.4) node [anchor=north west][inner sep=0.75pt]    {$\Sigma _{j_{2}}$};
% Text Node
\draw (281,107.4) node [anchor=north west][inner sep=0.75pt]    {$\Sigma _{j_{1}}$};
% Text Node
\draw (353,107.4) node [anchor=north west][inner sep=0.75pt]    {$C _{i}$};
% Text Node
%\draw (332,132.4) node [anchor=north west][inner sep=0.75pt]    {$C'_{j}$};
% Text Node
%\draw (372,132.4) node [anchor=north west][inner sep=0.75pt]    {$C''_{j}$};

\end{tikzpicture}
    \caption{Insert $C_i$.}
    \label{fig: insert}
\end{figure}

    The developing map $f$ is branched on $\Sigma_1,\cdots,\Sigma_k,C_1,\cdots,C_n$, with at most $2n$ folding curves. By Lemma \ref{lem-c1}, $A(f|_{C_i})\leq2\epsilon$. It follows that $A(f)\leq\sum_{i=1}^kA(f_i)+2n\epsilon$.
\end{proof}

It remains to find such a decomposition of $\Sigma$, and to take care of the area condition. To this end we use analytic techniques. In particular we relate the area to the energy, and borrow some ideas from the theory of harmonic maps.

\subsection{Area v.s. energy}

If $\Sigma$ is given a conformal structure $g$ we can define the energy of a developing $f$,
\[E_g(f)=\int_\Sigma\frac{1}{2} |df|_g^2\omega_g,\]
where $\omega_g$ is the volume element of the metric $g$. The expression $|df|_g^2\omega_g$ is a conformal invariant. Define
\[E(f)=\inf\{E_g(f):g\text{ is a conformal structure on }\Sigma\},\]
\[e_0(\rho)=\inf\{E(f):f\text{ is a developing map for }\rho\}.\]

\begin{lemma}
    \[E(f)=A(f).\]
\end{lemma}

\begin{proof}
    Clearly $E_g(f)\geq A(f)$ for any conformal structure $g$ on $\Sigma_g$. On the other hand, fix an arbitrary metric $g$ on $\Sigma$. Define a family of metrics $g_\epsilon:=f^*g_\H+\epsilon g\ (0<\epsilon\leq1)$ on $\Sigma$. Then
    \begin{align}\label{eq-E=A}
        \lim_{\epsilon\rightarrow0}E_{g_\epsilon}(f)=A(f).
    \end{align}
    Indeed, we have the pointwise inequalities
    \[|f^*\omega_\H|\leq\frac{1}{2} |df|_{g_\epsilon}^2\omega_{g_\epsilon}\leq\omega_{g_\epsilon}\leq\omega_{g_1}.\]
    Since $\lim_{\epsilon\rightarrow0}\omega_{g_\epsilon}=|f^*\omega_{\H}|$, we get
    \[\lim_{\epsilon\rightarrow0}\frac{1}{2} |df|_{g_\epsilon}^2\omega_{g_\epsilon}=|f^*\omega|.\]
    Equation (\ref{eq-E=A}) follows from the bounded convergence theorem.
\end{proof}

\begin{corollary}
    \[e_0(\rho)=a_0(\rho).\]
\end{corollary}

Note that
\[e_0(\rho)=\inf_{g,f}E_g(f).\]
If we fix $g$ and very $f$, $E_g(f)$ is minimized by the unique $\rho$-equivariant harmonic map $f_g$. So we can equally define
\[e_0(\rho)=\inf_g E_g(f_g).\]
Since conformal structures on $\Sigma$ are parametrized by the Teichm\"uller space $\T=\T(\Sigma)$, we may also write
\[e_0(\rho)=\inf_{X\in\T}E_X(f_X),\]
where $X\in\T$ are marked Riemann surfaces.

\subsection{Finding a decomposition}\label{subsec-ede-find decomp}

Take a sequence $X_n\in\T$ such that $\lim_{n\rightarrow\infty}E_{X_n}(f_n)=e_0(\rho)$, where $f_n=f_{X_n}$ is the unique $\rho$-equivariant harmonic map on $X_n$. Taking a subsequence, the projection of $X_n$ to $\M=\M(\Sigma)$ (the moduli space of Riemann surfaces homeomorphic to $\Sigma$) converges to a nodal surface $X_\infty\in\overline{\M}$. (The convergence need not respect the markings.) Let $Y_{1,\infty},\cdots,Y_{k,\infty}$ (surfaces with punctures) be the components of $X_\infty$. There are decompositions
\[X_n=Y_{1,n}\cup\cdots\cup Y_{k,n},\]
where $Y_{i,n}$ are the subsurfaces of $X_n$ with short geodesic boundary that converge to $Y_{i,\infty}$.

\begin{remark}
    Via the markings $\Sigma\cong X_n$, the decompositions of $X_n$ above give rise to decompositions of $\Sigma$. Since the convergence need not respect the markings, for different $n$ these give rise to different decompositions of $\Sigma$.
\end{remark}

These are not yet the decomposition we are looking for. The next step is to extract a subsequence so that $f_n|_{Y_{i,n}}$ converge to some developing map $f_{i,\infty}$ on $Y_{i,\infty}$. If $f_{i,\infty}$ is a constant, we will further decompose $Y_{i,n}$ into pairs of pants.

Since the total energy is bounded, we have uniform interior gradient and higher order estimates for $f_n$ in $Y_{i,n}$. Choose base points $P_{i,n}\in Y_{i,n}$ that convergess to some $P_{i,\infty}\in Y_{i,\infty}$. Let $\widetilde{P}_{i,n}\in \widetilde{X}_n$ be arbitrary lifts of $P_{i,n}$. Take $A_{i,n}\in\mathrm{PSL}_2(\R)$ such that $A_{i,n}(f_n(\widetilde{P}_{i,n}))=0\in\D$. Then, up to a subsequence, $A_{i,n}\circ f_n|_{\widetilde{Y}_{i,n}}$ converges to some harmonic map $f_{i,\infty}:\widetilde{Y}_{i,\infty}\rightarrow\H$. By arranging the indices we assume that $f_{i,\infty}$ is constant for $i>l$ and non-constant for $i\leq l$. For each $i>l$, we decompose $Y_{i,\infty}$ into pairs of pants, and take decompositions of $Y_{i,n}$ into pairs of pants that converge to the decomposition of $Y_{i,\infty}$. In this way we get refined decompositions
\[X_n=Y_{1,n}\cup\cdots\cup Y_{l,n}\cup Z_{1,n}\cup\cdots\cup Z_{m,n}\]
of $X_n$ into subsurfaces with geodesic boundaries such that
\begin{itemize}
    \item For each $1\leq i\leq l$, $\{Y_{i,n}\}_{n\in\N}$ converges to a component $Y_{i,\infty}$ of $X_\infty$ and $\left\{A_{i,n}\circ f_n|_{\widetilde{Y}_{i,n}}\right\}_{n\in\N}$ converges to a harmonic map $f_{i,\infty}:\widetilde{Y}_{i,\infty}\rightarrow\H$ which is non-constant (where $A_{i,n}\in\mathrm{PSL}_2(\R)$).
    \item For each $1\leq j\leq m$, $Z_{j,n}$ are pairs of pants such that $\left\{B_{j,n}\circ f_n|_{\widetilde{Z}_{j,n}}\right\}_{n\in\N}$ converges to a constant (where $B_{j,n}\in\mathrm{PSL}_2(\R)$).
\end{itemize}
We claim that for $n$ large enough, each $\rho|_{\pi_1(Y_{i,n})}$ and $\rho|_{\pi_1(Z_{j,n})}$ admits an $\epsilon$-good developing map. This will be proved in the next two subsections.

\subsection{Good developing maps on $Y_{i,n}$}

\begin{lemma}\label{lem-ede2}
    For every $1\leq i\leq l$ and every $\epsilon>0$, if $n$ is large enough, $\rho|_{\pi_1(Y_{i,n})}$ admits an $\epsilon$-good developing map with area $\leq A(f_n|_{Y_{i,n}})+\epsilon$.
\end{lemma}

\begin{proof}
    Since we are working with a fixed $i$, to simplify the notation we write $Y_n,Y_\infty,A_n,f_\infty$ for $Y_{i,n},Y_{i,\infty},A_{i,n},f_{i,\infty}$ respectively. Note that $A_n\circ f_n$ are developing maps for $A_n\circ\rho\circ A_n^{-1}$, which are conjugates of $\rho$. We may replace $f_n$ by $A_n\circ f_n$ and $\rho$ by $A_n\circ\rho\circ A_n^{-1}$. Then $f_n|_{\widetilde{Y}_n}$ converges to $f_\infty$.
    \begin{claim}
        $f_\infty$ is either holomorphic or anti-holomorphic.
    \end{claim}
    Since $\lim_{n\rightarrow\infty}E_{X_n}(f_n)=e_0(\rho)$, the Hopf differentials of $f_n|_{Y_n}$ converges to zero. So $f_\infty$ has vanishing Hopf differential, thus either holomorphic or anti-holomorphic. By symmetry we may assume that it is holomorphic.
    \begin{claim}
        $f_{\infty}$ is $\sigma$-equivariant for some representation $\sigma:\pi_1(Y_\infty)\rightarrow\mathrm{PSL}_2(\R)$ which sends the loop around each puncture to the identity or a parabolic element.
    \end{claim}
    By assumption $\rho(\gamma)$ is hyperbolic or parabolic for every boundary curve $\gamma$ of $Y_n$, taking limit we see that $f_{\infty}$ is $\sigma$-equivariant for some representation $\sigma$ such that $\sigma(\gamma)$ is non-elliptic for every loop $\gamma$ around a puncture. The crux of the claim is that $\sigma(\gamma)$ is non-hyperbolic. Suppose that $\sigma(\gamma)$ is hyperbolic. Then $\H/\rho(\gamma)$ is an annulus. This gives rise to a contradiction since every holomorphic map from a punctured disc to an annulus must be trivial on the fundamental group.

    Now we are ready to construct a developing map $\hat{f}_n$ on $Y_n$ with the desired properties. It will be obtained by perturbing $f_n$. Take disjoint small discs (or punctured discs) $\{D_{i,\infty}\}_{1\leq i\leq k}$ in $Y_\infty$ centered at the zeros of $df_\infty$ or the punctures of $Y_\infty$. Take discs or annuli $D_{i,n}\subset Y_n$ that converge to $D_{i,\infty}$. Since $f_n|_{Y_n\backslash\cup_{i=1}^kD_{i,n}}$ converges uniformly (up to any order) to $f_\infty|_{Y_\infty\backslash\cup_{i=1}^kD_{i,\infty}}$, and the latter one is orientation preserving, we conclude that $f_n|_{Y_n\backslash\cup_{i=1}^kD_{i,n}}$ is orientation preserving for $n$ large enough. We define $\hat{f}_n$ on $Y_n\backslash\cup_{i=1}^kD_{i,n}$ to be equal to $f_n|_{Y_n\backslash\cup_{i=1}^kD_{i,n}}$. It remains to extend $\hat{f}_n$ to each of these small discs or annuli. We analyze $f_n$ locally in each $D_{i,n}$. To simplify the notation we drop the index $i$.

    \paragraph{Case 1} $D_\infty$ is a disc centered at a zero of $df_\infty$.

    There is a local coordinate centered at the zero so that $f_\infty(z)=z^k\in\D$ for some $k\geq2$. By shrinking $D_\infty$ we may assume that $D_\infty=\{|z|<r\}$ for some $r>0$. We may assume $D_n=\{|z|<r\}$ as well, where we take local coordinates on $Y_n$ that converge to the above coordinate on $Y_\infty$. Let $\gamma_n=\partial D_n=\{|z|=r\}$. Then $f_n|_{\gamma_n}$ converges to the function $z\mapsto z^k$ up to any order. By Lemma \ref{lem-ede3} below, we can extend $f_n|_{\gamma_n}$ to a map $\hat{f}_n$ defined on $D_n$ that is orientation preserving except for one branch point.

    \paragraph{Case 2} $D_\infty$ is a disc centered at a puncture of $Y_\infty$, around which the monodromy of $f_\infty$ is parabolic.

    By an automorphism of $\H$, we may assume that the monodromy of $f_\infty$ around this puncture is $w\mapsto w+1$. So there is a local coordinate so that $f_\infty(z)=\frac{1}{2\pi i}\log z$. By shrinking $D_\infty$ we may assume that $D_\infty=\{|z|<r\}$ for some $r>0$. $D_n$ are annuli, whose outer boundary curves $\gamma_n$ converges to $\{|z|=r\}$. Since $f_n|_{\gamma_n}$ converges to $z\mapsto\frac{1}{2\pi i}\log z$, for large enough $n$, $f_n$ maps $\gamma_n$ to a simple closed curve $f(\gamma_n)$ in $\H/\rho(\gamma_n)$ which is disjoint from the axis of $\rho(\gamma_n)$ (since the axis of $\rho(\gamma_n)$ goes to infinity). We extend $\hat{f}_n$ to the annulus $D_n$ to be a homeomorphism from $D_n$ to the annulus in $\H/\rho(\gamma_n)$ bounded by $f(\gamma_n)$ and the axis of $\rho(\gamma_n)$.

    \paragraph{Case 3} $D_\infty$ is a disc centered at a puncture of $Y_\infty$, around which the monodromy of $f_\infty$ is the identity.

    In this case the puncture is a removable singularity of $f_\infty$. There is a local coordinate centered at the puncture so that $f_{\infty}(z)=z^k\in\D$ for some $k\in\N$. By shrinking $D_\infty$ we may assume that $D_\infty=\{|z|<r\}$ for some $r>0$. We may also assume that the outer boundaries of $D_n$ are $\gamma_n=\{|z|=r\}$, where we take local coordinates on $Y_n$ that converge to the above coordinate on $Y_\infty$. By assumption $\rho(\gamma_n)$ is a hyperbolic or parabolic element. By a further rotation of the coordinates, we may assume that the axis of $\rho(\gamma_n)$ is orthogonal to the real axis, and $\mathrm{Im}(\rho(\gamma_n)(0))<0$. (We use the model $\D$ here.)

    Let $\gamma_n$ be parametrized by $re^{i\theta}$, so that its universal covering $\widetilde{\gamma_n}$ is parametrized by $\theta\in\R$. By assumption $f_n|_{\widetilde{\gamma_n}}$ converges to the function $\theta\mapsto r^ke^{ik\theta}$, and $f_n(2\pi)=\rho(\gamma_n)\cdot f_n(0)$. We have $f_n(0)\approx f_n(2\pi)\approx r^k$ and $\arg(f_n(0)-f_n(2\pi))\approx\frac{\pi}{2}$ (recall the assumption on $\rho(\gamma_n)$ that we made at the end of the last paragraph). Let $\alpha_n$ be the straight line segment connecting $f_n(0)$ and $f_n(2\pi)$.

    Let $\beta_n=\{|z-\frac{1}{3}r|=\frac{2}{3}r\}$. The two circles $\gamma_n,\beta_n$ bound a simply connected, crescent-shaped region $U$ (see Figure \ref{fig: crescent}). By lemma \ref{lem-ede3} below, we can extend $\hat{f}_n$ to $U$ so that it maps $\beta_n$ homeomorphically to $\alpha_n$, and that it is orientation preserving on $U$ except for one branch point (Figure \ref{fig: crescent} indicates the situation for $k=2$). Finally we extend $\hat{f}_n$ to $D_n\backslash U$ as follows. $\alpha_n$ projects onto a simple closed curve in $\D/\rho(\gamma_n)$. We define $\hat{f}_n$ on $D_n\backslash U$ to be a homeomorphism onto an annulus in $\D/\rho(\gamma_n)$ bounded by $\alpha_n$ and an equidistant curve. This finishes the construction of $\hat{f}_n$ on $D_n$ in the last case.
    
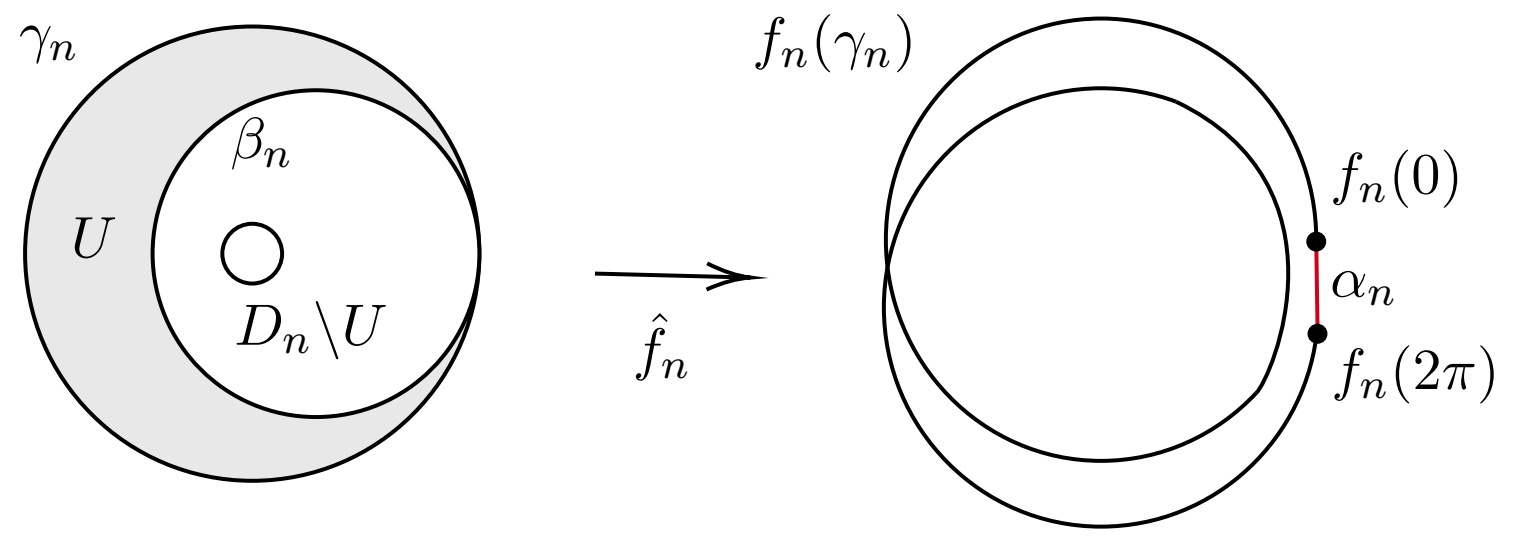
\begin{figure}
    \centering
    \tikzset{every picture/.style={line width=0.75pt}} %set default line width to 0.75pt        

\begin{tikzpicture}[x=0.75pt,y=0.75pt,yscale=-1,xscale=1]

%Shape: Arc [id:dp8036194979098606] 
\draw  [draw opacity=0] (394.95,155.95) .. controls (385.12,166.46) and (371.3,173) .. (356,173) .. controls (326.18,173) and (302,148.15) .. (302,117.5) .. controls (302,86.85) and (326.18,62) .. (356,62) .. controls (385.82,62) and (410,86.85) .. (410,117.5) .. controls (410,117.66) and (410,117.81) .. (410,117.97) -- (356,117.5) -- cycle ; \draw   (394.95,155.95) .. controls (385.12,166.46) and (371.3,173) .. (356,173) .. controls (326.18,173) and (302,148.15) .. (302,117.5) .. controls (302,86.85) and (326.18,62) .. (356,62) .. controls (385.82,62) and (410,86.85) .. (410,117.5) .. controls (410,117.66) and (410,117.81) .. (410,117.97) ;  
%Shape: Arc [id:dp5292275181105693] 
\draw  [draw opacity=0] (410.31,139.11) .. controls (407.99,167.33) and (384.56,189.5) .. (356,189.5) .. controls (325.9,189.5) and (301.5,164.88) .. (301.5,134.5) .. controls (301.5,104.12) and (325.9,79.5) .. (356,79.5) .. controls (362.61,79.5) and (368.94,80.69) .. (374.8,82.86) -- (356,134.5) -- cycle ; \draw   (410.31,139.11) .. controls (407.99,167.33) and (384.56,189.5) .. (356,189.5) .. controls (325.9,189.5) and (301.5,164.88) .. (301.5,134.5) .. controls (301.5,104.12) and (325.9,79.5) .. (356,79.5) .. controls (362.61,79.5) and (368.94,80.69) .. (374.8,82.86) ;  
%Curve Lines [id:da9599637434421745] 
\draw    (374.8,82.86) .. controls (419,104) and (399,152) .. (394.95,155.95) ;
%Straight Lines [id:da02271290156978578] 
\draw [color={rgb, 255:red, 208; green, 2; blue, 27 }  ,draw opacity=1 ][fill={rgb, 255:red, 255; green, 255; blue, 255 }  ,fill opacity=1 ]   (410,117.97) -- (410.31,139.11) ;
%Shape: Circle [id:dp6805221317081501] 
\draw  [fill={rgb, 255:red, 0; green, 0; blue, 0 }  ,fill opacity=1 ] (408,117.97) .. controls (408,116.87) and (408.89,115.97) .. (410,115.97) .. controls (411.1,115.97) and (412,116.87) .. (412,117.97) .. controls (412,119.08) and (411.1,119.97) .. (410,119.97) .. controls (408.89,119.97) and (408,119.08) .. (408,117.97) -- cycle ;
%Shape: Circle [id:dp8471259359191019] 
\draw  [fill={rgb, 255:red, 0; green, 0; blue, 0 }  ,fill opacity=1 ] (412.31,141.11) .. controls (412.31,140) and (411.42,139.11) .. (410.31,139.11) .. controls (409.21,139.11) and (408.31,140) .. (408.31,141.11) .. controls (408.31,142.21) and (409.21,143.11) .. (410.31,143.11) .. controls (411.42,143.11) and (412.31,142.21) .. (412.31,141.11) -- cycle ;
%Shape: Circle [id:dp35223764742704944] 
\draw  [fill={rgb, 255:red, 155; green, 155; blue, 155 }  ,fill opacity=0.23 ] (86,121) .. controls (86,89.52) and (111.52,64) .. (143,64) .. controls (174.48,64) and (200,89.52) .. (200,121) .. controls (200,152.48) and (174.48,178) .. (143,178) .. controls (111.52,178) and (86,152.48) .. (86,121) -- cycle ;
%Shape: Circle [id:dp04709887225534204] 
\draw  [fill={rgb, 255:red, 255; green, 255; blue, 255 }  ,fill opacity=1 ] (118,121) .. controls (118,98.36) and (136.36,80) .. (159,80) .. controls (181.64,80) and (200,98.36) .. (200,121) .. controls (200,143.64) and (181.64,162) .. (159,162) .. controls (136.36,162) and (118,143.64) .. (118,121) -- cycle ;
%Shape: Circle [id:dp23928968253460248] 
\draw   (135.5,121) .. controls (135.5,116.86) and (138.86,113.5) .. (143,113.5) .. controls (147.14,113.5) and (150.5,116.86) .. (150.5,121) .. controls (150.5,125.14) and (147.14,128.5) .. (143,128.5) .. controls (138.86,128.5) and (135.5,125.14) .. (135.5,121) -- cycle ;

%Straight Lines [id:da376547296392304] 
\draw    (229,126) -- (266,126.95) ;
\draw [shift={(268,127)}, rotate = 181.47] [color={rgb, 255:red, 0; green, 0; blue, 0 }  ][line width=0.75]    (10.93,-3.29) .. controls (6.95,-1.4) and (3.31,-0.3) .. (0,0) .. controls (3.31,0.3) and (6.95,1.4) .. (10.93,3.29)   ;

% Text Node
\draw (412,93.4) node [anchor=north west][inner sep=0.75pt]    {$f_{n}( 0)$};
% Text Node
\draw (412.31,142.51) node [anchor=north west][inner sep=0.75pt]    {$f_{n}( 2\pi )$};
% Text Node
\draw (412,124) node [anchor=north west][inner sep=0.75pt]    {$\alpha_n $};
% Text Node
\draw (96,110.4) node [anchor=north west][inner sep=0.75pt]    {$U$};
% Text Node
\draw (137,131.4) node [anchor=north west][inner sep=0.75pt]    {$D_{n} \backslash U$};
% Text Node
\draw (237,134.4) node [anchor=north west][inner sep=0.75pt]    {$\hat{f}_{n}$};
% Text Node
%\draw (222,199) node [anchor=north west][inner sep=0.75pt]   [align=left] {Figure.};
% Text Node
\draw (267,59.4) node [anchor=north west][inner sep=0.75pt]    {$f_{n}( \gamma_{n})$};
% Text Node
\draw (83,62.4) node [anchor=north west][inner sep=0.75pt]    {$\gamma_{n}$};
% Text Node
\draw (136,85.4) node [anchor=north west][inner sep=0.75pt]    {$\beta_{n}$};

\end{tikzpicture}
    \caption{Extend $\hat{f}_n$ to $U$.}
    \label{fig: crescent}
\end{figure}

    \begin{remark}
        The tip of the crescent-shaped region $U$ is also a branch point of $\hat{f}_n$.
    \end{remark}
    Finally it is easy to see that by taking the discs $D_{i,n}$ small enough we can make sure that $A(\hat{f}_n)\leq A(f_n|_{Y_{i,n}})+\epsilon$. This finishes the proof of Lemma \ref{lem-ede2}.
\end{proof}

\begin{lemma}\label{lem-ede3}
    Suppose that $\gamma:S^1\rightarrow\C\backslash\{0\}$ is piecewise $C^1$ with winding number $k$ about the origin and such that $\arg\gamma(\theta)$ is monotone increasing. Then $\gamma$ can be extended continuously to a piecewise $C^1$ map from $\overline{\D}$ to $\C$ that is orientation preserving except for a branch point of multiplicity $k-1$.
\end{lemma}

\begin{proof}
    Define $f(re^{i\theta})=r\gamma(\theta)$.
\end{proof}

\subsection{Good developing maps on $Z_{j,n}$}\label{subsec-ede-Z}

We first prove a general lemma constructing branched developing maps on a pair of pants.

\begin{lemma}\label{lem-ede-p1}
    Let $Z$ be a three-holed sphere. Let $\rho:\pi_1(Z)\rightarrow\mathrm{PSL}_2(\R)$ be a representation that sends every boundary curve to a hyperbolic or parabolic element. Then there exists a branched developing map $f$ for $\rho$ which maps every boundary curve homeomorpically onto an equidistant curve.
\end{lemma}

\begin{proof}
    Let $\gamma_1,\gamma_2,\gamma_3$ be the three boundary curves of $Z$. Pick a base point $P\in Z$. Take simple closed curves $\alpha_1,\alpha_2,\alpha_3$ based at $P$ homotopic to $\gamma_1,\gamma_2,\gamma_3$ respectively, such that $\alpha_i,\alpha_j\ (i\neq j)$ meet only at $P$. They cut $Z$ into four regions $C_1,C_2,C_3,T$, where $C_i$ are annuli with boundary curves $\alpha_i,\gamma_i$, and $T$ is a triangle with sides $\alpha_1,\alpha_2,\alpha_3$ (see Figure \ref{fig: pants1}). We will first construct the developing map $f$ on $T$, then extend it to the annuli $C_i$.

    \begin{figure}
        \centering
        \tikzset{every picture/.style={line width=0.75pt}} %set default line width to 0.75pt        

\begin{tikzpicture}[x=0.75pt,y=0.75pt,yscale=-1,xscale=1]

%Shape: Ellipse [id:dp1225520737379151] 
\draw   (239,52.5) .. controls (239,47.25) and (251.54,43) .. (267,43) .. controls (282.46,43) and (295,47.25) .. (295,52.5) .. controls (295,57.75) and (282.46,62) .. (267,62) .. controls (251.54,62) and (239,57.75) .. (239,52.5) -- cycle ;
%Shape: Ellipse [id:dp7106362167543303] 
\draw   (195.96,121.56) .. controls (198.19,119.04) and (209.52,125.44) .. (221.26,135.85) .. controls (232.99,146.26) and (240.7,156.74) .. (238.47,159.26) .. controls (236.24,161.78) and (224.91,155.38) .. (213.17,144.97) .. controls (201.43,134.55) and (193.73,124.07) .. (195.96,121.56) -- cycle ;
%Shape: Ellipse [id:dp9037731074299207] 
\draw   (305.27,155.72) .. controls (300.77,152.6) and (304.23,139.84) .. (312.99,127.21) .. controls (321.76,114.59) and (332.5,106.88) .. (337,110) .. controls (341.5,113.12) and (338.04,125.88) .. (329.27,138.51) .. controls (320.51,151.14) and (309.76,158.84) .. (305.27,155.72) -- cycle ;
%Curve Lines [id:da3558850514618388] 
\draw    (195.96,121.56) .. controls (233,123) and (236,67.5) .. (239,52.5) ;
%Curve Lines [id:da28254742111596287] 
\draw    (295,52.5) .. controls (293,72.5) and (313.72,103) .. (337,110) ;
%Shape: Arc [id:dp7330764656411166] 
\draw  [draw opacity=0][dash pattern={on 4.5pt off 4.5pt}] (237.79,69.63) .. controls (246.77,68.29) and (257.97,67.6) .. (270.11,67.82) .. controls (279.99,68) and (289.23,68.77) .. (297.13,69.94) -- (269.9,79.48) -- cycle ; \draw  [dash pattern={on 4.5pt off 4.5pt}] (237.79,69.63) .. controls (246.77,68.29) and (257.97,67.6) .. (270.11,67.82) .. controls (279.99,68) and (289.23,68.77) .. (297.13,69.94) ;  
%Shape: Arc [id:dp9286131662414808] 
\draw  [draw opacity=0][dash pattern={on 4.5pt off 4.5pt}] (252.26,152.53) .. controls (245.67,150.92) and (237.09,144.93) .. (229.6,136.1) .. controls (222.1,127.27) and (217.59,117.82) .. (217.08,111.06) -- (239.9,127.36) -- cycle ; \draw  [dash pattern={on 4.5pt off 4.5pt}] (252.26,152.53) .. controls (245.67,150.92) and (237.09,144.93) .. (229.6,136.1) .. controls (222.1,127.27) and (217.59,117.82) .. (217.08,111.06) ;  
%Shape: Arc [id:dp10516773632397292] 
\draw  [draw opacity=0][dash pattern={on 4.5pt off 4.5pt}] (319.65,101.04) .. controls (318.32,107.06) and (315.64,114.39) .. (311.81,121.93) .. controls (304.66,136) and (295.75,146.38) .. (290.03,147.91) -- (302.56,117.22) -- cycle ; \draw  [dash pattern={on 4.5pt off 4.5pt}] (319.65,101.04) .. controls (318.32,107.06) and (315.64,114.39) .. (311.81,121.93) .. controls (304.66,136) and (295.75,146.38) .. (290.03,147.91) ;  
%Shape: Arc [id:dp998189219660004] 
\draw  [draw opacity=0] (272.87,100.56) .. controls (266.47,98.94) and (258.76,94.68) .. (251.58,88.35) .. controls (244.76,82.34) and (239.74,75.63) .. (237.18,69.71) -- (261.17,77.46) -- cycle ; \draw   (272.87,100.56) .. controls (266.47,98.94) and (258.76,94.68) .. (251.58,88.35) .. controls (244.76,82.34) and (239.74,75.63) .. (237.18,69.71) ;  
%Shape: Arc [id:dp9823780326646456] 
\draw  [draw opacity=0] (297.13,69.94) .. controls (295.62,75.81) and (291.92,82.73) .. (286.43,89.2) .. controls (282.54,93.8) and (278.28,97.55) .. (274.15,100.21) -- (275.75,80.16) -- cycle ; \draw   (297.13,69.94) .. controls (295.62,75.81) and (291.92,82.73) .. (286.43,89.2) .. controls (282.54,93.8) and (278.28,97.55) .. (274.15,100.21) ;  
%Shape: Arc [id:dp9158642154793614] 
\draw  [draw opacity=0] (219.04,112.86) .. controls (222.57,108.95) and (232.69,104.36) .. (245.1,101.42) .. controls (258.08,98.34) and (269.61,97.99) .. (274.15,100.21) -- (246.95,109.2) -- cycle ; \draw   (219.04,112.86) .. controls (222.57,108.95) and (232.69,104.36) .. (245.1,101.42) .. controls (258.08,98.34) and (269.61,97.99) .. (274.15,100.21) ;  
%Shape: Arc [id:dp8358682287847108] 
\draw  [draw opacity=0] (272.87,100.56) .. controls (272.42,106.25) and (269.3,116.93) .. (264.43,128.71) .. controls (259.57,140.46) and (254.25,150.2) .. (250.55,154.56) -- (259.42,126.64) -- cycle ; \draw   (272.87,100.56) .. controls (272.42,106.25) and (269.3,116.93) .. (264.43,128.71) .. controls (259.57,140.46) and (254.25,150.2) .. (250.55,154.56) ;  
%Shape: Arc [id:dp28021817118762626] 
\draw  [draw opacity=0] (272.87,100.56) .. controls (278.98,99.59) and (286.19,99) .. (293.91,98.94) .. controls (303.71,98.87) and (312.69,99.66) .. (319.65,101.04) -- (293.98,108.74) -- cycle ; \draw   (272.87,100.56) .. controls (278.98,99.59) and (286.19,99) .. (293.91,98.94) .. controls (303.71,98.87) and (312.69,99.66) .. (319.65,101.04) ;  
%Shape: Arc [id:dp32812321081923246] 
\draw  [draw opacity=0] (290.27,150.09) .. controls (287.3,147.39) and (282.34,135.85) .. (278.23,121.59) .. controls (275.91,113.56) and (274.28,106.12) .. (273.49,100.46) -- (282.88,120.25) -- cycle ; \draw   (290.27,150.09) .. controls (287.3,147.39) and (282.34,135.85) .. (278.23,121.59) .. controls (275.91,113.56) and (274.28,106.12) .. (273.49,100.46) ;  
%Shape: Arc [id:dp43489808151534426] 
\draw  [draw opacity=0] (238.47,159.26) .. controls (245.61,154.36) and (256.29,150.68) .. (268.45,149.48) .. controls (282.99,148.05) and (296.14,150.51) .. (304.04,155.42) -- (270.47,169.93) -- cycle ; \draw   (238.47,159.26) .. controls (245.61,154.36) and (256.29,150.68) .. (268.45,149.48) .. controls (282.99,148.05) and (296.14,150.51) .. (304.04,155.42) ;

% Text Node
%\draw (129,176) node [anchor=north west][inner sep=0.75pt]   [align=left] {Figure. $\displaystyle Z$ is decomposed into four regions.};
% Text Node
\draw (257,100) node [anchor=north west][inner sep=0.75pt]   [align=left] {$\displaystyle P$};
% Text Node
\draw (261,23) node [anchor=north west][inner sep=0.75pt]   [align=left] {$\displaystyle \gamma _{1}$};
% Text Node
\draw (329,140) node [anchor=north west][inner sep=0.75pt]   [align=left] {$\displaystyle \gamma _{2}$};
% Text Node
\draw (200,150) node [anchor=north west][inner sep=0.75pt]   [align=left] {$\displaystyle \gamma _{3}$};
% Text Node
\draw (267,127) node [anchor=north west][inner sep=0.75pt]   [align=left] {$\displaystyle T$};
% Text Node
\draw (260,69) node [anchor=north west][inner sep=0.75pt]   [align=left] {$\displaystyle C _{1}$};
% Text Node
\draw (287,109) node [anchor=north west][inner sep=0.75pt]   [align=left] {$\displaystyle C _{2}$};
% Text Node
\draw (230,115) node [anchor=north west][inner sep=0.75pt]   [align=left] {$\displaystyle C _{3}$};

\end{tikzpicture}
        \caption{$Z$ is decomposed into four regions.}
        \label{fig: pants1}
    \end{figure}
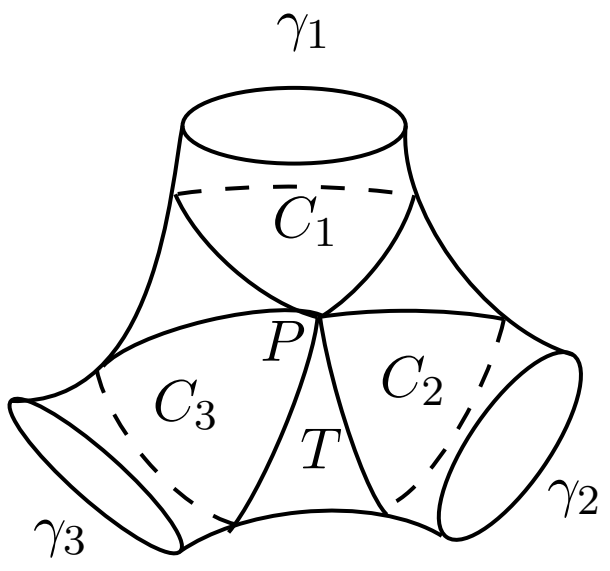

    Let $\pi:\widetilde{Z}\rightarrow Z$ be the universal covering. Then $T$ is a fundamental domain of $\pi^{-1}(T)$. To define $f:\pi^{-1}(T)\rightarrow\H$ that is $\rho$-equivariant, it suffices to define $f:T\rightarrow\H$ so that $f(P_1)=\rho(\alpha_2)\cdot f(P_3),\ f(P_2)=\rho(\alpha_1)^{-1}\cdot f(P_3)$, where $P_i$ is the vertex of $T$ opposite to side $\alpha_i\ (i=1,2,3)$. Take a point $Q_3\in\H$. Let $Q_1=\rho(\alpha_2)\cdot Q_3,\ Q_2=\rho(\alpha_1)^{-1}\cdot Q_3$. By perturbing $Q_3$, we may assume that $Q_1,Q_2,Q_3$ are not collinear (see Claim \ref{claim-ede} below). Let $U$ be the geodesic triangle in $\H$ with vertices $Q_1,Q_2,Q_3$. Let $\beta_i$ be the side of $U$ opposite to $Q_i$. We define $f:T\rightarrow\H$ to be a homeomorphism from $T$ to $U$, taking $P_i$ to $Q_i$, $\alpha_i$ to $\beta_i$ (see Figure \ref{fig: pants2}).

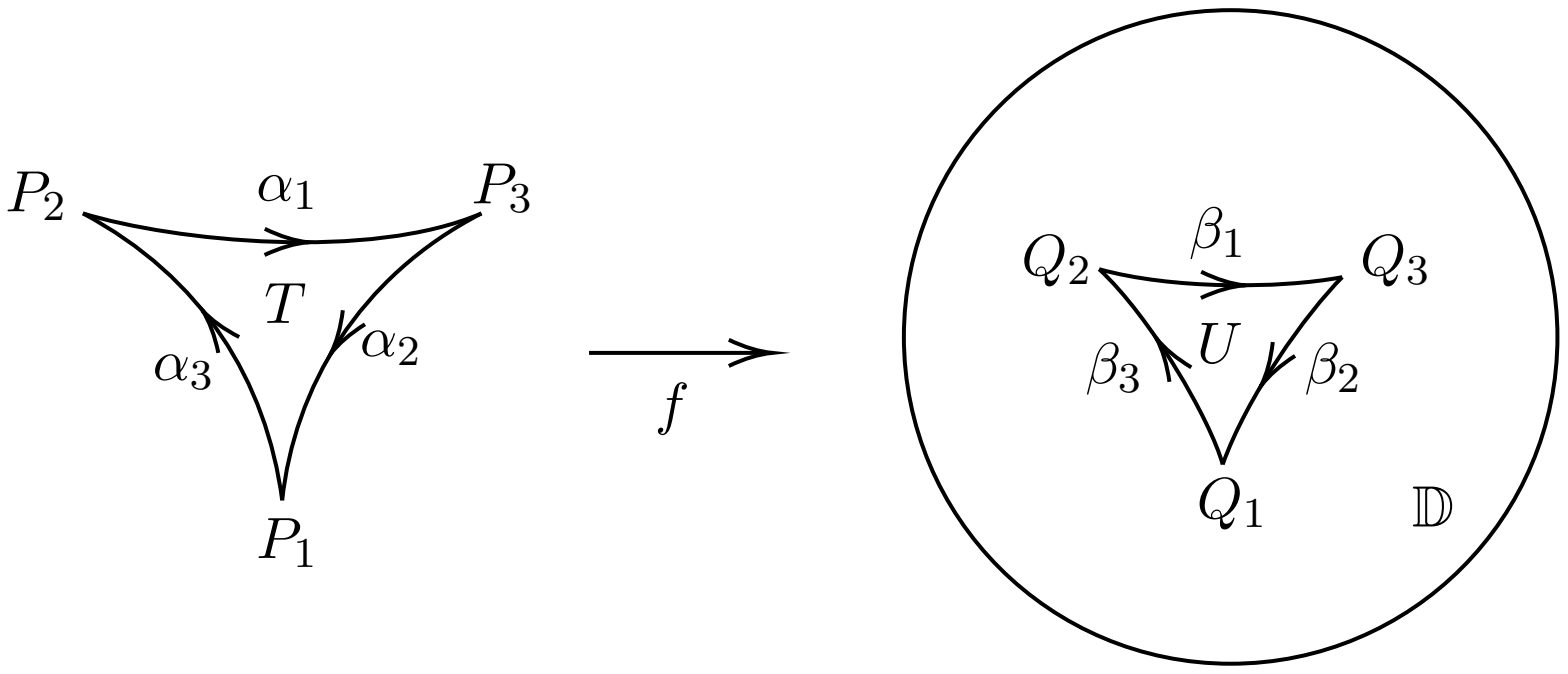
\begin{figure}
    \centering
    \tikzset{every picture/.style={line width=0.75pt}} %set default line width to 0.75pt        

\begin{tikzpicture}[x=0.75pt,y=0.75pt,yscale=-1,xscale=1]

%Curve Lines [id:da5310627799834882] 
\draw    (100,116) .. controls (123,123) and (172,128) .. (200,116) ;
\draw [shift={(156.53,123.24)}, rotate = 180.48] [color={rgb, 255:red, 0; green, 0; blue, 0 }  ][line width=0.75]    (10.93,-3.29) .. controls (6.95,-1.4) and (3.31,-0.3) .. (0,0) .. controls (3.31,0.3) and (6.95,1.4) .. (10.93,3.29)   ;
%Curve Lines [id:da8704606090686844] 
\draw    (100,116) .. controls (130,132) and (146,158) .. (150,188) ;
\draw [shift={(129.98,140.31)}, rotate = 52.75] [color={rgb, 255:red, 0; green, 0; blue, 0 }  ][line width=0.75]    (10.93,-3.29) .. controls (6.95,-1.4) and (3.31,-0.3) .. (0,0) .. controls (3.31,0.3) and (6.95,1.4) .. (10.93,3.29)   ;
%Curve Lines [id:da9491414783202362] 
\draw    (150,188) .. controls (153,160) and (169,132) .. (200,116) ;
\draw [shift={(162.13,151.29)}, rotate = 302.43] [color={rgb, 255:red, 0; green, 0; blue, 0 }  ][line width=0.75]    (10.93,-3.29) .. controls (6.95,-1.4) and (3.31,-0.3) .. (0,0) .. controls (3.31,0.3) and (6.95,1.4) .. (10.93,3.29)   ;
%Curve Lines [id:da7113146384777106] 
\draw    (355,130) .. controls (373,135) and (400,135) .. (416,132) ;
\draw [shift={(391.5,134.04)}, rotate = 180.59] [color={rgb, 255:red, 0; green, 0; blue, 0 }  ][line width=0.75]    (10.93,-3.29) .. controls (6.95,-1.4) and (3.31,-0.3) .. (0,0) .. controls (3.31,0.3) and (6.95,1.4) .. (10.93,3.29)   ;
%Curve Lines [id:da27311041753824994] 
\draw    (355,130) .. controls (368,142) and (383,168) .. (386,179) ;
\draw [shift={(369.35,147.33)}, rotate = 56.5] [color={rgb, 255:red, 0; green, 0; blue, 0 }  ][line width=0.75]    (10.93,-3.29) .. controls (6.95,-1.4) and (3.31,-0.3) .. (0,0) .. controls (3.31,0.3) and (6.95,1.4) .. (10.93,3.29)   ;
%Curve Lines [id:da21634736018663014] 
\draw    (416,132) .. controls (407,141) and (392,162) .. (386,179) ;
\draw [shift={(395.6,159.3)}, rotate = 301.83] [color={rgb, 255:red, 0; green, 0; blue, 0 }  ][line width=0.75]    (10.93,-3.29) .. controls (6.95,-1.4) and (3.31,-0.3) .. (0,0) .. controls (3.31,0.3) and (6.95,1.4) .. (10.93,3.29)   ;

%Shape: Circle [id:dp9275619679935689] 
\draw   (306,147) .. controls (306,101.71) and (342.71,65) .. (388,65) .. controls (433.29,65) and (470,101.71) .. (470,147) .. controls (470,192.29) and (433.29,229) .. (388,229) .. controls (342.71,229) and (306,192.29) .. (306,147) -- cycle ;

%Straight Lines [id:da631988869919302] 
\draw    (227,151) -- (271,151) ;
\draw [shift={(273,151)}, rotate = 180] [color={rgb, 255:red, 0; green, 0; blue, 0 }  ][line width=0.75]    (10.93,-3.29) .. controls (6.95,-1.4) and (3.31,-0.3) .. (0,0) .. controls (3.31,0.3) and (6.95,1.4) .. (10.93,3.29)   ;

% Text Node
\draw (79,104) node [anchor=north west][inner sep=0.75pt]   [align=left] {$\displaystyle P_{2}$};
% Text Node
\draw (196,102) node [anchor=north west][inner sep=0.75pt]   [align=left] {$\displaystyle P_{3}$};
% Text Node
\draw (142,191) node [anchor=north west][inner sep=0.75pt]   [align=left] {$\displaystyle P_{1}$};
% Text Node
\draw (144,132) node [anchor=north west][inner sep=0.75pt]   [align=left] {$\displaystyle T$};
% Text Node
\draw (142,105) node [anchor=north west][inner sep=0.75pt]   [align=left] {$\displaystyle \alpha _{1}$};
% Text Node
\draw (168,144) node [anchor=north west][inner sep=0.75pt]   [align=left] {$\displaystyle \alpha _{2}$};
% Text Node
\draw (116,150) node [anchor=north west][inner sep=0.75pt]   [align=left] {$\displaystyle \alpha _{3}$};
% Text Node
\draw (378,142) node [anchor=north west][inner sep=0.75pt]   [align=left] {$\displaystyle U$};
% Text Node
\draw (378,181) node [anchor=north west][inner sep=0.75pt]   [align=left] {$\displaystyle Q_{1}$};
% Text Node
\draw (334,120) node [anchor=north west][inner sep=0.75pt]   [align=left] {$\displaystyle Q_{2}$};
% Text Node
\draw (419,120) node [anchor=north west][inner sep=0.75pt]   [align=left] {$\displaystyle Q_{3}$};
% Text Node
\draw (376,113) node [anchor=north west][inner sep=0.75pt]   [align=left] {$\displaystyle \beta _{1}$};
% Text Node
\draw (405,147) node [anchor=north west][inner sep=0.75pt]   [align=left] {$\displaystyle \beta _{2}$};
% Text Node
\draw (350,147) node [anchor=north west][inner sep=0.75pt]   [align=left] {$\displaystyle \beta _{3}$};
% Text Node
\draw (242,157) node [anchor=north west][inner sep=0.75pt]   [align=left] {$\displaystyle f$};
% Text Node
%\draw (112,240) node [anchor=north west][inner sep=0.75pt]   [align=left] {Figure. $\displaystyle f$ maps $\displaystyle T$ homeomorphically to $\displaystyle U$.};
% Text Node
\draw (432,183) node [anchor=north west][inner sep=0.75pt]   [align=left] {$\displaystyle \mathbb{D}$};

\end{tikzpicture}
    \caption{$f$ maps $T$ homeomorphically to $U$.}
    \label{fig: pants2}
\end{figure}

    Now we extend $f$ to the annuli $C_i$. The projection of $\beta_i$ in $\H/\rho(\alpha_i)$ is a simple geodesic loop (Lemma \ref{lem-h7}). Take an equidistant curve $\eta_i\subset\H/\rho(\alpha_i)$ disjoint from $\beta_i$. Define $f:C_i\rightarrow\H/\rho(\alpha_i)$ to be a homeomorphism from $C_i$ to the annulus bounded by $\eta_i$ and $\beta_i$. We choose $\eta_i$ on the correct side of $\beta_i$ so that $f|_T$ and $f|_{C_i}$ have the same orientation. Therefore $f$ is a local homeomorphism away from $P$.

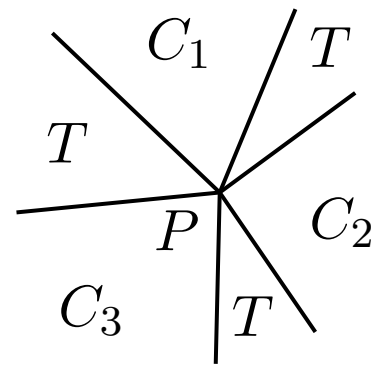
\begin{figure}
    \centering
    \tikzset{every picture/.style={line width=0.75pt}} %set default line width to 0.75pt        

\begin{tikzpicture}[x=0.75pt,y=0.75pt,yscale=-1,xscale=1]

%Straight Lines [id:da43131091775334685] 
\draw    (121,94) -- (163,134) ;
%Straight Lines [id:da14063662800022603] 
\draw    (182,88) -- (163,134) ;
%Straight Lines [id:da684868683638871] 
\draw    (163,134) -- (112,139) ;
%Straight Lines [id:da34998205650930714] 
\draw    (163,134) -- (162,177) ;
%Straight Lines [id:da362088545791335] 
\draw    (163,134) -- (187,169) ;
%Straight Lines [id:da1009386897039869] 
\draw    (163,134) -- (197,109) ;

% Text Node
\draw (145,137) node [anchor=north west][inner sep=0.75pt]   [align=left] {$\displaystyle P$};
% Text Node
\draw (164.5,158.5) node [anchor=north west][inner sep=0.75pt]   [align=left] {$\displaystyle T$};
% Text Node
\draw (184,91) node [anchor=north west][inner sep=0.75pt]   [align=left] {$\displaystyle T$};
% Text Node
\draw (118,115) node [anchor=north west][inner sep=0.75pt]   [align=left] {$\displaystyle T$};
% Text Node
\draw (143,89) node [anchor=north west][inner sep=0.75pt]   [align=left] {$\displaystyle C _{1}$};
% Text Node
\draw (184,134) node [anchor=north west][inner sep=0.75pt]   [align=left] {$\displaystyle C _{2}$};
% Text Node
\draw (121,156) node [anchor=north west][inner sep=0.75pt]   [align=left] {$\displaystyle C _{3}$};
% Text Node
%\draw (62,188) node [anchor=north west][inner sep=0.75pt]   [align=left] {Figure. Neighborhood of $\displaystyle P.$};

\end{tikzpicture}
    \caption{Neighborhood of $\displaystyle P$}
    \label{fig: pants3}
\end{figure}

    We analyze the local structure of $f$ around $P$. A small neighborhood of $P$ is partitioned into six sectorial regions (see Figure \ref{fig: pants3}). The map $f$ restricts to a homeomorphism from each of these sectorial region to a sectorial region in $\H$. Moreover the orientation of $f$ on different pieces are the same. We see that either $f$ is a local homeomorphism at $P$ or $P$ is a branch point of $f$. Therefore $f$ is a branched developing map.
\end{proof}

The proof of the following simple fact is left to the reader.

\begin{claim}\label{claim-ede}
    Let $A\neq B\in\mathrm{PSL}_2(\R)$ be non-identity elements. Then there exists $Q\in\H$ such that $Q,A(Q),B(Q)$ are not collinear. Moreover, such $Q$ are dense in $\H$.
\end{claim}

\begin{lemma}\label{lem-ede-p2}
    Notation as in Lemma \ref{lem-ede-p1}. Suppose that there exists $Q_3\in\H$ such that $Q_1:=\rho(\alpha_2)\cdot Q_3,\ Q_2:=\rho(\alpha_1)^{-1}\cdot Q_3$ satisfy $d(Q_3,Q_1),d(Q_3,Q_2)<\frac{1}{3}\epsilon$. Then there exists an $\epsilon$-good developing map for $\rho$ with area $<\epsilon$.
\end{lemma}

\begin{proof}
    Notation as in the proof of Lemma \ref{lem-ede-p1}. By assumption we can choose $Q_3$ so that $l_{\beta_i}<\frac{2}{3}\epsilon\ (i=1,2,3)$. Then, by Lemma \ref{lem-h7}, the equidistant curves $\eta_i$ can be chosen to have length $<\epsilon$. So $f$ is $\epsilon$-good. Moreover, by Lemma \ref{lem-h7}, the area of each annulus is $<\frac{1}{10}\epsilon$, while the area of the triangle $U$ is clearly $<\frac{1}{10}\epsilon$. Therefore the total area is $<\epsilon$.
\end{proof}

\begin{remark}
    In this case $\rho$ is of Euler class zero, and $P$ is a branch point of the developing map $f$.
\end{remark}

We now apply this to our situation. The fact that $f_n|_{Z_{j,n}}$ converges to a constant ensures the condition $d(Q_3,Q_1),d(Q_3,Q_2)<\frac{1}{3}\epsilon$.

\begin{lemma}\label{lem-ede4}
    For every $1\leq j\leq m$ and every $\epsilon>0$, if $n$ is large enough, $\rho|_{\pi_1(Z_{j,n})}$ admits an $\epsilon$-good developing map with area $<\epsilon$.
\end{lemma}

\begin{proof}
    By assumption $Z_{j,n}$ are pairs of pants such that $\left\{B_{j,n}\circ f_n|_{\widetilde{Z}_{j,n}}\right\}_{n\in\N}$ converges to a constant (where $B_{j,n}\in\mathrm{PSL}_2(\R)$). Take $P_{j,n}\in Z_{j,n}$ so that $B_{j,n}(f_n(\widetilde{P}_{j,n}))\in\H$ converges. Let $\alpha_1,\alpha_2,\alpha_3\in\pi_1(Z_{j,n},P_{j,n})$ homotopic to the three boundary curves, and let $Q=f_n(\widetilde{P}_{j,n})$. It follows that $d(Q,\rho(\alpha_i)Q)<\frac{1}{3}\epsilon\ (i=1,2,3)$ if $n$ is large enough. The proof is completed by Lemma \ref{lem-ede-p2}.
\end{proof}

Combining Lemma \ref{lem-ede2}, Lemma \ref{lem-ede4} and Lemma \ref{lem-ede1}, we see that for $n$ large enough, there exists a pseudo $\epsilon$-excellent developing map with area $\leq A(f_n)-4\chi(\Sigma)\epsilon$. But by assumption $\lim_{n\rightarrow\infty}E_{X_n}(f_n)=e_0(\rho)=a_0(\rho)$, while $A(f_n)\leq E_{X_n}(f_n)$. So for $n$ large enough, we have $A(f_n)\leq a_0(\rho)+\epsilon$. Since $\epsilon>0$ is arbitrary, this proves Proposition \ref{prop-ede}.

\section{Moving branch points: basic theory}\label{sec-m-basic}

In this section, we recall the well known surgery of moving branch points (see, for example, \cite{GBS14}), which start with a branched developing map, produces branched developing maps with the same monodromy but the position of whose branch points are different. This procedure can be intuitively think of as moving the branch points.

We prove some basic inequality concerning how various distance functions changes (the distance functions are given by integrating the pullback of $g_\H$ as usual). For example, by moving a branch point along a path of length $l$, its distance to the boundary of $\Sigma$ can change at most by $l$.

We then prove the important estimate, Lemma \ref{lem-m2}. It says that under the assumption that $\rho$ sends every simple closed curve to a non-elliptic element, $\mathrm{inj}_P$ cannot go to zero within finite time. This reflects the intuition that if no simple closed curve is sent to elliptic, the surface does not degenerate when moving branch points, as mentioned in the introduction.

\begin{definition}
Let $f$ be a developing map on $\Sigma$. Two simple curves $\sigma_1,\sigma_2:[0,l]\rightarrow\Sigma$ lying in the interior of $\Sigma$ are called $f$-associated if
\begin{itemize}
    \item $\sigma_1(0)=\sigma_2(0)$ and $\sigma_1,\sigma_2$ intersect only at $\sigma_1(0)$;
    \item $f\circ\widetilde{\sigma}_1=f\circ\widetilde{\sigma}_2$, where $\widetilde{\sigma}_1,\widetilde{\sigma}_2:[0,l]\rightarrow\widetilde{\Sigma}$ are lifts of $\sigma_1,\sigma_2$ such that $\widetilde{\sigma}_1(0)=\widetilde{\sigma}_2(0)$.
\end{itemize}
A collection $\sigma_1,\cdots,\sigma_m$ of simple curves are called associated if they are pairwise associated.
\end{definition}

\begin{figure}
    \centering

\tikzset{every picture/.style={line width=0.75pt}} %set default line width to 0.75pt        

\begin{tikzpicture}[x=0.75pt,y=0.75pt,yscale=-1,xscale=1]
%uncomment if require: \path (0,339); %set diagram left start at 0, and has height of 339

%Straight Lines [id:da8276289839936919] 
\draw    (138,79) -- (138,134) ;
\draw [shift={(138,99.5)}, rotate = 90] [color={rgb, 255:red, 0; green, 0; blue, 0 }  ][line width=0.75]    (10.93,-4.9) .. controls (6.95,-2.3) and (3.31,-0.67) .. (0,0) .. controls (3.31,0.67) and (6.95,2.3) .. (10.93,4.9)   ;
%Straight Lines [id:da3688238826502048] 
\draw    (138,134) -- (181,161) ;
\draw [shift={(164.58,150.69)}, rotate = 212.12] [color={rgb, 255:red, 0; green, 0; blue, 0 }  ][line width=0.75]    (10.93,-4.9) .. controls (6.95,-2.3) and (3.31,-0.67) .. (0,0) .. controls (3.31,0.67) and (6.95,2.3) .. (10.93,4.9)   ;
%Straight Lines [id:da4004548099537383] 
\draw    (138,134) -- (95,158) ;
\draw [shift={(111.26,148.92)}, rotate = 330.83] [color={rgb, 255:red, 0; green, 0; blue, 0 }  ][line width=0.75]    (10.93,-4.9) .. controls (6.95,-2.3) and (3.31,-0.67) .. (0,0) .. controls (3.31,0.67) and (6.95,2.3) .. (10.93,4.9)   ;
%Straight Lines [id:da8964962189221872] 
\draw    (209,127) -- (247,127) ;
\draw [shift={(249,127)}, rotate = 180] [color={rgb, 255:red, 0; green, 0; blue, 0 }  ][line width=0.75]    (10.93,-3.29) .. controls (6.95,-1.4) and (3.31,-0.3) .. (0,0) .. controls (3.31,0.3) and (6.95,1.4) .. (10.93,3.29)   ;
%Straight Lines [id:da25916822745627266] 
\draw    (297,105) -- (333,72) ;
\draw [shift={(319.42,84.45)}, rotate = 137.49] [color={rgb, 255:red, 0; green, 0; blue, 0 }  ][line width=0.75]    (10.93,-4.9) .. controls (6.95,-2.3) and (3.31,-0.67) .. (0,0) .. controls (3.31,0.67) and (6.95,2.3) .. (10.93,4.9)   ;
%Straight Lines [id:da731909631927951] 
\draw    (297,105) -- (298,152) ;
\draw [shift={(297.63,134.5)}, rotate = 268.78] [color={rgb, 255:red, 0; green, 0; blue, 0 }  ][line width=0.75]    (10.93,-4.9) .. controls (6.95,-2.3) and (3.31,-0.67) .. (0,0) .. controls (3.31,0.67) and (6.95,2.3) .. (10.93,4.9)   ;
%Straight Lines [id:da25746980223685845] 
\draw    (298,152) -- (338,186) ;
\draw [shift={(312.67,164.47)}, rotate = 40.36] [color={rgb, 255:red, 0; green, 0; blue, 0 }  ][line width=0.75]    (10.93,-4.9) .. controls (6.95,-2.3) and (3.31,-0.67) .. (0,0) .. controls (3.31,0.67) and (6.95,2.3) .. (10.93,4.9)   ;
%Straight Lines [id:da8128102442342527] 
\draw    (338,186) -- (373,154) ;
\draw [shift={(359.93,165.95)}, rotate = 137.56] [color={rgb, 255:red, 0; green, 0; blue, 0 }  ][line width=0.75]    (10.93,-4.9) .. controls (6.95,-2.3) and (3.31,-0.67) .. (0,0) .. controls (3.31,0.67) and (6.95,2.3) .. (10.93,4.9)   ;
%Straight Lines [id:da033529794741885666] 
\draw    (373,107) -- (373,154) ;
\draw [shift={(373,136.5)}, rotate = 270] [color={rgb, 255:red, 0; green, 0; blue, 0 }  ][line width=0.75]    (10.93,-4.9) .. controls (6.95,-2.3) and (3.31,-0.67) .. (0,0) .. controls (3.31,0.67) and (6.95,2.3) .. (10.93,4.9)   ;
%Straight Lines [id:da48996063354686015] 
\draw    (333,72) -- (373,107) ;
\draw [shift={(347.73,84.89)}, rotate = 41.19] [color={rgb, 255:red, 0; green, 0; blue, 0 }  ][line width=0.75]    (10.93,-4.9) .. controls (6.95,-2.3) and (3.31,-0.67) .. (0,0) .. controls (3.31,0.67) and (6.95,2.3) .. (10.93,4.9)   ;
%Straight Lines [id:da4049028745143822] 
\draw    (426,128) -- (466,128) ;
\draw [shift={(468,128)}, rotate = 180] [color={rgb, 255:red, 0; green, 0; blue, 0 }  ][line width=0.75]    (10.93,-3.29) .. controls (6.95,-1.4) and (3.31,-0.3) .. (0,0) .. controls (3.31,0.3) and (6.95,1.4) .. (10.93,3.29)   ;
%Straight Lines [id:da8655490481302773] 
\draw    (491,94) -- (533,124) ;
\draw [shift={(516.88,112.49)}, rotate = 215.54] [color={rgb, 255:red, 0; green, 0; blue, 0 }  ][line width=0.75]    (10.93,-4.9) .. controls (6.95,-2.3) and (3.31,-0.67) .. (0,0) .. controls (3.31,0.67) and (6.95,2.3) .. (10.93,4.9)   ;
%Straight Lines [id:da7450367995963774] 
\draw    (533,124) -- (576,90) ;
\draw [shift={(549.01,111.34)}, rotate = 321.67] [color={rgb, 255:red, 0; green, 0; blue, 0 }  ][line width=0.75]    (10.93,-4.9) .. controls (6.95,-2.3) and (3.31,-0.67) .. (0,0) .. controls (3.31,0.67) and (6.95,2.3) .. (10.93,4.9)   ;
%Straight Lines [id:da36155096283085886] 
\draw    (533,124) -- (533,176) ;
\draw [shift={(533,143)}, rotate = 90] [color={rgb, 255:red, 0; green, 0; blue, 0 }  ][line width=0.75]    (10.93,-4.9) .. controls (6.95,-2.3) and (3.31,-0.67) .. (0,0) .. controls (3.31,0.67) and (6.95,2.3) .. (10.93,4.9)   ;

% Text Node
\draw (166,161.4) node [anchor=north west][inner sep=0.75pt]    {$\sigma _{1}$};
% Text Node
\draw (144,76.4) node [anchor=north west][inner sep=0.75pt]    {$\sigma _{2}$};
% Text Node
\draw (92,159.4) node [anchor=north west][inner sep=0.75pt]    {$\sigma _{3}$};
% Text Node
\draw (215,101) node [anchor=north west][inner sep=0.75pt]   [align=left] {cut};
% Text Node
\draw (289,65.4) node [anchor=north west][inner sep=0.75pt]    {$\sigma _{2}^{+}$};
% Text Node
\draw (271,115.4) node [anchor=north west][inner sep=0.75pt]    {$\sigma _{3}^{-}$};
% Text Node
\draw (298,173.4) node [anchor=north west][inner sep=0.75pt]    {$\sigma _{3}^{+}$};
% Text Node
\draw (360,66.4) node [anchor=north west][inner sep=0.75pt]    {$\sigma _{2}^{-}$};
% Text Node
\draw (382,115.4) node [anchor=north west][inner sep=0.75pt]    {$\sigma _{1}^{+}$};
% Text Node
\draw (363,171.4) node [anchor=north west][inner sep=0.75pt]    {$\sigma _{1}^{-}$};
% Text Node
\draw (433,99) node [anchor=north west][inner sep=0.75pt]   [align=left] {glue};

\end{tikzpicture}

    \caption{Move along associated curves}
    \label{fig: Move along associated curves}
\end{figure}
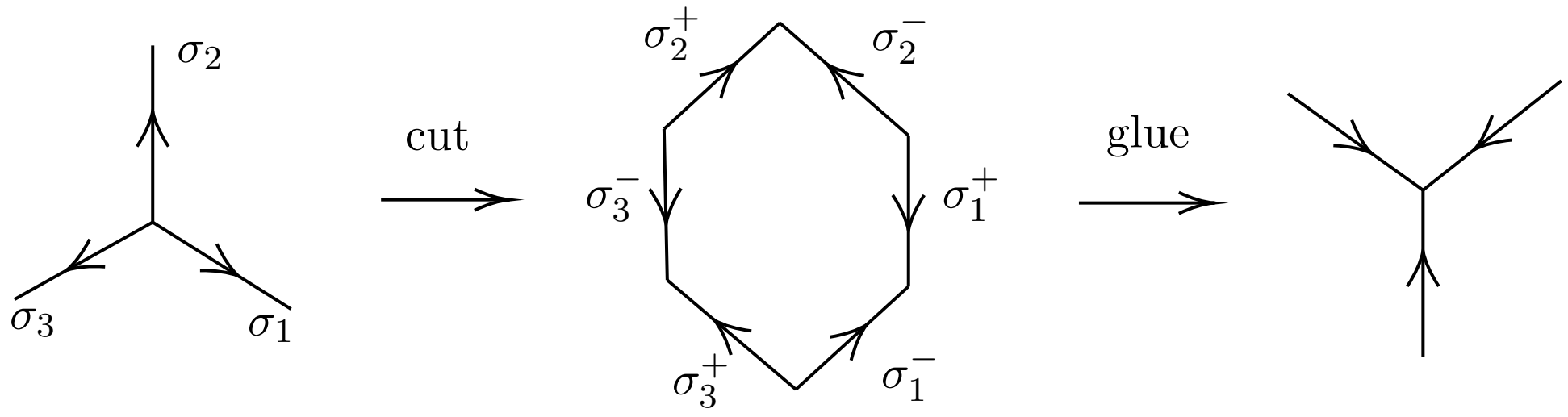

Let $\sigma_1,\cdots,\sigma_m$ be a collection of associated curves, with $\overrightarrow{\sigma_1},\cdots,\overrightarrow{\sigma_m}$ arranged in anticlockwise order. Cut $\Sigma$ along $\sigma_1,\cdots,\sigma_m$. The resulting surface $\Sigma^{cut}$ has a boundary curve composed of $\sigma_1^-,\sigma_1^+,\cdots,\sigma_m^-,\sigma_m^+$, where $\sigma_i^-,\sigma_i^+$ are copies of $\sigma_i$, all parametrized by $[0,l]$, such that $\sigma_i^+(0)=\sigma_{i+1}^-(0)\ (\forall i\in\Z/m\Z)$ (Figure \ref{fig: Move along associated curves} indicates the case $m=3$).

Now glue back by identifying, for every $t\in[0,l],i\in\Z/m\Z$, $\sigma_i^+(t)$ with $\sigma_{i+1}^-(t)$. The resulting surface $\Sigma_l$ is naturally homeomorphic to $\Sigma$. Since the gluing respects the developing map, $f$ defines a developing map $f_l$ on $\Sigma_l$, which has the same monodromy as $f$ if we identify $\Sigma_l$ with $\Sigma$ via a natural homeomorphism. We call $(\Sigma_l,f_l)$ the \textit{move of $(\Sigma,f)$ along the associated curves} $\sigma_1,\cdots,\sigma_m$.

If $f$ is a branched developing map, then $f_l$ is also branched. Moreover, $(\Sigma_l,f_l^*g_\H)\in[(\Sigma,g_f)]$, where $[S]$ is defined as follows.

\begin{definition}
    Let $f$ be a branched developing map on $\Sigma$. Define $[f]$ to be the set of branched developing maps that is homotopic to $f$ relative to $\partial\Sigma$.
\end{definition}

\begin{definition}
    Let $S=(\Sigma,g)$ be a branched hyperbolic surface. Let $f=\mathrm{dev}_g$. Define $[S]$ to be the set of pairs $S_1=(\Sigma_1,g_1)$ for which there exists $f_1\in[f]$ such that $S_1\cong(\Sigma,f_1^*g_{\H})$.
\end{definition}

\begin{remark}
    If $S_1\in[S]$, then there is an isometry $\partial S_1\cong\partial S$. Indeed, by definition $S_1\cong(\Sigma,f_1^*g_{\H})$ for some $f_1\in[f]$. Note that $f_1^*g_{\H}|_{\partial\Sigma}=f^*g_{\H}|_{\partial\Sigma}$. We get $\partial S_1\cong\partial S$. Therefore, for any $S_1\in[S]$, we always identify $\partial S_1$ with $\partial S$.
\end{remark}

Let $S=(\Sigma,g)$. A collection of simple curves $\sigma_1,\cdots,\sigma_m$ are called $g$-associated if they are $\mathrm{dev}_g$-associated. Suppose that $\sigma_1,\cdots,\sigma_m$ are $g$-associated, denote by $g_{cut}$ the pullback of $g$ to $\Sigma^{cut}$. Then $g^{cut}$ descends to a metric $g_l$ on $\Sigma_l$, which is precisely $f_l^*g_{\H}$. Therefore $S_l=(\Sigma_l,g_l)\in[S]$.

All branch points of $S$ away from $\sigma_1\cup\cdots\cup\sigma_m$ remain branch points of $S_l$. Other than these, assuming that $\sigma_1,\cdots,\sigma_m$ have no interior branch points, the only possible branch points of $S_l$ are the images of $\sigma_i^{\pm}(0),\sigma_i^{\pm}(l)$.

Let $P=\sigma_1(0)$. Then $P$ is a branch point of multiplicity at least $m-1$. The collection $\sigma_1,\cdots,\sigma_m$ is called \textit{complete} if each $\sigma_i$ does not contain interior branch points, and the multiplicity of $P$ is equal to $m-1$.

Now suppose that $\sigma_1,\cdots,\sigma_m$ is complete, intersecting at $P$. Let $P_l\in\Sigma_l$ be the point obtained by identifying $\sigma_1(l),\cdots,\sigma_m(l)$. The completeness implies that $\measuredangle(\overrightarrow{\sigma_i},\overrightarrow{\sigma_{i+1}})=2\pi\ (\forall i\in\Z/m\Z)$, thus the images of $\sigma_i^{\pm}(0)$ in $S_l$ are not branch points. We get

\begin{lemma}\label{lem-m0}
    Suppose that $\sigma_1,\cdots,\sigma_m$ is complete, then $P_l$ is the only branch point of $S_l$ other than the branch points of $S$ away from $\sigma_1\cup\cdots\cup\sigma_m$.
\end{lemma}

In particular, we have

\begin{lemma}\label{lem-m01}
    Suppose that $\sigma_1,\cdots,\sigma_m$ is complete, and $P$ is the only branch point of $S$ lying on $\sigma_1\cup\cdots\cup\sigma_m$. Then $\mathrm{Br}(S_l)=(\mathrm{Br}(S)\backslash\{P\})\cup\{P_l\}$.
\end{lemma}

Suppose that $\sigma_i$ have length $l$ with respect to $g$.

\begin{lemma}\label{lem-m11}
    Let $Q$ be any branch point of $S$ away from $\sigma_1\cup\cdots\cup\sigma_m$. Then
    \[d_S(P,Q)-l\leq d_{S_l}(P_l,Q)\leq d_S(P,Q)+l.\]
\end{lemma}

\begin{proof}
    Let $\sigma=\bigcup_i(\sigma_i^+\cup\sigma_i^-)$ be the boundary of $\Sigma^{cut}$ arising from the cut. Then we have
    \[d_S(P,Q)-l\leq d_{S^{cut}}(\sigma,Q)\leq d_S(P,Q),\]
    \[d_{S_l}(P_l,Q)-l\leq d_{S^{cut}}(\sigma,Q)\leq d_{S_l}(P_l,Q).\]
    Combining these two inequalities proves the lemma.
\end{proof}

Let $\gamma$ be a boundary curve of $\Sigma$. Recall the definition of $r_\gamma$ from Subsection \ref{subsec-c nbhd}. The following lemma is proved in the same way as Lemma \ref{lem-m11} and is left to the reader.

\begin{lemma}\label{lem-m12}
    For any boundary curve $\gamma$ of $\Sigma$, we have
    \[d_S(P,\gamma)-l\leq d_{S_l}(P_l,\gamma)\leq d_S(P,\gamma)+l,\]
    \[r_\gamma(S)-l\leq r_{\gamma}(S_l)\leq r_\gamma(S)+l.\]
\end{lemma}

Define
\[r_P(S)=d_S(P,\partial S\cup\mathrm{Br}(S)\backslash\{P\}).\]
Combining the previous two lemmas we get

\begin{lemma}\label{lem-m13}
    \[r_P(S)-l\leq r_{P_l}(S_l)\leq r_P(S)+l.\]
\end{lemma}

Define
\[s_P(S)=\min\{l_{\alpha}:\alpha\in\pi_1(\Sigma,P)\text{ non trivial}\},\]
\[\mathrm{inj}_P(S)=\min\{r_P(S),\frac{1}{2}s_P(S)\}.\]

\begin{lemma}\label{lem-m2}
    Suppose that the holonomy of $g$ around every simple closed curve is non-elliptic, then
    \[\mathrm{inj}_{P_l}\geq\min\{r_P-l,e^{-l}\mathrm{inj}_P\}.\]
\end{lemma}

\begin{proof}
    We may assume that $\sigma_i$ are parametrized by arc length. We have a family $(S_t,P_t)\ (t\in[0,l])$, obtained by the moves along $\sigma_i|_{[0,t]}$. Let $\lambda(t)=\mathrm{inj}_{P_t}$. Then what we need to prove is $\lambda(l)\geq\min\{r_P-l,e^{-l}\lambda(0)\}$. We first prove the following
    \begin{claim}
        If $\lambda(t)< r_P-t$ for some $t\in(0,l]$, then
        \[\limsup_{\Delta t\rightarrow0+}\frac{\log \lambda(t-\Delta t)-\log \lambda(t)}{\Delta t}\leq 1.\]
    \end{claim}
    By Lemma \ref{lem-m13} we have $r_{P_t}\geq r_P-t$. If $\lambda(t_0)< r_P-t_0$, then $\lambda(t_0)<r_{P_{t_0}}$. By definition $\lambda(t)=\mathrm{inj}_{P_t}=\min\{r_{P_t},\frac{1}{2}s_{P_t}\}$. So we have $\lambda(t_0)=\frac{1}{2}s_{P_{t_0}}<r_{P_{t_0}}$. Then there is a minimal geodesic loop $\alpha_{t_0}$ based at $P_{t_0}$ which does not meet $\partial\Sigma_{t_0}$, and without interior branch point. Then for $t$ sufficiently close to $t_0$ there is a geodesic loop $\alpha_{t}$ based at $P_{t}$ homotopic to $\alpha_{t_0}$. Let $A=\rho(\alpha_t)$, which is non-elliptic by assumption. Then $f_t(\widetilde{\alpha}_t)$ are the geodesic segments from $f_t(\widetilde{P}_t)$ to $A\cdot f_t(\widetilde{P}_t)$. Applying Lemma \ref{lem-h3} to the curve $t\mapsto f_t(\widetilde{P}_t)$, we get $\frac{d}{dt}\log l_{\alpha_t}\geq-1$. Since $\lambda(t)\leq \frac{1}{2}l_{\alpha_t},\ \lambda(t_0)=\frac{1}{2}l_{\alpha_{t_0}}$, the claim follows.

    If $\lambda(t)<r_P-t$ for all $t\in[0,1]$ then the claim implies $\lambda(l)\geq e^{-l}\lambda(0)$, the lemma follows. If not, define $t_0=\sup\{t\in[0,l]:\lambda(t)\geq r_P-t\}$. If $t_0=l$ then $\lambda(l)\geq r_P-l$, the lemma also follows. Otherwise suppose $t_0<l$, we have $\lambda(t_0)=r_P-t_0$. The claim implies $\lambda(l)\geq e^{t_0-l}\lambda(t_0)$. We show that $e^{t_0-l}\lambda(t_0)\geq\min\{r_P-l,e^{-l}\mathrm{inj}_P\}$.
    
    If $\lambda(t_0)\leq1$, then $e^{t_0-l}\lambda(t_0)\geq\lambda(t_0)+(e^{t_0-l}-1)\geq \lambda(t_0)+t_0-l=r_P-l$ and we are done. Otherwise $\lambda(t_0)>1$, then $e^{t_0-l}\lambda(t_0)=e^{-l}e^{t_0}\lambda(t_0)\geq e^{-l}(\lambda(t_0)+t_0)=e^{-l}r_P\geq e^{-l}\mathrm{inj}_P$ and we are also done.
\end{proof}

\section{Moving branch points: applications}\label{sec-m-app}

In this section we solve two problems: merging nearby branch points and moving branch points away from a boundary curve.

Throughout this section let $S=(\Sigma,g)$ be a branched hyperbolic surface, such that the holonomy of $g$ around every simple closed curve is non-elliptic. If $\gamma$ is a component of $\partial\Sigma$ we define
\[R_\gamma=\sup\{r_\gamma(S_1):S_1\in[S]\}.\]

\begin{definition}
    Let $P$ be a branch point of $g$. A pair $(S_1,P_1)$, where $P_1\in S_1$ is a branch point, is called an $l$-move of $(S,P)$ if $S_1\in[S]$, $\#\mathrm{Br}(S_1)\leq\#\mathrm{Br}(S)$, $r_{P_1}\geq r_P-l$ and for every boundary curve $\gamma$ of $\Sigma$, we have
    \[d(P,\gamma)-l\leq d(P_1,\gamma)\leq d(P,\gamma)+l,\quad r_{\gamma}(S_1)\geq r_\gamma(S)-l.\]
\end{definition}

\subsection{Merging nearby branch points}

Throughout this subsection let $P$ be a branch point of $g$. We show how to merge $P$ with nearby branch points. The principal result is Proposition \ref{prop-m1}, which intuitively says that by merging branch points near $P$, we end up with a developing map in which there is no branch point nearby $P$.

The strategy to merge branch points is intuitively clear. Let $Q$ be a branch point nearest to $P$. Let $\sigma$ be a minimizing geodesic from $P$ to $Q$. We move $P$ along $\sigma$ to collide it with $Q$. The move should be made in steps, each moving a distance less than $\mathrm{inj}_P$. Lemma \ref{lem-m2} makes sure that $\mathrm{inj}_P$ does not go to zero, so that we can eventually reach $P$. This process, merging with a nearest branch point, is the content of Lemma \ref{lem-ma3}. The branch point $Q$ will be implicit in the proof, but we are really thinking of approaching a branch point $Q$ inductively.

We then prove Proposition \ref{prop-m1}, by applying Lemma \ref{lem-ma3} iteratively, merging branch points nearby $P$ one after another, until there is no branch point nearby.

\begin{lemma}\label{lem-ma1}
    Suppose that $r_P<d(P,\partial S)$ and $\mathrm{inj}_P=r_P$. Then there is an $r_P$-move $(S_1,P_1)$ such that $\#\mathrm{Br}(S_1)<\#\mathrm{Br}(S)$.
\end{lemma}

\begin{proof}
    Since $r_P<d(P,\partial S)$, we have $r_P=d(P,\mathrm{Br}(S)\backslash\{P\})$. So there is a branch point $Q\neq P$ such that $r_P=d(P,Q)$. Let $\sigma:[0,l]\rightarrow\Sigma$ be a minimizing geodesic from $P$ to $Q$. Let $\sigma_0,\cdots,\sigma_m:[0,l]\rightarrow S$ be the geodesics such that $\sigma_k(0)=P$ and $\measuredangle(\overrightarrow{\sigma},\overrightarrow{\sigma_k})=2k\pi$, where $m=m(P)$.
    \begin{claim}
        $\sigma_0,\cdots,\sigma_m$ intersect only at $P$.
    \end{claim}
    Indeed, since $l=\mathrm{inj}_P$, if $\sigma_i,\sigma_j\ (i\neq j)$ intersect at some point other then $P$, it must be the case that $\sigma_i(l)=\sigma_j(l)$. But this would imply that $\rho(\sigma_i*\overline{\sigma}_j)$ is an elliptic element, contradicting our assumption.

    Therefore $\sigma_0,\cdots,\sigma_m$ is a complete collection of associated curves. Let $S_1$ be the move along $\sigma_0,\cdots,\sigma_m$, and let $P_1\in S_1$ be the point obtained by identifying $\sigma_0(l),\cdots,\sigma_m(l)$. By Lemma \ref{lem-m0}, all the branch points of $S_1$ other than $P_1$ are the branch points of $S$ away from $\sigma_1\cup\cdots\cup\sigma_m$. Note that $S$ has at least two branch points on $\sigma_1\cup\cdots\cup\sigma_m$, namely $P$ and $Q$. Therefore $\#\mathrm{Br}(S_1)<\#\mathrm{Br}(S)$. By Lemma \ref{lem-m12} and Lemma \ref{lem-m13} we deduce that $(S_1,P_1)$ is an $l$-move of $(S,P)$.
\end{proof}

\begin{lemma}\label{lem-ma2}
    Suppose that $r_P<d(P,\partial S)$. Let $l<\mathrm{inj}_P$. Then there is an $l$-move $(S_1,P_1)$ such that $r_{P_1}=r_P-l,\ \mathrm{inj}_{P_1}\geq\min\{r_{P_1},e^{-l}\mathrm{inj}_P\}$.
\end{lemma}

\begin{proof}
    Since $r_P<d(P,\partial S)$, there is a branch point $Q\neq P$ such that $r_P=d(P,Q)$. Let $\sigma:[0,L]\rightarrow S$ be a minimizing geodesic from $P$ to $Q$. Let $\sigma_0,\cdots,\sigma_m:[0,l]\rightarrow S$ be the (unit speed) geodesics such that $\sigma_k(0)=P$ and $\measuredangle(\overrightarrow{\sigma},\overrightarrow{\sigma_k})=2k\pi$, where $m=m(P)$. Then $\sigma_0,\cdots,\sigma_m$ is a complete collection of associated curves. Let $S_1$ be the move along $\sigma_0,\cdots,\sigma_m$, and let $P_1\in S_1$ be the point obtained by identifying $\sigma_0(l),\cdots,\sigma_m(l)$.
    \begin{claim}
        For $k\neq0$, $\sigma_k$ does not intersect $\sigma$
    \end{claim}
    $\sigma_k$ does not intersect $\sigma|_{[0,l]}$ because $l<\mathrm{inj}_P$, and $\sigma_k$ does not intersect $\sigma|_{(l,L]}$ because $\sigma$ is minimizing.

    Therefore $\sigma|_{[l,L]}$ is a geodesic on $S_1$ from $P_1$ to $Q$. This implies $d(P_1,Q)\leq L-l$. Since $r_{P_1}\leq d(P_1,Q)$, $L=r_P$, we get $r_{P_1}\leq r_P-l$. However, by Lemma \ref{lem-m13}, $r_{P_1}\geq r_P-l$. So we have $r_{P_1}= r_P-l$. The estimate on $\mathrm{inj}_{P_1}$ follows from \ref{lem-m2}.
\end{proof}

\begin{lemma}\label{lem-ma3}
    Lemma \ref{lem-ma1} holds without the assumption $\mathrm{inj}_P=r_P$.
\end{lemma}

\begin{proof}
    Let $l=\frac{r_P}{M}$ for some large enough integer $M$ such that $l<e^{-r_P}\mathrm{inj}_P$. Start with $(S_0,P_0)=(S,P)$, we construct $kl$-moves $(S_k,P_k)$ inductively in the way described in the next paragraph, satisfying $r_{P_k}=r_P-kl<d(P_k,\partial S),\ \mathrm{inj}_{P_k}\geq\min\{r_{P_k},e^{-kl}\mathrm{inj}_P\}$, until we reach some $k$ such that $\mathrm{inj}_{P_k}=r_{P_k}$.

    Suppose that $(S_k,P_k)$ is constructed, and $\mathrm{inj}_{P_k}<r_{P_k}$. Since $\mathrm{inj}_{P_k}\geq\min\{r_{P_k},e^{-kl}\mathrm{inj}_P\}$, we get $\mathrm{inj}_{P_k}\geq e^{-kl}\mathrm{inj}_P>l$. So we can apply Lemma \ref{lem-ma2} to find $(S_{k+1},P_{k+1})$, which is an $l$-move of $(S_k,P_k)$, with the desired properties.

    The process must stop at some $k$ with $k<M$. We then apply Lemma \ref{lem-ma1} to finish the proof.
\end{proof}

\begin{proposition}\label{prop-m1}
    Let $m$ be the number of branch points of $g$. Suppose that $d(P,\partial S)\geq 10^{m+1}l>0$. Then there exists $L\in[l,10^ml]$ and an $L$-move $(S',P')$ such that $r_{P'}\geq 9L$.
\end{proposition}

\begin{proof}
    Start with $(S_1,P_1)=(S,P)$, we construct $(S_k,P_k)$ with $S_k\in[S]$ inductively as follows.
    \begin{itemize}
        \item If $r_{P_k}<\min\{d(P_k,\partial S),10^kl\}$, apply Lemma \ref{lem-ma3} to construct a $10^kl$-move $(S_{k+1},P_{k+1})$ of $(S_k,P_k)$ such that $\#\mathrm{Br}(S_{k+1})<\#\mathrm{Br}(S_k)$.
        \item If $r_{P_k}\geq\min\{d(P_k,\partial S),10^kl\}$, stop.
    \end{itemize}
    Since there are $m$ branch points, and the number of branch points drops at each step, the procedure must stop at some $k\leq m$. We take $(S',P')=(S_k,P_k)$. Then $(S',P')$ is an $L$-move of $(S,P)$, where $L=\frac{10^k}{9}l>\sum_{i=1}^{k-1}10^il$, and we have $r_{P'}\geq\min\{d(P',\partial S),9L\}$. Note that $d(P',\partial S)\geq9L$. Indeed, by assumption $d(P,\partial S)\geq 10^{m+1}l\geq10L$, so we have $d(P',\partial S)\geq d(P,\partial S)-L\geq9L$. Therefore $r_{P'}\geq9L$.
\end{proof}

For good developing maps, the number of branch points can be bounded by the Euler characteristic of $\Sigma$.

\begin{lemma}\label{lem-number of branch point}
    Let $f$ be a good developing map on $\Sigma$. Then $f$ has at most $-\chi(\Sigma)$ branch points.
\end{lemma}

\begin{proof}
    Gauss-Bonnet theorem and Lemma \ref{lem1}.
\end{proof}

\subsection{Moving one branch point away from a boundary curve}

In this subsection we assume that a boundary curve $\gamma$ is an equidistant curve of length $\leq\frac{1}{100}$. We show how to move the nearest branch point away from $\gamma$. The principal result is the following

\begin{proposition}\label{prop-m2}
    Let $P$ be a branch point such that $d(P,\gamma)=r_\gamma\leq\frac{1}{2}$, $r_P\geq100l_\gamma$. Then
    \[R_\gamma\geq r_\gamma+\frac{1}{3}r_P.\]
\end{proposition}

The condition $d(P,\gamma)=r_\gamma$ means that $P$ is the nearest branch point to $\gamma$. In our applications $l_\gamma$ would be small, so the condition $r_P\geq100l_\gamma$ would be automatic.

The proof of the proposition is obtained by applying a finite number moves of $P$, each time moving away from $\gamma$ a little. The first lemma tells us that there is a correct direction to move.

\begin{lemma}\label{lem-m6}
    Let $P\in S$ (possibly a branch point). Let $\alpha,\beta$ be two simple orthogonal geodesics from $P$ to $\gamma$ of length $\leq1$. Suppose that $\alpha,\beta$ intersect only at $P$, and $|l_{\alpha}-l_{\beta}|\leq\frac{1}{100}$. Then the angle between $\alpha,\beta$ at $P$ is of the form $2k\pi+\theta$ for some $k\in\Z,|\theta|<\frac{1}{10}$.
\end{lemma}

\begin{proof}
    Let $Q_1,Q_2\in\gamma$ be the feet of $\alpha,\beta$. Let $\sigma\subset \gamma$ be an arc from $Q_1$ to $Q_2$ of length $\leq\frac{1}{2}l_\gamma$. Let $\eta$ be the simple closed curve $\alpha*\sigma*\overline{\beta}$. Let $\widetilde{\eta}=\widetilde{\alpha}*\widetilde{\sigma}*\widetilde{\bar{\beta}}$ be a lift of $\eta$ to $\widetilde{\Sigma}$. Let $\eta_1=\alpha_1*\sigma_1*\overline{\beta_1}=f(\widetilde{\eta})$, which is a curve in $\H$ from $P_1$ to $P_2$ (the red curve in Figure \ref{fig: lem-m6}). Let $Q_1'=f(\widetilde{Q}_1),Q_2=f(\widetilde{Q}_2)$. We have $\rho(\eta)\cdot P_1=P_2$, and the angle between $\rho(\eta)\cdot \alpha_1$ and $\beta_1$ at $P_2$ is equal to the angle between $\alpha,\beta$ at $P$ modulo $2\pi$. By assumption $\rho(\eta)$ is a hyperbolic or parabolic element. Let $\mathcal{P}:T_{P_1}\H\rightarrow T_{P_2}\H$ be the parallel transport along the geodesic connecting $P_1$ and $P_2$. Then by Lemma \ref{lem-h1}, the angle between $\mathcal{P}(\overrightarrow{\alpha_1})$ and $\rho(\eta)_*(\overrightarrow{\alpha_1})$ is less than $d(P_1,P_2)$. By Lemma \ref{lem-c31}, we have $d(P_1,P_2)\leq 8l_\sigma+|l_{\alpha}-l_{\beta}|$. Meanwhile, The geodesic quadrilateral with vertices $P_1,P_2,Q_2',Q_1'$ has area $\leq 8l_\sigma$. Thus, by the Gauss-Bonnet theorem and the fact that the angles at $Q_1'$ and $Q_2'$ are right angles, we deduce that $\angle P_1+\angle P_2\in[\pi-8l_\sigma,\pi]$. This shows that the angle between $\mathcal{P}(\overrightarrow{\alpha_1})$ and $\overrightarrow{\beta_1}$ is $\leq 8l_\sigma$. Combining these we see that the angle between $\rho(\eta)_*(\overrightarrow{\alpha_1})$ and $\overrightarrow{\beta_1}$ is less than $16l_\sigma+|l_{\alpha}-l_{\beta}|$. Since by assumption $l_\sigma\leq\frac{1}{2}l_\gamma\leq\frac{1}{200},\ |l_{\alpha}-l_{\beta}|\leq\frac{1}{100}$, we deduce that the angle between $\rho(\eta)\cdot \alpha_1$ and $\beta_1$ at $P_2$ is less than $\frac{1}{10}$, finishing the proof.
\end{proof}

\begin{figure}
    \centering
    \tikzset{every picture/.style={line width=0.75pt}} %set default line width to 0.75pt        

\begin{tikzpicture}[x=0.75pt,y=0.75pt,yscale=-1,xscale=1]

%Curve Lines [id:da7938046411269434] 
\draw    (269,204) .. controls (250,180) and (249,110) .. (267,83) ;
%Straight Lines [id:da9759204578684346] 
\draw [color={rgb, 255:red, 208; green, 2; blue, 27 }  ,draw opacity=1 ]   (170,126) -- (254,134) ;
%Straight Lines [id:da7633449986350773] 
\draw [color={rgb, 255:red, 208; green, 2; blue, 27 }  ,draw opacity=1 ]   (169,166) -- (255,158) ;
%Curve Lines [id:da17719085142115953] 
\draw    (169,166) .. controls (165,155) and (165,138) .. (170,126) ;
%Straight Lines [id:da2259097512096705] 
\draw [color={rgb, 255:red, 208; green, 2; blue, 27 }  ,draw opacity=1 ]   (254,134) -- (255,158) ;

% Text Node
\draw (271,59.4) node [anchor=north west][inner sep=0.75pt]    {$f(\widetilde{\gamma} )$};
% Text Node
\draw (256,140.4) node [anchor=north west][inner sep=0.75pt]    {$\sigma_1$};
% Text Node
\draw (258,113.4) node [anchor=north west][inner sep=0.75pt]    {$Q_{1}'$};
% Text Node
\draw (257,161.4) node [anchor=north west][inner sep=0.75pt]    {$Q_{2}'$};
% Text Node
\draw (308,210.4) node [anchor=north west][inner sep=0.75pt]    {$\mathbb{H}$};
% Text Node
\draw (200,115.4) node [anchor=north west][inner sep=0.75pt]    {$\alpha_1$};
% Text Node
\draw (200,167.4) node [anchor=north west][inner sep=0.75pt]    {$\beta_1$};
% Text Node
\draw (149,111.4) node [anchor=north west][inner sep=0.75pt]    {$P_{1}$};
% Text Node
\draw (152,167.4) node [anchor=north west][inner sep=0.75pt]    {$P_{2}$};
% Text Node
%\draw (216,240) node [anchor=north west][inner sep=0.75pt]   [align=left] {Figure.};

\end{tikzpicture}
    \caption{Lemma \ref{lem-m6}}
    \label{fig: lem-m6}
\end{figure}

\begin{lemma}\label{lem-m7}
    Suppose that $d(P,\gamma)=r_\gamma$. Let $l\leq\frac{1}{5}\mathrm{inj}_P$. If $r_\gamma+l\leq 1$, then there is an $l$-move $(S_1,P_1)$ such that
    \[d(P_1,\gamma)\geq r_\gamma+\frac{1}{2}l.\]
\end{lemma}

\begin{proof}
    Let $\alpha$ be a minimizing geodesic from $P$ to $\gamma$. Suppose that $P$ has multiplicity $m-1$. Let $\sigma_1,\cdots,\sigma_m:[0,l]\rightarrow S$ be the geodesics such that $\sigma_k(0)=P$ and $\measuredangle(\overrightarrow{\alpha},\overrightarrow{\sigma_k})=(2k-1)\pi$. Then $\sigma_1,\cdots,\sigma_m$ is a complete collection of associated curves. Let $(S1,P_1)$ be the move along $\sigma_1,\cdots,\sigma_m$. Let $\eta_k\ (k\in\Z/m\Z)$ be the curve on $S_1$ obtained by gluing $\sigma_k^-$ and $\sigma_{k-1}^+$ (see Figure \ref{fig: lem-m71}). Then $\beta_1:=\overline{\eta}_1*\alpha$ is a simple geodesic from $P_1$ to $\gamma$ of length $r_\gamma+l$. Let $\beta_2$ be a minimizing geodesic on $S_1$ from $P_1$ to $\gamma$. We need to prove that $l_\beta\geq r_\gamma+\frac{1}{2}l$. If $\beta_2=\beta_1$ then we are done. Now suppose $\beta_2\neq\beta_1$.
    
\begin{figure}
    \centering
    \tikzset{every picture/.style={line width=0.75pt}} %set default line width to 0.75pt        

\begin{tikzpicture}[x=0.75pt,y=0.75pt,yscale=-1,xscale=1]

%Curve Lines [id:da27265861811515646] 
\draw    (100,104) .. controls (117,97) and (136,73) .. (137,54) ;
%Straight Lines [id:da11367544250477968] 
\draw    (125,82) -- (200,139) ;
%Straight Lines [id:da9251477930776841] 
\draw    (201,100) -- (200,139) ;
%Straight Lines [id:da6097477645122852] 
\draw    (200,139) -- (234,123) ;
%Straight Lines [id:da6886881886867985] 
\draw    (200,139) -- (227,168) ;
%Straight Lines [id:da16603561514018017] 
\draw    (200,139) -- (162,144) ;
%Straight Lines [id:da7325712949612566] 
\draw    (200,139) -- (195,171) ;
%Straight Lines [id:da8960254755358535] 
\draw    (268,119.5) -- (310,119.5)(268,122.5) -- (310,122.5) ;
\draw [shift={(318,121)}, rotate = 180] [color={rgb, 255:red, 0; green, 0; blue, 0 }  ][line width=0.75]    (10.93,-3.29) .. controls (6.95,-1.4) and (3.31,-0.3) .. (0,0) .. controls (3.31,0.3) and (6.95,1.4) .. (10.93,3.29)   ;
%Curve Lines [id:da9237854644925165] 
\draw    (352,86) .. controls (373,81) and (391,67) .. (396,47) ;
%Straight Lines [id:da6675681225827035] 
\draw    (379,73) -- (461,149) ;
%Straight Lines [id:da34713932199609154] 
\draw    (461,149) -- (480,122) ;
%Straight Lines [id:da9006338786173894] 
\draw    (461,149) -- (495,158) ;
%Straight Lines [id:da8561826522932109] 
\draw    (461,149) -- (464,179) ;
%Straight Lines [id:da7725391500835564] 
\draw    (432,156) -- (461,149) ;
%Shape: Circle [id:dp9836369465037146] 
\draw  [fill={rgb, 255:red, 0; green, 0; blue, 0 }  ,fill opacity=1 ] (437,129) .. controls (437,127.9) and (437.9,127) .. (439,127) .. controls (440.1,127) and (441,127.9) .. (441,129) .. controls (441,130.1) and (440.1,131) .. (439,131) .. controls (437.9,131) and (437,130.1) .. (437,129) -- cycle ;

% Text Node
\draw (187,120.4) node [anchor=north west][inner sep=0.75pt]    {$\pi $};
% Text Node
\draw (176,128.4) node [anchor=north west][inner sep=0.75pt]    {$\pi $};
% Text Node
\draw (203,117.4) node [anchor=north west][inner sep=0.75pt]    {$2\pi $};
% Text Node
\draw (156,93.4) node [anchor=north west][inner sep=0.75pt]    {$\alpha $};
% Text Node
\draw (145,137.4) node [anchor=north west][inner sep=0.75pt]    {$\sigma _{1}$};
% Text Node
\draw (189,87.4) node [anchor=north west][inner sep=0.75pt]    {$\sigma _{m}$};
% Text Node
\draw (202,142.4) node [anchor=north west][inner sep=0.75pt]    {$P$};
% Text Node
\draw (118,39.4) node [anchor=north west][inner sep=0.75pt]    {$\gamma$};
% Text Node
%\draw (249,132) node [anchor=north west][inner sep=0.75pt]   [align=left] {simple move};
% Text Node
\draw (374,33.4) node [anchor=north west][inner sep=0.75pt]    {$\gamma$};
% Text Node
\draw (418,129.4) node [anchor=north west][inner sep=0.75pt]    {$P$};
% Text Node
\draw (445,151.4) node [anchor=north west][inner sep=0.75pt]    {$P_1$};
% Text Node
\draw (411,85.4) node [anchor=north west][inner sep=0.75pt]    {$\alpha $};
% Text Node
\draw (448,121.4) node [anchor=north west][inner sep=0.75pt]    {$\eta_{1}$};
% Text Node
%\draw (270,194) node [anchor=north west][inner sep=0.75pt]   [align=left] {Figure.};

\end{tikzpicture}
    \caption{Move $P$ away from $\gamma$.}
    \label{fig: lem-m71}
\end{figure}
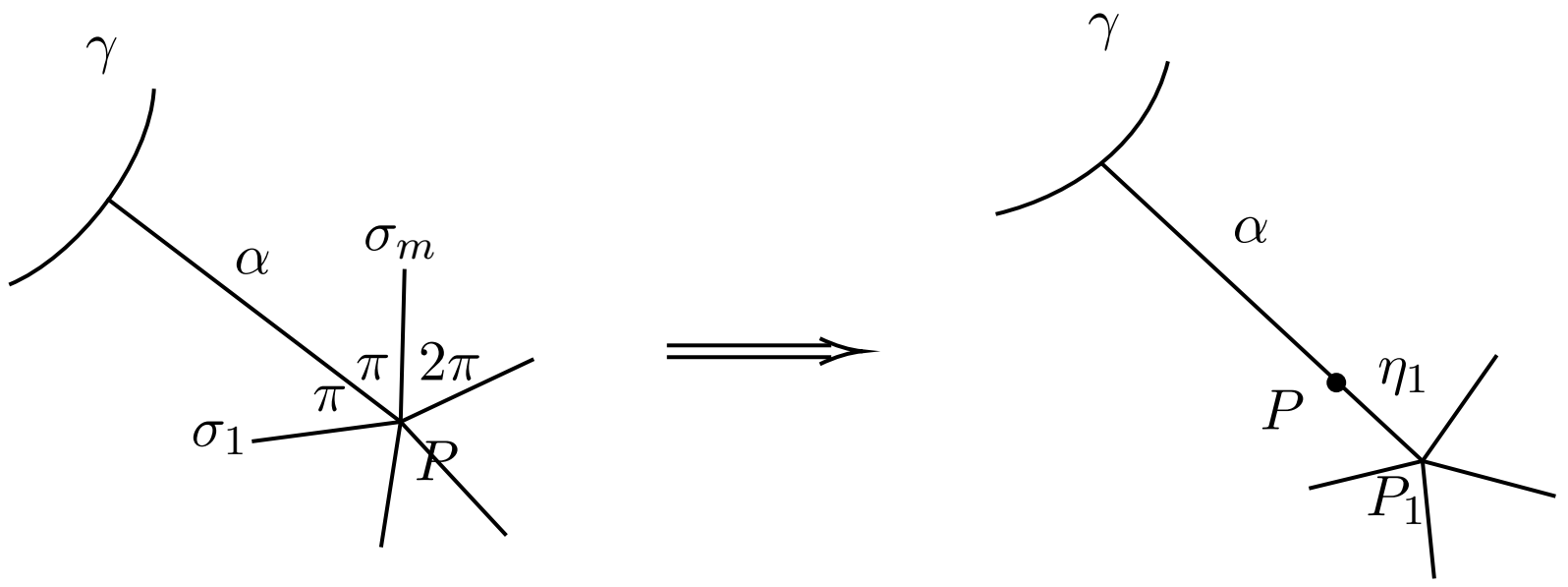

    \begin{claim}
        $\beta_1$ and $\beta_2$ intersect only at $P_1$.
    \end{claim}
    \begin{itemize}
        \item $\beta_1|_{(0,3l)}$ and $\beta_2|_{(0,3l)}$ do not intersect.

        This is because $\mathrm{inj}_{P_1}\geq 3l$, which follows from $\mathrm{inj}_P\geq5l$.
        \item $\beta_1|_{(0,3l)}$ and $\beta_2|_{[3l,l_{\beta_2}]}$ do not intersect.

        This is because $\beta_2$ is minimizing.
        \item $\beta_1|_{[3l,l_{\beta_1}]}$ and $\beta_2$ do not intersect.

        Since $\mathrm{inj}_\gamma(S_1)\geq\mathrm{inj}_\gamma-l=r_\gamma-l=l_{\beta_1}-2l$, $\gamma$ has a collar neighborhood of length $l_{\beta_1}-3l$ in $S_1$. $\beta_1|_{[3l,l_{\beta_1}]}$ is a minimizing geodesic between the two boundary curves of this collar neighborhood, so it can not intersect $\beta_2$.
    \end{itemize}
    These prove the claim.

    By Lemma \ref{lem-m6}, the angle formed by $\beta_1,\beta_2$ at $P_1$ is an integral multiple of $2\pi$ up to an error of $\frac{1}{10}$. So the angle between $\beta_2$ and $\eta_k$ is $\leq\frac{1}{10}$ for some $k$. This will imply that $r_\gamma\leq l_{\beta_2}-\frac{1}{2}l$. Indeed, let $Q$ be the point on $\beta_2$ such that $d(P_1,Q)=l$. One can join $P$ and $Q$ by a curve of length $\leq\frac{1}{2}l$ \footnote{Explicitly, such a curve can be written down by $t\mapsto\exp_{P_1}(l\cdot e^{it})\ (t\in[t_1,t_2])$, where $\exp_{P_1}(l\cdot e^{it_1})=P,\ \exp_{P_1}(l\cdot e^{it_2})=Q$.}. We get a curve from $P$ to $\gamma$ of length $\leq l_{\beta_2}-\frac{1}{2}l$ (the red curve in Figure \ref{fig: lem-m72}). Therefore $d(P_1,\gamma)=l_{\beta_2}\geq r_\gamma+\frac{1}{2}l$, as desired.
    
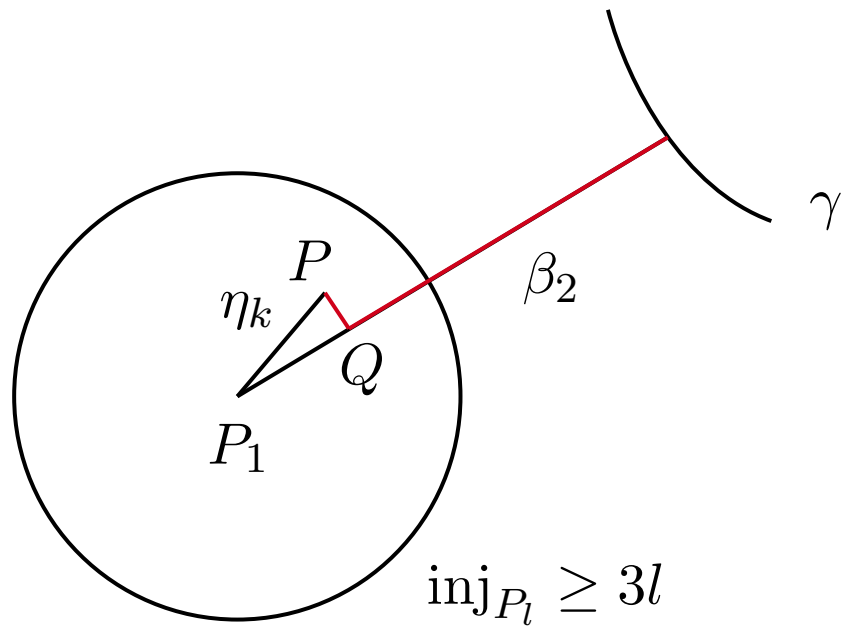
\begin{figure}
    \centering
    \tikzset{every picture/.style={line width=0.75pt}} %set default line width to 0.75pt        

\begin{tikzpicture}[x=0.75pt,y=0.75pt,yscale=-1,xscale=1]

%Shape: Circle [id:dp7221354333978164] 
\draw   (115,135) .. controls (115,104.07) and (140.07,79) .. (171,79) .. controls (201.93,79) and (227,104.07) .. (227,135) .. controls (227,165.93) and (201.93,191) .. (171,191) .. controls (140.07,191) and (115,165.93) .. (115,135) -- cycle ;
%Straight Lines [id:da7755357247926769] 
\draw    (171,135) -- (279,70) ;
%Curve Lines [id:da6991828852431058] 
\draw    (264,38) .. controls (269,57) and (283,83) .. (305,91) ;
%Straight Lines [id:da44985625932004725] 
\draw    (193,109) -- (171,135) ;
%Straight Lines [id:da7480040418424228] 
\draw [color={rgb, 255:red, 208; green, 2; blue, 27 }  ,draw opacity=1 ]   (193,109) -- (199,118) ;
%Straight Lines [id:da050635349318956635] 
\draw [color={rgb, 255:red, 208; green, 2; blue, 27 }  ,draw opacity=1 ]   (199,118) -- (279,70) ;

% Text Node
\draw (217,175.4) node [anchor=north west][inner sep=0.75pt]    {$\mathrm{inj} _{P_l}\geq 3l$};
% Text Node
\draw (313,82.4) node [anchor=north west][inner sep=0.75pt]    {$\gamma$};
% Text Node
\draw (162,140.4) node [anchor=north west][inner sep=0.75pt]    {$P_1$};
% Text Node
\draw (182,94.4) node [anchor=north west][inner sep=0.75pt]    {$P$};
% Text Node
\draw (195,120.4) node [anchor=north west][inner sep=0.75pt]    {$Q$};
% Text Node
\draw (165,107.4) node [anchor=north west][inner sep=0.75pt]    {$\eta_{k}$};
% Text Node
\draw (241,97.4) node [anchor=north west][inner sep=0.75pt]    {$\beta _{2}$};
% Text Node
%\draw (215,211) node [anchor=north west][inner sep=0.75pt]   [align=left] {Figure.};

\end{tikzpicture}
    \caption{A short cut from $P$ to $\gamma$.}
    \label{fig: lem-m72}
\end{figure}

\end{proof}

We next give a criterion for $d(P,\gamma)=r_\gamma$.

\begin{lemma}\label{lem-useful criterion}
    Suppose $r_\gamma\leq1$. Let $P$ be a branch point such that $d(P,\gamma)\leq r_\gamma+r_P-4 l_\gamma$. Then $d(P,\gamma)=r_\gamma$.
\end{lemma}

\begin{proof}
    Recall that we have $\mathrm{inj}_\gamma=r_\gamma$ by Lemma \ref{cor-c1}. Let $\beta$ be a path from $P$ to $\gamma$ of length $d(P,\gamma)$. Let $Q$ be the first intersection point of $\beta$ with $\exp_\gamma(\gamma\times[0,r_\gamma])$ (see Figure \ref{fig: lem-useful}). Then $Q$ must lie on $\gamma_1:=\exp_\gamma(\gamma\times\{r_\gamma\})$, and $d(P,\gamma)=d(P,Q)+r_\gamma$. By assumption $d(P,\gamma)\leq r_\gamma+r_P-4 l_\gamma$. So we have $d(P,Q)\leq r_P-4 l_\gamma$. By the definition of $r_\gamma$, $\gamma_1$ must intersect with $\mathrm{Br}(S)\cup(\partial S\backslash \gamma)$. Let $P_1$ be an intersection point. Then $d(P,P_1)\leq d(P,Q)+d(Q,P_1)\leq r_P-4 l_\gamma+\frac{1}{2}l_{\gamma_1}$. Since $l_{\gamma_1}<8 l_\gamma$ by Lemma \ref{lem-c31}, we get $d(P,P_1)<r_P$. So it must be the case that $P_1=P$. Thus $r_\gamma=d(P,\gamma)$.
\end{proof}

\begin{figure}
    \centering
    \tikzset{every picture/.style={line width=0.75pt}} %set default line width to 0.75pt        

\begin{tikzpicture}[x=0.75pt,y=0.75pt,yscale=-1,xscale=1]

%Curve Lines [id:da41585617303436495] 
\draw    (100,114) .. controls (111,126) and (118,150) .. (115,168) ;
%Curve Lines [id:da02843013044910625] 
\draw    (178,82) .. controls (196,104) and (207,138) .. (209,168) ;
%Curve Lines [id:da4684469064689779] 
\draw    (110,132) .. controls (150,102) and (208,101) .. (263,107) ;
%Shape: Circle [id:dp4020204670747063] 
\draw  [fill={rgb, 255:red, 0; green, 0; blue, 0 }  ,fill opacity=1 ] (204,144) .. controls (204,142.9) and (204.9,142) .. (206,142) .. controls (207.1,142) and (208,142.9) .. (208,144) .. controls (208,145.1) and (207.1,146) .. (206,146) .. controls (204.9,146) and (204,145.1) .. (204,144) -- cycle ;

% Text Node
\draw (260,112.4) node [anchor=north west][inner sep=0.75pt]    {$P$};
% Text Node
\draw (176,108.4) node [anchor=north west][inner sep=0.75pt]    {$Q$};
% Text Node
\draw (219,81.4) node [anchor=north west][inner sep=0.75pt]    {$\beta $};
% Text Node
\draw (111,176.4) node [anchor=north west][inner sep=0.75pt]    {$\gamma$};
% Text Node
\draw (194,177.4) node [anchor=north west][inner sep=0.75pt]    {$\gamma_{1}$};
% Text Node
\draw (210,135.4) node [anchor=north west][inner sep=0.75pt]    {$P_{1}$};
% Text Node
%\draw (168,211) node [anchor=north west][inner sep=0.75pt]   [align=left] {Figure.};

\end{tikzpicture}
    \caption{Lemma \ref{lem-useful criterion}}
    \label{fig: lem-useful}
\end{figure}

By Lemma \ref{lem-useful criterion}, we can give a sufficient condition to guarantee that the condition $d(P,\gamma)=r_\gamma$ remains true after a move as in Lemma \ref{lem-m7}.

\begin{corollary}\label{cor-m1}
    Notation as in Lemma \ref{lem-m7}. Suppose $r_P\geq 10l_\gamma$, then $d(P_1,\gamma)=r_{\gamma}(S_1)$.
\end{corollary}

\begin{proof}
    Combining $l_\gamma\leq\frac{1}{10}r_P$ and $l\leq\frac{1}{5}\mathrm{inj}_P\leq\frac{1}{5}r_P$, we get
    \begin{align}\label{eq-m1}
        3l+4l_\gamma\leq r_P.
    \end{align}
    Since $(S_1,P_1)$ is an $l$-move, we have $r_{\gamma}(S_1)\geq r_\gamma(S)-l,\ d(P,\gamma)\leq d(P_1,\gamma)+l,\ r_{P_1}\geq r_P-l$. Plugging these into (\ref{eq-m1}) we get $d(P_1,\gamma)\leq r_{\gamma}(S_1)+r_{P_1}-4 l_\gamma$. So the conditions of Lemma \ref{lem-useful criterion} are verified. We conclude that $d(P_1,\gamma)=r_{\gamma}(S_1)$.
\end{proof}

We can then apply Lemma \ref{lem-m7} inductively, to get the same result with the assumption $l\leq\frac{1}{5}\mathrm{inj}_P$ relaxed to $l\leq r_P-10l_\gamma$.

\begin{corollary}\label{cor-m2}
    Let $P$ be a branch point such that $d(P,\gamma)=r_\gamma$. Let $l=\min\{1-r_\gamma,r_P-10l_\gamma\}$. Then there is an $l$-move $(S',P')$ such that
    \[r_\gamma(S')\geq r_\gamma+\frac{1}{2}l.\]
\end{corollary}

\begin{proof}
    Let $l_0=\frac{l}{M}$ for some large enough integer $M$ so that $l_0<\min\{l_\gamma,\frac{1}{5}e^{-l}\mathrm{inj}_P\}$. Apply Lemma \ref{lem-m7} $M$ times, each time making an $l_0$-move, while increasing $r_\gamma$ by $\frac{1}{2}l_0$. The conditions of Lemma \ref{lem-m7} are justified by Lemma \ref{lem-m12} and Corollary \ref{cor-m1}. The corollary follows.
\end{proof}

Proposition \ref{prop-m2} follows immediately from Corollary \ref{cor-m2}. Indeed, $R_\gamma\geq r_{\gamma}(S')$ since $S'\in[S]$. Under the assumptions $r_\gamma\leq\frac{1}{2}$ and $r_P\geq 100l_\gamma$, we have
\[r_\gamma+\frac{1}{2}\min\{1-r_\gamma,r_P-10 l_\gamma\}\geq r_\gamma+\frac{1}{3}r_P.\]

\subsection{Moving more branch points away from a boundary curve}

We keep the same assumptions on $\gamma$ as in the previous subsection. We are going to combine the two methods above to give a lower bound of $R_\gamma$ in terms of $d(\gamma,\partial S\backslash \gamma)$. Proposition \ref{prop-m1} enables us to group the branch points near $\gamma$ into one branch point. We then use Proposition \ref{prop-m2} to move this branch point away from $\gamma$.

\begin{lemma}\label{lem-m9}
    Let $m$ be the number of branch points of $S$. Suppose that $10^{m+3}l_\gamma\leq r_\gamma\leq\frac{1}{10}$ and $d(\gamma,\partial S\backslash \gamma)\geq 2r_\gamma$. Then
    \[R_\gamma\geq (1+10^{-m-1})r_\gamma.\]
\end{lemma}

\begin{proof}
    Since $d(\gamma,\partial S\backslash \gamma)>r_\gamma$, there is a branch point $P$ such that $d(P,\gamma)=r_\gamma$. Since $d(\gamma,\partial S\backslash \gamma)\geq 2r_\gamma$, we have $d(P,\partial S\backslash \gamma)\geq r_\gamma$. Therefore $d(P,\partial S)=r_\gamma$. By Proposition \ref{prop-m1}, there exists $L\in[10^{-m-1} r_\gamma,10^{-1}r_\gamma]$ and $S'\in[S],\ P'\in\mathrm{Br}(S')$ such that
    \[r_{P'}\geq9L,\ d(P',\gamma)\leq r_\gamma+L,\ r_{\gamma}(S')\geq r_\gamma-L.\]
    By Lemma \ref{lem-useful criterion} one deduces $d(P',\gamma)=r_{\gamma}(S')$, so we can apply Proposition \ref{prop-m2} to get
    \[R_\gamma\geq r_{\gamma}(S')+\frac{1}{3}r_{P'}\geq (1+10^{-m-1})r_\gamma.\]
\end{proof}

The shortcoming of Lemma \ref{lem-m9} is that if $r_\gamma$ is small, it would not give an effective estimate on $R_\gamma$. To improve this, the trick is to first attach a collar at $\gamma$ so that $r_\gamma$ is not too small, then apply Lemma \ref{lem-m9}. The proof of the following lemma is left to the reader.

\begin{lemma}\label{lem-m8}
    Let $S_1$ be the surface obtained by attaching a collar of length $l$ at $\gamma$. Then $R_{\gamma_1}=R_\gamma+l$.
\end{lemma}

\begin{proposition}\label{prop-m3}
    Let $m$ be the number of branch points of $S$. Let $0<l\leq\frac{1}{10}$. Suppose that $l_\gamma\leq 10^{-m-4}l,\ d(\gamma,\partial S\backslash \gamma)\geq l$, then
    \[R_\gamma\geq10^{-m-2}l.\]
\end{proposition}

\begin{proof}
    If $r_\gamma>\frac{l}{10}$ then $R_\gamma\geq r_\gamma>\frac{l}{10}$ and we are done, so in the following we assume that $r_\gamma\leq\frac{l}{10}$. Let $S_1$ be the surface obtained by attaching a collar of length $\frac{l}{5}$ at $\gamma$. We have
    \begin{align*}
        r_{\gamma_1}=r_\gamma+\frac{l}{5},\ R_{\gamma_1}=R_\gamma+\frac{l}{5},\ d(\gamma_1,\partial S_1\backslash \gamma_1)=d(\gamma,\partial S\backslash \gamma)+\frac{l}{5}.
    \end{align*}
    Then $r_\gamma\leq\frac{l}{10}$ implies $r_{\gamma_1}\leq\frac{l}{2}$. Combining this with the assumption $d(\gamma,\partial S\backslash \gamma)\geq l$, we get $d(\gamma_1,\partial S_1\backslash \gamma_1)\geq 2r_{\gamma_1}$. Moreover, by Lemma \ref{lem-c3} we have $l_{\gamma_1}\leq 2l_\gamma$. So $r_{\gamma_1}\geq\frac{l}{5}\geq 10^{m+3}l_{\gamma_1}$. Therefore $S_1$ satisfies the assumptions of Lemma \ref{lem-m9}. We deduce that $R_{\gamma_1}\geq (1+10^{-m-1})r_{\gamma_1}\geq (1+10^{-m-1})\frac{l}{5}$, which implies $R_\gamma=R_{\gamma_1}-\frac{l}{5}\geq 10^{-m-1}\frac{l}{5}\geq10^{-m-2}l$.
\end{proof}

\section{Proof of Proposition \ref{prop2}}

Let $f$ be an $\epsilon$-excellent developing map on a closed surface $\Sigma$. Let $\gamma$ be a folding curve. Define
\[R_\gamma(f)=\min\{R_\gamma(S_+),R_\gamma(S_-)\},\]
where $S_+,S_-$ are the two components of $\Sigma\backslash\mathrm{Fold}(f)$ adjacent to $\gamma$, equipped with the metric $g_f$ (recall the definition of $R_\gamma(S_\pm)$ from the beginning of Section \ref{sec-m-app}). In parallel with Lemma \ref{lem-m8} we have the following lemma, whose proof is left to the reader.

\begin{lemma}\label{lem-mc1}
    Let $f_1$ be the developing map obtained by removing a pair of collars of length $l$ at $\gamma$ (which is $2\epsilon$-excellent by Lemma \ref{lem-c7}), then $R_{\gamma}(f_1)=R_\gamma(f)-l$.
\end{lemma}

\begin{lemma}\label{lem-mc2}
    Let $f$ be an $\epsilon$-excellent developing map with $n$ folding curves $(n\geq1)$. Then for any $\delta>0$, there exists a $2\epsilon$-excellent developing map $f_1$ with $n$ folding curves and the same monodromy as $f$, such that one of its folding curves $\gamma$ satisfies $R_{\gamma}(f_1)\leq \delta$.
\end{lemma}

\begin{proof}
     Let $\gamma$ be any folding curve of $f$. Let $\Sigma_\pm$ be the two components of $\mathrm{Fold}(f)$ adjacent to $\gamma$. By a homotopy of $f|_{\Sigma_\pm}$, we may assume that $r_\gamma>R_\gamma-\delta$. By Lemma \ref{lem-c7}, $\gamma$ has a pair of $2\epsilon$-thin collar neighborhoods of length $R_\gamma-\delta$. Let $f_1$ be the developing map obtained by removing this pair of collar neighborhoods. $f_1$ satisfies $R_{\gamma}(f_1)=\delta$ by Lemma \ref{lem-mc1}.
\end{proof}

We now prove Proposition \ref{prop2}. Let $f_0$ be an $\epsilon$-excellent developing map with $n$ folding curves ($\epsilon\leq\epsilon_0,n\geq1$). By Lemma \ref{lem-mc2}, there is a $2\epsilon$-excellent developing map $f$ with $n$ folding curves and the same monodromy as $f_0$, such that one of its folding curves, $\gamma_1$, satisfies $R_{\gamma_1}\leq 10^{\chi(\Sigma)-2}\delta_0$. By the definition of $R_\gamma$, there is a component $S_0$ of $\Sigma\backslash\mathrm{Fold}(f)$ adjacent to $\gamma_1$ such that $R_{\gamma_1}(S_0)\leq10^{\chi(\Sigma)-2}\delta_0$. By Proposition \ref{prop-m3}, we have $d_{S_0}(\gamma_1,\partial S_0\backslash \gamma_1)\leq \delta_0$, which means that there is an other component $\gamma_2$ of $\partial S_0$ such that $d_{S_0}(\gamma_1,\gamma_2)\leq \delta_0$. This finishes the proof of Proposition \ref{prop2}.

\section{Goldman's theorem on a one-holed torus}

In this section we recall a theorem of Goldman regarding representations of the fundamental groups of a one-holed torus. We then give some of its corollaries.

The following theorem is due to Goldman \cite{Goldman03} (Theorem 5.2.1), and in this form appears as Theorem 2.4 in \cite{MW16}.

\begin{theorem}\label{thm-goldman}
    Let $\Sigma$ be a one-holed torus, and let $\rho:\pi_1(\Sigma)\rightarrow\mathrm{PSL}_2(\R)$ be a representation of Euler class zero such that the boundary curve has trace in $[2,18]$. Then $\rho$ sends some simple closed curve to a non-hyperbolic element.
\end{theorem}

\begin{remark}
    The assumption that $\rho$ is of Euler class zero in the above theorem is superfluous. Indeed, the fact that the boundary curve has trace in $[2,+\infty)$ implies that $\rho$ is of Euler class zero.
\end{remark}

\begin{corollary}\label{cor-g0}
    Let $\Sigma$ be a one-holed torus and let $\alpha,\beta\subset\Sigma$ be simple closed curves with $i(\alpha,\beta)=1$. Let $\gamma$ be the boundary curve. Let $\rho:\pi_1(\Sigma)\rightarrow\mathrm{PSL}_2(\R)$ be a representation. Suppose that $|\mathrm{Tr}(\rho(\gamma))|\leq 18,\ |\mathrm{Tr}(\rho(\alpha))|,|\mathrm{Tr}(\rho(\beta))|\leq2\sqrt{2}$, then $\rho$ sends some simple closed curve to a non-hyperbolic element.
\end{corollary}

\begin{proof}
    A simple computation shows that if $g_1,g_2\in\mathrm{SL}_2(\R)$ satisfies $2<|\mathrm{Tr}(g_1)|,|\mathrm{Tr}(g_2)|\leq2\sqrt{2}$ then $\mathrm{Tr}([g_1,g_2])\geq-2$. Now apply Theorem \ref{thm-goldman}.
\end{proof}

\begin{corollary}\label{cor-g1}
    Let $f$ be a developing map on an arbitrary surface $\Sigma$ for some representation $\rho:\pi_1(\Sigma)\rightarrow\mathrm{PSL}_2(\R)$. Suppose that there are closed curves $\alpha,\beta$ on $\Sigma$ with $i(\alpha,\beta)=1$ such that $f(\alpha)$ and $f(\beta)$ both have length $\leq1$, then $\rho$ sends some simple closed curve to a non-hyperbolic element.
\end{corollary}

\begin{proof}
    We may assume that $\alpha,\beta$ are simple and intersect transversally at one point. A tubular neighborhood of $\alpha\cup\beta$ is a one-holed torus, which we denote by $\Sigma'$. Let $\gamma$ be the boundary curve of $\Sigma'$. The assumption implies that $|\mathrm{Tr}(\rho(\alpha))|,|\mathrm{Tr}(\rho(\beta))|\leq e^{\frac{1}{2}}+e^{-\frac{1}{2}}<2\sqrt{2}$, and that $f(\gamma)$ has length $\leq4$, so that $|\mathrm{Tr}(\rho(\gamma))|\leq e^2+e^{-2}<18$. Now apply Corollary \ref{cor-g0}.
\end{proof}

\begin{corollary}\label{cor-g2}
    Let $f$ be an excellent developing map, $\gamma_1,\gamma_2$ folding curves. Let $S_1,S_2$ be different components of $\Sigma\backslash\mathrm{Fold}(f)$, equipped with the metric $g_f$. Suppose that $d_{S_1}(\gamma_1,\gamma_2)\leq\frac{1}{2}$, then $d_{S_2}(\gamma_1,\gamma_2)>\frac{1}{2}$. (We set $d_{S_i}(\gamma_1,\gamma_2)=\infty$ if either $\gamma_1$ or $\gamma_2$ is not a component of $\partial S_i$.)
\end{corollary}

\begin{corollary}\label{cor-g3}
    Suppose that the holonomy of $S=(\Sigma,g)$ around every simple closed curve is hyperbolic. Let $\alpha\subset S$ be a simple closed curve of length $\leq1$. Then $\alpha$ disconnects $N_{\frac{1}{2}}(\alpha)$.
\end{corollary}

\section{Reduction to very standard developing maps}

The rest of the paper is devoted to the proof of Proposition \ref{prop3}. Recall that we are given two folding curves $\gamma_1,\gamma_2$ that are close to each other. Our goal is to cancel them or to merge them into a single folding curve. It is clear that we may assume that $\gamma_1,\gamma_2$ are boundary curves of a common component $\Sigma_0$ of $\Sigma\backslash\mathrm{Fold}(f)$ (which is in fact the case from the proof of Proposition \ref{prop2}). In this section we prove a dichotomy, showing that either $\Sigma_0$ is an annulus or we can modify the developing map so that $\gamma_1,\gamma_2$ are separated by a single branch point. The latter case is what we call standard. We then prove that if $\Sigma_0$ is an annulus then we can cancel $\gamma_1$ and $\gamma_2$. So it remains to consider standard developing maps. In particular we are reduced to Proposition \ref{prop-vs2}, which is proved in later sections.

\begin{definition}
    A surface $S=(\Sigma,g)$ is called $\epsilon$-good if the holonomy of $g$ around every simple closed curve is hyperbolic, and every boundary curve is an equidistant curve of length $\leq\epsilon$.
\end{definition}

\subsection{A criterion of collar}

\begin{lemma}\label{lem-d1}
    Let $S$ be an $\epsilon$-good surface. Let $\gamma_1,\gamma_2$ be two boundary curves of $S$. If $d(\gamma_1,\gamma_2)+10\epsilon\leq r_{\gamma_1}+r_{\gamma_2}\leq 1$, then $S$ is a collar.
\end{lemma}

\begin{proof}
    We have $\mathrm{inj}_{\gamma_i}=r_{\gamma_i}$ by Lemma \ref{cor-c1}. Let $l_i=r_{\gamma_i}-\epsilon$ and let $C_i$ be collar neighborhoods of $\gamma_i$ of length $l_i$. Let $\gamma_1'$ be the boundary curve of $C_1$ other than $\gamma_1$. Then $d(\gamma_1',\gamma_2)=d(\gamma_1,\gamma_2)-l_1\leq l_2-8\epsilon$. Since $l_{\gamma_1'}\leq 8\epsilon$ by Lemma \ref{lem-c31}, we get $\gamma_1'\subset C_2$. So $\gamma_1'$ is either homotopic to $\gamma_2$ or null homotopic. But $\gamma_1'$ is homotopic to $\gamma_1$. It follows that $\gamma_1$ is homotopic to $\gamma_2$. Therefore $S$ is an annulus, and thus a collar.
\end{proof}

\subsection{Standard triples}

\begin{definition}
    Let $S$ be an $\epsilon$-good suface. A standard triple $(P,\gamma_1,\gamma_2)$ of $S$ consists of a branch point $P$ and boundary curves $\gamma_1,\gamma_2$ such that $r_P\geq \delta_0,\ d(P,\gamma_i)=r_{\gamma_i}\leq \delta_1\ (i=1,2)$. We call $S$ standard if it has a standard triple.
\end{definition}

Proposition \ref{prop-m1} can be used to produce a standard triple.

\begin{lemma}\label{lem-d2}
    Let $S$ be an $\epsilon$-good surface $(\epsilon\leq100\epsilon_0)$. Suppose that there is a branch point $P$ and boundary curves $\gamma_1,\gamma_2$ such that $d(P,\partial S)\geq\frac{1}{10}\delta_1$, $d(P,\gamma_i)-2\delta_0\leq r_{\gamma_i}\leq\frac{1}{2}\delta_1$. Then there is $S'\in[S]$ which is standard.
\end{lemma}

\begin{proof}
    The assumption $d(P,\partial\Sigma)\geq\frac{1}{10}\delta_1$ enables us to apply Proposition \ref{prop-m1} with $l=\delta_0$ (recall that $\delta_0\ll\delta_1$). We get $S'\in[S],\ P'\in\mathrm{Br}(S')$ such that $r_{P'}\geq9L,\ d(P',\gamma_i)\leq d(P,\gamma_i)+L,\ r_{\gamma_i}(S')\geq r_{\gamma_i}-L$ for some $L\in[\delta_0,10^{-\chi(\Sigma)}\delta_0]$. In particular $r_{P'}\geq9L$ implies $r_{P'}\geq\delta_0$. By Lemma \ref{lem-useful criterion} and the assumption $d(P,\gamma_i)-2\delta_0\leq r_{\gamma_i}$, we get $d(P',\gamma_i)=r_{\gamma_i}(S')$. We deduce $r_{\gamma_i}(S')\leq\delta_1$ from the assumption $r_{\gamma_i}\leq\frac{1}{2}\delta_1$.
\end{proof}

\subsection{A dichotomy}\label{subsec-d}

\begin{definition}
    An $\epsilon$-excellent developing map $f$ is called standard if there is a component $S_0$ of $\Sigma\backslash\mathrm{Fold}(f)$ which is standard.
\end{definition}

In this subsection we prove the dichotomy, which says that if $\gamma_1,\gamma_2$ are boundaries of $S_0$ that are close to each other, then either $S_0$ is a collar or it is standard after a homotopy of $f|_{\Sigma_0}$.

The idea of proof is as follows. The condition that $S_0$ is not a collar essentially says that there is a branch point $P$ in between $\gamma_1$ and $\gamma_2$. We then use Proposition \ref{prop-m1} to merge $P$ with nearby branch points, producing a standard triple. In order to apply Proposition \ref{prop-m1}, we need to make room for the move. To do this, we attach collars at every boundary curve of $S_0$.

Thus we begin with two simple lemmas concerning attaching collars. The first one is essentially the triangle inequality, while the second one is an application of Lemma \ref{lem-d1}.

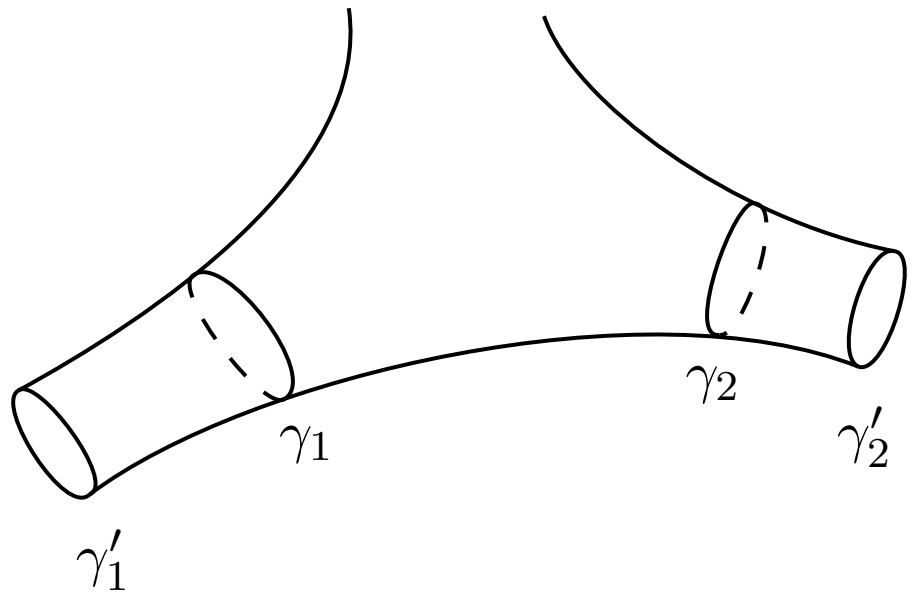
\begin{figure}
    \centering

\tikzset{every picture/.style={line width=0.75pt}} %set default line width to 0.75pt        

\begin{tikzpicture}[x=0.75pt,y=0.75pt,yscale=-1,xscale=1]
%uncomment if require: \path (0,339); %set diagram left start at 0, and has height of 339

%Curve Lines [id:da4153662845073185] 
\draw    (110,164) .. controls (143,146) and (199,109) .. (193,68) ;
%Curve Lines [id:da9387073277034541] 
\draw    (128,190) .. controls (168,160) and (268,137) .. (321,158) ;
%Curve Lines [id:da11099955952179874] 
\draw    (242,70) .. controls (249,91) and (289,121) .. (330,129) ;
%Curve Lines [id:da7068768402568921] 
\draw    (110,164) .. controls (104,168) and (122,196) .. (128,190) ;
%Curve Lines [id:da004744159818034777] 
\draw    (110,164) .. controls (115,160) and (135,187) .. (128,190) ;
%Curve Lines [id:da6961215582074012] 
\draw    (154,135) .. controls (162,129) and (186,160) .. (177,166) ;
%Curve Lines [id:da038452699641740584] 
\draw  [dash pattern={on 4.5pt off 4.5pt}]  (154,135) .. controls (148,139) and (172,170) .. (177,166) ;
%Curve Lines [id:da5048546333753102] 
\draw    (330,129) .. controls (323,127) and (314,155) .. (321,158) ;
%Curve Lines [id:da12234488512789721] 
\draw    (330,129) .. controls (337,130) and (328,160) .. (321,158) ;
%Curve Lines [id:da9418100775444542] 
\draw    (295,117) .. controls (290,115) and (277,151) .. (286,150) ;
%Curve Lines [id:da5019460776240675] 
\draw  [dash pattern={on 4.5pt off 4.5pt}]  (295,117) .. controls (303,118) and (292,151) .. (286,150) ;

% Text Node
\draw (174,171.4) node [anchor=north west][inner sep=0.75pt]    {$\gamma _{1}$};
% Text Node
\draw (276,156.4) node [anchor=north west][inner sep=0.75pt]    {$\gamma _{2}$};
% Text Node
\draw (123,197.4) node [anchor=north west][inner sep=0.75pt]    {$\gamma '_{1}$};
% Text Node
\draw (314,166.4) node [anchor=north west][inner sep=0.75pt]    {$\gamma '_{2}$};

\end{tikzpicture}

    \caption{Attach collars.}
    \label{fig: lem-d3}
\end{figure}

\begin{lemma}\label{lem-d3}
    Let $\gamma_1,\gamma_2$ be boundary curves of an $\epsilon$-good surface $S$. Suppose that $d(\gamma_1,\gamma_2)\leq \delta$. Let $S'$ be the surface obtained by attaching a collar of length $l\leq\frac{1}{10}$ at every boundary curve. Then for every branch point $P$ we have
    \[d'(P,\partial S')\geq l,\quad d'(P,\gamma_2')\leq d'(P,\gamma_1')+\delta,\]
    where $\gamma_i'$ is the boundary curve of $S'$ homotopic to $\gamma_i$, and we have written $d'$ for $d_{S'}$ to distinguish it with the distance function $d$ on $S$. (Figure \ref{fig: lem-d3})
\end{lemma}

\begin{proof}
    The first inequality is clear. To prove the second one, note that $d'(\gamma_1,\gamma_2)\leq d(\gamma_1,\gamma_2)\leq \delta$, $d'(P,\gamma_i')=d'(P,\gamma_i)+l$, while by the triangle inequality, $d'(P,\gamma_2)\leq d'(P,\gamma_1)+d'(\gamma_1,\gamma_2)$. We get $d'(P,\gamma_2')=d'(P,\gamma_2)+l\leq d'(P,\gamma_1)+d'(\gamma_1,\gamma_2)+l\leq d'(P,\gamma_1)+\delta+l=d'(P,\gamma_1')+\delta$.
\end{proof}

\begin{lemma}\label{lem-d4}
    Notation as in the previous lemma. Suppose that $S$ is not a collar and that $l\geq \delta\geq100\epsilon$. Then there is a branch point $P$ such that
    \[d'(P,\gamma_i')-2\delta\leq r_{\gamma_i'}\leq \delta+2l\quad(i=1,2).\]
\end{lemma}

\begin{proof}
    We have $d'(\gamma_1',\gamma_2')\leq d(\gamma_1,\gamma_2)+2l\leq \delta+2l$. Since $r_{\gamma_i'}\leq d'(\gamma_1',\gamma_2')$, we get $r_{\gamma_i'}\leq \delta+2l\ (i=1,2)$. This proves the second inequality. In particular $r_{\gamma_i'}\leq\frac{1}{2}$, so that $r_{\gamma_1'}+r_{\gamma_2'}\leq1$. By Lemma \ref{lem-d1}, if $S$ is not a collar, then $d'(\gamma_1',\gamma_2')+20\epsilon>r_{\gamma_1'}+r_{\gamma_2'}$. So $\min\{r_{\gamma_1'},r_{\gamma_2'}\}<\frac{1}{2}d'(\gamma_1',\gamma_2')+10\epsilon\leq\frac{1}{2}\delta+l+10\epsilon\leq \delta+l$. Without lost of generality we may assume that $r_{\gamma_1'}<\delta+l$. Since $d'(\gamma_1',\gamma)\geq 2l\geq \delta+l$ for any boundary curve $\gamma$ different from $\gamma_1'$, there is a branch point $P$ such that $d'(P,\gamma_1')=r_{\gamma_1'}$. By the previous lemma we have $d'(P,\gamma_2')\leq d'(P,\gamma_1')+\delta$. We get $d'(P,\gamma_2')\leq l+2\delta\leq r_{\gamma_2'}+2\delta$, proving the first inequality.
\end{proof}

Now we can prove the dichotomy.

\begin{proposition}\label{prop-d}
    Let $f$ be an $\epsilon$-excellent developing map $(\epsilon\leq2\epsilon_0)$. Suppose that a component $S_0=(\Sigma_0,g_f)$ of $\Sigma\backslash\mathrm{Fold}(f)$ has two boundary curves $\gamma_1,\gamma_2$ satisfying $d_{S_0}(\gamma_1,\gamma_2)\leq\delta_0$. Then either $S_0$ is a collar, or there is a standard $2\epsilon$-excellent developing map with the same monodromy and the same number of folding curves as $f$.
\end{proposition}

\begin{proof}
    Suppose that $S_0$ is not a collar. Let $f_1$ be the developing map obtained by attaching a pair of collars of length $\frac{1}{10}\delta_1$ along every boundary curve of $S_0$, which is $2\epsilon$-excellent by Lemma \ref{lem-c6} (note that $S_0$ has at most $-\chi(S_0)+2$ boundary curves). Let $S_1=(\Sigma_0,g_{f_1})$. By Lemma \ref{lem-d3} and Lemma \ref{lem-d4}, there is a branch point $P$ of $S_1$ such that $d_{S_1}(P,\partial S_1)\geq \frac{1}{10}\delta_1,\ d_{S_1}(P,\gamma_i)-2\delta_0\leq r_{\gamma_i}(S_1)\leq\frac{1}{2}\delta_1\ (i=1,2)$. Therefore $S_1$ satisfies the assumptions of Lemma \ref{lem-d2}. So $f_1$ becomes standard after a homotopy of $f_1|_{\Sigma_0}$.
\end{proof}

\subsection{The case of collar}\label{subsec-the case of annulus}

\begin{proposition}\label{prop-annulus}
    Let $\epsilon\leq2\epsilon_0$. Let $f$ be an $\epsilon$-excellent developing map with $n$ folding curves such that one component of $\Sigma\backslash\mathrm{Fold}(f)$ is a collar of length $\leq\delta_0$. Then there exists a $2\epsilon$-excellent developing map with $n-2$ folding curves and the same monodromy as $f$.
\end{proposition}

\begin{figure}
    \centering
    \tikzset{every picture/.style={line width=0.75pt}} %set default line width to 0.75pt        

\begin{tikzpicture}[x=0.75pt,y=0.75pt,yscale=-1,xscale=1]

%Curve Lines [id:da9112661562474601] 
\draw    (74,97) .. controls (97,105) and (120,107) .. (138,97) ;
%Curve Lines [id:da9350144473156888] 
\draw    (74,137) .. controls (94,127) and (120,127) .. (138,136) ;
%Curve Lines [id:da5054781057582207] 
\draw    (138,97) .. controls (131,106) and (131,126) .. (138,136) ;
%Curve Lines [id:da6217962509777276] 
\draw  [dash pattern={on 4.5pt off 4.5pt}]  (138,97) .. controls (145,108) and (143,124) .. (138,136) ;
%Curve Lines [id:da375118944552153] 
\draw    (138,97) .. controls (155,104) and (177,108) .. (193,98) ;
%Curve Lines [id:da41306226555782] 
\draw    (138,136) .. controls (152,128) and (183,129) .. (193,136) ;
%Curve Lines [id:da4462676057646341] 
\draw    (193,98) .. controls (209,105) and (250,114) .. (276,83) ;
%Curve Lines [id:da8730665986492127] 
\draw    (193,98) .. controls (187,109) and (187,124) .. (193,136) ;
%Curve Lines [id:da206006176537474] 
\draw  [dash pattern={on 4.5pt off 4.5pt}]  (193,98) .. controls (198,110) and (199,125) .. (193,136) ;
%Curve Lines [id:da669729971746779] 
\draw    (193,136) .. controls (213,127) and (246,126) .. (275,155) ;
%Curve Lines [id:da8582036554073554] 
\draw    (249,101) .. controls (242,112) and (242,125) .. (248,136) ;
%Curve Lines [id:da4218202553814909] 
\draw  [dash pattern={on 4.5pt off 4.5pt}]  (249,101) .. controls (255,109) and (255,126) .. (248,136) ;
%Straight Lines [id:da7543835855585607] 
\draw    (306,125.5) -- (322,125.5)(306,128.5) -- (322,128.5) ;
\draw [shift={(330,127)}, rotate = 180] [color={rgb, 255:red, 0; green, 0; blue, 0 }  ][line width=0.75]    (10.93,-3.29) .. controls (6.95,-1.4) and (3.31,-0.3) .. (0,0) .. controls (3.31,0.3) and (6.95,1.4) .. (10.93,3.29)   ;
%Curve Lines [id:da8109572261553177] 
\draw    (363,104) .. controls (391,115) and (431,119) .. (467,86) ;
%Curve Lines [id:da5717166759222067] 
\draw    (366,146) .. controls (388,128) and (440,132) .. (468,169) ;
%Curve Lines [id:da6625738695679272] 
\draw    (440,105) .. controls (433,117) and (434,132) .. (439,145) ;
%Curve Lines [id:da767290075680199] 
\draw  [dash pattern={on 4.5pt off 4.5pt}]  (440,105) .. controls (448,112) and (447,137) .. (439,145) ;

% Text Node
\draw (129,139.4) node [anchor=north west][inner sep=0.75pt]    {$\gamma_2$};
% Text Node
\draw (185,137.4) node [anchor=north west][inner sep=0.75pt]    {$\gamma_1$};
% Text Node
\draw (234,141.4) node [anchor=north west][inner sep=0.75pt]    {$\gamma_2'$};
% Text Node
\draw (158,109.4) node [anchor=north west][inner sep=0.75pt]    {$C$};
% Text Node
\draw (212,109.4) node [anchor=north west][inner sep=0.75pt]    {$C_1$};
% Text Node
\draw (270,93.4) node [anchor=north west][inner sep=0.75pt]    {$\Sigma _{1}$};
% Text Node
\draw (452,112.4) node [anchor=north west][inner sep=0.75pt]    {$\Sigma _{1} \backslash C _{1}$};
% Text Node
\draw (94,107.4) node [anchor=north west][inner sep=0.75pt]    {$\Sigma _{2}$};
% Text Node
\draw (386,114.4) node [anchor=north west][inner sep=0.75pt]    {$\Sigma _{2}$};
% Text Node
%\draw (196,188) node [anchor=north west][inner sep=0.75pt]   [align=left] {Figure. Cancel $\displaystyle \Pi $ and $\displaystyle \Pi _{1}$.};

\end{tikzpicture}
    \caption{Cancel $C$ with $C_1$.}
    \label{fig: case of annulus}
\end{figure}
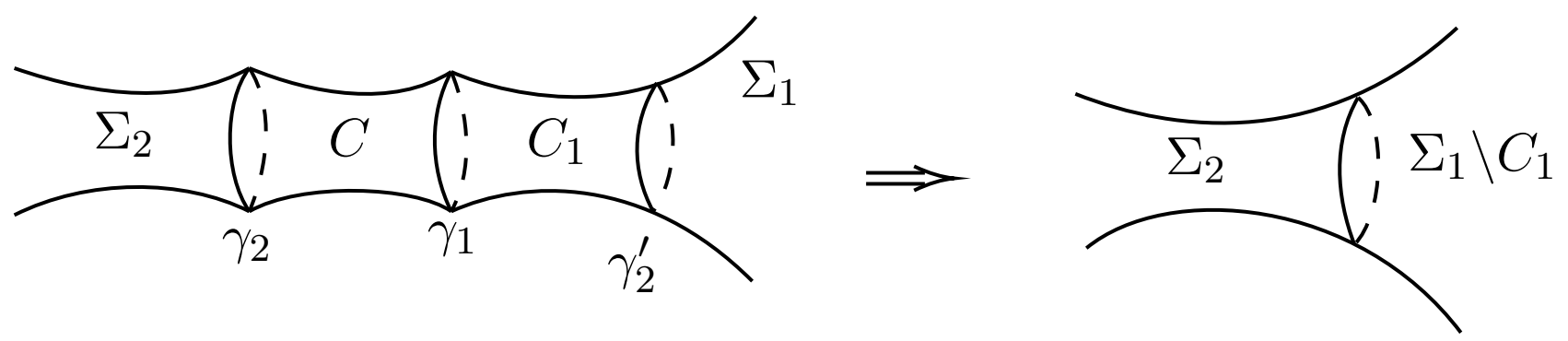

\begin{proof}
    Let $C$ be a component of $\Sigma\backslash\mathrm{Fold}(f)$ which is a collar of length $l\leq \delta_0$. Let $\gamma_1,\gamma_2$ be the boundary curves curves of $C$. Let $S_1=(\Sigma_1,g_f)$ be the component of $\Sigma\backslash\mathrm{Fold}(f)$ adjacent to $\gamma_1$ other than $C$. Replace $f$ by the $2\epsilon$-excellent developing map obtained by attaching a pair of collars of length $10^{\chi(\Sigma)}$ along every boundary curve of $\Sigma_1$ different from $\gamma_1,\gamma_2$ (Lemma \ref{lem-c6}). Then $d_{S_1}(\gamma_1,\partial S_1\backslash (\gamma_1\cup \gamma_2))\geq 10^{\chi(\Sigma)}$. Combining this with Corollary \ref{cor-g2} we get $d_{S_1}(\gamma_1,\partial S_1\backslash \gamma_1)\geq 10^{\chi(\Sigma)}$. Then Proposition \ref{prop-m3} implies $R_{\gamma_1}(S_1)\geq 10^{2\chi(\Sigma)-2}>\delta_0$. Therefore, by a homotopy of $f|_{\Sigma_1}$, we may assume that $r_{\gamma_1}(S_1)>\delta_0$. By assumption $C$ has length $l\leq\delta$. So there is a collar neighborhood $C_1\subset S_1$ of $\gamma_1$ of length $l$. We can `cancel' $C$ with $C_1$ (Figure \ref{fig: case of annulus}). Indeed, $f$ maps $C$ and $C_1$ to the same collar in $\H/\rho(\gamma_1)$. Let $\gamma_2'$ be the boundary curve of $C_1$ other than $\gamma_1$. Then $f$ maps $\gamma_2$ and $\gamma_2'$ homeomorphically to the same equidistant curve in $\H/\rho(\gamma_1)$, giving rise to an identification $\gamma_2\cong\gamma_2'$. We remove $C\cup C_1$ from $\Sigma$ and identify $\gamma_2$ with $\gamma_2'$ via the above identification. The resulting surface is naturally homeomorphic to $\Sigma$, and $f$ restricts to an $2\epsilon$-excellent developing map with $n-2$ folding cures (the developing map is a local homeomorphism along the curve obtained by gluing $\gamma_2$ and $\gamma_2'$).
\end{proof}

\subsection{Very standard developing maps}

We need a slightly stronger version of standard developing maps.

\begin{definition}
    An excellent developing map $f$ is called very standard if there is a component $S_0$ of $\Sigma\backslash\mathrm{Fold}(f)$ with a standard triple $(P,\gamma_1,\gamma_2)$ such that $r_{\gamma_i}(S_i)\geq 10\delta_1\ (i=1,2)$, where $S_i$ is the component of $\Sigma\backslash\mathrm{Fold}(f)$ other than $S_0$ which is adjacent to $\gamma_i$.
\end{definition}

The point is that we need long enough collar neighborhoods in $S_1,S_2$, so that later on when we perform surgery, there is enough room to cancel with collars in $S_0$. In this subsection we show that standard developing maps can always be deformed to very standard ones.

\begin{lemma}\label{lem-vs1}
    Let $\gamma_1,\gamma_2$ be boundary curves of an $\epsilon$-good surface $S$. Suppose that $r_{\gamma_1}>10\delta_1$, and $d(\gamma_2,\partial S\backslash \gamma_2)\geq 10^{-\chi(S)+4}\delta_1$. Then there exists $S'\in[S]$ such that $r_{\gamma_i}(S')\geq 10\delta_1\ (i=1,2)$.
\end{lemma}

\begin{proof}
    The idea is to apply Proposition \ref{prop-m3} to $\gamma_2$ after `freezing' a collar neighborhood of $\gamma_1$. Since $r_{\gamma_1}>10\delta_1$, there is a collar neighborhood $C$ of $\gamma_1$ of length $10\delta_1$. Let $S_1=S\backslash C$. Then $d(\gamma_2,\partial S_1\backslash \gamma_2)\geq d(\gamma_2,\partial S\backslash \gamma_2)-10\delta_1>10^{-\chi(S)+3}\delta_1$. By Proposition \ref{prop-m3}, we get $R_{\gamma_2}(S_1)>10\delta_1$. Therefore there is $S'_1\in[S_1]$ such that $r_{\gamma_2}(S_1')\geq 10\delta_1$. Let $S'$ be the surface obtained by attaching $C$ to $S_1'$ along $\gamma_1$. Then $S'\in[S]$, and $r_{\gamma_1}(S')\geq 10\delta_1$.
\end{proof}

\begin{proposition}\label{prop-vs1}
    Let $f$ be a standard $\epsilon$-excellent developing map with $n$ folding curves $(\epsilon\leq100\epsilon_0)$. Then there exists a very standard $4\epsilon$-excellent developing map with $n$ folding curves and the same monodromy as $f$.
\end{proposition}

\begin{proof}
    Let $\Sigma_0$ be a component of $\Sigma\backslash\mathrm{Fold}(f)$ with a standard triple $(P,\gamma_1,\gamma_2)$. Let $\Sigma_i$ be the component of $\Sigma\backslash\mathrm{Fold}(f)$ other than $\Sigma_0$ which is adjacent to $\gamma_i\ (i=1,2)$.
    
    \paragraph{Case 1} $\Sigma_1\neq\Sigma_2$.

    Let $f_1$ be the $2\epsilon$-excellent developing map obtained by attaching a collar of length $10^{-\chi(\Sigma)+4}\delta_1$ along every boundary curve of $\Sigma_1,\Sigma_2$ other than $\gamma_1,\gamma_2$. Let $S_i=(\Sigma_i,g_{f_1})$. Then $d_{S_i}(\gamma_i,\partial S_i\backslash \gamma_i)\geq 10^{-\chi(\Sigma)+4}\delta_1\ (i=1,2)$. By Proposition \ref{prop-m3} we get $R_{\gamma_i}(S_i)>10\delta_1$. Therefore, after a homotopy of $f_1|_{\Sigma_1}$ and $f_1|_{\Sigma_2}$, we have $r_{\gamma_i}(S_i)\geq 10\delta_1$, giving the desired very standard developing map.

    \paragraph{Case 2} $\Sigma_1=\Sigma_2$.

    Let $f_1$ be the $2\epsilon$-excellent developing map obtained by attaching a collar of length $10^{-\chi(\Sigma)+4}\delta_1$ at every boundary curve of $\Sigma_1$ other than $\gamma_1,\gamma_2$. Let $S_1=(\Sigma_1,g_{f_1})$. Then $d_{S_1}(\gamma_1,\partial S_1\backslash (\gamma_1\cup \gamma_2))\geq 10^{-\chi(\Sigma)+4}\delta_1$. Combining this with Corollary \ref{cor-g2} we get $d_{S_1}(\gamma_1,\partial S_1\backslash \gamma_1)\geq 10^{-\chi(\Sigma)+4}\delta_1$. Therefore Proposition \ref{prop-m3} implies $R_{\gamma_1}(S_1)>10\delta_1$. By a homotopy of $f_1|_{\Sigma_1}$, we may assume that $r_{\gamma_1}(S_1)>10\delta_1$. Now let $f_2$ be the $4\epsilon$-excellent developing map obtained from $f_1$ by attaching a collar of length $10^{-\chi(\Sigma)+4}\delta_1$ along every boundary curve of $\Sigma_1$ other than $\gamma_1,\gamma_2$. Let $S_2=(\Sigma_1,g_{f_2})$. By the same reasoning as above we have $d_{S_2}(\gamma_2,\partial S_2\backslash \gamma_2)\geq 10^{-\chi(\Sigma)+4}\delta_1$. We see that $S_2$ satisfies the assumptions of Lemma \ref{lem-vs1}. Therefore after a homotopy of $f_2|_{\Sigma_1}$, we have $r_{\gamma_i}(S_2)\geq 10\delta_1\ (i=1,2)$. This gives the desired very standard developing map.
\end{proof}

The rest of the paper is devoted to the proof of the following proposition:

\begin{proposition}\label{prop-vs2}
    Let $f$ be a very standard $\epsilon$-excellent developing map with $n$ folding curves $(\epsilon\leq100\epsilon_0)$. Then there exists a $5\epsilon$-excellent developing map with $n-1$ folding curves and the same monodromy as $f$.
\end{proposition}

Assuming Proposition \ref{prop-vs2}, we finish the proof of Proposition \ref{prop3}.

\paragraph{Proof of Proposition \ref{prop3}}

Let $f$ be an $\epsilon$-excellent developing map with $n$ folding curves $(\epsilon\leq2\epsilon_0)$, such that two folding curves $\gamma_1,\gamma_2$ satisfy $d_f(\gamma_1,\gamma_2)\leq\delta_0$. We may assume that $\gamma_1,\gamma_2$ are boundary curves of a component $\Sigma_0$ of $\Sigma\backslash\mathrm{Fold}(f)$. By Proposition \ref{prop-d}, either $\Sigma_0$ is a collar, or there is a standard $2\epsilon$-excellent developing map with $n$ folding curves and the same monodromy as $f$. In the first case, by Proposition \ref{prop-annulus}, there exists a $2\epsilon$-excellent developing map with $n-2$ folding curves and the same monodromy as $f$. In the second case, by Proposition \ref{prop-vs1}, there exists a very standard $8\epsilon$-excellent developing map with $n$ folding curves and the same monodromy as $f$. Then by Proposition \ref{prop-vs2}, there exists a $40\epsilon$-excellent developing map with $n-1$ folding curves and the same monodromy as $f$, as desired.

\section{Flower neighborhoods}

We pause a little and see what geometric structure we can get from a standard developing map. In this section we introduce flower neighborhoods of a boundary curve. These are variants of collar neighborhoods, and one can think of them as collar neighborhoods with a branch point on the boundary. They will be used to perform surgery in later sections.

Throughout this section let $S=(\Sigma,g)$ be a branched hyperbolic surface whose holonomy around every simple closed curve is non-elliptic. Denote the distance function by $d$. Suppose that a boundary curve $\gamma$ is an equidistant curve. Let $N\gamma$ denote the normal bundle of $\gamma$ in $\Sigma$. We trivialize $N\gamma$ by the inward pointing unit normal vector field, giving rise to an isomorphism $N\gamma\cong \gamma\times\R$. The exponential map
\[\exp_\gamma:\gamma\times[0,+\infty)\supset\mathcal{U}_\gamma\rightarrow S\]
is a local diffeomorphism. If $l\leq d(\gamma,S\backslash\gamma)$, then $\mathcal{U}_\gamma\cap(\gamma\times[0,l])$ is of the form $\bigcup_{x\in\gamma}(\{x\}\times[0,u(x)])$, where $u:\gamma\rightarrow(0,l]$ is some function such that $u(x)=l$ for all but finitely many $x$. If $u(x)<l$ for some $x\in\gamma$, then $\exp_\gamma(x,u(x))$ is a branch point.

\subsection{Star-shaped neighborhoods}

Star-shaped neighborhoods can be think of as generalizations of collar neighborhoods.

\begin{definition}
    A star-shaped neighborhood of $\gamma$ is an embedded annulus $C\subset S$ which is the diffeomorphic image under $\exp_\gamma$ of an open subset $U\subset\gamma\times[0,+\infty)$ of the form $U=\{(x,t)\in\gamma\times[0,+\infty):0\leq t<\psi(x)\}$ for some continuous function $\psi:\gamma\rightarrow(0,+\infty)$. We call $\psi$ the height function. The length of $C$ is defined as $\sup_\gamma\psi$.
\end{definition}

Note that $\exp_\gamma^*g$ is a hyperbolic metric on $\mathcal{U}_\gamma\subset\gamma\times\R$, which can be uniquely extended to a complete hyperbolic metric on $\gamma\times\R$, making all the fibers $\R\rightarrow\{x\}\times\R\ (x\in\gamma)$ unit speed geodesics. We call this the induced hyperbolic metric on $\gamma\times\R$.

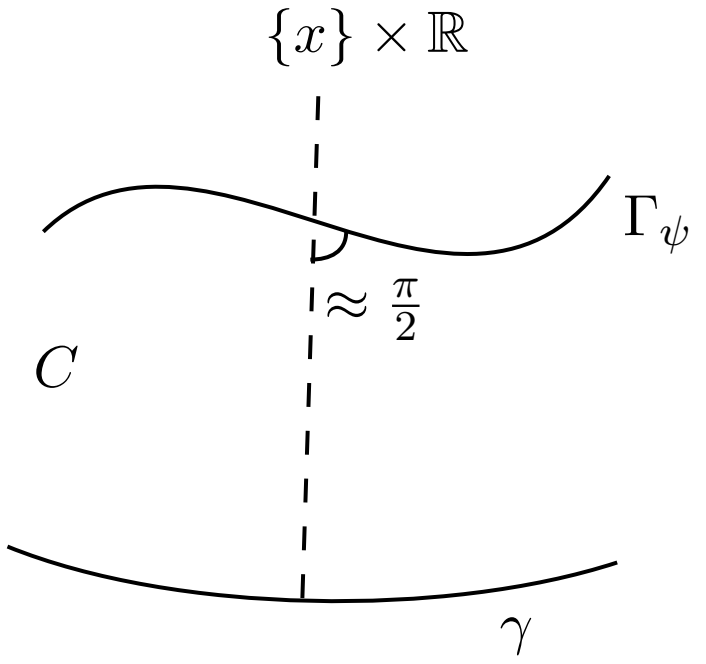
\begin{figure}
    \centering
    \tikzset{every picture/.style={line width=0.75pt}} %set default line width to 0.75pt        

\begin{tikzpicture}[x=0.75pt,y=0.75pt,yscale=-1,xscale=1]

%Curve Lines [id:da3203657085151589] 
\draw    (78,192) .. controls (120,209) and (188,210) .. (231,196) ;
%Straight Lines [id:da8487763429756797] 
\draw  [dash pattern={on 4.5pt off 4.5pt}]  (156,79) -- (152,205) ;
%Curve Lines [id:da11792840192033116] 
\draw    (154,120) .. controls (159,120) and (163,118) .. (163,113) ;
%Curve Lines [id:da5902656950437865] 
\draw    (87,113) .. controls (127,74) and (192,153) .. (229,99) ;

% Text Node
\draw (156,123.4) node [anchor=north west][inner sep=0.75pt]    {$\approx \frac{\pi }{2}$};
% Text Node
\draw (200,208.4) node [anchor=north west][inner sep=0.75pt]    {$\gamma $};
% Text Node
\draw (231,102.4) node [anchor=north west][inner sep=0.75pt]    {$\Gamma _{\psi }$};
% Text Node
\draw (83,140.4) node [anchor=north west][inner sep=0.75pt]    {$C$};
% Text Node
\draw (141,55.4) node [anchor=north west][inner sep=0.75pt]    {$\{x\}\times \mathbb{R}$};
% Text Node
%\draw (93,236) node [anchor=north west][inner sep=0.75pt]   [align=left] {Figure. $\displaystyle \theta $-flatness.};

\end{tikzpicture}
    \caption{$\theta$-flatness}
    \label{fig: theta-flatness}
\end{figure}

\begin{definition}
    Let $\theta\in(0,1]$. A function $\psi:\gamma\rightarrow(0,+\infty)$ is called $\theta$-flat if it is Lipschitz and for every $x\in\gamma$, the angle between $\Gamma_\psi$ (the graph of $\psi$) and the fiber $\{x\}\times\R\subset\gamma\times\R$ at $(x,\psi(x))$ is in the range $[\frac{\pi}{2}-\theta,\frac{\pi}{2}+\theta]$. (Here the angle is measured by the induced hyperbolic metric on $\gamma\times\R$.) A star-shaped neighborhood of $\gamma$ is called $\theta$-flat if its height function $\psi$ is $\theta$-flat (see Figure \ref{fig: theta-flatness}).
\end{definition}

The following lemma says that away from the branch points, a $1$-flat star-shaped neighborhood extends to an embedding on its boundary. The proof is similar in spirit to that of Lemma \ref{lem-c4}.

\begin{lemma}\label{lem-fs1}
    Suppose that $l_\gamma\leq\frac{1}{100}$. Let $C\subset S$ be a $1$-flat star-shaped neighborhood of $\gamma$ of length $\leq\frac{1}{100}$. Let $B\subset\Gamma_\psi$ be the set of points that are mapped to branch points under $\exp_\gamma$. Then $\exp_\gamma$ restricts to an embedding on $\Gamma_\psi\backslash B$.
\end{lemma}

\begin{proof}
    Suppose not. Then there are $x_1\neq x_2\in\gamma$ such that $\exp_\gamma(x_1,\psi(x_1))=\exp_\gamma(x_2,\psi(x_2))=Q\in S$ is not a branch point. Let $\alpha_i=\exp_\gamma(\{x_i\}\times[0,\psi(x_i)])\ (i=1,2)$. Then $\alpha_1,\alpha_2$ are orthogonal geodesics from $Q$ to $\gamma$. By Lemma \ref{lem-m6}, the angle formed by $\alpha_1$ and $\alpha_2$ at $Q$ is less than $\frac{1}{10}$. This contradicts $1$-flatness (see Figure \ref{lem-fs1}).
\end{proof}

\begin{figure}
    \centering
    \tikzset{every picture/.style={line width=0.75pt}} %set default line width to 0.75pt        

\begin{tikzpicture}[x=0.75pt,y=0.75pt,yscale=-1,xscale=1]

%Curve Lines [id:da6116645549170693] 
\draw    (119,114) .. controls (159,103) and (213,131) .. (222,155) ;
%Curve Lines [id:da3926623875485382] 
\draw    (162,90) .. controls (173,132) and (237,145) .. (262,141) ;
%Straight Lines [id:da20807621613995242] 
\draw  [dash pattern={on 4.5pt off 4.5pt}]  (193,127) -- (187,193) ;
%Straight Lines [id:da2523067407172275] 
\draw  [dash pattern={on 4.5pt off 4.5pt}]  (193,127) -- (271,95) ;

% Text Node
\draw (164,158.4) node [anchor=north west][inner sep=0.75pt]    {$\alpha _{1}$};
% Text Node
\draw (260,101.4) node [anchor=north west][inner sep=0.75pt]    {$\alpha _{2}$};
% Text Node
\draw (186,102.4) node [anchor=north west][inner sep=0.75pt]    {$Q$};
% Text Node
\draw (107,89.4) node [anchor=north west][inner sep=0.75pt]    {$\Gamma _{\psi }$};
% Text Node
\draw (255,147.4) node [anchor=north west][inner sep=0.75pt]    {$\Gamma _{\psi }$};
% Text Node
\draw (119,146.4) node [anchor=north west][inner sep=0.75pt]    {$C$};
% Text Node
\draw (212,55.4) node [anchor=north west][inner sep=0.75pt]    {$C$};
% Text Node
%\draw (154,207) node [anchor=north west][inner sep=0.75pt]   [align=left] {Figure.};

\end{tikzpicture}
    \caption{If the boundary touches itself, the flatness condition would imply that the angle between the two orthogonal geodesics cannot be small, a contradiction.}
    \label{fig: lem-fs1}
\end{figure}
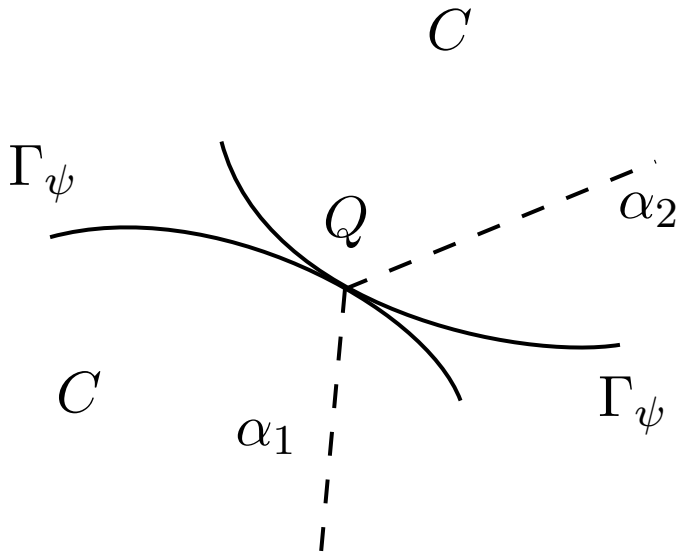

The following Lemma can be think of as a generalization of Lemma \ref{cor-c2}.

\begin{lemma}\label{cor-fs1}
    Suppose that $l_\gamma\leq\frac{1}{100}$. Let $U\subset\mathcal{U}_\gamma$ be a $1$-flat star-shaped neighborhood of $\gamma\times\{0\}$ of length $\leq\frac{1}{100}$ ($\mathcal{U}_\gamma$ is given the pullback hyperbolic metric). Then $\exp_\gamma|_U$ is an embedding. Therefore $\exp_\gamma (U)$ is star-shaped neighborhood of $\gamma$.
\end{lemma}

\begin{proof}
    Let $\psi$ be the height function of $U$. Replacing $\psi$ by $\psi-\epsilon$ for an arbitrarily small $\epsilon$ we may assume that $\exp_\gamma(\overline{U})$ does not contain any branch point.
    
    For $t>0$ define $U_t=U\cap(\gamma\times[0,t))$. Then $U_t$ is $1$-flat for every $t$. For $t$ small enough $\exp_\gamma|_{U_t}$ is an embedding. Suppose that $\exp_\gamma|_U$ is not an embedding. Let $t_0=\sup\{t:\exp_\gamma|_{U_t}\text{ is an embedding}\}$. Then $\exp_\gamma|_{U_{t_0}}$ is an embedding. By Lemma \ref{lem-fs1} we deduce that $\exp_\gamma|_{\overline{U}_{t_0}}$ is an embedding. Since $\exp_\gamma$ is a local diffeomorphism, we deduce that $\exp_\gamma|_{U_{t'}}$ is an embedding for some $t'>t_0$, a contradiction.
\end{proof}

We prove a lemma for later use.

\begin{lemma}\label{lem-fs2}
    Let $C$ be a $\frac{1}{2}$-flat star-shaped neighborhood of $\gamma$ and let $\alpha\subset\partial C$ be a simple geodesic arc (or loop) such that $d(\alpha,\gamma)\geq 2l_\alpha,\ l_{\alpha}\leq\frac{1}{100}$. Then for any geodesic triangle $T$ with sides $\alpha_1,\alpha_2,\alpha_3$, such that $l_{\alpha_2}\leq l_{\alpha_1}=l_\alpha$, and the angle between $\alpha_1,\alpha_2$ is $\leq1$, there exists an isometric embedding $T\hookrightarrow C$ which maps $\alpha_1$ to $\alpha$ (respecting any preferred orientations on $\alpha_1$ and $\alpha$), and maps $\alpha_2,\alpha_3$ to simple geodesic arcs in $C$.
\end{lemma}

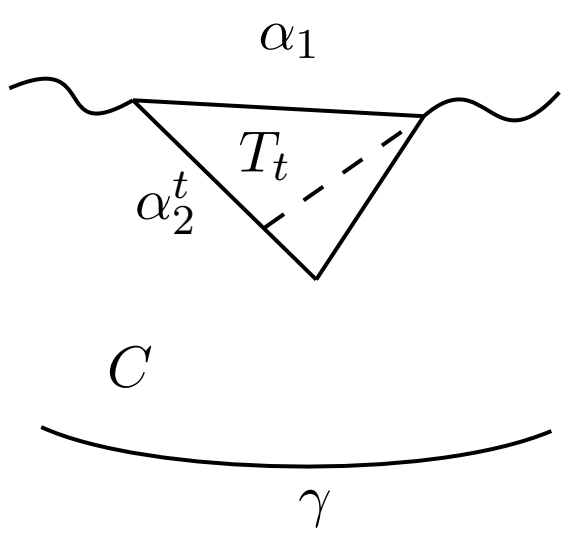
\begin{figure}
    \centering

\tikzset{every picture/.style={line width=0.75pt}} %set default line width to 0.75pt        

\begin{tikzpicture}[x=0.75pt,y=0.75pt,yscale=-1,xscale=1]
%uncomment if require: \path (0,339); %set diagram left start at 0, and has height of 339

%Curve Lines [id:da5553138267921144] 
\draw    (196,221) .. controls (222,233) and (293,235) .. (324,222) ;
%Straight Lines [id:da5274053545184917] 
\draw    (219,139) -- (292,143) ;
%Straight Lines [id:da12726510434395966] 
\draw    (219,139) -- (265,184) ;
%Straight Lines [id:da3158250590241649] 
\draw    (265,184) -- (292,143) ;
%Curve Lines [id:da6661582219685488] 
\draw    (188,136) .. controls (211,126) and (198,151) .. (219,139) ;
%Curve Lines [id:da9546391531207] 
\draw    (292,143) .. controls (308,129) and (309,156) .. (326,137) ;
%Straight Lines [id:da8698971744005968] 
\draw  [dash pattern={on 4.5pt off 4.5pt}]  (252,171) -- (292,143) ;

% Text Node
\draw (259,235.4) node [anchor=north west][inner sep=0.75pt]    {$\gamma $};
% Text Node
\draw (211,199.4) node [anchor=north west][inner sep=0.75pt]    {$C$};
% Text Node
\draw (249,118.4) node [anchor=north west][inner sep=0.75pt]    {$\alpha _{1}$};
% Text Node
\draw (218,155.4) node [anchor=north west][inner sep=0.75pt]    {$\alpha _{2}^{t}$};
% Text Node
\draw (244,145.4) node [anchor=north west][inner sep=0.75pt]    {$T_{t}$};

\end{tikzpicture}

    \caption{Try to show that $\iota$ is an embedding, by considering $\sup t$ for which $\iota|_{T_t}$ is an embedding.}
    \label{fig: lem-fs2}
\end{figure}

\begin{proof}
    We may assume that $C\subset\gamma\times[0,+\infty)$ and $\gamma=\gamma\times\{0\}$. Then $\partial C=\Gamma_\psi$. Since $d(\alpha,\gamma)\geq 2l_\alpha$, there is a local isometric embedding $\iota:T\rightarrow\gamma\times[0,+\infty)$ that maps $\alpha_1$ to $\alpha$ (respecting preferred orientations on $\alpha_1$ and $\alpha$), and maps $T$ to the same side of $\alpha$ as $C$. For $t\in(0,1]$, let $\alpha_2^t\subset\alpha_2$ be the subsegment of length $t\cdot l_{\alpha_2}$ which shares a common vertex with $\alpha_1$. Let $T_t\subset T$ be the geodesic triangle spanned by $\alpha_1$ and $\alpha_2^t$ (see Figure \ref{fig: lem-fs2}). Define
    \[t_0=\sup\{t:\iota|_{T_t}\text{ is an embedding}\}.\]
    Then $\iota|_{T_{t_0}}$ is an embedding. The following claim follows from Lemma \ref{cor-h1} and the $\frac{1}{2}$-flatness of $C$.
    \begin{claim}
        $\iota$ maps the two sides of $T_{t_0}$ other than $\alpha_1$ to simple geodesic arcs contained in $C$.
    \end{claim}
    Therefore the images of the three sides of $T_{t_0}$ meet only at vertices. Since $\iota$ is a local isometric embedding, if $t_0<1$, then there is $t>t_0$ such that $\iota|_{T_t}$ is still an embedding, a contradiction. Therefore $t=1$, so $T_t=T$. The claim implies that $\iota(T)\subset C$.
\end{proof}

\subsection{Flower neighborhoods}

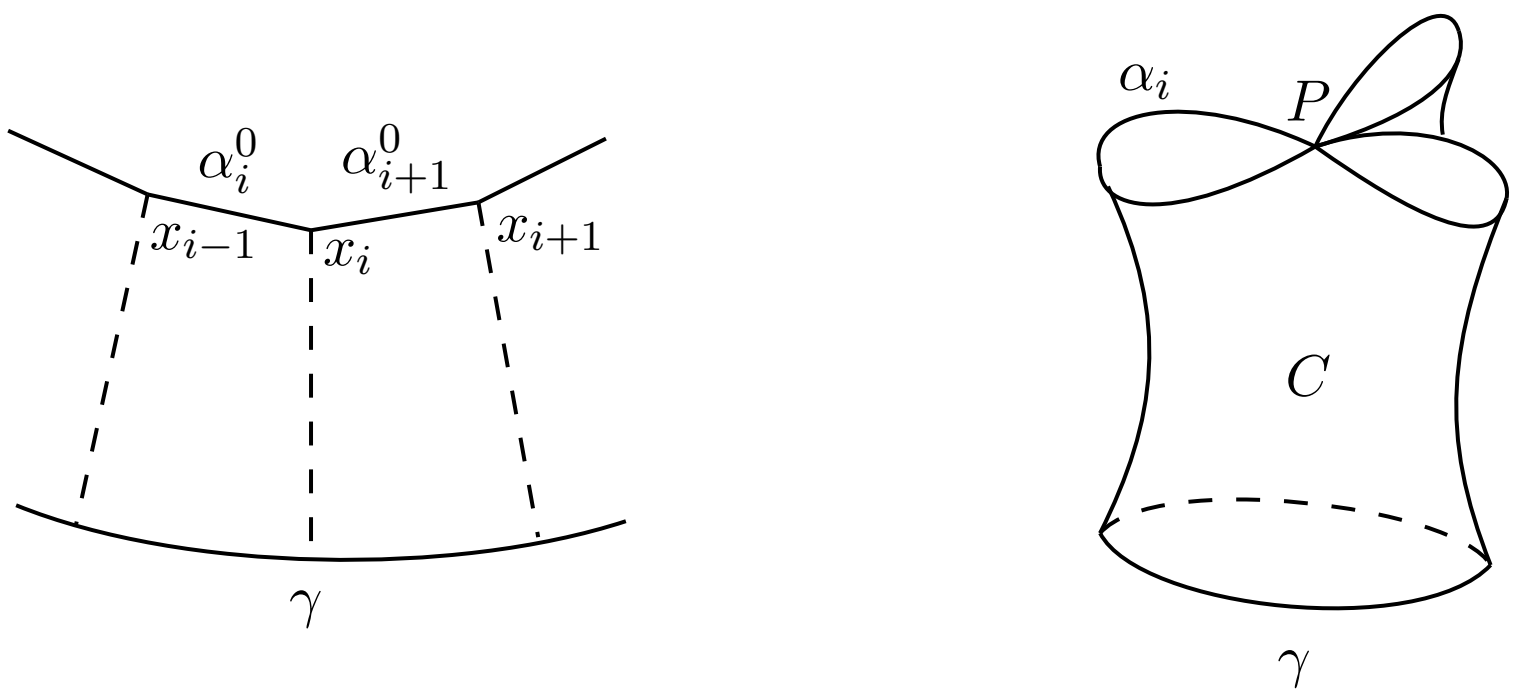
\begin{figure}
    \centering
    \tikzset{every picture/.style={line width=0.75pt}} %set default line width to 0.75pt        

\begin{tikzpicture}[x=0.75pt,y=0.75pt,yscale=-1,xscale=1]

%Curve Lines [id:da3203657085151589] 
\draw    (78,192) .. controls (120,209) and (188,210) .. (231,196) ;
%Straight Lines [id:da06587995078968467] 
\draw    (76,98) -- (111,114) ;
%Straight Lines [id:da6447717906476725] 
\draw    (111,114) -- (152,123) ;
%Straight Lines [id:da43691349550493364] 
\draw    (152,123) -- (194,116) ;
%Straight Lines [id:da8781828428513039] 
\draw    (194,116) -- (226,100) ;
%Straight Lines [id:da8487763429756797] 
\draw  [dash pattern={on 4.5pt off 4.5pt}]  (152,123) -- (152,205) ;
%Straight Lines [id:da2563181364422976] 
\draw  [dash pattern={on 4.5pt off 4.5pt}]  (194,116) -- (209,200) ;
%Straight Lines [id:da3385026549635789] 
\draw  [dash pattern={on 4.5pt off 4.5pt}]  (111,114) -- (93,197) ;
%Curve Lines [id:da9281824323296846] 
\draw    (350,107) .. controls (346,93) and (372,87) .. (404,102) ;
%Curve Lines [id:da11543856247009121] 
\draw    (350,107) .. controls (349,117) and (366,124) .. (404,102) ;
%Curve Lines [id:da8585528972464899] 
\draw    (404,102) .. controls (431,93) and (454,104) .. (452,115) ;
%Curve Lines [id:da7151916330482812] 
\draw    (404,102) .. controls (428,119) and (449,130) .. (452,115) ;
%Curve Lines [id:da044972869870847254] 
\draw    (404,102) .. controls (434,93) and (443,82) .. (440,73) ;
%Curve Lines [id:da09211209100462836] 
\draw    (404,102) .. controls (413,83) and (436,60) .. (440,73) ;
%Curve Lines [id:da8399280511609293] 
\draw    (352,112) .. controls (368,146) and (364,171) .. (350,199) ;
%Curve Lines [id:da48108932880514566] 
\draw    (452,115) .. controls (437,152) and (435,174) .. (448,207) ;
%Curve Lines [id:da5536295738630528] 
\draw    (350,199) .. controls (360,218) and (430,226) .. (448,207) ;
%Curve Lines [id:da3627150928629319] 
\draw  [dash pattern={on 4.5pt off 4.5pt}]  (350,199) .. controls (364,183) and (439,192) .. (448,207) ;
%Curve Lines [id:da2919704639671187] 
\draw    (440,80) .. controls (437,88) and (435,93) .. (436,99) ;

% Text Node
\draw (122,96) node [anchor=north west][inner sep=0.75pt]    {$\alpha^0 _{i}$};
% Text Node
%\draw (118,187) node [anchor=north west][inner sep=0.75pt]    {$\gamma _{i}$};
% Text Node
\draw (145,212) node [anchor=north west][inner sep=0.75pt]    {$\gamma$};
% Text Node
%\draw (80,203) node [anchor=north west][inner sep=0.75pt]    {$x_{i-1}$};
% Text Node
%\draw (202,209) node [anchor=north west][inner sep=0.75pt]    {$x_{i+1}$};
% Text Node
%\draw (168,189) node [anchor=north west][inner sep=0.75pt]    {$\gamma _{i+1}$};
% Text Node
\draw (158,95) node [anchor=north west][inner sep=0.75pt]    {$\alpha^0 _{i+1}$};
% Text Node
\draw (110,120) node [anchor=north west][inner sep=0.75pt]    {$x_{i-1}$};
% Text Node
\draw (153.5,124) node [anchor=north west][inner sep=0.75pt]    {$x_i$};
% Text Node
\draw (197,118) node [anchor=north west][inner sep=0.75pt]    {$x_{i+1}$};
% Text Node
\draw (395,153) node [anchor=north west][inner sep=0.75pt]    {$C$};
% Text Node
\draw (393,227) node [anchor=north west][inner sep=0.75pt]    {$\gamma $};
% Text Node
\draw (353,80) node [anchor=north west][inner sep=0.75pt]    {$\alpha _{i}$};
% Text Node
\draw (395,84) node [anchor=north west][inner sep=0.75pt]    {$P$};
% Text Node
%\draw (167,249) node [anchor=north west][inner sep=0.75pt]   [align=left] {Figure. Flower neighborhood.};

\end{tikzpicture}
    \caption{Flower neighborhood}
    \label{fig: flower neighborhood}
\end{figure}

\begin{definition}
    Let $\theta\in(0,\frac{1}{10}]$. Let $P$ be a branch point. A $\theta$-flower neighborhood of $\gamma$ with pistil $P$ is a $\theta$-flat star-shaped neighborhood $C$ with the following property. There are points $x_1,\cdots,x_n\in\Gamma_\psi$ ($\psi$ is the height function) dividing $\Gamma_\psi$ into arcs $\alpha^0_1,\cdots,\alpha^0_n$, such that for every $i$, $\exp_\gamma(x_i)=P$ and $\exp_\gamma$ maps $\alpha^0_i$ to a simple geodesic loop $\alpha_i$ based at $P$. We call each $\alpha_i$ a petal. (Figure \ref{fig: flower neighborhood})
\end{definition}

\begin{remark}
    The assumption that $P$ is a branch point is not really necessary, but this is the case for all our applications. Moreover, $P$ must be a branch point if there are at least two petals.
\end{remark}

Flower neighborhood exists in our setting.

\begin{lemma}\label{lem-ff1}
    Suppose that $l_\gamma\leq\epsilon$. Let $P$ be a branch point such that $d(P,\gamma)=r_\gamma\leq\frac{1}{200},\ r_P\geq 10\epsilon$. Then there exist an $\epsilon$-flower neighborhood with pistil $P$.
\end{lemma}

\begin{proof}
    By Lemma \ref{lem-fs1} and Lemma \ref{cor-fs1}, it suffices to show that there exists an $\epsilon$-flat star shaped neighborhood $U\subset\mathcal{U}_\gamma$ of $\gamma\times\{0\}$ of length $\leq\frac{1}{100}$ such that $\partial U$ is a piecewise geodesic, whose vertices are mapped to $P$ by $\exp_\gamma$.
    
    Let $l=r_\gamma+5\epsilon$. Since $r_P\geq 10\epsilon$, we have $d(\gamma,\partial S\backslash\gamma)\geq l$ \footnote{Follow a similar argument to that of Lemma \ref{lem-useful criterion}. $P$ lies on $\gamma_1:=\exp_\gamma(\gamma\times\{r_\gamma\})$. A curve from $\gamma$ to $\partial S\backslash\gamma$ would intersect with $\gamma_1$ at some point $Q$. We have $d(P,Q)\leq\epsilon,\ d(Q,\partial S\backslash\gamma)\leq d(\gamma,\partial S\backslash\gamma)-r_\gamma$.}. We have $\mathcal{U}_\gamma\cap(\gamma\times[0,l])=(\gamma\times[0,l])\backslash\bigcup_{i=1}^n(\{x_i\}\times(u(x_i),l])$ for finitely many points $\{(x_i,u(x_i))\}_{1\leq i\leq n}\subset\gamma\times\R$. Moreover, $\exp(x_i,u(x_i))=P$ for every $i$, and we have $r_\gamma=\min_iu(x_i)$. Without loss of generality, assume $u(x_1)=\min_iu(x_i)$. Let $\gamma\times\{a\}$ be the core geodesic of the hyperbolic annulus $\gamma\times\R$, where $a\in\R$.

    \paragraph{Case 1} $u(x_1)\geq a$.

    Let $\alpha$ be the simple geodesic loop based at $(x_1,u(x_1))$ (see Figure \ref{fig: lem-ff1-case1}). Then $\alpha\subset\gamma\times[0,u(x_1)]\subset\mathcal{U}_\gamma$. The annulus $U$ bounded by $\alpha$ and $\gamma\times\{0\}$ is the desired $\epsilon$-flat star-shaped neighborhood. (The flatness condition follows from Lemma \ref{lem-h11}.)

    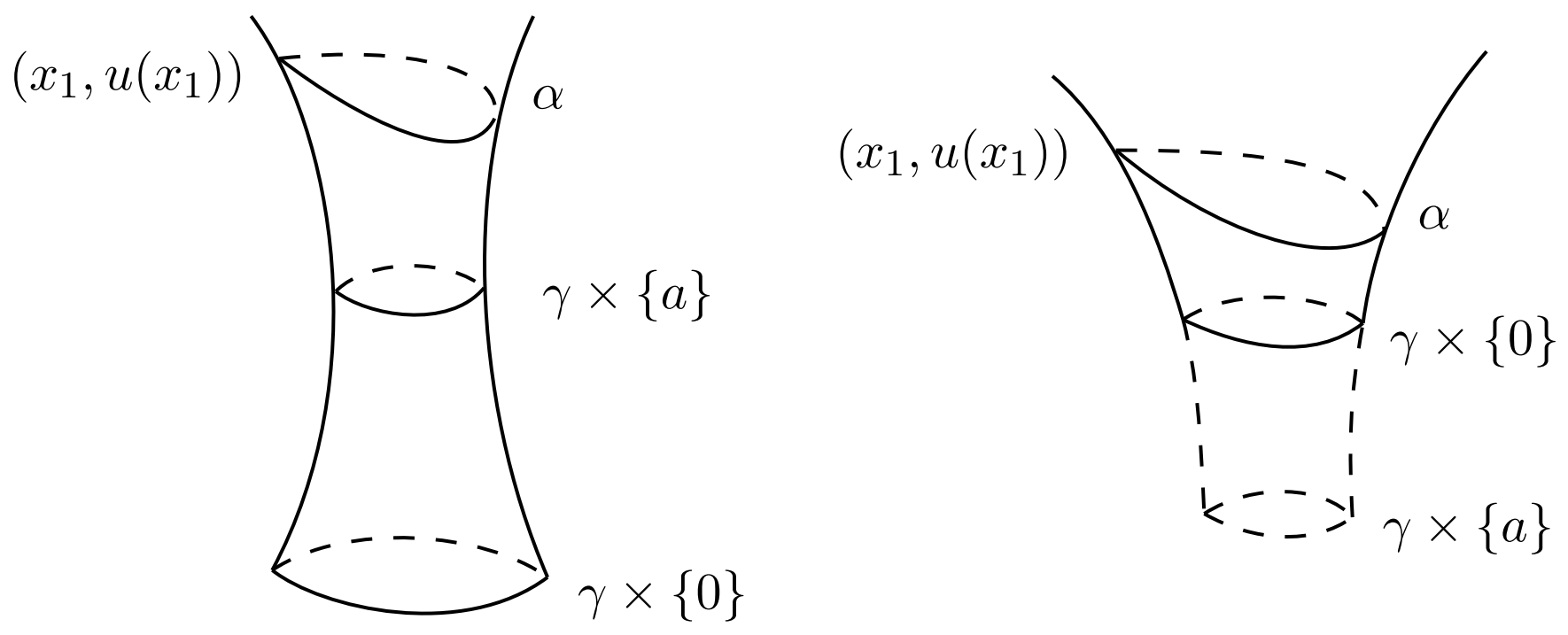
\begin{figure}
        \centering
        \tikzset{every picture/.style={line width=0.75pt}} %set default line width to 0.75pt        

\begin{tikzpicture}[x=0.75pt,y=0.75pt,yscale=-1,xscale=1]

%Curve Lines [id:da5926439381721786] 
\draw    (98,65) .. controls (122,97) and (133,168) .. (104,222) ;
%Curve Lines [id:da06896703223362877] 
\draw    (182,224) .. controls (160,173) and (158,107) .. (178,65) ;
%Curve Lines [id:da540319193887979] 
\draw    (106,77) .. controls (128,94) and (159,110) .. (167,94) ;
%Curve Lines [id:da20244784063449772] 
\draw  [dash pattern={on 4.5pt off 4.5pt}]  (106,77) .. controls (143,73) and (170,80) .. (167,94) ;
%Curve Lines [id:da8938953077509532] 
\draw    (122,143) .. controls (131,150) and (154,154) .. (164,142) ;
%Curve Lines [id:da08632083611661501] 
\draw  [dash pattern={on 4.5pt off 4.5pt}]  (122,143) .. controls (132,133) and (155,134) .. (164,142) ;
%Curve Lines [id:da8286221900651853] 
\draw    (104,222) .. controls (118,233) and (158,242) .. (182,224) ;
%Curve Lines [id:da7225597647300915] 
\draw  [dash pattern={on 4.5pt off 4.5pt}]  (104,222) .. controls (125,208) and (164,211) .. (182,224) ;
%Curve Lines [id:da6832619539314017] 
\draw    (325,82) .. controls (344,98) and (355,127) .. (362,151) ;
%Curve Lines [id:da9366104182292021] 
\draw    (413,152) .. controls (416,128) and (429,97) .. (448,75) ;
%Curve Lines [id:da9690563964937001] 
\draw  [dash pattern={on 4.5pt off 4.5pt}]  (362,151) .. controls (366,166) and (367,184) .. (368,206) ;
%Curve Lines [id:da9020570914131997] 
\draw  [dash pattern={on 4.5pt off 4.5pt}]  (410,207) .. controls (409,189) and (409,173) .. (413,152) ;
%Curve Lines [id:da24045723610819547] 
\draw    (362,151) .. controls (379,159) and (399,163) .. (413,152) ;
%Curve Lines [id:da3345496006353694] 
\draw  [dash pattern={on 4.5pt off 4.5pt}]  (362,151) .. controls (377,141) and (404,144) .. (413,152) ;
%Curve Lines [id:da5097249565563298] 
\draw  [dash pattern={on 4.5pt off 4.5pt}]  (368,206) .. controls (382,197) and (402,198) .. (410,207) ;
%Curve Lines [id:da058639388648107915] 
\draw  [dash pattern={on 4.5pt off 4.5pt}]  (368,206) .. controls (382,214) and (399,215) .. (410,207) ;
%Curve Lines [id:da15554204579402686] 
\draw    (343,103) .. controls (370,126) and (404,138) .. (419,126) ;
%Curve Lines [id:da583942477358431] 
\draw  [dash pattern={on 4.5pt off 4.5pt}]  (343,103) .. controls (383,103) and (418,106) .. (419,126) ;

% Text Node
\draw (189,220.4) node [anchor=north west][inner sep=0.75pt]    {$\gamma \times \{0\}$};
% Text Node
\draw (179,135.4) node [anchor=north west][inner sep=0.75pt]    {$\gamma \times \{a\}$};
% Text Node
\draw (176,83.4) node [anchor=north west][inner sep=0.75pt]    {$\alpha $};
% Text Node
\draw (28,73.4) node [anchor=north west][inner sep=0.75pt]    {$(x_{1} ,u(x_{1} ))$};
% Text Node
\draw (427,117.4) node [anchor=north west][inner sep=0.75pt]    {$\alpha $};
% Text Node
\draw (262,95.4) node [anchor=north west][inner sep=0.75pt]    {$(x_{1} ,u(x_{1} ))$};
% Text Node
\draw (419,148.4) node [anchor=north west][inner sep=0.75pt]    {$\gamma \times \{0\}$};
% Text Node
\draw (417,201.4) node [anchor=north west][inner sep=0.75pt]    {$\gamma \times \{a\}$};
% Text Node
%\draw (137,257) node [anchor=north west][inner sep=0.75pt]   [align=left] {Figure. Case 1: $(x_{1} ,u(x_{1} ))$ is above the axis.};

\end{tikzpicture}
        \caption{Case 1: $(x_{1} ,u(x_{1} ))$ is above the axis.}
        \label{fig: lem-ff1-case1}
    \end{figure}
    
    \paragraph{Case 2} $u(x_1)<a$.

    Let $K$ be the convex hull of $(\gamma\times[a+\infty))\bigcup\{(x_i,u(x_i))\}_{1\leq i\leq n}$. Then $\partial K$ is a piecewise geodesic with vertices in $\{(x_i,u(x_i))\}_{1\leq i\leq n}$. Moreover, $\partial K$ is $\epsilon$-flat, as one can easily deduce from Lemma \ref{lem-h11}. Therefore The annulus $U$ bounded by $\partial K$ and $\gamma\times\{0\}$ is the desired $\epsilon$-flat star-shaped neighborhood.

    \begin{figure}
        \centering
        \tikzset{every picture/.style={line width=0.75pt}} %set default line width to 0.75pt        

\begin{tikzpicture}[x=0.5pt,y=0.5pt,yscale=-1,xscale=1]
%uncomment if require: \path (0,339); %set diagram left start at 0, and has height of 339

%Curve Lines [id:da16457785274803627] 
\draw    (196,51) .. controls (212,112) and (212,177) .. (152,241) ;
%Curve Lines [id:da5636833040096483] 
\draw    (269,50) .. controls (255,104) and (266,185) .. (317,246) ;
%Curve Lines [id:da4316417577384474] 
\draw    (152,241) .. controls (187,275) and (277,276) .. (317,246) ;
%Curve Lines [id:da09461113181773384] 
\draw  [dash pattern={on 4.5pt off 4.5pt}]  (152,241) .. controls (192,211) and (293,220) .. (317,246) ;
%Curve Lines [id:da15379861724042965] 
\draw    (196,172) .. controls (207,166) and (229,176) .. (233,187) ;
%Curve Lines [id:da000391739402249236] 
\draw    (233,187) .. controls (239,176) and (256,175) .. (264,181) ;
%Curve Lines [id:da5505776026484825] 
\draw    (264,181) .. controls (265,173) and (269,169) .. (275,169) ;

% Text Node
%\draw (22,302) node [anchor=north west][inner sep=0.75pt]   [align=left] {Figure. Case 2: $(x_{1} ,u(x_{1} ))$ is below the axis.};
% Text Node
\draw (227,100.4) node [anchor=north west][inner sep=0.75pt]    {$K$};
% Text Node
\draw (284,156.4) node [anchor=north west][inner sep=0.75pt]    {$\partial K$};
% Text Node
\draw (294,258.4) node [anchor=north west][inner sep=0.75pt]    {$\gamma \times \{0\}$};

\end{tikzpicture}
        \caption{Case 2: $(x_{1} ,u(x_{1} ))$ is below the axis.}
        \label{fig: lem-ff1-case2}
    \end{figure}
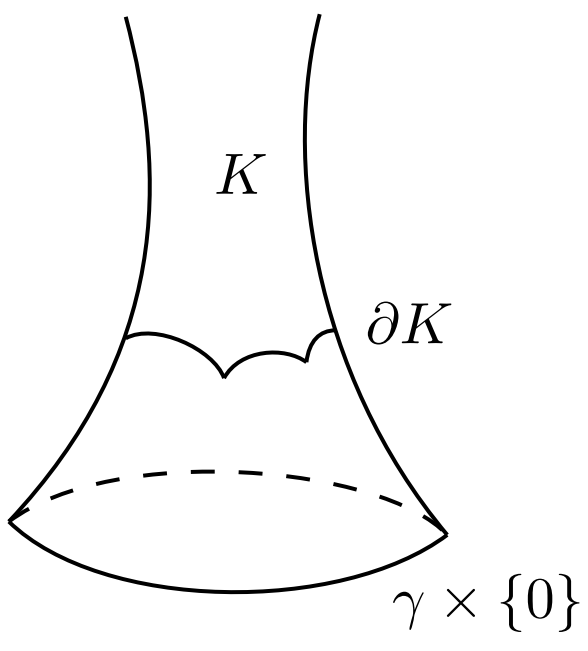
\end{proof}

\begin{lemma}\label{lem-ff2}
    Suppose that $l_\gamma\leq\epsilon,\ r_\gamma\leq\frac{1}{10}$. Let $C$ be a flower neighborhood of $\gamma$. Then every pistil of $C$ has length $\leq3\epsilon$
\end{lemma}

\begin{proof}
    Follows from Lemma \ref{lem-c3}.
\end{proof}

\subsection{$\theta$-extensions}

\begin{figure}
    \centering
    \tikzset{every picture/.style={line width=0.75pt}} %set default line width to 0.75pt        

\begin{tikzpicture}[x=0.75pt,y=0.75pt,yscale=-1,xscale=1]
%uncomment if require: \path (0,339); %set diagram left start at 0, and has height of 339

%Curve Lines [id:da18913243152613401] 
\draw    (164,197) .. controls (213,189) and (269,190) .. (311,197) ;
%Curve Lines [id:da6611201185035477] 
\draw    (164,197) .. controls (191,178) and (220,158) .. (238,138) ;
%Curve Lines [id:da6747058265091662] 
\draw    (238,138) .. controls (254,156) and (288,181) .. (311,197) ;
%Curve Lines [id:da8614122600397857] 
\draw    (180,186) .. controls (182,190) and (183,190) .. (182,195) ;
%Curve Lines [id:da4715015400530308] 
\draw    (291,194) .. controls (291,188) and (291,190) .. (295,186) ;

% Text Node
\draw (234,198) node [anchor=north west][inner sep=0.75pt]    {$\alpha $};
% Text Node
\draw (194,175.4) node [anchor=north west][inner sep=0.75pt]    {$\theta $};
% Text Node
\draw (269,175.4) node [anchor=north west][inner sep=0.75pt]    {$\theta $};
% Text Node
%\draw (179,226) node [anchor=north west][inner sep=0.75pt]   [align=left] {Figure. A $\displaystyle \theta $-cap.};

\end{tikzpicture}
    \caption{A $\theta $-cap}
    \label{fig: cap}
\end{figure}
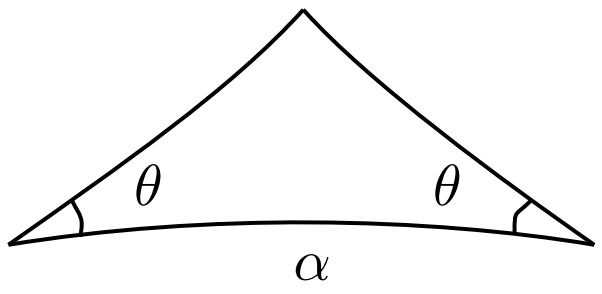

\begin{definition}
    Let $\alpha$ be a geodesic segment. A $\theta$-cap based on $\alpha$ is a geodesic triangle with $\alpha$ being one of its sides such that the angles at the vertices of $\alpha$ are both equal to $\theta$. (Figure \ref{fig: cap})
\end{definition}

\begin{convention}\label{con-interor}
    A triangle refers to its interior.
\end{convention}

\begin{definition}
    Let $C\subset S$ be a subset. Let $\alpha\subset\partial C$ be a simple (open) geodesic arc. A $\theta$-cap of $C$ at $\alpha$ is an embedded $\theta$-cap $T\subset S$ based on $\alpha$ which is disjoint from $C$. We call $C\cup\alpha\cup T$ a $\theta$-extension of $C$ at $\alpha$.
\end{definition}

\begin{remark}
    As in Convention \ref{con-interor}, a $\theta$-cap refers to its interior, so that we allow $\partial T$ to have self-intersection.
\end{remark}

If $C$ is a $1$-petal flower, with petal $\alpha$, by a $\theta$-extension of $C$ we mean a $\theta$-extension of $C$ at $\alpha$.

A $\theta$-extension makes room for performing certain surgeries. It is designed mainly for Model 3 in the next section.

\begin{lemma}\label{lem-ft1}
    Suppose that $l_\gamma\leq\epsilon$. Let $P$ be a branch point such that $d(P,\gamma)=r_\gamma\leq\frac{1}{100},\ r_P\geq 10\epsilon$, and let $C$ be a $\theta$-flower neighborhood of $\gamma$ with pistil $P$ and $n$ petals. Let $\alpha$ be a petal of $C$. Let $\theta_1>0$ satisfies $\theta_2=\theta+\theta_1+8\epsilon<\frac{1}{10}$. Then either $C$ admits a $\theta_1$-extension at $\alpha$, or there exists a $\theta_2$-flower neighborhood with $n+1$ petals.
\end{lemma}

We present the idea before giving a formal proof. Suppose that $C$ does not admit a $\theta_1$-extension at $\alpha$. For small enough $\phi$, $C$ admits a $\phi$-extension at $\alpha$. We gradually increase $\phi$, until the $\phi$-cap meet some obstruction. By Lemma \ref{lem-fs1}, the obstruction cannot be a self intersection or another petal of $C$. It must be the case that the upper boundary of the cap hits some branch point. This branch point must be $P$ (thanks to the condition $r_P\geq 10\epsilon$). In this way we get a flower neighborhood with one more petal, by taking the union of $C$ with a triangle based on $\alpha$, with tip vertex $P$.

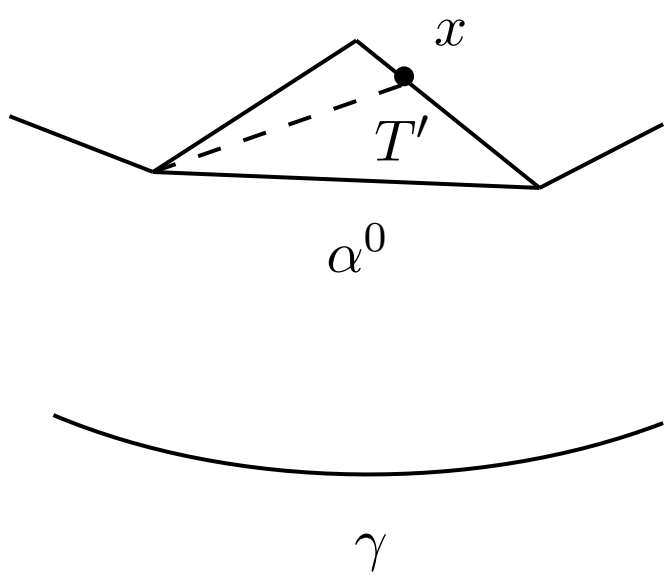
\begin{figure}
    \centering

\tikzset{every picture/.style={line width=0.75pt}} %set default line width to 0.75pt        

\begin{tikzpicture}[x=0.75pt,y=0.75pt,yscale=-1,xscale=1]
%uncomment if require: \path (0,339); %set diagram left start at 0, and has height of 339

%Curve Lines [id:da7943948701463031] 
\draw    (197,214) .. controls (242,233) and (303,234) .. (350,216) ;
%Straight Lines [id:da9169137334492355] 
\draw    (222,153) -- (319,157) ;
%Straight Lines [id:da4876988457187199] 
\draw    (186,139) -- (222,153) ;
%Straight Lines [id:da3808590824339422] 
\draw    (319,157) -- (350,141) ;
%Straight Lines [id:da31797546686472633] 
\draw    (222,153) -- (273,120) ;
%Straight Lines [id:da9641786427910588] 
\draw    (273,120) -- (319,157) ;
%Straight Lines [id:da17049025108809868] 
\draw  [dash pattern={on 4.5pt off 4.5pt}]  (222,153) -- (285,131) ;
%Shape: Circle [id:dp1854028719260391] 
\draw  [fill={rgb, 255:red, 0; green, 0; blue, 0 }  ,fill opacity=1 ] (283,129) .. controls (283,127.9) and (283.9,127) .. (285,127) .. controls (286.1,127) and (287,127.9) .. (287,129) .. controls (287,130.1) and (286.1,131) .. (285,131) .. controls (283.9,131) and (283,130.1) .. (283,129) -- cycle ;

% Text Node
\draw (271,242.4) node [anchor=north west][inner sep=0.75pt]    {$\gamma $};
% Text Node
\draw (264,164.4) node [anchor=north west][inner sep=0.75pt]    {$\alpha ^{0}$};
% Text Node
\draw (291,113.4) node [anchor=north west][inner sep=0.75pt]    {$x$};
% Text Node
\draw (276,137.4) node [anchor=north west][inner sep=0.75pt]    {$T'$};

\end{tikzpicture}

    \caption{Gradually increase $\phi$, until the $\phi$-cap meet a branch point.}
    \label{fig: lem-ft1}
\end{figure}

\begin{proof}
    Suppose that $C=\exp_\gamma(U)$. Let $\alpha^0\subset\Gamma_\psi$ be the (open) arc that is mapped to $\alpha$. Let $T_\phi\subset\gamma\times\R$ be the $\phi$-cap of $U$ based on $\alpha$ (recall that $\gamma\times\R$ is given the induced hyperbolic metric). Since $\theta+\theta_1+8\epsilon<\frac{1}{10}$, we deduce by Lemma \ref{lem-h6} that $E_\phi:=U\cup\alpha^0\cup T_{\phi}$ is a $\frac{1}{10}$-flat star-shaped neighborhood of $\gamma\times\{0\}$ for every $\phi\leq\theta_1$. By Lemma \ref{cor-fs1}, if $E_\phi\subset\mathcal{U}_\gamma\ (\phi\leq\theta_1)$ then $\exp_\gamma|_{E_\phi}$ is an embedding, so $\exp_\gamma(E_\phi)$ is a $\phi$-extension of $C$ at $\alpha$.

    Let $\theta_0=\sup\{\phi\in(0,\theta_1]:E_\phi\subset\mathcal{U}_\gamma\}$. Then $E_{\theta_0}\subset\mathcal{U}_\gamma$. If $\theta_0=\theta_1$ then $C$ admits a $\theta_1$-extension at $\alpha$ and we are done. If $\theta_0<\theta_1$, then there is $x\in\partial T_{\theta_0}\backslash(\alpha^0\cup\{x_0,x_1\})$ (where $x_0,x_1$ are the endpoints of $\alpha^0$) such that $\exp_\gamma(x)$ is a branch point. The condition $r_P\geq 10\epsilon$ shows that this branch point must be $P$. Let $T'\subset T_{\theta_0}$ be the geodesic triangle based on $\alpha^0$ with tip vertex $x$ (see Figure \ref{fig: lem-ft1}). Then $\exp_\gamma(U\cup\alpha^0\cup T')$ is a $\theta_2$-flower neighborhood with $n+1$ petals.
\end{proof}

\subsection{Flower neighborhoods do not intersect}

\begin{lemma}\label{lem-fd1}
    Suppose that two boundary curves $\gamma_1,\gamma_2$ are equidistant curves of length $\leq\epsilon$. Let $C_1,C_2$ be flower neighborhoods of $\gamma_1,\gamma_2$ of length $\leq\frac{1}{10}$. Then $C_1\cap C_2=\varnothing$.
\end{lemma}

\begin{figure}
    \centering

\tikzset{every picture/.style={line width=0.75pt}} %set default line width to 0.75pt        

\begin{tikzpicture}[x=0.75pt,y=0.75pt,yscale=-1,xscale=1]
%uncomment if require: \path (0,339); %set diagram left start at 0, and has height of 339

%Curve Lines [id:da6637623004297187] 
\draw    (109,200) .. controls (137,213) and (216,218) .. (248,204) ;
%Curve Lines [id:da15404001626307184] 
\draw    (150,94) .. controls (173,97) and (194,97) .. (214,96) ;
%Curve Lines [id:da2576093289275001] 
\draw    (97,122) .. controls (114,118) and (135,105) .. (150,94) ;
%Curve Lines [id:da2797012635704539] 
\draw    (214,96) .. controls (233,110) and (252,126) .. (273,140) ;
%Curve Lines [id:da7902942691945625] 
\draw    (120,113) .. controls (160,126) and (220,131) .. (251,125) ;
%Straight Lines [id:da7668984609590694] 
\draw    (191,108) -- (189,212) ;

%Curve Lines [id:da5984054771305141] 
\draw    (368,130) .. controls (403,135) and (435,138) .. (468,131) ;
%Curve Lines [id:da6720630121082287] 
\draw    (318,176) .. controls (338,168) and (357,147) .. (368,130) ;
%Curve Lines [id:da8710916090331056] 
\draw    (468,131) .. controls (481,155) and (503,182) .. (528,199) ;
%Curve Lines [id:da8122772676749406] 
\draw    (337,165) .. controls (380,173) and (460,185) .. (508,183) ;
%Curve Lines [id:da14613675551675165] 
\draw    (337,165) .. controls (378,160) and (435,147) .. (468,131) ;

% Text Node
\draw (136,75.4) node [anchor=north west][inner sep=0.75pt]    {$P$};
% Text Node
\draw (218,78.4) node [anchor=north west][inner sep=0.75pt]    {$P$};
% Text Node
\draw (197,101.4) node [anchor=north west][inner sep=0.75pt]    {$Q$};
% Text Node
\draw (195,129.4) node [anchor=north west][inner sep=0.75pt]    {$Q_{1}$};
% Text Node
\draw (216,219.4) node [anchor=north west][inner sep=0.75pt]    {$\gamma $};
% Text Node
\draw (139,124.4) node [anchor=north west][inner sep=0.75pt]    {$\beta _{0}$};
% Text Node
\draw (421,187.4) node [anchor=north west][inner sep=0.75pt]    {$\beta _{0}$};
% Text Node
\draw (443,154.4) node [anchor=north west][inner sep=0.75pt]    {$T$};
% Text Node
\draw (477,116.4) node [anchor=north west][inner sep=0.75pt]    {$P$};
% Text Node
\draw (350,114.4) node [anchor=north west][inner sep=0.75pt]    {$P$};

\end{tikzpicture}

    \caption{If two Flower neighborhoods $C_1,C_2$ intersect, then one petal $\beta$ of $C_2$ would intersect $C_1$, from with one deduce that the pistil $P$ of $C_1$ lies in the interior of $C_2$, a contradiction.}
    \label{fig: lem-fd1}
\end{figure}
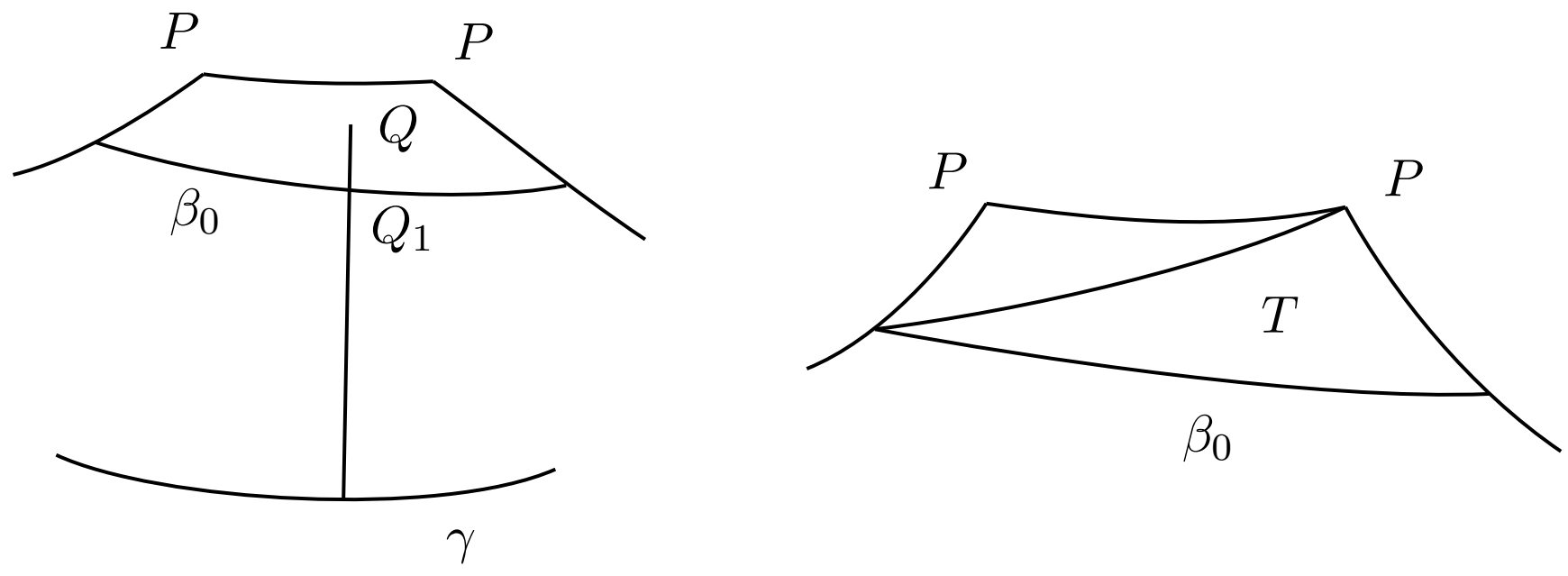

\begin{proof}
    We argue by contradiction. Suppose that $Q\in C_1\cap C_2$. Let $\eta$ be the orthogonal geodesic from $Q$ to $\gamma_1$ contained in $C_1$. Let $Q_1$ be the last intersection point of $\eta$ with $\overline{C_2}$ (see Figure \ref{fig: lem-fd1}). Then $Q_1$ lies on a petal of $C_2$, which we denote by $\beta$. Let $\beta_0\subset\beta$ be the component of $C_1\cap\beta$ containing $Q_1$. Then $\beta_0$ is an (open) geodesic segment contained in $C_1$ with both endpoints lying on $\partial C_1$. We denote by $U$ the component of $C_1\backslash\beta_0$ not adjacent to $\gamma_1$. There exists a geodesic triangle $T\subset U$ based on $\beta_0$ with tip vertex $P$ (the pistil of $C_1$). By Lemma \ref{lem-fd4} below, the sides of $T$ are all $\frac{1}{5}$-flat. Therefore the base angles of $T$ are less than $\frac{1}{2}$. However, by Lemma \ref{lem-fd2} below (applied to $C_2$), there should be an embedded $1$-cap $T'\subset C_2$ based on $\beta$. The triangles $T,T'$ are on the same side of $\beta$, the base of $T$ is contained in that $T'$ and the base angles of $T$ are strictly less than that of $T'$. So the tip of $T$ is in the interior of $T'$. But the tip of $T$ is a branch point. This yields a contradiction.
\end{proof}

\begin{lemma}\label{lem-fd4}
    Suppose that $l_\gamma\leq\epsilon$ and let $C$ be a $\frac{1}{10}$-flat star-shaped neighborhood of $\gamma$ of length $\leq\frac{1}{10}$. Let $\alpha\subset\overline{C}$ be a simple geodesic arc such that both of its end points lie on $\partial C$. Then $\alpha$ is $\frac{1}{5}$-flat (i.e., it is the graph of a $\frac{1}{5}$-flat function).
\end{lemma}

\begin{proof}
    By Lemma \ref{lem-h6}.
\end{proof}

\begin{lemma}\label{lem-fd2}
    Suppose that $l_\gamma\leq\epsilon$ and let $C$ be a flower neighborhood of $\gamma$ of length $\leq\frac{1}{10}$. Then for every petal $\alpha$ of $C$ there exists an embedded $1$-cap $T\subset C$ based on $\alpha$.
\end{lemma}

\begin{proof}
    Follows from Lemma \ref{lem-fs2}.
\end{proof}

\section{Three models}\label{sec-three models}

This section can serve as a black box. We describe three models, on which specific surgeries are performed to merge folding curves. In the next and last section we will reduce standard developing maps to these models, thereby proving Proposition \ref{prop-vs2}. The reader who is eager to see the strategy of the proof of Proposition \ref{prop-vs2} may proceed to the next section, and refer to the specific models when mentioned.

Let $S$ be a good surface, and let $\gamma_1,\gamma_2$ be two boundary curves. We call $(S,\gamma_1,\gamma_2)$ a model if it satisfies one of the three conditions described in Model 1, 2, 3 below.

\paragraph{Common assumptions:}\textit{For all three models we assume that $r_{\gamma_i}\leq 3\delta_1\ (i=1,2)$, and there are flower neighborhoods $C_1,C_2\subset S$ of $\gamma_1,\gamma_2$ respectively, with the same pistil $P$.}

\paragraph{Model 1:}\textit{$C_i$ has $i$ petal(s), and the petal of $C_1$ is a petal of $C_2$.}

\paragraph{Model 2:}\textit{Both $C_1,C_2$ are $1$-petal flowers, and the two petals are sides of an embedded geodesic triangle in which the angle formed by these two petals is $\leq\frac{\pi}{2}$, and the third side is a simple geodesic loop based at $P$.}

\paragraph{Model 3:}\textit{Both $C_1,C_2$ are $1$-petal flowers, and they admit $5\epsilon$-extensions that are disjoint.}

Now let $f$ be an $\epsilon$-excellent developing map with $n$ folding curves. Let $S_0$ be a component of $\Sigma\backslash\mathrm{Fold}(f)$. Let $\gamma_1,\gamma_2$ be two boundary curves of $S_0$ such that $r_{\gamma_i}(S_i)\geq 10\delta_1\ (i=1,2)$, where $S_i$ is the component of $\Sigma\backslash\mathrm{Fold}(f)$ different from $S_0$ that is adjacent to $\gamma_i$. In this section we prove the following proposition:

\begin{proposition}\label{prop-model}
    If $(S_0,\gamma_1,\gamma_2)$ is a model, then there exists a $5\epsilon$-excellent developing map with $n-1$ folding curves and the same monodromy as $f$.
\end{proposition}

Before discussing the models, we make precise the surgeries that we are going to perform.

\subsection{Some surgeries}\label{subsec-3m-some sur}

Let $f$ be a developing map on a closed surface $\Sigma$, for some representation $\rho$.

\paragraph{Removing a pair of annuli}

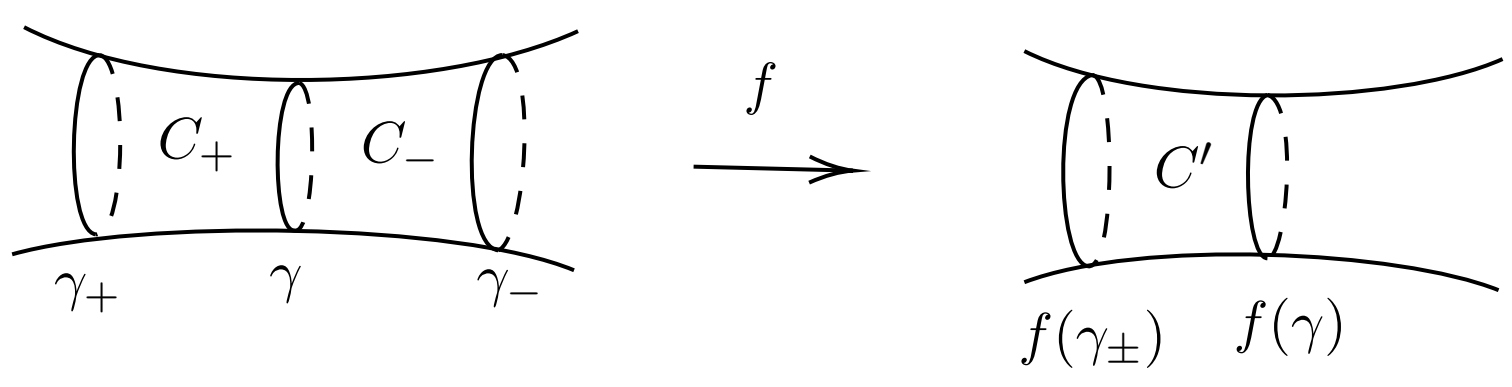
\begin{figure}
    \centering

\tikzset{every picture/.style={line width=0.75pt}} %set default line width to 0.75pt        

\begin{tikzpicture}[x=0.75pt,y=0.75pt,yscale=-1,xscale=1]
%uncomment if require: \path (0,339); %set diagram left start at 0, and has height of 339

%Curve Lines [id:da7932902687682782] 
\draw    (78,87) .. controls (113,105) and (182,104) .. (217,88) ;
%Curve Lines [id:da1758852584290176] 
\draw    (75,144) .. controls (107,135) and (187,136) .. (216,148) ;
%Curve Lines [id:da7327100132719255] 
\draw    (147,101) .. controls (140,101) and (140,139) .. (146,138) ;
%Curve Lines [id:da0446116656537584] 
\draw  [dash pattern={on 4.5pt off 4.5pt}]  (147,101) .. controls (152,102) and (151,139) .. (146,138) ;
%Curve Lines [id:da13867108431680575] 
\draw    (97,94) .. controls (89,94) and (88,138) .. (96,139) ;
%Curve Lines [id:da8539176354980044] 
\draw  [dash pattern={on 4.5pt off 4.5pt}]  (97,94) .. controls (104,95) and (104,139) .. (96,139) ;
%Curve Lines [id:da06885483949967763] 
\draw    (198,94) .. controls (189,94) and (187,142) .. (197,143) ;
%Curve Lines [id:da5287941147189459] 
\draw  [dash pattern={on 4.5pt off 4.5pt}]  (198,94) .. controls (207,95) and (204,140) .. (197,143) ;

%Straight Lines [id:da08217781899309529] 
\draw    (246,122) -- (284,122.95) ;
\draw [shift={(286,123)}, rotate = 181.43] [color={rgb, 255:red, 0; green, 0; blue, 0 }  ][line width=0.75]    (10.93,-3.29) .. controls (6.95,-1.4) and (3.31,-0.3) .. (0,0) .. controls (3.31,0.3) and (6.95,1.4) .. (10.93,3.29)   ;
%Curve Lines [id:da29238484427173284] 
\draw    (329,93) .. controls (359,108) and (422,107) .. (449,95) ;
%Curve Lines [id:da0028032539160133263] 
\draw    (329,151) .. controls (358,140) and (423,143) .. (448,153) ;
%Curve Lines [id:da3750936154492661] 
\draw    (390,104) .. controls (383,106) and (384,145) .. (390,145) ;
%Curve Lines [id:da7867183800895464] 
\draw  [dash pattern={on 4.5pt off 4.5pt}]  (390,104) .. controls (398,106) and (395,145) .. (390,145) ;
%Curve Lines [id:da8126649724461994] 
\draw    (345,147) .. controls (337,146) and (336,101) .. (346,99) ;
%Curve Lines [id:da9705286454090859] 
\draw  [dash pattern={on 4.5pt off 4.5pt}]  (346,99) .. controls (352,99) and (352,148) .. (345,147) ;

% Text Node
\draw (257,94.4) node [anchor=north west][inner sep=0.75pt]    {$f$};
% Text Node
\draw (84,147.4) node [anchor=north west][inner sep=0.75pt]    {$\gamma _{+}$};
% Text Node
\draw (190,146.4) node [anchor=north west][inner sep=0.75pt]    {$\gamma _{-}$};
% Text Node
\draw (138,145.4) node [anchor=north west][inner sep=0.75pt]    {$\gamma $};
% Text Node
\draw (110,108.4) node [anchor=north west][inner sep=0.75pt]    {$C_{+}$};
% Text Node
\draw (161,109.4) node [anchor=north west][inner sep=0.75pt]    {$C_{-}$};
% Text Node
\draw (360,114.4) node [anchor=north west][inner sep=0.75pt]    {$C'$};
% Text Node
\draw (380,153.4) node [anchor=north west][inner sep=0.75pt]    {$f( \gamma )$};
% Text Node
\draw (326,156.4) node [anchor=north west][inner sep=0.75pt]    {$f( \gamma _{\pm })$};

\end{tikzpicture}

    \caption{Removing annuli. $C_{\pm}$ are mapped to the same annulus $C'$. Cut $C_+\cup C_-$ out of $\Sigma$ and glue its boundary curves $\gamma_+,\gamma_-$.}
    \label{fig: sur-remove annuli}
\end{figure}

The surgery of removing annuli is a straight forward generalization of removing collars. Let $C_+,C_-\subset\Sigma$ be closed annuli with disjoint interior, sharing a common boundary curve $\gamma=C_+\cap C_-$. Suppose that $\rho(\gamma)$ is hyperbolic. Suppose that $f|_{C_+}$ (resp. $f|_{C_-}$) is orientation preserving (resp. orientation reversing), and that they are both embeddings into $\H/\rho(\gamma)$, with the same image $C'\subset\H/\rho(\gamma)$. Let $\gamma_{\pm}$ be the boundary curve of $C_{\pm}$ other than $\gamma$. Then the isomorphisms $C_+\cong C'\cong C_-$ induced by $f$ give rise to an identification $\gamma_+\cong\gamma_-$. Remove the open annulus $(C_+\backslash\gamma_+)\cup (C_-\backslash\gamma_-)$ from $\Sigma$ and identify its boundary curves $\gamma_+$ with $\gamma_-$ via the above identification (see Figure \ref{fig: sur-remove annuli}). The resulting surface $\Sigma_1$ is naturally homeomorphic to $\Sigma$, and $f$ restricts to a developing map $f_1$ on $\Sigma_1$ with the same monodromy as $f$. Let $\gamma_1$ denote the curve on $\Sigma_1$ obtained by gluing $\gamma_{\pm}$. If $f$ is branched on a neighborhood of $\gamma_+\cup\gamma_-$, then $\gamma_1$ is a folding curve of $f_1$.

\paragraph{Removing a pair of discs}

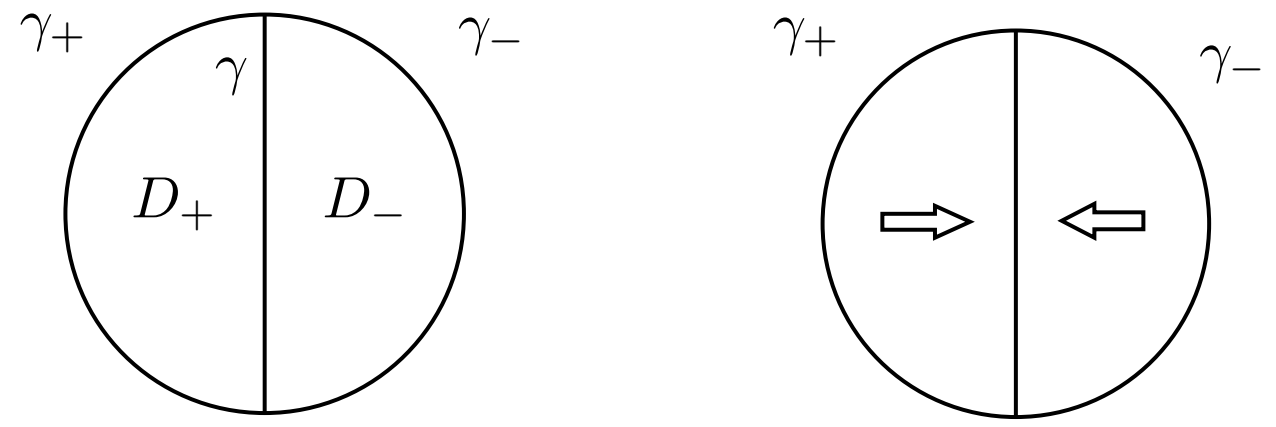
\begin{figure}
    \centering

\tikzset{every picture/.style={line width=0.75pt}} %set default line width to 0.75pt        

\begin{tikzpicture}[x=0.75pt,y=0.75pt,yscale=-1,xscale=1]
%uncomment if require: \path (0,339); %set diagram left start at 0, and has height of 339

%Shape: Circle [id:dp6262396400104863] 
\draw   (100,158) .. controls (100,130.39) and (122.39,108) .. (150,108) .. controls (177.61,108) and (200,130.39) .. (200,158) .. controls (200,185.61) and (177.61,208) .. (150,208) .. controls (122.39,208) and (100,185.61) .. (100,158) -- cycle ;
%Straight Lines [id:da19048628217929908] 
\draw    (150,108) -- (150,208) ;
%Shape: Circle [id:dp3560440404996157] 
\draw   (290,160.5) .. controls (290,133.71) and (311.71,112) .. (338.5,112) .. controls (365.29,112) and (387,133.71) .. (387,160.5) .. controls (387,187.29) and (365.29,209) .. (338.5,209) .. controls (311.71,209) and (290,187.29) .. (290,160.5) -- cycle ;
%Straight Lines [id:da4895374657958015] 
\draw    (338.5,112) -- (338.5,209) ;
%Right Arrow [id:dp004758436795505916] 
\draw   (305,158) -- (318.2,158) -- (318.2,156) -- (327,160) -- (318.2,164) -- (318.2,162) -- (305,162) -- cycle ;
%Left Arrow [id:dp3431925993489333] 
\draw   (350,159.75) -- (358.2,155.5) -- (358.2,157.63) -- (370.5,157.63) -- (370.5,161.88) -- (358.2,161.88) -- (358.2,164) -- cycle ;

% Text Node
\draw (115,147.4) node [anchor=north west][inner sep=0.75pt]    {$D_{+}$};
% Text Node
\draw (136,117.4) node [anchor=north west][inner sep=0.75pt]    {$\gamma $};
% Text Node
\draw (163,147.4) node [anchor=north west][inner sep=0.75pt]    {$D_{-}$};
% Text Node
\draw (87,106.4) node [anchor=north west][inner sep=0.75pt]    {$\gamma _{+}$};
% Text Node
\draw (197,107.4) node [anchor=north west][inner sep=0.75pt]    {$\gamma _{-}$};
% Text Node
\draw (276,107.4) node [anchor=north west][inner sep=0.75pt]    {$\gamma _{+}$};
% Text Node
\draw (383,114.4) node [anchor=north west][inner sep=0.75pt]    {$\gamma _{-}$};

\end{tikzpicture}

    \caption{Removing discs. Remove the disc $D_+\cup D_-$, and identify $\gamma_+$ with $\gamma_-$.}
    \label{fig: sur-remove discs}
\end{figure}

The above surgery can in fact be decomposed into performing a more elementary surgery twice, namely removing discs. Let $D_+,D_-\subset\Sigma$ be embedded closed disc with disjoint interior, such that $\gamma=D_+\cap D_-$ is a connected arc. Let $\gamma_{\pm}\subset\partial D_{\pm}$ be the closure of $\partial D_{\pm}\backslash\gamma$. Suppose that $f|_{D_+}$ (resp. $f|_{D_-}$) is orientation preserving (resp. orientation reversing), and that they are both embeddings into $\H$, with the same image $D'\subset\H$. Then the isomorphisms $D_+\cong D'\cong D_-$ induced by $f$ give rise to an identification $\gamma_+\cong\gamma_-$. Remove the open disc $(D_+\backslash\gamma_+)\cup (D_-\backslash\gamma_-)$ from $\Sigma$ and identify $\gamma_+$ with $\gamma_-$ via the above identification (see Figure \ref{fig: sur-remove discs}). The resulting surface $\Sigma_1$ is naturally homeomorphic to $\Sigma$, and $f$ restricts to a developing map on $\Sigma_1$ with the same monodromy as $f$.

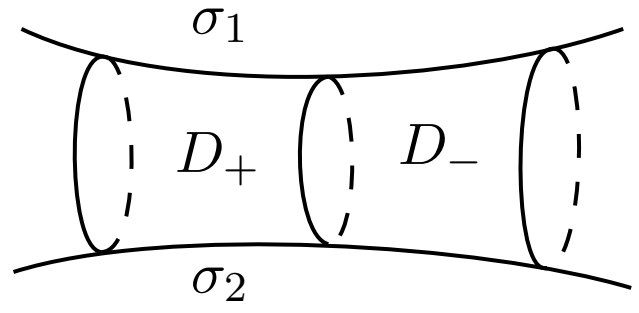
\begin{figure}
    \centering

\tikzset{every picture/.style={line width=0.75pt}} %set default line width to 0.75pt        

\begin{tikzpicture}[x=0.75pt,y=0.75pt,yscale=-1,xscale=1]
%uncomment if require: \path (0,339); %set diagram left start at 0, and has height of 339

%Curve Lines [id:da7709153542646819] 
\draw    (96,116) .. controls (134,134) and (207,130) .. (247,116) ;
%Curve Lines [id:da8420980488734031] 
\draw    (94,177) .. controls (132,165) and (212,170) .. (249,181) ;
%Curve Lines [id:da8462215668744201] 
\draw    (173,128) .. controls (163,129) and (164,168) .. (173,170) ;
%Curve Lines [id:da45009794182673646] 
\draw  [dash pattern={on 4.5pt off 4.5pt}]  (173,128) .. controls (181,130) and (181,170) .. (173,170) ;
%Curve Lines [id:da5100501551144773] 
\draw    (117,123) .. controls (107,121) and (107,171) .. (116,172) ;
%Curve Lines [id:da32944989682634673] 
\draw  [dash pattern={on 4.5pt off 4.5pt}]  (117,123) .. controls (125,124) and (127,172) .. (116,172) ;
%Curve Lines [id:da9239432240394818] 
\draw    (229,121) .. controls (219,122) and (219,175) .. (227,176) ;
%Curve Lines [id:da13748821118270782] 
\draw  [dash pattern={on 4.5pt off 4.5pt}]  (229,121) .. controls (239,119) and (238,178) .. (227,176) ;

% Text Node
\draw (133,140.4) node [anchor=north west][inner sep=0.75pt]    {$D_{+}$};
% Text Node
\draw (137,109.4) node [anchor=north west][inner sep=0.75pt]    {$\sigma_1$};
% Text Node
\draw (137,174.4) node [anchor=north west][inner sep=0.75pt]    {$\sigma_2$};
% Text Node
\draw (189,138.4) node [anchor=north west][inner sep=0.75pt]    {$D_{-}$};

\end{tikzpicture}

    \caption{$\sigma_1,\sigma_2$ divide $C_+$ into two discs $D_+,D_+'$. In the figure, $D_+$ lies in the front, while $D_+'$ lies in the back.}
    \label{fig: sur-remove discs2}
\end{figure}

We show that removing a pair of annuli can be decomposed into a two-step surgery, removing a pair of discs at each step. Indeed, let $C_+,C_-$ be closed annuli as in the setting of removing annuli. Take disjoint simple arcs $\sigma_1,\sigma_2\subset C_+$, each joining the two boundary components of $C_+$ (see Figure \ref{fig: sur-remove discs2}). Denote by $D_+,D_+'$ the closure of the two components of $C_+\backslash(\alpha_+\cup\beta_+)$. Let $D_-,D_-'\subset C_-$ be the images of $D_+,D_+'$ under the isomorphism $C_+\cong C_-$. Then removing $C_{\pm}$ has the same effect as first removing $D_{\pm}$, then removing $D_{\pm}'$.

\paragraph{A variant of removing annuli}

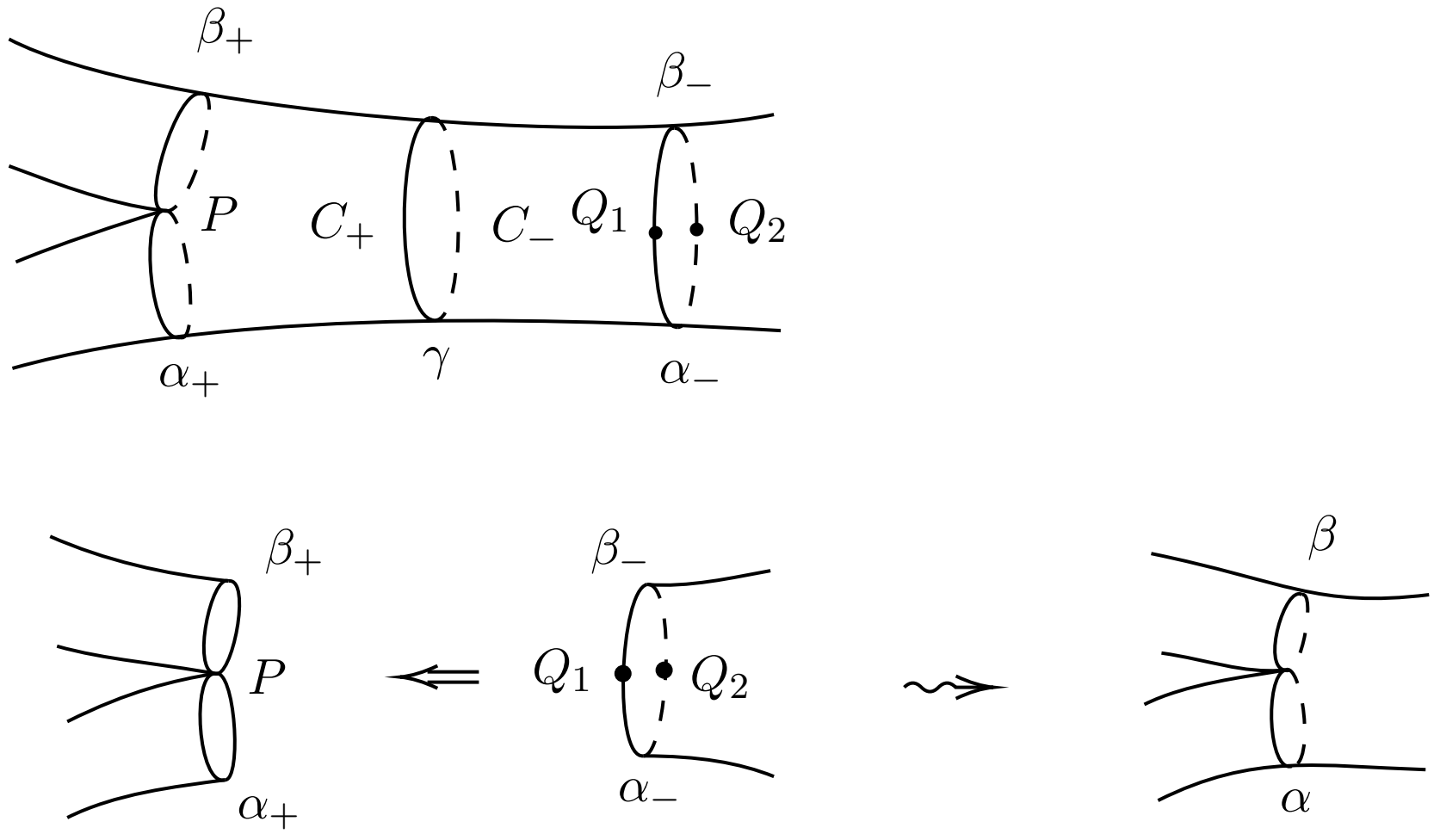
\begin{figure}
    \centering

\tikzset{every picture/.style={line width=0.75pt}} %set default line width to 0.75pt        

\begin{tikzpicture}[x=0.75pt,y=0.75pt,yscale=-1,xscale=1]
%uncomment if require: \path (0,339); %set diagram left start at 0, and has height of 339

%Curve Lines [id:da5901925258682795] 
\draw    (107,41) .. controls (147,62) and (282,73) .. (330,63) ;
%Curve Lines [id:da24561255632088397] 
\draw    (108,137) .. controls (176,118) and (279,123) .. (332,126) ;
%Curve Lines [id:da3938759175768053] 
\draw    (301,67) .. controls (294,65) and (292,126) .. (302,125) ;
%Curve Lines [id:da6094098231587965] 
\draw  [dash pattern={on 4.5pt off 4.5pt}]  (301,67) .. controls (310,66) and (309,126) .. (302,125) ;
%Shape: Circle [id:dp5160023970443959] 
\draw  [fill={rgb, 255:red, 0; green, 0; blue, 0 }  ,fill opacity=1 ] (294,97.5) .. controls (294,96.67) and (294.67,96) .. (295.5,96) .. controls (296.33,96) and (297,96.67) .. (297,97.5) .. controls (297,98.33) and (296.33,99) .. (295.5,99) .. controls (294.67,99) and (294,98.33) .. (294,97.5) -- cycle ;
%Shape: Circle [id:dp7647988811319245] 
\draw  [fill={rgb, 255:red, 0; green, 0; blue, 0 }  ,fill opacity=1 ] (306,96.5) .. controls (306,95.67) and (306.67,95) .. (307.5,95) .. controls (308.33,95) and (309,95.67) .. (309,96.5) .. controls (309,97.33) and (308.33,98) .. (307.5,98) .. controls (306.67,98) and (306,97.33) .. (306,96.5) -- cycle ;
%Curve Lines [id:da9604296366245649] 
\draw    (107,78) .. controls (118,82) and (135,89) .. (152,91) ;
%Curve Lines [id:da20680239751168594] 
\draw    (109,106) .. controls (121,101) and (139,95) .. (152,91) ;
%Curve Lines [id:da9814029804856508] 
\draw    (152,91) .. controls (145,91) and (156,55) .. (163,57) ;
%Curve Lines [id:da5072866978972067] 
\draw  [dash pattern={on 4.5pt off 4.5pt}]  (163,57) .. controls (170,55) and (160,94) .. (152,91) ;
%Curve Lines [id:da6706297733243359] 
\draw    (152,91) .. controls (145,92) and (148,130) .. (157,128) ;
%Curve Lines [id:da906121915806987] 
\draw  [dash pattern={on 4.5pt off 4.5pt}]  (152,91) .. controls (159,89) and (163,130) .. (157,128) ;
%Curve Lines [id:da628646590491332] 
\draw    (230,64) .. controls (220,64) and (219,122) .. (231,123) ;
%Curve Lines [id:da4160555052690267] 
\draw  [dash pattern={on 4.5pt off 4.5pt}]  (230,64) .. controls (240,62) and (241,124) .. (231,123) ;

%Curve Lines [id:da2185908715989997] 
\draw    (119,186) .. controls (135,193) and (152,197) .. (171,199) ;
%Curve Lines [id:da01781741798894043] 
\draw    (124,268) .. controls (138,261) and (155,259) .. (170,257) ;
%Curve Lines [id:da6915567777910329] 
\draw    (124,240) .. controls (142,231) and (155,228) .. (167,226) ;
%Curve Lines [id:da7510314822309843] 
\draw    (121,218) .. controls (135,222) and (151,223) .. (167,226) ;
%Curve Lines [id:da07656007423021205] 
\draw    (171,199) .. controls (166,197) and (160,225) .. (167,226) ;
%Curve Lines [id:da3694087224831851] 
\draw    (171,199) .. controls (178,198) and (172,228) .. (167,226) ;
%Curve Lines [id:da13118837473901623] 
\draw    (167,226) .. controls (160,226) and (162,259) .. (170,257) ;
%Curve Lines [id:da6214098662991753] 
\draw    (167,226) .. controls (173,224) and (175,257) .. (170,257) ;
%Curve Lines [id:da4901672837983313] 
\draw    (293,200) .. controls (305,201) and (318,198) .. (329,196) ;
%Curve Lines [id:da8573565879187492] 
\draw    (292,250) .. controls (308,250) and (320,252) .. (330,256) ;
%Curve Lines [id:da5480598758419549] 
\draw    (293,200) .. controls (285,202) and (283,251) .. (292,250) ;
%Curve Lines [id:da9828188202109746] 
\draw  [dash pattern={on 4.5pt off 4.5pt}]  (293,200) .. controls (302,199) and (299,250) .. (292,250) ;
%Shape: Circle [id:dp6027117313448836] 
\draw  [fill={rgb, 255:red, 0; green, 0; blue, 0 }  ,fill opacity=1 ] (284,226) .. controls (284,224.9) and (284.9,224) .. (286,224) .. controls (287.1,224) and (288,224.9) .. (288,226) .. controls (288,227.1) and (287.1,228) .. (286,228) .. controls (284.9,228) and (284,227.1) .. (284,226) -- cycle ;
%Shape: Circle [id:dp05424244404151579] 
\draw  [fill={rgb, 255:red, 0; green, 0; blue, 0 }  ,fill opacity=1 ] (296,225) .. controls (296,223.9) and (296.9,223) .. (298,223) .. controls (299.1,223) and (300,223.9) .. (300,225) .. controls (300,226.1) and (299.1,227) .. (298,227) .. controls (296.9,227) and (296,226.1) .. (296,225) -- cycle ;
%Straight Lines [id:da02391062187116888] 
\draw    (244,228.5) -- (229,228.5)(244,225.5) -- (229,225.5) ;
\draw [shift={(221,227)}, rotate = 360] [color={rgb, 255:red, 0; green, 0; blue, 0 }  ][line width=0.75]    (10.93,-3.29) .. controls (6.95,-1.4) and (3.31,-0.3) .. (0,0) .. controls (3.31,0.3) and (6.95,1.4) .. (10.93,3.29)   ;
%Straight Lines [id:da4943040861291037] 
\draw    (368,230) .. controls (369.67,228.33) and (371.33,228.33) .. (373,230) .. controls (374.67,231.67) and (376.33,231.67) .. (378,230) .. controls (379.67,228.33) and (381.33,228.33) .. (383,230) -- (384,230) -- (392,230) ;
\draw [shift={(394,230)}, rotate = 180] [color={rgb, 255:red, 0; green, 0; blue, 0 }  ][line width=0.75]    (10.93,-3.29) .. controls (6.95,-1.4) and (3.31,-0.3) .. (0,0) .. controls (3.31,0.3) and (6.95,1.4) .. (10.93,3.29)   ;
%Curve Lines [id:da7831338926885625] 
\draw    (440,191) .. controls (479,199) and (488,207) .. (521,203) ;
%Curve Lines [id:da3516551734766692] 
\draw    (442,263) .. controls (470,249) and (489,253) .. (520,254) ;
%Curve Lines [id:da9201277842610998] 
\draw    (443,220) .. controls (457,222) and (466,226) .. (479,225) ;
%Curve Lines [id:da0917425731209146] 
\draw    (438,235) .. controls (450,229) and (467,227) .. (479,225) ;
%Curve Lines [id:da7692443222806324] 
\draw    (479,225) .. controls (472,224) and (479,200) .. (484,203) ;
%Curve Lines [id:da7529705847084259] 
\draw  [dash pattern={on 4.5pt off 4.5pt}]  (484,203) .. controls (488,202) and (484,227) .. (479,225) ;
%Curve Lines [id:da33310003509556096] 
\draw    (479,225) .. controls (473,227) and (474,255) .. (481,253) ;
%Curve Lines [id:da4196643310297604] 
\draw  [dash pattern={on 4.5pt off 4.5pt}]  (479,225) .. controls (484,222) and (488,253) .. (481,253) ;

% Text Node
\draw (269,83.4) node [anchor=north west][inner sep=0.75pt]    {$Q_{1}$};
% Text Node
\draw (315,86.4) node [anchor=north west][inner sep=0.75pt]    {$Q_{2}$};
% Text Node
\draw (295,133.4) node [anchor=north west][inner sep=0.75pt]    {$\alpha _{-}$};
% Text Node
\draw (294,43.4) node [anchor=north west][inner sep=0.75pt]    {$\beta _{-}$};
% Text Node
\draw (149,134.4) node [anchor=north west][inner sep=0.75pt]    {$\alpha _{+}$};
% Text Node
\draw (160,30.4) node [anchor=north west][inner sep=0.75pt]    {$\beta _{+}$};
% Text Node
\draw (226,129.4) node [anchor=north west][inner sep=0.75pt]    {$\gamma $};
% Text Node
\draw (193,87.4) node [anchor=north west][inner sep=0.75pt]    {$C_{+}$};
% Text Node
\draw (246,88.4) node [anchor=north west][inner sep=0.75pt]    {$C_{-}$};
% Text Node
\draw (161,85.4) node [anchor=north west][inner sep=0.75pt]    {$P$};
% Text Node
\draw (258,217.4) node [anchor=north west][inner sep=0.75pt]    {$Q_{1}$};
% Text Node
\draw (304,219.4) node [anchor=north west][inner sep=0.75pt]    {$Q_{2}$};
% Text Node
\draw (172,260.4) node [anchor=north west][inner sep=0.75pt]    {$\alpha _{+}$};
% Text Node
\draw (283,255.4) node [anchor=north west][inner sep=0.75pt]    {$\alpha _{-}$};
% Text Node
\draw (180,182.4) node [anchor=north west][inner sep=0.75pt]    {$\beta _{+}$};
% Text Node
\draw (175,220.4) node [anchor=north west][inner sep=0.75pt]    {$P$};
% Text Node
\draw (275,182.4) node [anchor=north west][inner sep=0.75pt]    {$\beta _{-}$};
% Text Node
\draw (484,178.4) node [anchor=north west][inner sep=0.75pt]    {$\beta $};
% Text Node
\draw (476,258.4) node [anchor=north west][inner sep=0.75pt]    {$\alpha $};

\end{tikzpicture}

    \caption{A variant of removing annuli. Cut $C_+\cup C_-$ out of $\Sigma$ and glue $\alpha_+$ with $\alpha_-$, $\beta_+$ with $\beta_-$. $P,Q_1,Q_2$ are identified to one point.}
    \label{fig: sur-variant remove annuli1}
\end{figure}

The point of introducing the surgery of removing discs is to make precise the following variant of removing annuli. Let $C_-\subset\Sigma$ be a closed annulus, with boundary curves $\gamma,\gamma_-$. Suppose that $f$ is orientation reversing on $C_-$, and embeds $C_-$ into $\H/\rho(\gamma)$. Let $C_+$ be a half open annulus, such that $\overline{C_+}\cap C_-=C_+\cap C_-=\gamma$. Suppose that $f$ is orientation preserving on $C_+$, and embeds $C_+$ into $\H/\rho(\gamma)$, such that $f(C_+)=f(C_-\backslash\gamma_-)$. This gives rise to an orientation reversing isomorphism $C_+\cong C_-\backslash\gamma_-$. Suppose that $\overline{C_+}\backslash C_+=\alpha_+\cup\beta_+$, where $\alpha_+,\beta_+$ are simple closed curves, such that $\alpha_+\cap\beta_+=\{P\}$ (see Figure \ref{fig: sur-variant remove annuli1}). Suppose further that the isomorphism $C_-\backslash\gamma_-\cong C_+$ extends to a continuous map $\iota:C_-\rightarrow\overline{C_+}$, which is bijective except for $\iota(Q_1)=\iota(Q_2)=P$, where $Q_1,Q_2\in\gamma_-$. Let $\alpha_-=\iota^{-1}(\alpha_+),\ \beta_-=\iota^{-1}(\beta_+)$. Intuitively, we can remove $C_+\cup (C_-\backslash\gamma_-)$ from $\Sigma$ and identify $\alpha_+$ with $\alpha_-$, $\beta_+$ with $\beta_-$ via $\iota$ to get a surface naturally homeomorphic to $\Sigma$, on which $f$ restricts to a developing map with the same monodromy as $f$. This can be made more precise by a two-step surgery which removes a pair of discs at each step, as follows.

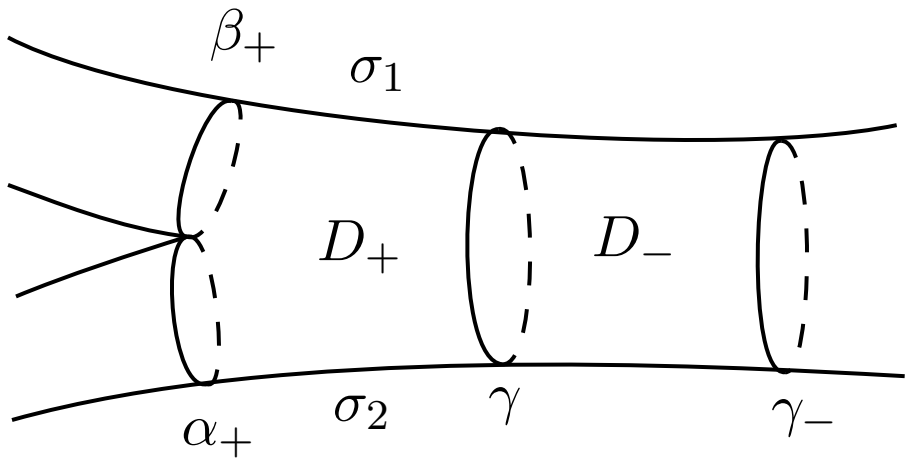
\begin{figure}
    \centering

\tikzset{every picture/.style={line width=0.75pt}} %set default line width to 0.75pt        

\begin{tikzpicture}[x=0.75pt,y=0.75pt,yscale=-1,xscale=1]
%uncomment if require: \path (0,339); %set diagram left start at 0, and has height of 339

%Curve Lines [id:da6462082802279246] 
\draw    (136,67) .. controls (176,88) and (311,99) .. (359,89) ;
%Curve Lines [id:da5357128086525302] 
\draw    (137,163) .. controls (205,144) and (308,149) .. (361,152) ;
%Curve Lines [id:da9666047290177009] 
\draw    (330,93) .. controls (323,91) and (321,152) .. (331,151) ;
%Curve Lines [id:da6775309924232138] 
\draw  [dash pattern={on 4.5pt off 4.5pt}]  (330,93) .. controls (339,92) and (338,152) .. (331,151) ;
%Curve Lines [id:da11407066258490828] 
\draw    (136,104) .. controls (147,108) and (164,115) .. (181,117) ;
%Curve Lines [id:da23206646687763277] 
\draw    (138,132) .. controls (150,127) and (168,121) .. (181,117) ;
%Curve Lines [id:da9499988436497597] 
\draw    (181,117) .. controls (174,117) and (185,81) .. (192,83) ;
%Curve Lines [id:da09245042757169719] 
\draw  [dash pattern={on 4.5pt off 4.5pt}]  (192,83) .. controls (199,81) and (189,120) .. (181,117) ;
%Curve Lines [id:da049364024767532766] 
\draw    (181,117) .. controls (174,118) and (177,156) .. (186,154) ;
%Curve Lines [id:da5482455707313327] 
\draw  [dash pattern={on 4.5pt off 4.5pt}]  (181,117) .. controls (188,115) and (192,156) .. (186,154) ;
%Curve Lines [id:da14194999706053757] 
\draw    (259,90) .. controls (249,90) and (248,148) .. (260,149) ;
%Curve Lines [id:da7629357447889238] 
\draw  [dash pattern={on 4.5pt off 4.5pt}]  (259,90) .. controls (269,88) and (270,150) .. (260,149) ;

% Text Node
\draw (178,161.4) node [anchor=north west][inner sep=0.75pt]    {$\alpha _{+}$};
% Text Node
\draw (185,58.4) node [anchor=north west][inner sep=0.75pt]    {$\beta _{+}$};
% Text Node
\draw (255,153.4) node [anchor=north west][inner sep=0.75pt]    {$\gamma $};
% Text Node
\draw (212,111.4) node [anchor=north west][inner sep=0.75pt]    {$D_{+}$};
% Text Node
\draw (281,110.4) node [anchor=north west][inner sep=0.75pt]    {$D_{-}$};
% Text Node
\draw (220,70.4) node [anchor=north west][inner sep=0.75pt]    {$\sigma _{1}$};
% Text Node
\draw (216,155.4) node [anchor=north west][inner sep=0.75pt]    {$\sigma _{2}$};
% Text Node
\draw (326,156.4) node [anchor=north west][inner sep=0.75pt]    {$\gamma _{-}$};

\end{tikzpicture}

    \caption{$\sigma_1,\sigma_2$ divide $C_+$ into two discs $D_+,D_+'$. In the figure, $D_+$ lies in the front, while $D_+'$ lies in the back.}
    \label{fig: sur-variant remove annuli2}
\end{figure}

Let $\sigma_1,\sigma_2$ be disjoint simple arcs in $\overline{C_+}$, such that $\sigma_1$ joins $\gamma$ with $\beta_+\backslash\{P\}$, while $\sigma_2$ joins $\gamma$ with $\alpha_+\backslash\{P\}$ (see Figure \ref{fig: sur-variant remove annuli2}). Let $D_+,D_+'$ be the closure of the two components of $C_+\backslash(\sigma_1\cup\sigma_2)$, and let $D_-=\iota^{-1}(D_+),D_-'=\iota^{-1}(D_+')$. It is readily verified that $D_{\pm}$ satisfies the conditions of removing discs, so we can remove this pair to obtain a surface $\Sigma_1$ with a developing map $f_1$. Then one can verify that the image of $D_{\pm}'$ in $\Sigma_1$ satisfies the condition of removing discs, so we can further remove this pair to obtain a surface $\Sigma_2$ with a developing map $f_2$. Let $Q$ be the image of $P$ in $\Sigma_2$, and let $\alpha,\beta$ be the image of $\alpha_+,\beta_+$ in $\Sigma_2$ respectively. By a local analysis, it is easy to prove the following lemma:

\begin{lemma}\label{lem-3m1}
    Notation as above. Suppose that $f$ is a local homeomorphism on a neighborhood of $\gamma_-$ and $\beta_+\backslash\{P\}$, while $\alpha_+$ is a folding curve. Then $f_2$ is a local homeomorphism on a neighborhood of $\alpha\backslash\{Q\}$, while $\beta$ is a folding curve of $f_2$.
\end{lemma}

\paragraph{Attaching a pair of annuli}

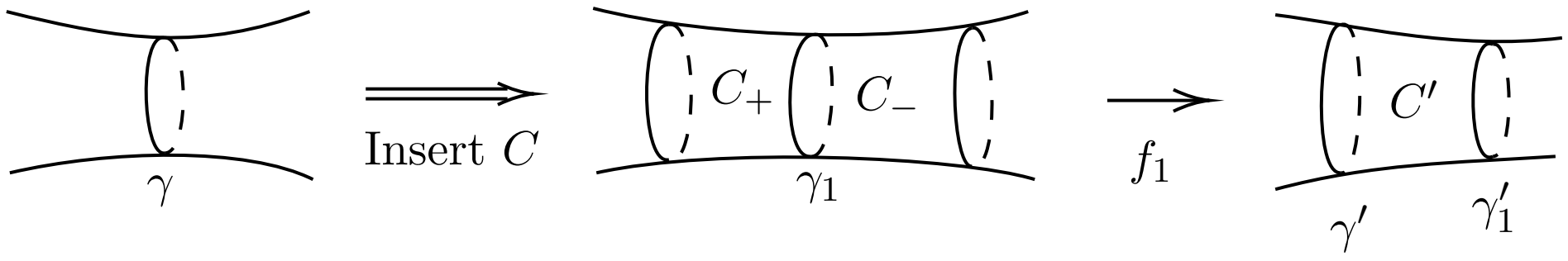
\begin{figure}
    \centering

\tikzset{every picture/.style={line width=0.75pt}} %set default line width to 0.75pt        

\begin{tikzpicture}[x=0.75pt,y=0.75pt,yscale=-1,xscale=1]
%uncomment if require: \path (0,339); %set diagram left start at 0, and has height of 339

%Curve Lines [id:da8269150543559334] 
\draw    (78,100) .. controls (114,109) and (137,112) .. (170,100) ;
%Curve Lines [id:da561335177386727] 
\draw    (79,149) .. controls (109,142) and (151,141) .. (171,151) ;
%Curve Lines [id:da9206825289995298] 
\draw    (126,108) .. controls (118,106) and (119,144) .. (126,143) ;
%Curve Lines [id:da538809965370062] 
\draw  [dash pattern={on 4.5pt off 4.5pt}]  (126,108) .. controls (133,107) and (134,143) .. (126,143) ;
%Straight Lines [id:da5970834249240549] 
\draw    (187,123.5) -- (230,123.5)(187,126.5) -- (230,126.5) ;
\draw [shift={(238,125)}, rotate = 180] [color={rgb, 255:red, 0; green, 0; blue, 0 }  ][line width=0.75]    (10.93,-3.29) .. controls (6.95,-1.4) and (3.31,-0.3) .. (0,0) .. controls (3.31,0.3) and (6.95,1.4) .. (10.93,3.29)   ;

%Curve Lines [id:da18843867154561] 
\draw    (256,99) .. controls (291,108) and (355,111) .. (388,100) ;
%Curve Lines [id:da5945834482951167] 
\draw    (257,149) .. controls (293,141) and (367,144) .. (390,151) ;
%Curve Lines [id:da5826974788328192] 
\draw    (323,107) .. controls (314,107) and (313,145) .. (322,144) ;
%Curve Lines [id:da3945521123132347] 
\draw  [dash pattern={on 4.5pt off 4.5pt}]  (323,107) .. controls (331,105) and (330,145) .. (322,144) ;
%Curve Lines [id:da8538741441354442] 
\draw    (279,104) .. controls (270,103) and (270,147) .. (279,145) ;
%Curve Lines [id:da4729306592733451] 
\draw  [dash pattern={on 4.5pt off 4.5pt}]  (279,104) .. controls (288,102) and (288,148) .. (279,145) ;
%Curve Lines [id:da07859516067198824] 
\draw    (371,105) .. controls (363,106) and (365,148) .. (371,147) ;
%Curve Lines [id:da3921089841398592] 
\draw  [dash pattern={on 4.5pt off 4.5pt}]  (371,105) .. controls (378,103) and (380,149) .. (371,147) ;

%Curve Lines [id:da886169484292377] 
\draw    (463,101) .. controls (487,104) and (517,112) .. (550,108) ;
%Curve Lines [id:da5740328475045804] 
\draw    (463,154) .. controls (493,146) and (517,146) .. (548,144) ;
%Curve Lines [id:da6128350342472482] 
\draw    (529,110) .. controls (521,106) and (521,148) .. (529,144) ;
%Curve Lines [id:da5006707299219293] 
\draw  [dash pattern={on 4.5pt off 4.5pt}]  (529,110) .. controls (535,108) and (536,144) .. (529,144) ;
%Curve Lines [id:da029456072298989944] 
\draw    (483,104) .. controls (476,103) and (474,150) .. (483,149) ;
%Curve Lines [id:da9825983349834561] 
\draw  [dash pattern={on 4.5pt off 4.5pt}]  (483,104) .. controls (491,101) and (490,150) .. (483,149) ;

%Straight Lines [id:da04552457445544145] 
\draw    (412,127) -- (441,127) ;
\draw [shift={(443,127)}, rotate = 180] [color={rgb, 255:red, 0; green, 0; blue, 0 }  ][line width=0.75]    (10.93,-3.29) .. controls (6.95,-1.4) and (3.31,-0.3) .. (0,0) .. controls (3.31,0.3) and (6.95,1.4) .. (10.93,3.29)   ;

% Text Node
\draw (119,148.4) node [anchor=north west][inner sep=0.75pt]    {$\gamma $};
% Text Node
\draw (316,147.4) node [anchor=north west][inner sep=0.75pt]    {$\gamma _{1}$};
% Text Node
\draw (290,116.4) node [anchor=north west][inner sep=0.75pt]    {$C_{+}$};
% Text Node
\draw (334,118.4) node [anchor=north west][inner sep=0.75pt]    {$C_{-}$};
% Text Node
\draw (185,135) node [anchor=north west][inner sep=0.75pt]   [align=left] {Insert $\displaystyle C$};
% Text Node
\draw (417,137.75) node [anchor=north west][inner sep=0.75pt]    {$f_{1}$};
% Text Node
\draw (521,150.4) node [anchor=north west][inner sep=0.75pt]    {$\gamma '_{1}$};
% Text Node
\draw (496,119.4) node [anchor=north west][inner sep=0.75pt]    {$C'$};
% Text Node
\draw (478,157.4) node [anchor=north west][inner sep=0.75pt]    {$\gamma '$};

\end{tikzpicture}

    \caption{Attach annuli. Insert $C$ and extend $f$.}
    \label{fig: sur-attach annuli}
\end{figure}

The inverse surgery of removing annuli is attaching annuli, which is a straight forward generalization of attaching collars. Let $\gamma\subset\Sigma$ be a simple closed curve. Suppose that $\rho(\gamma)$ is hyperbolic and $f$ maps $\gamma$ homeomorphically onto a simple closed curve $\gamma'$ in $\H/\rho(\gamma)$. Let $\gamma_1'\subset\H/\rho(\gamma)$ be a core curve disjoint from $\gamma'$. Then $\gamma',\gamma_1'$ bound a (closed) annulus $C'$ (see Figure \ref{fig: sur-attach annuli}). Cut $\Sigma$ along $\gamma$ and insert an annulus $C$. Let $\gamma_1\subset C$ be a core curve which divides $C$ into two annuli $C_+,C_-$. Extend $f$ to $C$ by mapping both $C_+$ and $C_-$ homeomorphically onto $C'$, $\gamma_1$ homeomorhically onto $\gamma_1'$. Then $\gamma_1$ is a folding curve of $f_1$.

Suppose that $\gamma$ is a folding curve of $f$. Let $N$ be a neighborhood of $\gamma$ such that $N\backslash\gamma= N_+\sqcup N_-$ and $f$ is orientation preserving (resp. reversing) on $N_+$ (resp. $N_-$). Orient $\gamma$ so that $N_+$ lies on the left. Suppose that $C'$ lies on the right of $\gamma'$. Then $f_1$ is branched on a neighborhood of the two boundary curves of the inserted annulus $C$.

Our main application of attaching annuli is the following example, which turns a simple geodesic loop folding curve into an equidistant folding curve.

\paragraph{Example} Let $f$ be a developing map of area $A(f)\leq a_0(\rho)+4\epsilon$. Suppose that $f$ is branched away from $n-1$ folding curves, and all folding curves except one, $\gamma$, are mapped homeomorphically onto equidistant curves of length $\leq\epsilon$, while $\gamma$ is mapped homeomorphically onto a simple geodesic loop $\gamma'$ in $\H/\rho(\gamma)$ of length $\leq4\epsilon$. By Lemma \ref{lem-h7}, there is an equidistant curve $\gamma_1'\subset\H/\rho(\gamma)$ of length $\leq5\epsilon$, disjoint from $\gamma'$, such that the annulus $C'$ bounded by $\gamma'$ and $\gamma'_1$ has area $\leq\frac{1}{2}\epsilon$. We can choose $\gamma_1'$ on the correct side of $\gamma'$ so that the developing map $f_1$ obtained by attaching a pair of annuli $C'$ is branched on a neighborhood of the boundary curves of the inserted annulus. Then $f_1$ is a $5\epsilon$-excellent developing map with $n-1$ folding curves.

\subsection{Model 1 and 2}

Since $r_{\gamma_i}(S_i)\geq 10\delta_1$, $\gamma_i$ has a collar neighborhood of length $5\delta_1$ in $S_i$, which we denote by $N_i$. We have $N_1\cap N_2=\varnothing$.

Denote by $\exp_{\gamma_i}^+:\gamma_i\times[0,+\infty)\supset\mathcal{U}_{\gamma_i}^+\rightarrow S_0$ the exponential map of $\gamma_i$ in $S_0$, and $\exp_{\gamma_i}^-:\gamma_i\times[0,+\infty)\supset\mathcal{U}_{\gamma_i}^-\rightarrow S_i$ the exponential map of $\gamma_i$ in $S_i$.

\paragraph{Common assumptions:}\textit{For all three models we assume that $r_{\gamma_i}(S_0)\leq 3\delta_1\ (i=1,2)$, and there are flower neighborhoods $C_1,C_2\subset S_0$ of $\gamma_1,\gamma_2$ respectively, with the same pistil $P$.}

We have $C_i=\exp_{\gamma_i}^+(C_i^0)$ for some star-shaped neighborhood $C_i^0\subset\gamma_i\times[0,+\infty)$ of $\gamma_i\times\{0\}$. Since $r_{\gamma_i}(S_0)\leq 3\delta_1$, $C_i$ has length $\leq 3\delta_1+3\epsilon$. Therefore $\exp_{\gamma_i}^-|_{C_i^0}$ is an embedding into $N_i$. Its image $C_i':=\exp_{\gamma_i}^-(C_i^0)\subset N_i$ is a star-shaped neighborhood, and we have an orientation reversing isomorphism $C_i'\cong C_i$, which extends to a continuous map $\iota_i:\overline{C_i'}\rightarrow\overline{C_i}$. The restriction $\iota_i:\overline{C_i'}\backslash\iota_i^{-1}(P)\rightarrow\overline{C_i}\backslash\{P\}$ is bijective.

\paragraph{Model 1:}\textit{$C_i$ has $i$ petal(s), and the petal of $C_1$ is a petal of $C_2$.}

Let $\alpha$ be the common petal of $C_1,C_2$ and let $\beta$ be the other petal of $C_2$. let $\alpha_1'=\iota_1^{-1}(\alpha),\alpha_2'=\iota_2^{-1}(\alpha),\beta'=\iota_2^{-1}(\beta)$. In brief, the surgery is as follows: cut $C_1'\cup C_1\cup C_2\cup C_2'$ out of $\Sigma$, and glue $\alpha_1'$ with $\alpha_2'$, $\beta$ with $\beta'$ (see Figure \ref{fig: model 1}).

This surgery can be decomposed into two steps. First remove the pair of annuli $C_1,C_1'$ and identify $\alpha_1'$ with $\alpha$ via $\iota_1$. $f$ restricts to a developing map $f_1$. We still denote by $\alpha$ the curve obtained by gluing $\alpha$ with $\alpha_1'$, which is a folding curve of $f_1$. Then we remove the pair of annuli $C_2,C_2'$, and identify $\alpha_2'$ with $\alpha$, $\beta'$ with $\beta$ via $\iota_2$. $f_1$ restricts to a developing map $f_2$. We still denote by $\alpha,\beta$ their images in the new surface. By Lemma \ref{lem-3m1}, $\beta$ is a folding curve of $f_2$, while $f_2$ is a local homeomorphism on a neighborhood of $\alpha\backslash\beta$. Therefore $f_2$ has $n-1$ folding curves. $f_2$ maps $\beta$ homeomorphically onto a simple geodesic loop in $\H/\rho(\beta)$, which is of length $\leq3\epsilon$ by Lemma \ref{lem-ff2}. By the example at the end of the last subsection, we can attach a pair of annuli at $\beta$ to turn it into an equidistant curve of length $\leq5\epsilon$.

\begin{figure}
    \centering

\tikzset{every picture/.style={line width=0.75pt}} %set default line width to 0.75pt        

\begin{tikzpicture}[x=0.75pt,y=0.75pt,yscale=-1,xscale=1]
%uncomment if require: \path (0,339); %set diagram left start at 0, and has height of 339

%Curve Lines [id:da7641335973309826] 
\draw    (179,39) .. controls (214,63) and (235,64) .. (273,69) ;
%Curve Lines [id:da08728233426935939] 
\draw    (171,129) .. controls (199,121) and (237,113) .. (271,110) ;
%Curve Lines [id:da4943193979575362] 
\draw    (164,72) .. controls (178,76) and (192,81) .. (202,87) ;
%Curve Lines [id:da12146088751766726] 
\draw    (202,87) .. controls (199,78) and (209,59) .. (216,58) ;
%Curve Lines [id:da7159902195999998] 
\draw  [dash pattern={on 4.5pt off 4.5pt}]  (202,87) .. controls (209,86) and (220,65) .. (216,58) ;
%Curve Lines [id:da3745256317948673] 
\draw    (202,87) .. controls (196,92) and (201,117) .. (208,119) ;
%Curve Lines [id:da36825118868766127] 
\draw  [dash pattern={on 4.5pt off 4.5pt}]  (202,87) .. controls (208,88) and (213,115) .. (208,119) ;
%Curve Lines [id:da03472265634575433] 
\draw    (273,69) .. controls (263,73) and (262,106) .. (271,110) ;
%Curve Lines [id:da4606518064086128] 
\draw  [dash pattern={on 4.5pt off 4.5pt}]  (273,69) .. controls (281,73) and (281,109) .. (271,110) ;
%Curve Lines [id:da6640181388806906] 
\draw    (273,69) .. controls (293,66) and (339,51) .. (354,44) ;
%Curve Lines [id:da6286479414078368] 
\draw    (271,110) .. controls (289,113) and (339,130) .. (353,139) ;
%Curve Lines [id:da33078533358808937] 
\draw    (339,51) .. controls (332,54) and (326,74) .. (323,88) ;
%Curve Lines [id:da8755339388364067] 
\draw    (323,88) .. controls (326,103) and (330,129) .. (338,130) ;
%Curve Lines [id:da05286688169328102] 
\draw  [dash pattern={on 4.5pt off 4.5pt}]  (339,51) .. controls (346,46) and (352,75) .. (356,88) ;
%Curve Lines [id:da16561526420305872] 
\draw  [dash pattern={on 4.5pt off 4.5pt}]  (338,130) .. controls (345,128) and (353,106) .. (356,88) ;
%Curve Lines [id:da6583123569593298] 
\draw    (167,96) .. controls (179,92) and (191,89) .. (202,87) ;
%Curve Lines [id:da023458736279594605] 
\draw    (171,129) .. controls (164,128) and (159,100) .. (167,96) ;
%Curve Lines [id:da1740387336228213] 
\draw  [dash pattern={on 4.5pt off 4.5pt}]  (167,96) .. controls (173,99) and (176,125) .. (171,129) ;
%Curve Lines [id:da7445569895630973] 
\draw    (115,96) .. controls (132,99) and (158,99) .. (167,96) ;
%Curve Lines [id:da34681329263435456] 
\draw    (119,138) .. controls (132,135) and (159,129) .. (171,129) ;
%Curve Lines [id:da9443713183618625] 
\draw    (126,98) .. controls (119,100) and (125,132) .. (132,135) ;
%Curve Lines [id:da804866216834738] 
\draw  [dash pattern={on 4.5pt off 4.5pt}]  (126,98) .. controls (133,100) and (138,131) .. (132,135) ;

%Curve Lines [id:da2993899104509542] 
\draw    (157,198) .. controls (168,205) and (178,210) .. (192,211) ;
%Curve Lines [id:da7037600329132584] 
\draw    (144,230) .. controls (156,229) and (173,233) .. (184,239) ;
%Curve Lines [id:da3753944435840858] 
\draw    (192,211) .. controls (185,214) and (178,232) .. (184,239) ;
%Curve Lines [id:da02032031558567271] 
\draw    (192,211) .. controls (198,216) and (190,238) .. (184,239) ;
%Curve Lines [id:da6096001483115285] 
\draw    (89,262) .. controls (99,266) and (115,265) .. (123,261) ;
%Curve Lines [id:da6204844477982024] 
\draw    (91,300) .. controls (102,296) and (118,293) .. (128,293) ;
%Curve Lines [id:da29519465567259595] 
\draw    (123,261) .. controls (116,267) and (120,293) .. (128,293) ;
%Curve Lines [id:da32204154330303825] 
\draw    (123,261) .. controls (130,262) and (135,289) .. (128,293) ;
%Curve Lines [id:da28794117121590435] 
\draw    (300,213) .. controls (318,208) and (335,199) .. (345,193) ;
%Curve Lines [id:da717024951178017] 
\draw    (300,277) .. controls (316,282) and (332,291) .. (340,298) ;
%Curve Lines [id:da5232678061029455] 
\draw    (300,213) .. controls (294,214) and (289,227) .. (284,243) ;
%Curve Lines [id:da01100331424626333] 
\draw    (284,243) .. controls (289,260) and (292,277) .. (300,277) ;
%Curve Lines [id:da9847046638954964] 
\draw  [dash pattern={on 4.5pt off 4.5pt}]  (300,213) .. controls (309,213) and (314,229) .. (318,244) ;
%Curve Lines [id:da7146029238334308] 
\draw  [dash pattern={on 4.5pt off 4.5pt}]  (300,277) .. controls (309,277) and (313,265) .. (318,244) ;
%Straight Lines [id:da9390556655173036] 
\draw    (153.93,271.5) -- (165.94,270.9)(154.07,274.5) -- (166.08,273.9) ;
\draw [shift={(174,272)}, rotate = 177.14] [color={rgb, 255:red, 0; green, 0; blue, 0 }  ][line width=0.75]    (10.93,-3.29) .. controls (6.95,-1.4) and (3.31,-0.3) .. (0,0) .. controls (3.31,0.3) and (6.95,1.4) .. (10.93,3.29)   ;
%Straight Lines [id:da820924237327796] 
\draw    (249,252.5) -- (239,252.5)(249,249.5) -- (239,249.5) ;
\draw [shift={(231,251)}, rotate = 360] [color={rgb, 255:red, 0; green, 0; blue, 0 }  ][line width=0.75]    (10.93,-3.29) .. controls (6.95,-1.4) and (3.31,-0.3) .. (0,0) .. controls (3.31,0.3) and (6.95,1.4) .. (10.93,3.29)   ;
%Curve Lines [id:da21099304806953223] 
\draw    (452,210) .. controls (463,216) and (474,222) .. (487,219) ;
%Curve Lines [id:da690335803572314] 
\draw    (448,246) .. controls (462,242) and (474,243) .. (484,244) ;
%Curve Lines [id:da5068631552644637] 
\draw    (487,219) .. controls (481,221) and (477,243) .. (484,244) ;
%Curve Lines [id:da3550915276044847] 
\draw  [dash pattern={on 4.5pt off 4.5pt}]  (487,219) .. controls (493,222) and (489,242) .. (484,244) ;
%Curve Lines [id:da8110123268965953] 
\draw    (487,219) .. controls (500,217) and (518,213) .. (529,202) ;
%Curve Lines [id:da8308293059705617] 
\draw    (452,282) .. controls (480,266) and (512,274) .. (529,283) ;
%Curve Lines [id:da4972218727879546] 
\draw    (445,256) .. controls (456,250) and (477,245) .. (484,244) ;
%Curve Lines [id:da30594621998144544] 
\draw    (484,244) .. controls (475,246) and (477,272) .. (484,273) ;
%Curve Lines [id:da18002702851520025] 
\draw  [dash pattern={on 4.5pt off 4.5pt}]  (484,244) .. controls (489,245) and (490,270) .. (484,273) ;

%Straight Lines [id:da4123477324255962] 
\draw    (382,248) .. controls (383.67,246.33) and (385.33,246.33) .. (387,248) .. controls (388.67,249.67) and (390.33,249.67) .. (392,248) .. controls (393.67,246.33) and (395.33,246.33) .. (397,248) .. controls (398.67,249.67) and (400.33,249.67) .. (402,248) -- (410,248) ;
\draw [shift={(412,248)}, rotate = 180] [color={rgb, 255:red, 0; green, 0; blue, 0 }  ][line width=0.75]    (10.93,-3.29) .. controls (6.95,-1.4) and (3.31,-0.3) .. (0,0) .. controls (3.31,0.3) and (6.95,1.4) .. (10.93,3.29)   ;

% Text Node
\draw (196,197.4) node [anchor=north west][inner sep=0.75pt]    {$\beta $};
% Text Node
\draw (358,56.4) node [anchor=north west][inner sep=0.75pt]    {$\beta '$};
% Text Node
\draw (356,102.4) node [anchor=north west][inner sep=0.75pt]    {$\alpha _{2}'$};
% Text Node
\draw (265,119.4) node [anchor=north west][inner sep=0.75pt]    {$\gamma _{2}$};
% Text Node
\draw (216,36.4) node [anchor=north west][inner sep=0.75pt]    {$\beta $};
% Text Node
\draw (205,125.4) node [anchor=north west][inner sep=0.75pt]    {$\alpha$};
% Text Node
\draw (160,137.4) node [anchor=north west][inner sep=0.75pt]    {$\gamma _{1}$};
% Text Node
\draw (100,107.4) node [anchor=north west][inner sep=0.75pt]    {$\alpha _{1}'$};
% Text Node
\draw (290,82.4) node [anchor=north west][inner sep=0.75pt]    {$C_{2}'$};
% Text Node
\draw (227,84.4) node [anchor=north west][inner sep=0.75pt]    {$C_{2}$};
% Text Node
\draw (178,101.4) node [anchor=north west][inner sep=0.75pt]    {$C_{1}$};
% Text Node
\draw (136,106.4) node [anchor=north west][inner sep=0.75pt]    {$C_{1}'$};
% Text Node
\draw (271,216.4) node [anchor=north west][inner sep=0.75pt]    {$\beta '$};
% Text Node
\draw (265,260.4) node [anchor=north west][inner sep=0.75pt]    {$\alpha _{2}'$};
% Text Node
\draw (137,279.4) node [anchor=north west][inner sep=0.75pt]    {$\alpha _{1}'$};
% Text Node
\draw (157,208.4) node [anchor=north west][inner sep=0.75pt]    {$+$};
% Text Node
\draw (95,270.4) node [anchor=north west][inner sep=0.75pt]    {$-$};
% Text Node
\draw (327,244.4) node [anchor=north west][inner sep=0.75pt]    {$-$};
% Text Node
\draw (483,193.4) node [anchor=north west][inner sep=0.75pt]    {$\beta $};
% Text Node
\draw (475,279.4) node [anchor=north west][inner sep=0.75pt]    {$\alpha$};
% Text Node
\draw (454,219.4) node [anchor=north west][inner sep=0.75pt]    {$+$};
% Text Node
\draw (456,254.4) node [anchor=north west][inner sep=0.75pt]    {$-$};
% Text Node
\draw (514,235.4) node [anchor=north west][inner sep=0.75pt]    {$-$};

\end{tikzpicture}

    \caption{Model 1. Cut $C_1'\cup C_1\cup C_2\cup C_2'$ out of $\Sigma$, and glue $\alpha_1'$ with $\alpha_2'$, $\beta$ with $\beta'$. `$+$' and `$-$' stand for orientation.}
    \label{fig: model 1}
\end{figure}

\paragraph{Model 2:}\textit{Both $C_1,C_2$ are $1$-petal flowers, and the two petals are sides of an embedded geodesic triangle in which the angle formed by these two petals is $\leq\frac{\pi}{2}$, and the third side is a simple geodesic loop based at $P$.}

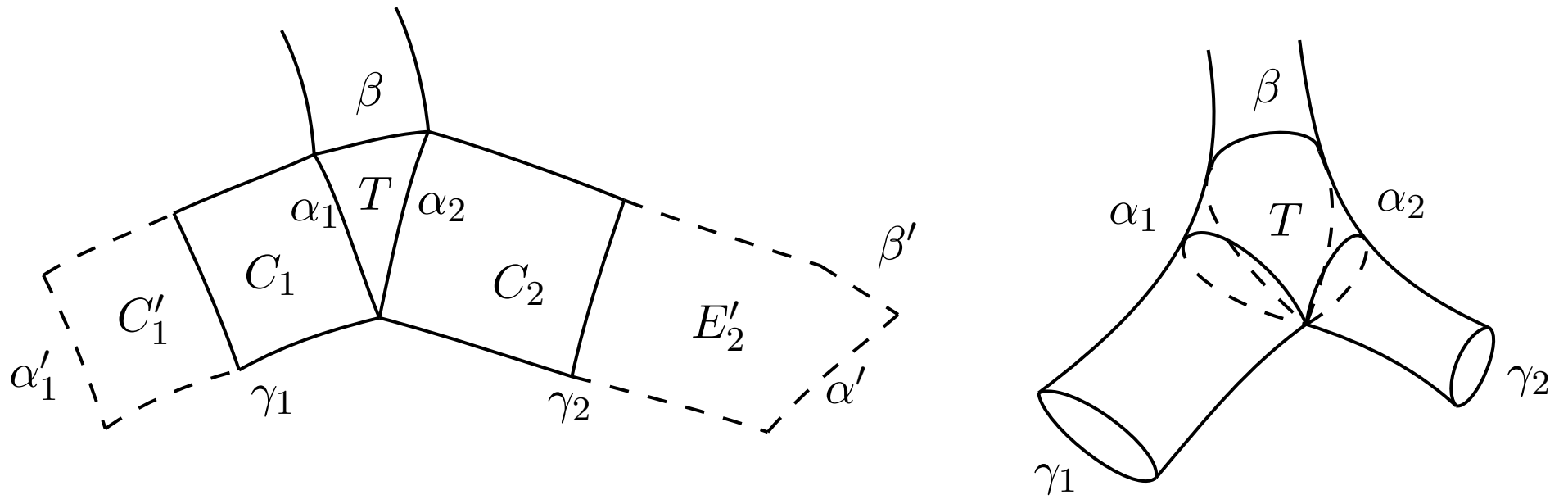
\begin{figure}
    \centering
    
\tikzset{every picture/.style={line width=0.75pt}} %set default line width to 0.75pt        

\begin{tikzpicture}[x=0.75pt,y=0.75pt,yscale=-1,xscale=1]
%uncomment if require: \path (0,339); %set diagram left start at 0, and has height of 339

%Curve Lines [id:da5039567731588581] 
\draw    (150,110) .. controls (163,107) and (172,104) .. (185,103) ;
%Curve Lines [id:da3441236899616442] 
\draw    (150,110) .. controls (157,121) and (164,146) .. (170,160) ;
%Curve Lines [id:da6475668594366567] 
\draw    (185,103) .. controls (178,120) and (174,141) .. (170,160) ;
%Curve Lines [id:da5736752307151006] 
\draw    (185,103) .. controls (203,108) and (233,119) .. (245,124) ;
%Curve Lines [id:da26317994157902613] 
\draw    (170,160) .. controls (189,165) and (213,173) .. (229,178) ;
%Curve Lines [id:da5643559606683708] 
\draw    (245,124) .. controls (240,139) and (233,159) .. (229,178) ;
%Curve Lines [id:da6836542165767561] 
\draw  [dash pattern={on 4.5pt off 4.5pt}]  (245,124) .. controls (269,133) and (289,139) .. (305,144) ;
%Curve Lines [id:da06501055233739983] 
\draw  [dash pattern={on 4.5pt off 4.5pt}]  (229,178) .. controls (252,184) and (271,189) .. (289,195) ;
%Curve Lines [id:da6268024946223831] 
\draw  [dash pattern={on 4.5pt off 4.5pt}]  (305,144) .. controls (313,149) and (321,154) .. (329,159) ;
%Curve Lines [id:da5012056756849272] 
\draw  [dash pattern={on 4.5pt off 4.5pt}]  (289,195) .. controls (301,181) and (316,170) .. (329,159) ;
%Curve Lines [id:da4498665316409859] 
\draw    (107,128) .. controls (122,121) and (137,116) .. (150,110) ;
%Curve Lines [id:da13918390562496674] 
\draw    (127,176) .. controls (143,167) and (158,163) .. (170,160) ;
%Curve Lines [id:da47540602967294066] 
\draw    (107,128) .. controls (114,143) and (122,161) .. (127,176) ;
%Curve Lines [id:da18005314655661864] 
\draw  [dash pattern={on 4.5pt off 4.5pt}]  (67,147) .. controls (79,139) and (95,134) .. (107,128) ;
%Curve Lines [id:da5165749756964981] 
\draw  [dash pattern={on 4.5pt off 4.5pt}]  (86,194) .. controls (100,184) and (116,179) .. (127,176) ;
%Curve Lines [id:da09625181592939658] 
\draw  [dash pattern={on 4.5pt off 4.5pt}]  (67,147) .. controls (75,163) and (81,178) .. (86,194) ;
%Curve Lines [id:da5359619276515958] 
\draw    (140,72) .. controls (146,84) and (149,96) .. (150,110) ;
%Curve Lines [id:da6766289183096943] 
\draw    (175,65) .. controls (180,75) and (184,91) .. (185,103) ;

%Curve Lines [id:da1152350370937909] 
\draw    (372,183) .. controls (412,153) and (432,125) .. (424,78) ;
%Curve Lines [id:da585586908907715] 
\draw    (408,209) .. controls (424,190) and (436,175) .. (454,162) ;
%Curve Lines [id:da15647282205586222] 
\draw    (452,75) .. controls (458,121) and (468,144) .. (510,163) ;
%Curve Lines [id:da5942637510045666] 
\draw    (454,162) .. controls (473,168) and (488,176) .. (500,187) ;
%Curve Lines [id:da2625677302112719] 
\draw    (372,183) .. controls (382,179) and (410,200) .. (408,209) ;
%Curve Lines [id:da5124461523281999] 
\draw    (372,183) .. controls (373,193) and (401,215) .. (408,209) ;
%Curve Lines [id:da5105375007721469] 
\draw    (510,163) .. controls (503,163) and (495,181) .. (500,187) ;
%Curve Lines [id:da08695195525947397] 
\draw    (500,187) .. controls (506,188) and (515,166) .. (510,163) ;
%Curve Lines [id:da5432286160418404] 
\draw    (417,136) .. controls (424,128) and (447,145) .. (454,162) ;
%Curve Lines [id:da6254470027172131] 
\draw  [dash pattern={on 4.5pt off 4.5pt}]  (417,136) .. controls (412,145) and (434,158) .. (454,162) ;
%Curve Lines [id:da16566454542047582] 
\draw    (454,162) .. controls (456,148) and (467,130) .. (472,136) ;
%Curve Lines [id:da26774207096698066] 
\draw  [dash pattern={on 4.5pt off 4.5pt}]  (454,162) .. controls (465,156) and (475,143) .. (472,136) ;
%Curve Lines [id:da672122706851039] 
\draw    (425,114) .. controls (426,103) and (457,99) .. (458,109) ;
%Curve Lines [id:da742974194977921] 
\draw  [dash pattern={on 4.5pt off 4.5pt}]  (425,114) .. controls (419,123) and (430,141) .. (454,162) ;
%Curve Lines [id:da10571413867629997] 
\draw  [dash pattern={on 4.5pt off 4.5pt}]  (458,109) .. controls (463,119) and (465,134) .. (454,162) ;

% Text Node
\draw (180,120.4) node [anchor=north west][inner sep=0.75pt]    {$\alpha _{2}$};
% Text Node
\draw (203,142.4) node [anchor=north west][inner sep=0.75pt]    {$C_{2}$};
% Text Node
\draw (264,153.4) node [anchor=north west][inner sep=0.75pt]    {$E_{2}'$};
% Text Node
\draw (220,181.4) node [anchor=north west][inner sep=0.75pt]    {$\gamma _{2}$};
% Text Node
\draw (55,168.4) node [anchor=north west][inner sep=0.75pt]    {$\alpha _{1}'$};
% Text Node
\draw (88,151.4) node [anchor=north west][inner sep=0.75pt]    {$C_{1}'$};
% Text Node
\draw (129,179.4) node [anchor=north west][inner sep=0.75pt]    {$\gamma _{1}$};
% Text Node
\draw (127,139.4) node [anchor=north west][inner sep=0.75pt]    {$C_{1}$};
% Text Node
\draw (141,122.4) node [anchor=north west][inner sep=0.75pt]    {$\alpha _{1}$};
% Text Node
\draw (161,83.4) node [anchor=north west][inner sep=0.75pt]    {$\beta $};
% Text Node
\draw (162,115.4) node [anchor=north west][inner sep=0.75pt]    {$T$};
% Text Node
\draw (305,173.4) node [anchor=north west][inner sep=0.75pt]    {$\alpha '$};
% Text Node
\draw (321,129.4) node [anchor=north west][inner sep=0.75pt]    {$\beta '$};
% Text Node
\draw (369,203.4) node [anchor=north west][inner sep=0.75pt]    {$\gamma _{1}$};
% Text Node
\draw (514,173.4) node [anchor=north west][inner sep=0.75pt]    {$\gamma _{2}$};
% Text Node
\draw (436,82.4) node [anchor=north west][inner sep=0.75pt]    {$\beta $};
% Text Node
\draw (392,123.4) node [anchor=north west][inner sep=0.75pt]    {$\alpha _{1}$};
% Text Node
\draw (474,119.4) node [anchor=north west][inner sep=0.75pt]    {$\alpha _{2}$};
% Text Node
%\draw (230,250) node [anchor=north west][inner sep=0.75pt]   [align=left] {Figure. Model 2.};
% Text Node
\draw (441,123.4) node [anchor=north west][inner sep=0.75pt]    {$T$};

\end{tikzpicture}
    \caption{Model 2. $C_2\cong C_2'$ extends to an isomorphism $C_2\cup\alpha_2\cup T\cong E_2'$.}
    \label{fig: model 2}
\end{figure}

Denote by $\alpha_i$ the petal of $C_i$. Without loss of generality, we may assume that $l_{\alpha_1}\leq l_{\alpha_2}$. Let $T$ be the embedded geodesic triangle in the assumption, and denote its third side by $\beta$. Let $E_2=C_2\cup\alpha_2\cup T$. Let $\alpha_2'=\iota_2^{-1}(\alpha_2)$. By Lemma \ref{lem-h5}, the isomorphism $\iota_2^{-1}:C_2\cong C_2'$ can be extended to an isomorphism $\hat\iota_2^{-1}:E_2\cong E_2'=C_2'\cup\alpha_2'\cup T'$, for some embedded geodesic triangle $T'\subset N_2$ isomorphic to $T$ (see Figure \ref{fig: model 2}). Moreover, $\overline{E_2'}$ is a closed annulus. $\hat\iota_2$ can be extended to a continuous map $\hat\iota_2:\overline{E_2'}\rightarrow\overline{E_2}$, which is bijective except for $\hat\iota_2(P')=\hat\iota_2(Q)=P$, where $P'$ is the pistil of $C_2'$, and $Q$ is the vertex of $T'$ opposite to $\alpha_2'$. Let $\alpha'=\hat\iota_2^{-1}(\alpha_1),\ \beta'=\hat\iota_2^{-1}(\beta),\ \alpha_1'=\iota_1^{-1}(\alpha_1)$. We can then proceed in exactly the same way as in model 1: Cut $C_1'\cup C_1\cup E_2\cup E_2'$ out of $\Sigma$, and glue $\alpha_1'$ with $\alpha'$, $\beta$ with $\beta$.

\subsection{More surgeries}

To discuss Model 3, we need to introduce one more surgery, namely \textit{merging folding curves}.

\paragraph{Attaching a pair of discs}

\begin{figure}
    \centering

\tikzset{every picture/.style={line width=0.75pt}} %set default line width to 0.75pt        

\begin{tikzpicture}[x=0.75pt,y=0.75pt,yscale=-1,xscale=1]
%uncomment if require: \path (0,339); %set diagram left start at 0, and has height of 339

%Curve Lines [id:da7609210017041207] 
\draw    (96,79) .. controls (108,106) and (85,136) .. (97,155) ;
%Straight Lines [id:da12210365892377484] 
\draw    (132,118) -- (179,118) ;
\draw [shift={(181,118)}, rotate = 180] [color={rgb, 255:red, 0; green, 0; blue, 0 }  ][line width=0.75]    (10.93,-3.29) .. controls (6.95,-1.4) and (3.31,-0.3) .. (0,0) .. controls (3.31,0.3) and (6.95,1.4) .. (10.93,3.29)   ;
%Curve Lines [id:da41087049287148525] 
\draw    (236,85) .. controls (212,93) and (212,147) .. (237,153) ;
%Curve Lines [id:da8805143673064597] 
\draw    (236,85) .. controls (261,95) and (261,147) .. (237,153) ;
%Curve Lines [id:da2152699111433214] 
\draw    (218,199) .. controls (181,218) and (175,264) .. (218,281) ;
%Straight Lines [id:da4200135949473708] 
\draw    (218,199) -- (218,281) ;
%Curve Lines [id:da45019982596954833] 
\draw    (300,204) .. controls (262,217) and (260,267) .. (302,281) ;
%Straight Lines [id:da31431394078025476] 
\draw    (300,204) -- (302,281) ;
%Straight Lines [id:da7274849911000724] 
\draw    (301,120) -- (349,120) ;
\draw [shift={(351,120)}, rotate = 180] [color={rgb, 255:red, 0; green, 0; blue, 0 }  ][line width=0.75]    (10.93,-3.29) .. controls (6.95,-1.4) and (3.31,-0.3) .. (0,0) .. controls (3.31,0.3) and (6.95,1.4) .. (10.93,3.29)   ;
%Straight Lines [id:da38545226349607076] 
\draw    (430,75) -- (430,154) ;
%Curve Lines [id:da6532767646007066] 
\draw    (430,154) .. controls (384,144) and (379,87) .. (430,75) ;
%Curve Lines [id:da5048533161914851] 
\draw    (430,75) .. controls (484,89) and (473,149) .. (430,154) ;
%Straight Lines [id:da04516684705378382] 
\draw    (332,202) -- (369.51,168.34) ;
\draw [shift={(371,167)}, rotate = 138.09] [color={rgb, 255:red, 0; green, 0; blue, 0 }  ][line width=0.75]    (10.93,-3.29) .. controls (6.95,-1.4) and (3.31,-0.3) .. (0,0) .. controls (3.31,0.3) and (6.95,1.4) .. (10.93,3.29)   ;

% Text Node
\draw (143,92) node [anchor=north west][inner sep=0.75pt]   [align=left] {cut};
% Text Node
\draw (202,75.4) node [anchor=north west][inner sep=0.75pt]    {$\gamma _{+}$};
% Text Node
\draw (256,75.4) node [anchor=north west][inner sep=0.75pt]    {$\gamma _{-}$};
% Text Node
\draw (69,75.4) node [anchor=north west][inner sep=0.75pt]    {$\gamma $};
% Text Node
\draw (324,140) node [anchor=north west][inner sep=0.75pt]   [align=left] {glue};
% Text Node
\draw (403,109.4) node [anchor=north west][inner sep=0.75pt]    {$D_{+}$};
% Text Node
\draw (440,110.4) node [anchor=north west][inner sep=0.75pt]    {$D_{-}$};
% Text Node
\draw (434,86.4) node [anchor=north west][inner sep=0.75pt]    {$\beta $};
% Text Node
\draw (377,72.4) node [anchor=north west][inner sep=0.75pt]    {$\gamma _{+}$};
% Text Node
\draw (461,70.4) node [anchor=north west][inner sep=0.75pt]    {$\gamma _{-}$};
% Text Node
\draw (171,205.4) node [anchor=north west][inner sep=0.75pt]    {$\gamma '_{+}$};
% Text Node
\draw (224,242.4) node [anchor=north west][inner sep=0.75pt]    {$\beta _{+}$};
% Text Node
\draw (193,232.4) node [anchor=north west][inner sep=0.75pt]    {$D_{+}$};
% Text Node
\draw (277,233.4) node [anchor=north west][inner sep=0.75pt]    {$D_{-}$};
% Text Node
\draw (251,211.4) node [anchor=north west][inner sep=0.75pt]    {$\gamma '_{-}$};
% Text Node
\draw (308,246.4) node [anchor=north west][inner sep=0.75pt]    {$\beta _{-}$};

\end{tikzpicture}

    \caption{Attach discs. Cut along $\gamma$ and glue in $D_{\pm}$.}
    \label{fig: attach discs}
\end{figure}
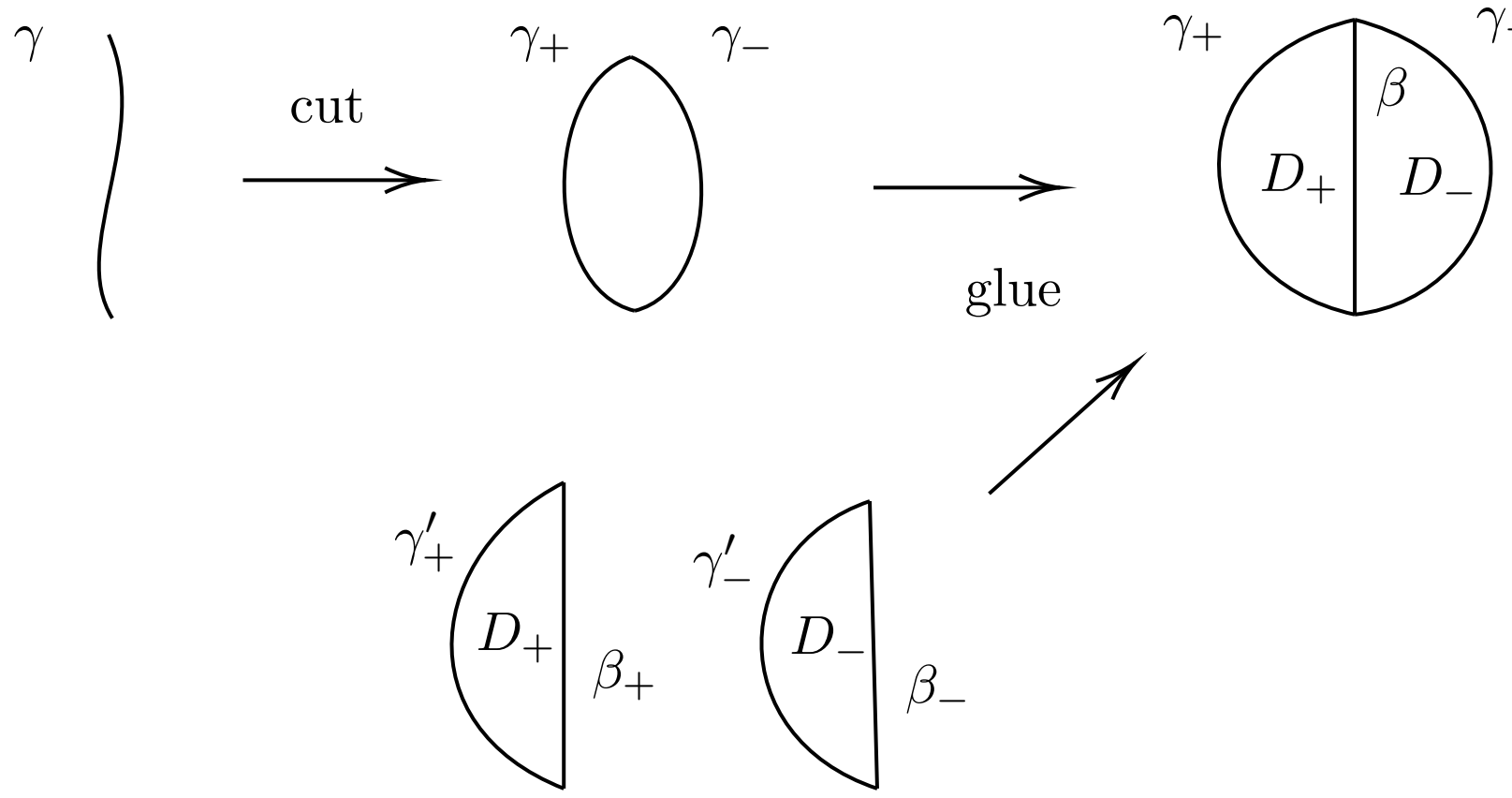

The inverse surgery of removing discs is attaching discs, as follows. Let $\gamma\subset\Sigma$ be a simple arc. Suppose that $\gamma'=f(\widetilde{\gamma})$ is a simple arc in $\H$, where $\widetilde{\gamma}\subset\widetilde{\Sigma}$ is a lift of $\gamma$. Let $D\subset\H$ be a closed disc such that $\gamma'\subset\partial D$. Let $\beta$ be the closure of $\partial D\backslash\gamma'$. Let $D_+,D_-$ be two copies of $D$, and denote by $\gamma'_{\pm},\beta_{\pm}$ the corresponding arcs in $\partial D_{\pm}$. Cut $\Sigma$ along $\gamma$. The resulting surface $\Sigma^{cut}$ has a boundary curve composed of two copies of $\gamma$, denoted $\gamma_+,\gamma_-$, such that $\Sigma^{cut}$ lies on the left of $\gamma_+$. Now let $\Sigma_1=\Sigma^{cut}\cup D_+\cup D_-$, identifying $\gamma_{\pm}$ with $\gamma_{\pm}'$ via $f$, and identifying $\beta_+$ with $\beta_-$ via the identity (see Figure \ref{fig: attach discs}). Extend $f$ to $\Sigma_1$ by mapping $D_{\pm}$ to $D$ via the identity.

\paragraph{Merging folding curves}

Let $\alpha,\beta\subset\Sigma$ be simple closed curves such that $\alpha\cap\beta=\{P\},\ i(\alpha,\beta)=0$. Let $N$ be a small neighborhood of $\alpha\cup\beta$. We have $N\backslash(\alpha\cup\beta)=N_0\sqcup N_1\sqcup N_2$, where $\overline{N_0}$ contains $\alpha\cup\beta$, $\overline{N_1}$ contains $\alpha$, $\overline{N_2}$ contains $\beta$. Orient $\alpha,\beta$ so that $N_0$ lies on the left. Let $\gamma=\alpha*\beta$. Let $\widetilde{\gamma}=\widetilde{\alpha}*\widetilde{\beta}\subset\widetilde{\Sigma}$ be a lift of $\gamma$. Let $f$ be a developing map on $\Sigma$. Suppose that $\gamma'=f(\widetilde{\gamma})$ is a simple curve in $\H$. Let $D\subset\H$ be a closed disc such that $\gamma'\subset\partial D$, and $D$ lies on the right of $\gamma'$. Let $\eta$ be the closure of $\partial D\backslash\gamma'$. Let $D_+,D_-$ be two copies of $D$. Denote by $\gamma'_{\pm},\eta_{\pm}$ the corresponding arcs in $\partial D_{\pm}$.

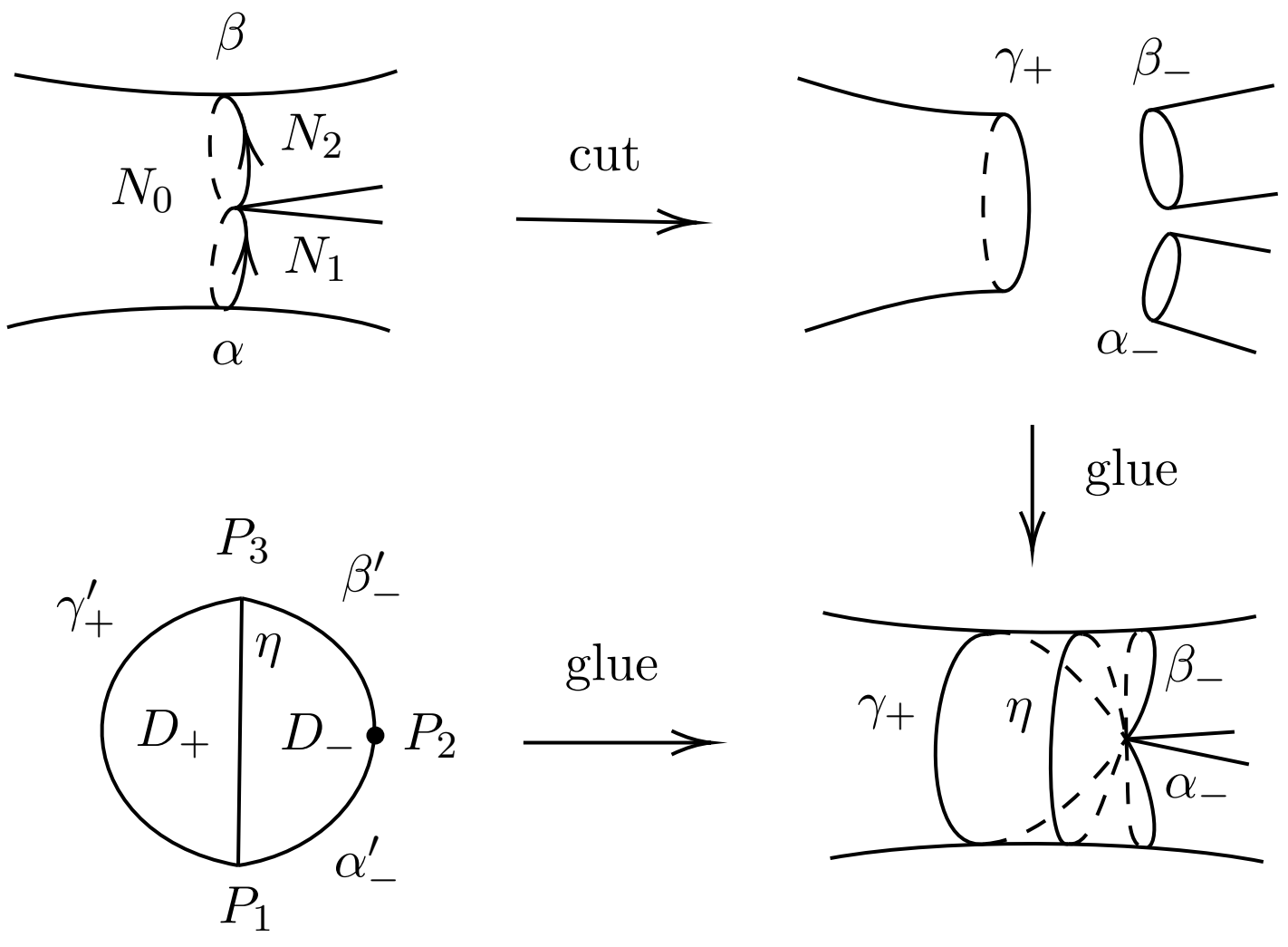
\begin{figure}
    \centering

\tikzset{every picture/.style={line width=0.75pt}} %set default line width to 0.75pt        

\begin{tikzpicture}[x=0.75pt,y=0.75pt,yscale=-1,xscale=1]
%uncomment if require: \path (0,339); %set diagram left start at 0, and has height of 339

%Curve Lines [id:da6932778110059908] 
\draw    (88,84) .. controls (120,90) and (165,93) .. (194,83) ;
%Curve Lines [id:da9341425165804076] 
\draw    (86,154) .. controls (111,147) and (167,146) .. (192,155) ;
%Straight Lines [id:da005165414945121571] 
\draw    (149,121) -- (190,115) ;
%Straight Lines [id:da9790747971890299] 
\draw    (149,121) -- (190,125) ;
%Curve Lines [id:da5534484862852777] 
\draw    (149,121) .. controls (156,120) and (153,90) .. (146,90) ;
\draw [shift={(151.72,98.91)}, rotate = 80.49] [color={rgb, 255:red, 0; green, 0; blue, 0 }  ][line width=0.75]    (10.93,-3.29) .. controls (6.95,-1.4) and (3.31,-0.3) .. (0,0) .. controls (3.31,0.3) and (6.95,1.4) .. (10.93,3.29)   ;
%Curve Lines [id:da45372193528190397] 
\draw  [dash pattern={on 4.5pt off 4.5pt}]  (146,90) .. controls (140,91) and (141,123) .. (149,121) ;
%Curve Lines [id:da3986661736251248] 
\draw    (149,121) .. controls (155,120) and (152,150) .. (146,149) ;
\draw [shift={(152.24,128.49)}, rotate = 91.88] [color={rgb, 255:red, 0; green, 0; blue, 0 }  ][line width=0.75]    (10.93,-3.29) .. controls (6.95,-1.4) and (3.31,-0.3) .. (0,0) .. controls (3.31,0.3) and (6.95,1.4) .. (10.93,3.29)   ;
%Curve Lines [id:da06880578343655175] 
\draw  [dash pattern={on 4.5pt off 4.5pt}]  (149,121) .. controls (144,120) and (139,149) .. (146,149) ;
%Straight Lines [id:da6929620790406648] 
\draw    (227,124) -- (275,124.96) ;
\draw [shift={(277,125)}, rotate = 181.15] [color={rgb, 255:red, 0; green, 0; blue, 0 }  ][line width=0.75]    (10.93,-3.29) .. controls (6.95,-1.4) and (3.31,-0.3) .. (0,0) .. controls (3.31,0.3) and (6.95,1.4) .. (10.93,3.29)   ;
%Curve Lines [id:da23319206524027747] 
\draw    (305,85) .. controls (321,90) and (338,95) .. (362,95) ;
%Curve Lines [id:da5089574321393866] 
\draw    (307,155) .. controls (323,150) and (340,144) .. (362,144) ;
%Curve Lines [id:da5675421619426837] 
\draw    (362,95) .. controls (371,95) and (372,144) .. (362,144) ;
%Curve Lines [id:da8931917013390764] 
\draw  [dash pattern={on 4.5pt off 4.5pt}]  (362,95) .. controls (355,96) and (354,143) .. (362,144) ;
%Straight Lines [id:da3673776265580254] 
\draw    (402,94) -- (437,87) ;
%Straight Lines [id:da018525266562273268] 
\draw    (408,121) -- (438,117) ;
%Curve Lines [id:da626059966764401] 
\draw    (402,94) .. controls (410,91) and (415,119) .. (408,121) ;
%Curve Lines [id:da33804769013633407] 
\draw    (402,94) .. controls (398,97) and (401,122) .. (408,121) ;
%Straight Lines [id:da20175583052016488] 
\draw    (408,128) -- (436,133) ;
%Straight Lines [id:da6813701375423575] 
\draw    (403,152) -- (432,161) ;
%Curve Lines [id:da9212076185731439] 
\draw    (408,128) .. controls (415,130) and (407,154) .. (403,152) ;
%Curve Lines [id:da3044613671863112] 
\draw    (403,152) .. controls (397,148) and (405,129) .. (408,128) ;
%Straight Lines [id:da7898330453873852] 
\draw    (151,229) -- (150,303) ;
%Curve Lines [id:da05820334495059909] 
\draw    (151,229) .. controls (99,237) and (100,294) .. (150,303) ;
%Curve Lines [id:da5827440449776496] 
\draw    (151,229) .. controls (205,241) and (195,296) .. (150,303) ;
%Shape: Circle [id:dp5422891661204344] 
\draw  [fill={rgb, 255:red, 0; green, 0; blue, 0 }  ,fill opacity=1 ] (190,267) .. controls (190,265.9) and (189.1,265) .. (188,265) .. controls (186.9,265) and (186,265.9) .. (186,267) .. controls (186,268.1) and (186.9,269) .. (188,269) .. controls (189.1,269) and (190,268.1) .. (190,267) -- cycle ;

%Straight Lines [id:da6344530945797584] 
\draw    (229,269) -- (279,269) ;
\draw [shift={(281,269)}, rotate = 180] [color={rgb, 255:red, 0; green, 0; blue, 0 }  ][line width=0.75]    (10.93,-3.29) .. controls (6.95,-1.4) and (3.31,-0.3) .. (0,0) .. controls (3.31,0.3) and (6.95,1.4) .. (10.93,3.29)   ;
%Curve Lines [id:da06143153267812518] 
\draw    (312,233) .. controls (337,239) and (397,241) .. (432,234) ;
%Curve Lines [id:da8804087122455566] 
\draw    (314,301) .. controls (347,295) and (404,296) .. (434,302) ;
%Curve Lines [id:da1761163734729425] 
\draw    (358,239) .. controls (341,237) and (336,300) .. (357,297) ;
%Curve Lines [id:da4826873283785117] 
\draw  [dash pattern={on 4.5pt off 4.5pt}]  (358,239) .. controls (371,235) and (389,257) .. (396,268) ;
%Curve Lines [id:da3625334765863173] 
\draw  [dash pattern={on 4.5pt off 4.5pt}]  (357,297) .. controls (371,296) and (386,282) .. (396,268) ;
%Curve Lines [id:da3715654038259646] 
\draw    (384,239) .. controls (374,235) and (372,300) .. (380,297) ;
%Curve Lines [id:da1919685448478372] 
\draw    (396,268) .. controls (404,259) and (407,235) .. (400,238) ;
%Curve Lines [id:da9967496576982645] 
\draw  [dash pattern={on 4.5pt off 4.5pt}]  (396,268) .. controls (396,254) and (395,239) .. (400,238) ;
%Curve Lines [id:da25740858751257034] 
\draw  [dash pattern={on 4.5pt off 4.5pt}]  (384,239) .. controls (390,239) and (394,256) .. (396,268) ;
%Curve Lines [id:da42181829064511833] 
\draw  [dash pattern={on 4.5pt off 4.5pt}]  (380,297) .. controls (387,298) and (390,282) .. (396,268) ;
%Curve Lines [id:da3970000943946511] 
\draw    (402,298) .. controls (406,296) and (403,279) .. (396,268) ;
%Curve Lines [id:da6675131279588005] 
\draw  [dash pattern={on 4.5pt off 4.5pt}]  (396,268) .. controls (397,284) and (395,299) .. (402,298) ;
%Straight Lines [id:da11494925391558974] 
\draw    (396,268) -- (426,266) ;
%Straight Lines [id:da6377368863788953] 
\draw    (396,268) -- (430,275) ;
%Straight Lines [id:da11326037727365068] 
\draw    (370,181) -- (370,214) ;
\draw [shift={(370,216)}, rotate = 270] [color={rgb, 255:red, 0; green, 0; blue, 0 }  ][line width=0.75]    (10.93,-3.29) .. controls (6.95,-1.4) and (3.31,-0.3) .. (0,0) .. controls (3.31,0.3) and (6.95,1.4) .. (10.93,3.29)   ;

% Text Node
\draw (113,108.4) node [anchor=north west][inner sep=0.75pt]    {$N_{0}$};
% Text Node
\draw (161,127.4) node [anchor=north west][inner sep=0.75pt]    {$N_{1}$};
% Text Node
\draw (160,93.4) node [anchor=north west][inner sep=0.75pt]    {$N_{2}$};
% Text Node
\draw (141,156.4) node [anchor=north west][inner sep=0.75pt]    {$\alpha $};
% Text Node
\draw (142,65.4) node [anchor=north west][inner sep=0.75pt]    {$\beta $};
% Text Node
\draw (240,100) node [anchor=north west][inner sep=0.75pt]   [align=left] {cut};
% Text Node
\draw (358,75.4) node [anchor=north west][inner sep=0.75pt]    {$\gamma _{+}$};
% Text Node
\draw (396,72.4) node [anchor=north west][inner sep=0.75pt]    {$\beta _{-}$};
% Text Node
\draw (386,153.4) node [anchor=north west][inner sep=0.75pt]    {$\alpha _{-}$};
% Text Node
\draw (239,240) node [anchor=north west][inner sep=0.75pt]   [align=left] {glue};
% Text Node
\draw (361,255.4) node [anchor=north west][inner sep=0.75pt]    {$\eta $};
% Text Node
\draw (320,254.4) node [anchor=north west][inner sep=0.75pt]    {$\gamma _{+}$};
% Text Node
\draw (405,240.4) node [anchor=north west][inner sep=0.75pt]    {$\beta _{-}$};
% Text Node
\draw (405,276.4) node [anchor=north west][inner sep=0.75pt]    {$\alpha _{-}$};
% Text Node
\draw (383,186) node [anchor=north west][inner sep=0.75pt]   [align=left] {glue};
% Text Node
\draw (120,258.4) node [anchor=north west][inner sep=0.75pt]    {$D_{+}$};
% Text Node
\draw (160,259.4) node [anchor=north west][inner sep=0.75pt]    {$D_{-}$};
% Text Node
\draw (153,237.4) node [anchor=north west][inner sep=0.75pt]    {$\eta $};
% Text Node
\draw (98,222.4) node [anchor=north west][inner sep=0.75pt]    {$\gamma '_{+}$};
% Text Node
\draw (175,293.4) node [anchor=north west][inner sep=0.75pt]    {$\alpha '_{-}$};
% Text Node
\draw (177,214.4) node [anchor=north west][inner sep=0.75pt]    {$\beta '_{-}$};
% Text Node
\draw (143,307.4) node [anchor=north west][inner sep=0.75pt]    {$P_{1}$};
% Text Node
\draw (194,259.4) node [anchor=north west][inner sep=0.75pt]    {$P_{2}$};
% Text Node
\draw (142,205.4) node [anchor=north west][inner sep=0.75pt]    {$P_{3}$};

\end{tikzpicture}

    \caption{Merging folding curves. Identify $\gamma_+$ with $\gamma_+'$, $\alpha_-$ with $\alpha_-'$, $\beta_-$ with $\beta_-'$ (we have $\gamma_-'=\alpha_-'*\beta_-'$). The three points $P_1,P_2,P_3$ are identified.}
    \label{fig: sur-merge}
\end{figure}

Cut $\Sigma$ along $\alpha,\beta$. The resulting surface $\Sigma^{cut}$ has three boundary curves $\gamma_+,\alpha_-,\beta_-$, where $\gamma_+$ is a boundary of $N_0$, $\alpha_-$ is a boundary of $N_1$, and $\beta_-$ is a boundary of $N_2$. We still denote by $f$ its restriction on $\Sigma^{cut}$ (the universal covering of $\Sigma^{cut}$ is a subsurface of $\widetilde{\Sigma}$, so we can restrict $f$). Let $\Sigma_1=\Sigma^{cut}\cup D_+\cup D_-$, identifying $\gamma_+$ with $\gamma_+'$ via $f$, identifying $\eta_+$ with $\eta_-$ via the identity, and attaching $\alpha_-,\beta_-$ to $\gamma_-'$ via $f$ (see Figure \ref{fig: sur-merge}).

\begin{lemma}\label{lem-3m2}
    $f$ can be extended to a developing map $f_1$ on $\Sigma_1$ by mapping $D_{\pm}$ to $D$ via the identity. Moreover, there is a homeomorphism between $\Sigma_1$ and $\Sigma$ such that $f_1$ has the same monodromy as $f$ under this identification of $\Sigma_1$ and $\Sigma$.
\end{lemma}

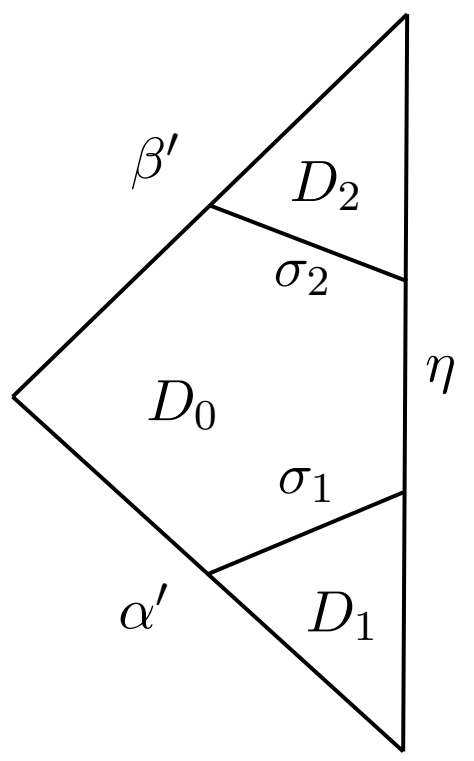
\begin{figure}
    \centering

\tikzset{every picture/.style={line width=0.75pt}} %set default line width to 0.75pt        

\begin{tikzpicture}[x=0.75pt,y=0.75pt,yscale=-1,xscale=1]
%uncomment if require: \path (0,339); %set diagram left start at 0, and has height of 339

%Straight Lines [id:da33555034549439156] 
\draw    (352,79) -- (351,264) ;
%Straight Lines [id:da5131468995619916] 
\draw    (253,175) -- (352,79) ;
%Straight Lines [id:da3078087322548926] 
\draw    (253,175) -- (351,264) ;
%Straight Lines [id:da22692147304419652] 
\draw    (302.5,127) -- (352,146) ;
%Straight Lines [id:da8259132625789573] 
\draw    (302,219.5) -- (351,199) ;

% Text Node
\draw (278,220.4) node [anchor=north west][inner sep=0.75pt]    {$\alpha '$};
% Text Node
\draw (281,107.4) node [anchor=north west][inner sep=0.75pt]    {$\beta '$};
% Text Node
\draw (355,163.4) node [anchor=north west][inner sep=0.75pt]    {$\eta $};
% Text Node
\draw (318,191.4) node [anchor=north west][inner sep=0.75pt]    {$\sigma _{1}$};
% Text Node
\draw (317,139.4) node [anchor=north west][inner sep=0.75pt]    {$\sigma _{2}$};
% Text Node
\draw (285,169.4) node [anchor=north west][inner sep=0.75pt]    {$D_{0}$};
% Text Node
\draw (325,222.4) node [anchor=north west][inner sep=0.75pt]    {$D_{1}$};
% Text Node
\draw (321,114.4) node [anchor=north west][inner sep=0.75pt]    {$D_{2}$};

\end{tikzpicture}

    \caption{$\sigma_1,\sigma_2$ divide $D$ into three discs $D_0,D_1,D_2$. Attach a pair of $D_0,D_1,D_2$ in order.}
    \label{fig: sur-merge2}
\end{figure}

\begin{proof}
    This is the result of attaching three pairs of discs, one after another. Let $\alpha'=f(\widetilde{\alpha}),\beta'=f(\widetilde{\beta})$. Take two disjoint simple curves $\sigma_1,\sigma_2$ in $D$, such that $\sigma_1$ joins an interior point of $\alpha
    '$ with an interior point of $\eta$, while $\sigma_2$ joins an interior point of $\beta'$ with an interior point of $\eta$. $\sigma_1,\sigma_2$ divide $D$ into three discs $D_0,D_1,D_2$ (see Figure \ref{fig: sur-merge2}). Then $\Sigma_1,f_1$ is obtained by first attaching a pair of $D_0$, then attaching a pair of $D_1$, and finally attaching a pair of $D_2$.
\end{proof}

By a local analysis one proves the following lemma.

\begin{lemma}\label{lem-3m3}
    Suppose that $f$ is orientation preserving on $N_0$, while orientation reversing on $N_1,N_2$. Then $\eta$ is a folding curve of $f_1$, and there is a neighborhood $N$ of $D_+\cup D_-\subset\Sigma_1$ such that $f$ is branched on $N\backslash\eta$.
\end{lemma}

\subsection{Model 3}

\paragraph{Model 3:}\textit{Both $C_1,C_2$ are $1$-petal flowers, and they admit $5\epsilon$-extensions that are disjoint.}

Model 3 splits into three cases, according to the angle between $\alpha_1$ and $\alpha_2$, where $\alpha_i$ is the petal of $C_i$. The surgeries for different cases are similar. In fact the three surgeries are topologically the same, but we need to adjust the geometric positions of the cut and paste loci. We will see that only the last case actually uses of the extensions of $C_i$ (which provides room for the surgery).

Without loss of generality we may assume that $f$ is orientation preserving on $\Sigma_0$. Orient $\alpha_i$ so that $C_i$ lies on the right of $\alpha_i$.

\paragraph{Case 1} $\measuredangle(\overleftarrow{\alpha_1},\overrightarrow{\alpha_2})=2k\pi+\theta$ or $\measuredangle(\overleftarrow{\alpha_2},\overrightarrow{\alpha_1})=2k\pi+\theta$, for some $k\in\Z,\ \theta\in(0,\pi)$.

Without loss of generality we may assume that $\measuredangle(\overleftarrow{\alpha_1},\overrightarrow{\alpha_2})=2k\pi+\theta,\ k\in\Z,\ \theta\in(0,\pi)$. Let $\alpha=\alpha_1*\alpha_2$ and let $\widetilde{\alpha}=\widetilde{\alpha}_1*\widetilde{\alpha}_2\subset\widetilde{\Sigma}$ be a lift of $\alpha$. Let $\alpha_1'=f(\widetilde{\alpha}_1),\alpha_2'=f(\widetilde{\alpha}_2)$. Then $\measuredangle(\overleftarrow{\alpha_1'},\overrightarrow{\alpha_2'})=\theta$.
Let $T$ be the geodesic triangle spanned by $\alpha_1',\alpha_2'$, and denote its third side by $\beta$. Then $T$ lies on the right of $\alpha_1',\alpha_2'$.

The surgery goes as follows. First remove the pair of annuli $C_1,C_1'$ and the pair of annuli $C_2,C_2'$ (recall that $C_i'=\exp_{\gamma_i}^-(C_i^0)\subset\Sigma_i$). Denote the resulting surface by $\Sigma'$, and denote the image of $\alpha_i$ in $\Sigma'$ still by $\alpha_i$. We have $\alpha_1\cap\alpha_2=\{P\}$ (recall that $C_1\cap C_2=\varnothing$), and $i(\alpha_1,\alpha_2)=0$ (since $\alpha_i\sim\gamma_i$ and $i(\gamma_1,\gamma_2)=0$). So we are in the right position to apply the surgery of merging folding curves. Cut $\Sigma'$ along $\alpha_1,\alpha_2$ and attach a pair of $T$. Denote the resulting surface by $\Sigma''$. By Lemma \ref{lem-3m2}, $f$ can be extended to a developing map $f_1$ on $\Sigma''$, which has the same monodromy as $f$, after identifying $\Sigma''$ with $\Sigma$ via some homeomorphism. Note that we can take $N_0\subset\Sigma_0\backslash(C_1\cup C_2),\ N_1\subset\Sigma_1\backslash C_1',\ N_2\subset\Sigma_2\backslash C_2'$, in particular $f$ is orientation preserving on $N_0$, while orientation reversing on $N_1,N_2$. Therefore by Lemma \ref{lem-3m3}, $f_1$ is branched away from $n-1$ folding curves, with $\beta$ being one of these folding curves. The area of $f_1$ satisfies $A(f_1)\leq A(f)+2\mathrm{Area}(T)\leq 4\epsilon$. Therefore, by the example at the end of Subsection \ref{subsec-3m-some sur}, by attaching a pair of annuli at $\beta$, we get a $5\epsilon$-excellent developing map with $n-1$ folding curves. This finishes the discussion of Case 1.

\begin{claim}
    $\measuredangle(\overleftarrow{\alpha_1},\overrightarrow{\alpha_2})+\measuredangle(\overleftarrow{\alpha_2},\overrightarrow{\alpha_1})\equiv\theta\mod2\pi$ for some $|\theta|<4\epsilon$.
\end{claim}

\begin{proof}
    $\measuredangle(\overleftarrow{\alpha_1},\overrightarrow{\alpha_2})+\measuredangle(\overleftarrow{\alpha_2},\overrightarrow{\alpha_1})\equiv\measuredangle(\overleftarrow{\alpha_1},\overrightarrow{\alpha_1})+\measuredangle(\overleftarrow{\alpha_2},\overrightarrow{\alpha_2})\mod2\pi$. By Lemma \ref{lem-h1}, we have $\measuredangle(\overleftarrow{\alpha_i},\overrightarrow{\alpha_i})\equiv\pi+\theta_i\mod2\pi$ for some $|\theta_i|<2\epsilon\ (i=1,2)$. This proves the claim.
\end{proof}

By the claim, if neither $\measuredangle(\overleftarrow{\alpha_1},\overrightarrow{\alpha_2})$ nor $\measuredangle(\overleftarrow{\alpha_2},\overrightarrow{\alpha_1})$ is of the form $2k\pi+\theta$ for some $k\in\Z,\ \theta\in(0,\pi)$, then
\begin{itemize}
    \item either $\measuredangle(\overleftarrow{\alpha_1},\overrightarrow{\alpha_2})\equiv\pi+\theta_1,\ \measuredangle(\overleftarrow{\alpha_2},\overrightarrow{\alpha_1})\equiv\pi+\theta_2\mod2\pi$ for some $\theta_1,\theta_2\in[0,4\epsilon)$;
    \item or $\measuredangle(\overleftarrow{\alpha_1},\overrightarrow{\alpha_2})\equiv-\theta_1,\ \measuredangle(\overleftarrow{\alpha_2},\overrightarrow{\alpha_1})\equiv-\theta_2\mod2\pi$ for some $\theta_1,\theta_2\in[0,4\epsilon)$.
\end{itemize}

\paragraph{Case 2} $\measuredangle(\overleftarrow{\alpha_1},\overrightarrow{\alpha_2})\equiv\pi+\theta\mod2\pi$ for some $\theta\in[0,4\epsilon)$.

The solution for this case is almost the same as in Case 1. The only difference is that instead of attaching a pair of geodesic triangles, we attach a pair of concave quadrilaterals.

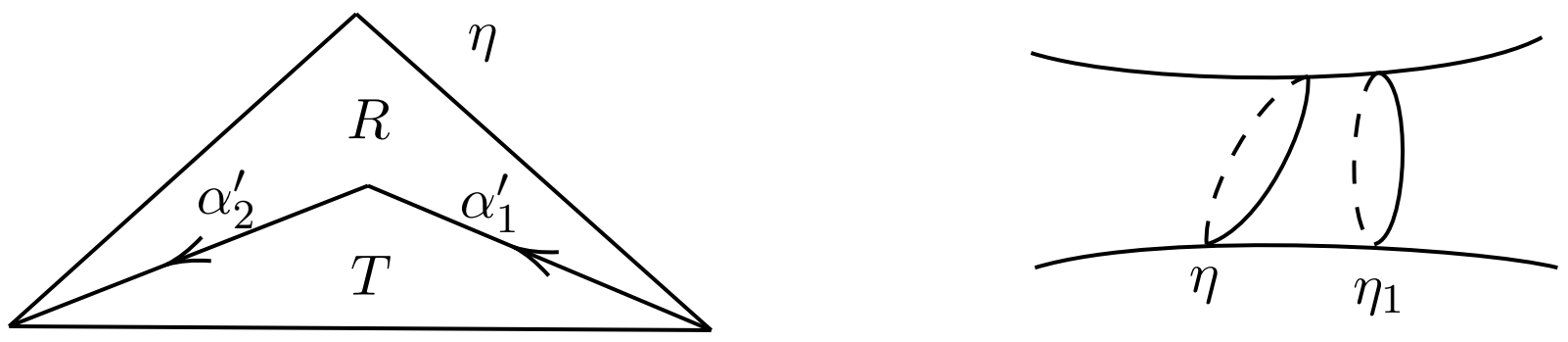
\begin{figure}
    \centering

\tikzset{every picture/.style={line width=0.75pt}} %set default line width to 0.75pt        

\begin{tikzpicture}[x=0.75pt,y=0.75pt,yscale=-1,xscale=1]
%uncomment if require: \path (0,339); %set diagram left start at 0, and has height of 339

%Straight Lines [id:da07490395525151616] 
\draw    (181,143) -- (269,180) ;
\draw [shift={(218.55,158.79)}, rotate = 22.8] [color={rgb, 255:red, 0; green, 0; blue, 0 }  ][line width=0.75]    (10.93,-3.29) .. controls (6.95,-1.4) and (3.31,-0.3) .. (0,0) .. controls (3.31,0.3) and (6.95,1.4) .. (10.93,3.29)   ;
%Straight Lines [id:da7176099823189929] 
\draw    (181,143) -- (89,179) ;
\draw [shift={(129.41,163.19)}, rotate = 338.63] [color={rgb, 255:red, 0; green, 0; blue, 0 }  ][line width=0.75]    (10.93,-3.29) .. controls (6.95,-1.4) and (3.31,-0.3) .. (0,0) .. controls (3.31,0.3) and (6.95,1.4) .. (10.93,3.29)   ;
%Straight Lines [id:da6549025780249812] 
\draw    (89,179) -- (269,180) ;
%Straight Lines [id:da7298747507203596] 
\draw    (178,99) -- (89,179) ;
%Straight Lines [id:da07479202800545937] 
\draw    (178,99) -- (269,180) ;

%Curve Lines [id:da6676937582041809] 
\draw    (351,109) .. controls (380,118) and (458,118) .. (482,105) ;
%Curve Lines [id:da3510186699366217] 
\draw    (352,164) .. controls (385,154) and (466,159) .. (486,164) ;
%Curve Lines [id:da37309842995380593] 
\draw    (422,115) .. controls (423,127) and (410,154) .. (396,158) ;
%Curve Lines [id:da12429252381988198] 
\draw  [dash pattern={on 4.5pt off 4.5pt}]  (396,158) .. controls (395,145) and (410,118) .. (422,115) ;
%Curve Lines [id:da5612787777365491] 
\draw    (440,114) .. controls (449,116) and (448,156) .. (439,158) ;
%Curve Lines [id:da4715697436070636] 
\draw  [dash pattern={on 4.5pt off 4.5pt}]  (440,114) .. controls (433,117) and (431,154) .. (439,158) ;

% Text Node
\draw (175,159.4) node [anchor=north west][inner sep=0.75pt]    {$T$};
% Text Node
\draw (135.5,137.4) node [anchor=north west][inner sep=0.75pt]    {$\alpha' _{2}$};
% Text Node
\draw (203,138.4) node [anchor=north west][inner sep=0.75pt]    {$\alpha' _{1}$};
% Text Node
\draw (174,119.4) node [anchor=north west][inner sep=0.75pt]    {$R$};
% Text Node
\draw (205,100.4) node [anchor=north west][inner sep=0.75pt]    {$\eta $};
% Text Node
\draw (390,162.4) node [anchor=north west][inner sep=0.75pt]    {$\eta $};
% Text Node
\draw (432,165.4) node [anchor=north west][inner sep=0.75pt]    {$\eta _{1}$};

\end{tikzpicture}

    \caption{Attach a pair of $R$, and then further attach a pair of annuli to turn $\eta$ into an equidistant curve.}
    \label{fig: m3case2}
\end{figure}

Let $\alpha_1',\alpha_2'$ be as in Case 1. We have $\measuredangle(\overleftarrow{\alpha_1'},\overrightarrow{\alpha_2'})=\pi+\theta$. Let $T$ be the (possibly degenerate) geodesic triangle spanned by $\alpha_1',\alpha_2'$, and let $\beta$ be the third side of $T$. Now $T$ is either degenerate (in which case we set $T=\varnothing$) or lies on the left of $\alpha_1',\alpha_2'$. The angle of $T$ opposite to $\beta$ is equal to $\pi-\theta$. Therefore the angles at the vertices of $\beta$ are less than $4\epsilon$. Let $T_1$ be the $4\epsilon$-cap based on $\beta$, lying on the same side of $\beta$ as $T$. Then $T\subset T_1$. Let $R=\overline{T_1\backslash T}$ (see Figure \ref{fig: m3case2}). Then $R$ satisfies the conditions for the surgery of merging folding curves. We can then proceed in exactly the same way as in Case 1, with $T$ replaced by $R$. Let $\eta=\overline{\partial R\backslash(\alpha'_1\cup\alpha'_2)}$. Then we arrive at a developing map $f_1$, branched away from $n-1$ folding curves, with $\eta$ being one of its folding curves. Note that $\eta$ projects onto a simple closed curve in $\H/\rho(\eta)$ (by Lemma \ref{lem-h5}). Therefore we can attach a pair of annulus to turn it into an equidistant curve of length $\leq5\epsilon$.

\paragraph{Case 3} $\measuredangle(\overleftarrow{\alpha_1},\overrightarrow{\alpha_2})\equiv-\theta_1,\ \measuredangle(\overleftarrow{\alpha_2},\overrightarrow{\alpha_1})\equiv-\theta_2\mod2\pi$ for some $\theta_1,\theta_2\in[0,4\epsilon)$.

By symmetry we may assume $l_{\alpha_1}\geq l_{\alpha_2}$. For the moment assume $\theta_1\neq0$. Let $\alpha_1',\alpha_2'$ be as in Case 1. We have $\measuredangle(\overrightarrow{\alpha_2'},\overleftarrow{\alpha_1'})=\theta_1$. Let $T$ be the geodesic triangle spanned by $\alpha_1',\alpha_2'$.

\begin{claim}\label{claim-mod3}
    The angles of $T$ at the two vertices of $\alpha_1'$ are both less than $4\epsilon$.
\end{claim}

\begin{proof}
    Let $\widetilde{P}\in\widetilde{\Sigma}$ be the common endpoint of $\widetilde{\alpha}_1$ and $\widetilde{\alpha}_2$. Let $\widetilde{\alpha}'_1$ be the lift of $\alpha_1$ that starts at $\widetilde{P}$, and let $\widetilde{\alpha}'_2$ be the lift of $\alpha_2$ that terminates at $\widetilde{P}$. Let $\alpha_1''=f(\widetilde{\alpha}'_1),\alpha_2''=f(\widetilde{\alpha}'_2)$ (see Figure \ref{fig: m3case31}). Then $\measuredangle(\overrightarrow{\alpha_1''},\overleftarrow{\alpha_2''})=\theta_2$. The claim follows from Lemma \ref{lem-h4}.
\end{proof}

\begin{figure}
    \centering

\tikzset{every picture/.style={line width=0.75pt}} %set default line width to 0.75pt        

\begin{tikzpicture}[x=0.75pt,y=0.75pt,yscale=-1,xscale=1]
%uncomment if require: \path (0,339); %set diagram left start at 0, and has height of 339

%Straight Lines [id:da47988428970891517] 
\draw    (181,91) -- (253,115) ;
\draw [shift={(210.36,100.79)}, rotate = 18.43] [color={rgb, 255:red, 0; green, 0; blue, 0 }  ][line width=0.75]    (10.93,-3.29) .. controls (6.95,-1.4) and (3.31,-0.3) .. (0,0) .. controls (3.31,0.3) and (6.95,1.4) .. (10.93,3.29)   ;
%Straight Lines [id:da5745381944191227] 
\draw    (132,129) -- (253,115) ;
\draw [shift={(198.46,121.31)}, rotate = 173.4] [color={rgb, 255:red, 0; green, 0; blue, 0 }  ][line width=0.75]    (10.93,-3.29) .. controls (6.95,-1.4) and (3.31,-0.3) .. (0,0) .. controls (3.31,0.3) and (6.95,1.4) .. (10.93,3.29)   ;
%Straight Lines [id:da8316459934532194] 
\draw    (253,115) -- (323,101) ;
\draw [shift={(281.14,109.37)}, rotate = 348.69] [color={rgb, 255:red, 0; green, 0; blue, 0 }  ][line width=0.75]    (10.93,-3.29) .. controls (6.95,-1.4) and (3.31,-0.3) .. (0,0) .. controls (3.31,0.3) and (6.95,1.4) .. (10.93,3.29)   ;
%Straight Lines [id:da7444243453159844] 
\draw    (253,115) -- (357,136) ;
\draw [shift={(310.88,126.69)}, rotate = 191.42] [color={rgb, 255:red, 0; green, 0; blue, 0 }  ][line width=0.75]    (10.93,-3.29) .. controls (6.95,-1.4) and (3.31,-0.3) .. (0,0) .. controls (3.31,0.3) and (6.95,1.4) .. (10.93,3.29)   ;
%Straight Lines [id:da4188102027461982] 
\draw    (181,91) -- (132,129) ;

% Text Node
\draw (177,100.4) node [anchor=north west][inner sep=0.75pt]    {$T$};
% Text Node
\draw (291,130.4) node [anchor=north west][inner sep=0.75pt]    {$\alpha ''_{1}$};
% Text Node
\draw (284,83.4) node [anchor=north west][inner sep=0.75pt]    {$\alpha ''_{2}$};
% Text Node
\draw (191,128.4) node [anchor=north west][inner sep=0.75pt]    {$\alpha '_{1}$};
% Text Node
\draw (207,80.4) node [anchor=north west][inner sep=0.75pt]    {$\alpha '_{2}$};

\end{tikzpicture}

    \caption{Claim \ref{claim-mod3}}
    \label{fig: m3case31}
\end{figure}

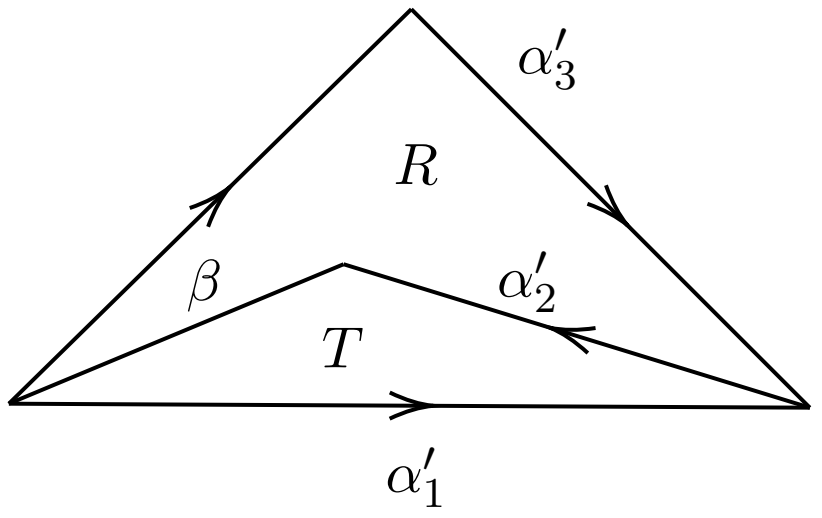
\begin{figure}
    \centering

\tikzset{every picture/.style={line width=0.75pt}} %set default line width to 0.75pt        

\begin{tikzpicture}[x=0.75pt,y=0.75pt,yscale=-1,xscale=1]
%uncomment if require: \path (0,339); %set diagram left start at 0, and has height of 339

%Straight Lines [id:da15170791387682525] 
\draw    (281,91) -- (381,191) ;
\draw [shift={(335.24,145.24)}, rotate = 225] [color={rgb, 255:red, 0; green, 0; blue, 0 }  ][line width=0.75]    (10.93,-3.29) .. controls (6.95,-1.4) and (3.31,-0.3) .. (0,0) .. controls (3.31,0.3) and (6.95,1.4) .. (10.93,3.29)   ;
%Straight Lines [id:da44640199087124066] 
\draw    (281,91) -- (180,190) ;
\draw [shift={(235.5,135.6)}, rotate = 135.57] [color={rgb, 255:red, 0; green, 0; blue, 0 }  ][line width=0.75]    (10.93,-3.29) .. controls (6.95,-1.4) and (3.31,-0.3) .. (0,0) .. controls (3.31,0.3) and (6.95,1.4) .. (10.93,3.29)   ;
%Straight Lines [id:da357709635926557] 
\draw    (180,190) -- (381,191) ;
\draw [shift={(286.5,190.53)}, rotate = 180.29] [color={rgb, 255:red, 0; green, 0; blue, 0 }  ][line width=0.75]    (10.93,-3.29) .. controls (6.95,-1.4) and (3.31,-0.3) .. (0,0) .. controls (3.31,0.3) and (6.95,1.4) .. (10.93,3.29)   ;
%Straight Lines [id:da43142726858666014] 
\draw    (264,155) -- (180,190) ;
%Straight Lines [id:da7682073499769343] 
\draw    (264,155) -- (381,191) ;
\draw [shift={(315.81,170.94)}, rotate = 17.1] [color={rgb, 255:red, 0; green, 0; blue, 0 }  ][line width=0.75]    (10.93,-3.29) .. controls (6.95,-1.4) and (3.31,-0.3) .. (0,0) .. controls (3.31,0.3) and (6.95,1.4) .. (10.93,3.29)   ;

% Text Node
\draw (275,123.4) node [anchor=north west][inner sep=0.75pt]    {$R$};
% Text Node
\draw (223,152.4) node [anchor=north west][inner sep=0.75pt]    {$\beta $};
% Text Node
\draw (273,199.4) node [anchor=north west][inner sep=0.75pt]    {$\alpha '_{1}$};
% Text Node
\draw (301,150.4) node [anchor=north west][inner sep=0.75pt]    {$\alpha '_{2}$};
% Text Node
\draw (306,94.4) node [anchor=north west][inner sep=0.75pt]    {$\alpha '_{3}$};
% Text Node
\draw (257,169.4) node [anchor=north west][inner sep=0.75pt]    {$T$};

\end{tikzpicture}

    \caption{$R$ lies on the right of $\alpha_3'*\alpha_2'$.}
    \label{fig: m3case32}
\end{figure}

Recall that $C_1$ admits a $5\epsilon$-extension (at $\alpha_1$) that is disjoint from $C_2$. Let $E$ be the $4\epsilon$-extension of $C_1$, and let $T_1=E\backslash(C_1\cup\alpha_1)$ be the $4\epsilon$-cap of $C_1$ at $\alpha_1$. Then $\alpha_3:=\partial E$ is a simple loop based at $P$, and we have $\alpha_2\cap\alpha_3={P}$. Orient $\alpha_3$ so that it is homotopic to $\alpha_1$. Since $\alpha_3\sim\alpha_1$, we have $i(\alpha_2,\alpha_3)=0$. Extend the lift $\widetilde{\alpha}_1$ continuously to a lift $\widetilde{\bar{T}}_1$ and let $T_1'=f(\widetilde{T}_1),\alpha_3'=f(\widetilde{\alpha}_3)$. Then $T_1'$ is the $4\epsilon$-cap based on $\alpha_1'$, lying on the same side of $\alpha_1'$ as $T$, and $\alpha_3'=\overline{\partial T_1'\backslash\alpha_1'}$. By Claim \ref{claim-mod3}, we have $T\subset T_1'$. Let $R=\overline{T_1'\backslash T}$. Then $\alpha_3,\alpha_2,R$ satisfies the conditions for the surgery of merging folding curves (see Figure \ref{fig: m3case32}).

Let $E'$ be the $4\epsilon$-extension of $C_1'$ (recall that $C_1'=\exp_{\gamma_1}^-(C_1^0)\subset\Sigma_1$). Then $f$ embeds both $\overline{E}$ and $\overline{E'}$ into $\H/\rho(\gamma_1)$, with the same image. Remove the pair of annuli $E,E'$ and identify $\partial E'$ with $\partial E=\alpha_3$. Then remove the pair of annuli $C_2,C_2'$ and identify $\partial C_2'$ with $\partial C_2=\alpha_2$. Denote the images of $\alpha_3,\alpha_2$ in the new surface still by $\alpha_3,\alpha_2$. Now cut along $\alpha_3,\alpha_2$, and attach a pair of $R$. By Lemma \ref{lem-3m2} and Lemma \ref{lem-3m3}, $f$ can be extended to a developing map $f_1$ on the new surface, which is branched away from $n-1$ folding curve. The new folding curve is $\beta=\overline{\partial R\backslash(\alpha_2'\cup\alpha_3')}$, which is a simple geodesic loop. We then further attach a pair of annuli at $\beta$ to turn it into an equidistant curve of length $\leq5\epsilon$.

If $\theta_1=0$, then $l_{\alpha_1}\neq l_{\alpha_2}$, otherwise $\rho(\alpha_1*\alpha_2)$ would be an elliptic element or the identity. Therefore $l_{\alpha_1}>l_{\alpha_2}$. The above discussion goes through without change, with $T$ replaced by $\varnothing$.

\section{Nice triples and simple geodesic loops}

The remaining task is to prove Proposition \ref{prop-vs2}. The branch point $P$ in the standard triple may has high multiplicity, so the local structure near $P$ may still be complicated. The strategy is to reduce the multiplicity of $P$, possibly with the cost of reducing $r_P$ as well, with the hope that the local structure gets simpler. We show that either the multiplicity of $P$ can be further reduced, or we get a model. In the course of the reduction procedure we may leave the class of standard developing maps, so we introduce a notion slightly more general than standard.

\begin{definition}
    Let $S$ be a good surface. A nice triple $(P,\gamma_1,\gamma_2)$ of $S$ consists of a branch point $P$ and two boundary curves $\gamma_1,\gamma_2$ such that $d(P,\gamma_i)=r_{\gamma_i}\ (i=1,2)$.
\end{definition}

\begin{definition}
    A nice triple $(P,\gamma_1,\gamma_2)$ is called reducible if there exists $S'\in[S]$ and a nice triple $(P',\gamma_1,\gamma_2)$ of $S'$ such that $r_{P'}\geq\frac{1}{10}r_P,\ r_{\gamma_i}(S')\leq r_{\gamma_i}+10\epsilon\ (i=1,2),\ m(P')<m(P)$, where $m(P)$ denotes the multiplicity of the branch point $P$.
\end{definition}

Throughout this and the next section let $\epsilon\leq100\epsilon_0$. Let $\epsilon_1=1000\epsilon$. Proposition \ref{prop-vs2} follows from the proposition below.

\begin{proposition}\label{prop-rm}
    Let $S$ be an $\epsilon$-good surface with an irreducible nice triple $(P,\gamma_1,\gamma_2)$ such that $r_P\geq \epsilon_1$ and $r_{\gamma_i}\leq2\delta_1\ (i=1,2)$, then there is a model $S'\in[S]$.
\end{proposition}

\paragraph{Proof of Proposition \ref{prop-vs2}}

Let $f$ be a very standard $\epsilon$-excellent developing map with $n$ folding curves. Then there is a component $S_0=(\Sigma_0,g_f)$ of $\Sigma\backslash\mathrm{Fold}(f)$ with a standard triple $(P,\gamma_1,\gamma_2)$ such that $r_{\gamma_i}(S_i')\geq 10\delta_1\ (i=1,2)$, where $S_i'$ is the component of $\Sigma\backslash\mathrm{Fold}(f)$ other than $S_0$ which is adjacent to $\gamma_i$. If $(P,\gamma_1,\gamma_2)$ is reducible, then by a homotopy of $f|_{\Sigma_0}$, we get a developing map $f_1$ with a nice triple $(P_1,\gamma_1,\gamma_2)$ for $S_1=(\Sigma_0,g_{f_1})$ such that $r_{P_1}(S_1)\geq\frac{1}{10}r_P(S_0),\ r_{\gamma_i}(S_1)\leq r_{\gamma_i}(S)+10\epsilon\ (i=1,2)$ and $m(P_1)<m(P)$. If $(P_1,\gamma_1,\gamma_2)$ is reducible we continue this process, until we reach some $f_k$ with an irreducible nice triple $(P_k,\gamma_1,\gamma_2)$. We have $k\leq m(P)\leq-\chi(\Sigma)$. Then $r_{P_k}(S_k)\geq10^{\chi(\Sigma)}r_P(S_0)\geq 10^{\chi(\Sigma)}\delta_0,\ r_{\gamma_i}(S_k)\leq r_{\gamma_i}(S_0)-10\chi(\Sigma)\epsilon\leq 2\delta_1\ (i=1,2)$. Therefore $S_k$ satisfies the conditions of Proposition \ref{prop-rm}, so by a further homotopy of $f_k|_{\Sigma_0}$ we get a model. Applying Proposition \ref{prop-model}, we finish the proof.

The rest of the paper is devoted to the proof of Proposition \ref{prop-rm}.

This section serves as a preparation for the proof. We study simple geodesic loops based at $P$, a useful tool for analyzing the local geometry around $P$. If we cut along a simple geodesic loop, we get geodesic loop boundaries. We first introduce a surgery that transform such boundary curves to equidistant curves, so that the theory developed in the previous sections can be applied.

\subsection{Surgery: attaching a $1$-petal flower}

An (abstract) $1$-petal flower is a `cornered' version of a collar.

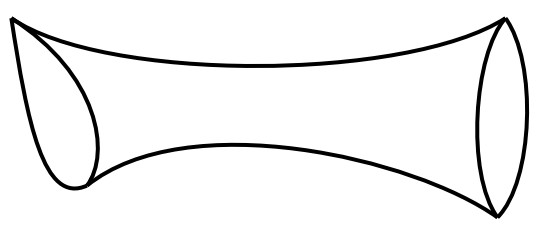
\begin{figure}
    \centering
    \tikzset{every picture/.style={line width=0.75pt}} %set default line width to 0.75pt        

\begin{tikzpicture}[x=0.75pt,y=0.75pt,yscale=-1,xscale=1]

%Curve Lines [id:da02349344286974364] 
\draw    (99,107) .. controls (122,123) and (198,123) .. (223,107) ;
%Curve Lines [id:da5772501920691919] 
\draw    (118,149) .. controls (143,129) and (200,142) .. (221,157) ;
%Curve Lines [id:da19735996251728438] 
\draw    (223,107) .. controls (215,116) and (213,146) .. (221,157) ;
%Curve Lines [id:da025672822158196773] 
\draw    (223,107) .. controls (231,118) and (230,148) .. (221,157) ;
%Curve Lines [id:da7858910836479619] 
\draw    (99,107) .. controls (121,121) and (124,141) .. (118,149) ;
%Curve Lines [id:da9446287107794905] 
\draw    (99,107) .. controls (103,133) and (107,154) .. (118,149) ;

% Text Node
%\draw (96,180) node [anchor=north west][inner sep=0.75pt]   [align=left] {Figure. A $1$-petal flower.};

\end{tikzpicture}
    \caption{A $1$-petal flower}
    \label{fig: 1-petal flower}
\end{figure}

\begin{definition}
    A $1$-petal flower is a (non-complete) hyperbolic annulus such that one of its boundary curves is an equidistant curve, while the other is a simple geodesic loop (Figure 
    \ref{fig: 1-petal flower}). The simple geodesic loop is called the petal.
\end{definition}

The length of a $1$-petal flower is defined to be its length as a flower neighborhood of the equidistant curve.

The surgery of attaching $1$-petal flowers is a simple variant of attaching collars. Let $S$ be a branched hyperbolic surface and let $\gamma$ be a geodesic loop boundary curve, with base point $P$, such that $\rho(\gamma)$ is hyperbolic. Then $\mathrm{dev}_S$ maps $\gamma$ homeomorphically to a simple geodesic loop $\gamma_0$ in $\H/\rho(\gamma)$ (Lemma \ref{lem-h7}). Let $C$ be an annulus in $\H/\rho(\gamma)$ bounded by $\gamma_0$ and an equidistant curve $\gamma_1$, such that outward pointing vectors of $S$ along $\gamma$ are mapped (by $\mathrm{dev}_S$) to inward pointing vectors of $C$ along $\gamma_0$. Then we can attach $C$ to $S$, identifying $\gamma$ with $\gamma_0$ via $\mathrm{dev}_S$. We call this \textit{attaching a $1$-petal flower along $\gamma$}. (We have already encountered this construction in the proof of Lemma \ref{lem-ede-p1}.)

\begin{figure}
    \centering

\tikzset{every picture/.style={line width=0.75pt}} %set default line width to 0.75pt        

\begin{tikzpicture}[x=0.75pt,y=0.75pt,yscale=-1,xscale=1]
%uncomment if require: \path (0,339); %set diagram left start at 0, and has height of 339

%Curve Lines [id:da8946110677670374] 
\draw    (100,107) .. controls (125.55,102) and (138.55,99) .. (159.55,85) ;
%Curve Lines [id:da6030693413019931] 
\draw    (84,145) .. controls (115.55,151) and (138.55,158) .. (157.55,169) ;
%Curve Lines [id:da6599482715053148] 
\draw    (100,107) .. controls (99.55,129) and (89.55,148) .. (84,145) ;
%Curve Lines [id:da15404042650251204] 
\draw    (84,145) .. controls (77.55,144) and (83.55,124) .. (100,107) ;
%Curve Lines [id:da9638582503158045] 
\draw    (215,118) .. controls (244.55,116) and (268.55,114) .. (292.55,96) ;
%Curve Lines [id:da8139182429508349] 
\draw    (215,150) .. controls (244.55,151) and (273.55,160) .. (291.55,172) ;
%Curve Lines [id:da1930132193749634] 
\draw    (215,150) .. controls (207.55,150) and (208.55,120) .. (215,118) ;
%Curve Lines [id:da2400068683797334] 
\draw    (215,150) .. controls (221.55,149) and (222.55,119) .. (215,118) ;
%Curve Lines [id:da2789569738903819] 
\draw    (231.55,151) .. controls (237.55,154) and (248.55,140) .. (252,115) ;
%Curve Lines [id:da07285160886524578] 
\draw  [dash pattern={on 4.5pt off 4.5pt}]  (231.55,151) .. controls (225.55,147) and (235.55,131) .. (252,115) ;
%Straight Lines [id:da5890596901318677] 
\draw    (167,128.5) -- (174.55,128.5)(167,131.5) -- (174.55,131.5) ;
\draw [shift={(182.55,130)}, rotate = 180] [color={rgb, 255:red, 0; green, 0; blue, 0 }  ][line width=0.75]    (10.93,-3.29) .. controls (6.95,-1.4) and (3.31,-0.3) .. (0,0) .. controls (3.31,0.3) and (6.95,1.4) .. (10.93,3.29)   ;
%Curve Lines [id:da08906873023234474] 
\draw    (357,100) .. controls (391.55,95) and (409.55,99) .. (436.55,85) ;
%Curve Lines [id:da23040103237176834] 
\draw    (358,129) .. controls (371.55,127) and (386.55,129) .. (395.55,130) ;
%Curve Lines [id:da7206393887959359] 
\draw    (385,160) .. controls (403.55,159) and (423.55,166) .. (439.55,169) ;
%Curve Lines [id:da33936141191746805] 
\draw    (385,160) .. controls (391.55,161) and (396.55,152) .. (395.55,130) ;
%Curve Lines [id:da54245262295717] 
\draw    (385,160) .. controls (380.55,158) and (385.55,140) .. (395.55,130) ;
%Curve Lines [id:da2480790116383863] 
\draw    (495,100) .. controls (527.55,99) and (554.55,104) .. (581.55,86) ;
%Curve Lines [id:da5706122205107393] 
\draw    (497.55,132) .. controls (507.55,128) and (525.1,132) .. (532.55,132) ;
%Curve Lines [id:da1047400600288323] 
\draw    (503.55,163) .. controls (526.1,158) and (559.55,166) .. (585.55,169) ;
%Curve Lines [id:da7774912110255832] 
\draw    (501.55,138) .. controls (509.55,136) and (523.55,134) .. (532.55,132) ;
%Curve Lines [id:da08343690494886036] 
\draw    (501.55,138) .. controls (507.55,137) and (509.55,163) .. (503.55,163) ;
%Curve Lines [id:da07419689767935644] 
\draw    (503.55,163) .. controls (498.55,162) and (497.55,141) .. (501.55,138) ;
%Curve Lines [id:da7622368998113929] 
\draw    (522.55,162) .. controls (528.55,163) and (532.55,153) .. (532.55,132) ;
%Curve Lines [id:da9549158617953807] 
\draw  [dash pattern={on 4.5pt off 4.5pt}]  (522.55,162) .. controls (516.55,162) and (522.55,145) .. (532.55,132) ;
%Straight Lines [id:da73294982955271] 
\draw    (454.55,132.5) -- (462.55,132.5)(454.55,135.5) -- (462.55,135.5) ;
\draw [shift={(470.55,134)}, rotate = 180] [color={rgb, 255:red, 0; green, 0; blue, 0 }  ][line width=0.75]    (10.93,-3.29) .. controls (6.95,-1.4) and (3.31,-0.3) .. (0,0) .. controls (3.31,0.3) and (6.95,1.4) .. (10.93,3.29)   ;

% Text Node
\draw (76,148.4) node [anchor=north west][inner sep=0.75pt]    {$\gamma $};
% Text Node
\draw (206,158.4) node [anchor=north west][inner sep=0.75pt]    {$\gamma_1$};
% Text Node
\draw (91,85.4) node [anchor=north west][inner sep=0.75pt]    {$P$};
% Text Node
\draw (246,91.4) node [anchor=north west][inner sep=0.75pt]    {$P$};
% Text Node
\draw (404,116.4) node [anchor=north west][inner sep=0.75pt]    {$P$};
% Text Node
\draw (373,165.4) node [anchor=north west][inner sep=0.75pt]    {$\gamma $};
% Text Node
\draw (493,165.4) node [anchor=north west][inner sep=0.75pt]    {$\gamma_1$};
% Text Node
\draw (538,121.4) node [anchor=north west][inner sep=0.75pt]    {$P$};

\end{tikzpicture}

    \caption{Attach a $1$-petal flower.}
    \label{fig: add 1-petal flower}
\end{figure}
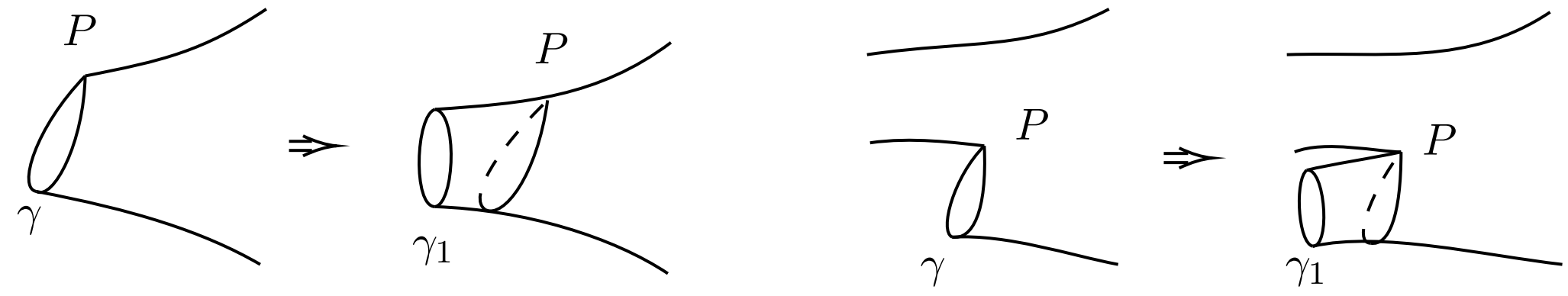

If the angle between $\overrightarrow{\gamma}$ and $\overleftarrow{\gamma}$ at $P$ interior to $S$ is less than $2\pi$, then $P$ is not a branch point on the resulting surface. Otherwise $P$ is a branch point on the resulting surface (see Figure \ref{fig: add 1-petal flower}).

\subsection{Free and trapped geodesic loops}

Let $S$ be a branched hyperbolic surface. Let $\alpha\subset S$ be a simple geodesic loop based at $P$.

\begin{definition}
    We say that $\alpha$ is \textit{free} if either $\measuredangle(\overleftarrow{\alpha},\overrightarrow{\alpha})<2\pi$ or $\measuredangle(\overrightarrow{\alpha},\overleftarrow{\alpha})<2\pi$. Otherwise $\alpha$ is called \textit{trapped}.
\end{definition}

Suppose that the holonomy of $S$ around every simple closed curve is hyperbolic, and that $l_\alpha\leq\frac{1}{100}$. By Corollary \ref{cor-g3}, $\alpha$ disconnects $N_{\frac{1}{2}}(P)$. Let $N$ be one component of $N_{\frac{1}{2}}(P)\backslash\alpha$ \footnote{Approximate $N$ by a subsurface with piecewise smooth boundary so that the theory so developed applies, or observe that the proofs go through in this situation without change.}. If the angle formed by $\overrightarrow\alpha$ and $\overleftarrow\alpha$ at $P$ interior to $N$ is less than $2\pi$, we say that $\alpha$ is free in $N$. Otherwise we say that $\alpha$ is trapped in $N$.

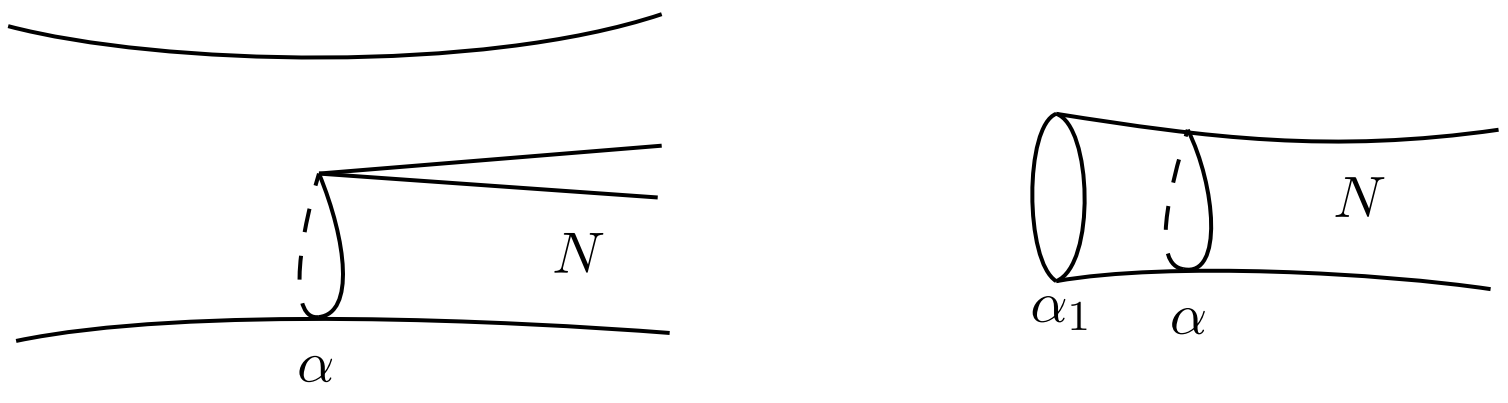
\begin{figure}
    \centering

\tikzset{every picture/.style={line width=0.75pt}} %set default line width to 0.75pt        

\begin{tikzpicture}[x=0.75pt,y=0.75pt,yscale=-1,xscale=1]
%uncomment if require: \path (0,339); %set diagram left start at 0, and has height of 339

%Curve Lines [id:da9597837501635287] 
\draw    (127,83) .. controls (165,93) and (246,95) .. (291,80) ;
%Curve Lines [id:da2298430890134714] 
\draw    (129,162) .. controls (172,153) and (254,157) .. (293,160) ;
%Straight Lines [id:da8882404306772747] 
\draw    (205,120) -- (291,113) ;
%Straight Lines [id:da5502122390526561] 
\draw    (205,120) -- (290,126) ;
%Curve Lines [id:da9924029286014756] 
\draw    (205,156) .. controls (214,155) and (212,137) .. (205,120) ;
%Curve Lines [id:da7594279389865072] 
\draw  [dash pattern={on 4.5pt off 4.5pt}]  (205,156) .. controls (198,157) and (199,140) .. (205,120) ;
%Curve Lines [id:da4000702479383026] 
\draw    (390,105) .. controls (428,111) and (459,115) .. (501,109) ;
%Curve Lines [id:da4300978415456074] 
\draw    (390,147) .. controls (418,142) and (472,145) .. (499,149) ;
%Curve Lines [id:da32452422589929686] 
\draw    (390,105) .. controls (399,108) and (400,143) .. (390,147) ;
%Curve Lines [id:da4820650418364919] 
\draw    (390,147) .. controls (382,142) and (382,108) .. (390,105) ;
%Curve Lines [id:da19565941525810837] 
\draw    (423,109) .. controls (430,123) and (432,146) .. (422,144) ;
%Curve Lines [id:da05773184679792154] 
\draw  [dash pattern={on 4.5pt off 4.5pt}]  (422,144) .. controls (414,143) and (418,125) .. (423,109) ;

% Text Node
\draw (198,164.4) node [anchor=north west][inner sep=0.75pt]    {$\alpha $};
% Text Node
\draw (262,133.4) node [anchor=north west][inner sep=0.75pt]    {$N$};
% Text Node
\draw (417,152.4) node [anchor=north west][inner sep=0.75pt]    {$\alpha $};
% Text Node
\draw (382,149.4) node [anchor=north west][inner sep=0.75pt]    {$\alpha_1$};
% Text Node
\draw (458,119.4) node [anchor=north west][inner sep=0.75pt]    {$N$};

\end{tikzpicture}

    \caption{If $\alpha$ is free, then it is the petal of a $1$-petal flower.}
    \label{fig: lem-free}
\end{figure}

\begin{lemma}\label{lem-free}
    Let $l_\alpha<l\leq\frac{1}{2}$. Suppose that $\alpha$ is free in $N$, and that $r_P>l$. Then $\alpha$ is the petal of a $1$-petal flower of length $\geq l-l_\alpha$ in $N$ (see Figure \ref{fig: lem-free}).
\end{lemma}

\begin{proof}
    Let $N_1$ be the surface obtained from $N$ by attaching a $1$-petal flower of length $l_\alpha$ along $\alpha$. Let $\alpha_1$ be the new boundary curve of $N_1$, which is an equidistant curve. The assumption that the angle formed by $\overrightarrow\alpha$ and $\overleftarrow\alpha$ at $P$ is less than $2\pi$ implies that $P$ is not a branch point of $N_1$. So $r_{\alpha_1}\geq r_P(S)>l$. By Lemma \ref{cor-c2}, we get $\mathrm{inj}_{\alpha_1}>l$. So $\alpha_1$ has a collar neighborhood of length $l$. Removing the short $1$-petal flower that we attached, we get the desired $1$-petal flower in $N$.
\end{proof}

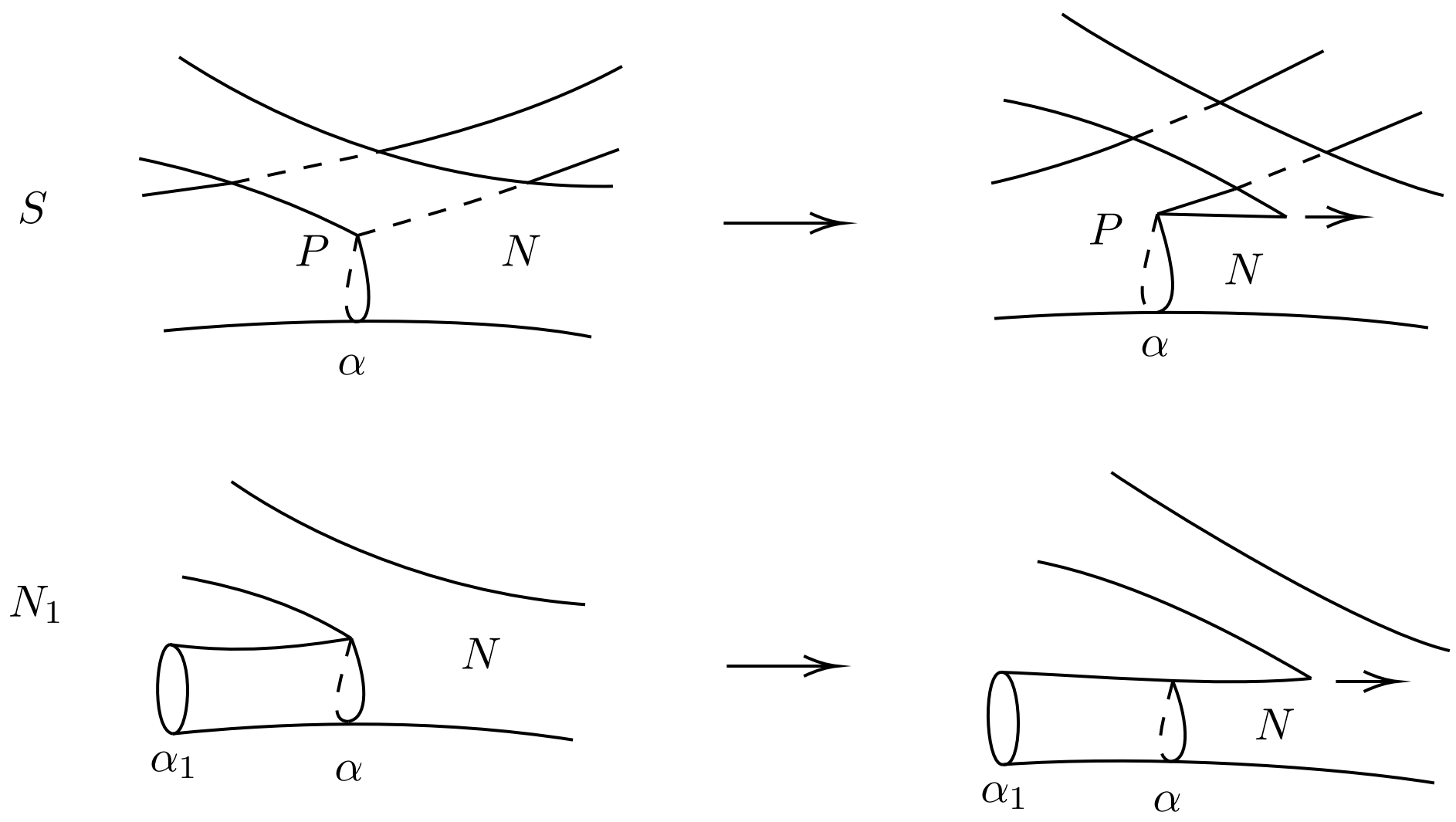
\begin{figure}
    \centering

\tikzset{every picture/.style={line width=0.75pt}} %set default line width to 0.75pt        

\begin{tikzpicture}[x=0.75pt,y=0.75pt,yscale=-1,xscale=1]
%uncomment if require: \path (0,339); %set diagram left start at 0, and has height of 339

%Curve Lines [id:da6345476543381797] 
\draw    (90,136) .. controls (131,132) and (194,131) .. (229,138) ;
%Curve Lines [id:da6370642660489599] 
\draw    (95,47) .. controls (127,68) and (177,90) .. (236,89) ;
%Curve Lines [id:da6118674280909501] 
\draw    (82,80) .. controls (106,85) and (130,93) .. (153,105) ;
%Straight Lines [id:da2508932760217387] 
\draw    (83,92) -- (112,88) ;
%Straight Lines [id:da7700693452253352] 
\draw  [dash pattern={on 4.5pt off 4.5pt}]  (112,88) -- (159,78) ;
%Curve Lines [id:da5737931747611248] 
\draw    (159,78) .. controls (189,71) and (215,63) .. (239,50) ;
%Curve Lines [id:da4813226080039297] 
\draw  [dash pattern={on 4.5pt off 4.5pt}]  (153,105) .. controls (170,100) and (192,94) .. (208,88) ;
%Straight Lines [id:da9061394565227702] 
\draw    (208,88) -- (238,77) ;
%Curve Lines [id:da831846204774454] 
\draw    (153,105) .. controls (157,118) and (159,134) .. (152,133) ;
%Curve Lines [id:da09919680612845316] 
\draw  [dash pattern={on 4.5pt off 4.5pt}]  (152,133) .. controls (147,130) and (150,120) .. (153,105) ;
%Straight Lines [id:da2843989098181139] 
\draw    (272,101) -- (309,101) ;
\draw [shift={(311,101)}, rotate = 180] [color={rgb, 255:red, 0; green, 0; blue, 0 }  ][line width=0.75]    (10.93,-3.29) .. controls (6.95,-1.4) and (3.31,-0.3) .. (0,0) .. controls (3.31,0.3) and (6.95,1.4) .. (10.93,3.29)   ;
%Curve Lines [id:da22221231082686554] 
\draw    (360,132) .. controls (400,129) and (461,129) .. (501,135) ;
%Straight Lines [id:da19431372276099335] 
\draw    (413,98) -- (455,99) ;
%Straight Lines [id:da5648695745042586] 
\draw    (413,98) -- (438,90) ;
%Curve Lines [id:da31629771081212577] 
\draw    (363,61) .. controls (400,68) and (428,83) .. (455,99) ;
%Curve Lines [id:da6397883686266633] 
\draw    (382,33) .. controls (410,52) and (480,86) .. (506,92) ;
%Straight Lines [id:da1836115821697638] 
\draw  [dash pattern={on 4.5pt off 4.5pt}]  (438,90) -- (468,78) ;
%Straight Lines [id:da6146533249820817] 
\draw    (468,78) -- (499,65) ;
%Curve Lines [id:da12578979177822158] 
\draw    (359,88) .. controls (373,85) and (395,78) .. (406,73) ;
%Straight Lines [id:da507987990142928] 
\draw  [dash pattern={on 4.5pt off 4.5pt}]  (406,73) -- (433,62) ;
%Straight Lines [id:da47489284186796465] 
\draw    (433,62) -- (467,45) ;
%Curve Lines [id:da9109829277112558] 
\draw    (413,98) .. controls (418,114) and (421,128) .. (413,130) ;
%Curve Lines [id:da062103365709558545] 
\draw  [dash pattern={on 4.5pt off 4.5pt}]  (413,98) .. controls (408,118) and (405,129) .. (413,130) ;
%Straight Lines [id:da8870844465465015] 
\draw    (461,99) -- (477,99) ;
\draw [shift={(479,99)}, rotate = 180] [color={rgb, 255:red, 0; green, 0; blue, 0 }  ][line width=0.75]    (10.93,-3.29) .. controls (6.95,-1.4) and (3.31,-0.3) .. (0,0) .. controls (3.31,0.3) and (6.95,1.4) .. (10.93,3.29)   ;
%Curve Lines [id:da3077919223517286] 
\draw    (93,267) .. controls (132,263) and (178,262) .. (223,269) ;
%Curve Lines [id:da46648796320857455] 
\draw    (96,216) .. controls (118,220) and (135,226) .. (151,236) ;
%Curve Lines [id:da44002869998762395] 
\draw    (112,185) .. controls (132,199) and (174,221) .. (227,225) ;
%Curve Lines [id:da57223715658821] 
\draw    (92,238) .. controls (114,241) and (133,239) .. (151,236) ;
%Curve Lines [id:da846471218898484] 
\draw    (92,238) .. controls (100,239) and (99,267) .. (93,267) ;
%Curve Lines [id:da5865085817117067] 
\draw    (92,238) .. controls (87,239) and (86,266) .. (93,267) ;
%Curve Lines [id:da5137936298285727] 
\draw    (151,236) .. controls (157,252) and (156,262) .. (150,263) ;
%Curve Lines [id:da3213907095943985] 
\draw  [dash pattern={on 4.5pt off 4.5pt}]  (151,236) .. controls (147,251) and (143,263) .. (150,263) ;

%Straight Lines [id:da9463948499705496] 
\draw    (273,245) -- (307,245) ;
\draw [shift={(309,245)}, rotate = 180] [color={rgb, 255:red, 0; green, 0; blue, 0 }  ][line width=0.75]    (10.93,-3.29) .. controls (6.95,-1.4) and (3.31,-0.3) .. (0,0) .. controls (3.31,0.3) and (6.95,1.4) .. (10.93,3.29)   ;
%Curve Lines [id:da9859437886473987] 
\draw    (363,277) .. controls (405,274) and (465,277) .. (503,283) ;
%Curve Lines [id:da520951346356917] 
\draw    (362,247) .. controls (391,248) and (432,252) .. (463,249) ;
%Curve Lines [id:da3723565638916352] 
\draw    (374,211) .. controls (404,217) and (436,233) .. (463,249) ;
%Curve Lines [id:da3887111286067595] 
\draw    (398,182) .. controls (411,191) and (478,233) .. (508,240) ;
%Curve Lines [id:da050873805571555186] 
\draw    (418,250) .. controls (423,262) and (424,276) .. (417,276) ;
%Curve Lines [id:da36125622796185064] 
\draw  [dash pattern={on 4.5pt off 4.5pt}]  (418,250) .. controls (415,262) and (411,274) .. (417,276) ;
%Curve Lines [id:da4165742441736642] 
\draw    (362,247) .. controls (369,246) and (370,278) .. (363,277) ;
%Curve Lines [id:da9021719091053005] 
\draw    (363,277) .. controls (357,278) and (356,248) .. (362,247) ;
%Straight Lines [id:da4597462569583458] 
\draw    (471,250) -- (489,250) ;
\draw [shift={(491,250)}, rotate = 180] [color={rgb, 255:red, 0; green, 0; blue, 0 }  ][line width=0.75]    (10.93,-3.29) .. controls (6.95,-1.4) and (3.31,-0.3) .. (0,0) .. controls (3.31,0.3) and (6.95,1.4) .. (10.93,3.29)   ;

% Text Node
\draw (145,142.4) node [anchor=north west][inner sep=0.75pt]    {$\alpha $};
% Text Node
\draw (198,103.4) node [anchor=north west][inner sep=0.75pt]    {$N$};
% Text Node
\draw (131,103.4) node [anchor=north west][inner sep=0.75pt]    {$P$};
% Text Node
\draw (41,89.4) node [anchor=north west][inner sep=0.75pt]    {$S$};
% Text Node
\draw (433,109.4) node [anchor=north west][inner sep=0.75pt]    {$N$};
% Text Node
\draw (406,136.4) node [anchor=north west][inner sep=0.75pt]    {$\alpha $};
% Text Node
\draw (389,96.4) node [anchor=north west][inner sep=0.75pt]    {$P$};
% Text Node
\draw (38,217.4) node [anchor=north west][inner sep=0.75pt]    {$N_1$};
% Text Node
\draw (144,274.4) node [anchor=north west][inner sep=0.75pt]    {$\alpha $};
% Text Node
\draw (84,271.4) node [anchor=north west][inner sep=0.75pt]    {$\alpha_1$};
% Text Node
\draw (185,234.4) node [anchor=north west][inner sep=0.75pt]    {$N$};
% Text Node
\draw (354,281.4) node [anchor=north west][inner sep=0.75pt]    {$\alpha_1$};
% Text Node
\draw (410,284.4) node [anchor=north west][inner sep=0.75pt]    {$\alpha $};
% Text Node
\draw (443,257.4) node [anchor=north west][inner sep=0.75pt]    {$N$};

\end{tikzpicture}

    \caption{If $\alpha$ is trapped, then we can `move some multiplicity of $P$ away'.}
    \label{fig: lem-trap}
\end{figure}

\begin{lemma}\label{lem-trap}
    Let $200l_\alpha\leq l\leq\frac{1}{10}$. Suppose that $\alpha$ is trapped in $N$, and that $r_P\geq l$. Then there is $N'\in[N]$ such that $\alpha$ is free in $N'$ and $r_\alpha(N')\geq\frac{1}{4}l$ (see Figure \ref{fig: lem-trap}).
\end{lemma}

\begin{proof}
     Let $N_1$ be the surface obtained from $N$ by attaching a $1$-petal flower $C$ of length $l+l_\alpha$ along $\alpha$. Then $r_P(N_1)\geq l$ (since $\alpha$ is trapped in $N$, $P$ is a branch point of $N_1$). By Proposition \ref{prop-m2}, we have $R_{\alpha_1}\geq r_{\alpha_1}+\frac{1}{3}r_P(N_1)> l-l_{\alpha}+\frac{1}{3}l$. Therefore there is $N_1'\in[N_1]$ such that $r_{\alpha_1}(N_1')> l-l_{\alpha}+\frac{1}{3}l$. In particular, by Lemma \ref{cor-c2}, $\alpha_1$ has a collar neiborhood of length $l-l_{\alpha}+\frac{1}{3}l$ in $N_1'$, so $C\subset N_1'$ naturally. Let $N'=N_1'\backslash C$. Then $N'\in[N]$ and $r_\alpha(N')\geq\frac{1}{3}l-l_\alpha>\frac{1}{4}l$, as desired.
\end{proof}

\subsection{Dividing and non-dividing loops}

Let $(P,\gamma_1,\gamma_2)$ be a nice triple such that $r_{\gamma_i}\leq\frac{1}{10}\ (i=1,2)$. Let $\alpha$ be a simple geodesic loop based at $P$ of length $\leq5\epsilon$. Since $d(P,\gamma_1),d(P,\gamma_2)\leq \frac{1}{10}$, we have $\gamma_1,\gamma_2\subset N_{\frac{1}{2}}(P)$.

\begin{definition}
    We call $\alpha$ \textit{dividing} if $\gamma_1,\gamma_2$ lie in different components of $N_{\frac{1}{2}}(P)\backslash\alpha$. Otherwise $\alpha$ is called \textit{non-dividing}.
\end{definition}

\begin{lemma}\label{lem-noncross trap}
    Suppose that $\alpha$ is non-dividing and trapped, and that $r_P\geq\epsilon_1$. Then there exists $S'\in[S]$ with a nice triple $(P',\gamma_1,\gamma_2)$ such that $r_{P'}\geq\frac{1}{4}r_P,\ r_{\gamma_i}(S')\leq r_{\gamma_i}+5\epsilon\ (i=1,2),\ m(P')<m(P)$.
\end{lemma}

\begin{proof}
    Let $N$ be the components of $N_{\frac{1}{2}}(P)\backslash\alpha$ not containing $\gamma_1,\gamma_2$. By Lemma \ref{lem-trap}, there is $N'\in[N]$ such that $\alpha$ is free in $N'$ and $r_{\alpha}(N')\geq \frac{1}{4}r_P$. Let $S'=(S\backslash N)\cup N'$, attaching $\partial(S\backslash N)\cap\partial N$ to $\partial N'\cong\partial N$ via the identity. Then $S'\in[S]$. Let $P'$ be the image of $P$ in $S'$. Since $\alpha$ is trapped in $N$, while free in $N'$, the cone angle at $P'$ is less than the cone angle at $P$. Therefore $m(P')<m(P)$. For any branch or boundary point $Q$ of $S\backslash N$ not lying on $\alpha$, we have $d(\alpha,Q)\geq r_P-l_{\alpha}>\frac{1}{4}r_P$. Combining this with $r_{\alpha}(N')\geq \frac{1}{4}r_P$ we get $r_{P'}\geq \frac{1}{4}r_P$.
    \begin{claim}
        $r_{\gamma_i}-l_\alpha\leq r_{\gamma_i}(S')\leq d_{S'}(\gamma_i,P')\leq r_{\gamma_i}+l_\alpha\ (i=1,2).$
    \end{claim}
    Indeed, $r_{\gamma_i}(S')\geq r_{\gamma_i}(S\backslash N)\geq r_{\gamma_i}-l_\alpha$, while $d_{S'}(\gamma_i,P')\leq d_{S\backslash N}(\gamma_i,P)\leq d_{S\backslash N}(\gamma_i,\alpha)+l_{\alpha}\leq d_S(\gamma_i,P)+l_\alpha=r_{\gamma_i}+l_\alpha$, proving the claim.

    Since $l_\alpha\leq5\epsilon$, the claim implies $r_{\gamma_i}(S')\leq r_{\gamma_i}+5\epsilon$. Moreover, we deduce by Lemma \ref{lem-useful criterion} that $d_{S'}(\gamma_i,P')=r_{\gamma_i}(S')$, so $(P',\gamma_1,\gamma_2)$ is a nice triple of $S'$.
\end{proof}

\begin{lemma}\label{lem-noncross trap2}
    Suppose that $\alpha$ is non-dividing, and that $r_P\geq\epsilon_1$. Let $N$ be the component of $N_{\frac{1}{2}}(P)\backslash\alpha$ containing $\gamma_1,\gamma_2$. Then $\alpha$ is trapped in $N$.
\end{lemma}

\begin{proof}
    We argue by contradiction. Suppose that $\alpha$ is free in $N$. By Lemma \ref{lem-free}, $\alpha$ is the petal of a $1$-petal flower $C\subset N$ of length $r_P-10\epsilon>500\epsilon$. Let $\beta$ be the other boundary component of $C$. Then $d(\gamma_i,\beta)\leq d(\gamma_i,P)-500\epsilon\ (i=1,2)$. Since $d(\gamma_i,P)=r_{\gamma_i}$, we get $d(\gamma_1,\gamma_2)\leq d(\gamma_1,\beta)+l_{\beta}+d(\beta,\gamma_2)<r_{\gamma_1}+r_{\gamma_2}-10\epsilon$. We deduce by Lemma \ref{lem-d1} that $S$ is a collar, a contradiction.
\end{proof}

\subsection{Reducibility from two non-dividing loops}

We show that a nice triple is reducible if there are two non-dividing geodesic loops based at $P$ such that the angle between them is in certain range.

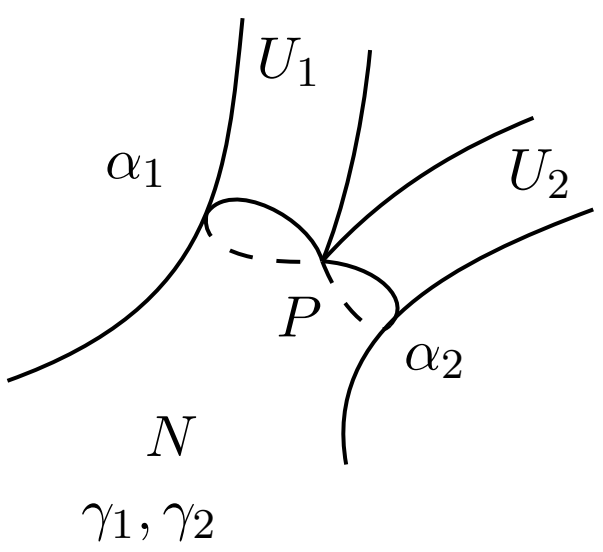
\begin{figure}
    \centering

\tikzset{every picture/.style={line width=0.75pt}} %set default line width to 0.75pt        

\begin{tikzpicture}[x=0.75pt,y=0.75pt,yscale=-1,xscale=1]
%uncomment if require: \path (0,339); %set diagram left start at 0, and has height of 339

%Curve Lines [id:da9913848333697405] 
\draw    (131,129) .. controls (184,110) and (186,78) .. (190,38) ;
%Curve Lines [id:da9738056403148035] 
\draw    (210,99) .. controls (215,87) and (220,67) .. (222,46) ;
%Curve Lines [id:da07192249995446542] 
\draw    (210,99) .. controls (227,80) and (246,70) .. (263,63) ;
%Curve Lines [id:da8432112454123474] 
\draw    (216,150) .. controls (210,116) and (243,99) .. (278,86) ;
%Curve Lines [id:da9531661107653834] 
\draw    (181,87) .. controls (187,78) and (207,88) .. (210,99) ;
%Curve Lines [id:da8247247959394096] 
\draw  [dash pattern={on 4.5pt off 4.5pt}]  (181,87) .. controls (179,96) and (192,100) .. (210,99) ;
%Curve Lines [id:da8960969584159836] 
\draw    (210,99) .. controls (224,100) and (233,109) .. (227,115) ;
%Curve Lines [id:da0017561031369816016] 
\draw  [dash pattern={on 4.5pt off 4.5pt}]  (210,99) .. controls (215,112) and (224,119) .. (227,115) ;

% Text Node
\draw (197,106.4) node [anchor=north west][inner sep=0.75pt]    {$P$};
% Text Node
\draw (154,70.4) node [anchor=north west][inner sep=0.75pt]    {$\alpha _{1}$};
% Text Node
\draw (229,118.4) node [anchor=north west][inner sep=0.75pt]    {$\alpha _{2}$};
% Text Node
\draw (192,41.4) node [anchor=north west][inner sep=0.75pt]    {$U_{1}$};
% Text Node
\draw (255,69.4) node [anchor=north west][inner sep=0.75pt]    {$U_{2}$};
% Text Node
\draw (164,136.4) node [anchor=north west][inner sep=0.75pt]    {$N$};
% Text Node
\draw (148,158.4) node [anchor=north west][inner sep=0.75pt]    {$\gamma _{1} ,\gamma _{2}$};

\end{tikzpicture}

    \caption{Two non-dividing loops}
    \label{fig: two noncross}
\end{figure}

Throughout this subsection let $(P,\gamma_1,\gamma_2)$ be a nice triple such that $r_P\geq \epsilon_1,\ r_{\gamma_i}\leq10\delta_0\ (i=1,2)$. Let $\alpha_1,\alpha_2$ be two non-dividing simple geodesic loops of length $\leq5\epsilon$ based at $P$. Let $N_i,N_i'$ be the two components of $N_{\frac{1}{2}}(P)\backslash\alpha_i\ (i=1,2)$, where $N_i'$ is the one containing $\gamma_1,\gamma_2$. By Lemma \ref{lem-noncross trap} and Lemma \ref{lem-noncross trap2}, we may assume that $\alpha_i$ is free in $N_i$ (otherwise the triple is reducible and we are done). By Lemma \ref{lem-free}, $\alpha_i$ is the petal of a $1$-petal flower $U_i\subset N_i$ of length $r_P-5\epsilon$. Since $\overline{U_1}$ does not contain any geodesic loops based at $P$ other than $\alpha_1$, we deduce that $\mathring\alpha_2$ is contained in $N_1'$. Therefore $N_2\subset N_1'$. It follows that $\overline{U_1}\cap\overline{U_2}=\{P\}$. Let $N=N_{\frac{1}{2}}(P)\backslash(N_1\cup N_2)$ (see Figure \ref{fig: two noncross}). Orient $\alpha_i$ so that $U_i$ lie on the left.

\begin{lemma}\label{lem-two noncross 1}
    Suppose that $\measuredangle(\overleftarrow{\alpha_1},\overrightarrow{\alpha_2})\in[\frac{1}{2}\pi,\frac{3}{2}\pi]$. Then $(P,\gamma_1,\gamma_2)$ is reducible.
\end{lemma}

\begin{figure}
    \centering

\tikzset{every picture/.style={line width=0.75pt}} %set default line width to 0.75pt        

\begin{tikzpicture}[x=0.75pt,y=0.75pt,yscale=-1,xscale=1]
%uncomment if require: \path (0,339); %set diagram left start at 0, and has height of 339

%Curve Lines [id:da9182584329186294] 
\draw    (147,170) .. controls (159,168) and (182,169) .. (192,172) ;
\draw [shift={(175.79,169.46)}, rotate = 183.39] [color={rgb, 255:red, 0; green, 0; blue, 0 }  ][line width=0.75]    (10.93,-4.9) .. controls (6.95,-2.3) and (3.31,-0.67) .. (0,0) .. controls (3.31,0.67) and (6.95,2.3) .. (10.93,4.9)   ;
%Curve Lines [id:da8109374157016175] 
\draw    (192,172) .. controls (205,179) and (218,192) .. (226,204) ;
\draw [shift={(215.15,190.45)}, rotate = 225.29] [color={rgb, 255:red, 0; green, 0; blue, 0 }  ][line width=0.75]    (10.93,-4.9) .. controls (6.95,-2.3) and (3.31,-0.67) .. (0,0) .. controls (3.31,0.67) and (6.95,2.3) .. (10.93,4.9)   ;
%Curve Lines [id:da48269246231889007] 
\draw    (192,172) .. controls (206,139) and (231,94) .. (258,66) ;
%Curve Lines [id:da8060280926621474] 
\draw    (226,204) .. controls (237,176) and (260,125) .. (284,100) ;
%Curve Lines [id:da7229606146676705] 
\draw    (258,66) .. controls (270,76) and (281,91) .. (284,100) ;
%Curve Lines [id:da3695972332410483] 
\draw    (216,65) .. controls (230,63) and (248,63) .. (258,66) ;
%Curve Lines [id:da981536270672176] 
\draw    (147,170) .. controls (164,134) and (188,90) .. (216,65) ;

% Text Node
\draw (236,94.4) node [anchor=north west][inner sep=0.75pt]    {$\sigma $};
% Text Node
\draw (221,139.4) node [anchor=north west][inner sep=0.75pt]    {$V'_{2}$};
% Text Node
\draw (178,121.4) node [anchor=north west][inner sep=0.75pt]    {$V'_{1}$};
% Text Node
\draw (163,177.4) node [anchor=north west][inner sep=0.75pt]    {$\alpha _{1}$};
% Text Node
\draw (197,192.4) node [anchor=north west][inner sep=0.75pt]    {$\alpha _{2}$};
% Text Node
\draw (257,44.4) node [anchor=north west][inner sep=0.75pt]    {$Q$};

\end{tikzpicture}

    \caption{Identify $\sigma_1^+$ with $\sigma_2^-$, $\sigma_1^-$ with $\sigma_2^+$. Intuitively, $Q$ carries one multiplicity of $P$ and is moved away. This is somewhat analogous to particle fission.}
    \label{fig: two noncross 1}
\end{figure}
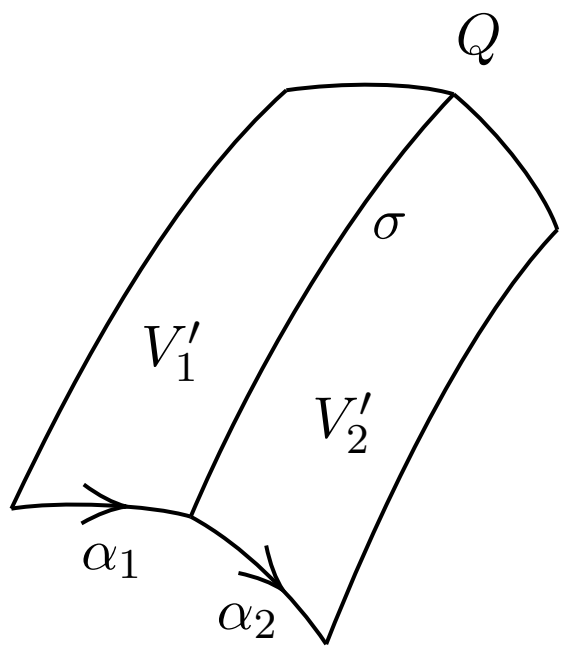

\begin{proof}
    Let $\theta=\measuredangle(\overleftarrow{\alpha_1},\overrightarrow{\alpha_2})\in[\frac{1}{2}\pi,\frac{3}{2}\pi]$. Then $\pi-\frac{1}{2}\theta\in[\frac{1}{4}\pi,\frac{3}{4}\pi]$. Let $l=\frac{1}{2}r_P$. By Lemma \ref{cor-h1}, there is a simple geodesic $\sigma_1:[0,l]\rightarrow\overline{U_1}$ of length $l$ such that $\sigma_1(0)=P$, $\measuredangle(\overrightarrow{\sigma_1},\overleftarrow{\alpha_1})=\pi-\frac{1}{2}\theta$. Similarly there is a simple geodesic $\sigma_2:[0,l]\rightarrow\overline{U_2}$ of length $l$ such that $\sigma_2(0)=P$, $\measuredangle(\overrightarrow{\alpha_2},\overrightarrow{\sigma_2})=\pi-\frac{1}{2}\theta$. We have $\measuredangle(\overrightarrow{\sigma_1},\overrightarrow{\sigma_2})=2\pi$. Thus $\sigma_1,\sigma_2$ is an associated pair. Let $S_1$ be the move of $S$ along $\sigma_1,\sigma_2$. Thus $S_1$ is obtained from $S$ by cutting along $\sigma_1,\sigma_2$ and identify $\sigma_1^+$ with $\sigma_2^-$, $\sigma_1^-$ with $\sigma_2^+$ (see Figure \ref{fig: two noncross 1}). Let $Q$ be the branch point on $S_1$ obtained by identifying $\sigma_1(l)$ and $\sigma_2(l)$, and let $P_1$ be the image of $\sigma_1^-(0)=\sigma_2^+(0)$ in $S_1$. We have $m(P_1)=m(P)-1$.
    \begin{claim}\label{claim-two noncross 1-1}
        $(P_1,\gamma_1,\gamma_2)$ is a nice triple of $S_1$.
    \end{claim}
    Indeed, $r_{\gamma_i}(S_1)\geq r_{\gamma_i}(N)\geq r_{\gamma_i}-5\epsilon$, while $d_{S_1}(P_1,\gamma_i)\leq d_N(P_1,\gamma_i)\leq d_S(P,\gamma_i)+5\epsilon=r_{\gamma_i}(S)+5\epsilon$. It follows that $d_{S_1}(P_1,\gamma_i)\leq r_{\gamma_i}(S_1)+10\epsilon$. Claim \ref{claim-two noncross 1-1} then follows from Lemma \ref{lem-useful criterion}. Moreover, we have $r_{\gamma_i}(S_1)\leq r_{\gamma_i}(S)+5\epsilon$.
    \begin{claim}\label{claim-two noncross 1-2}
        $r_{P_1}\geq\frac{1}{3}l$.
    \end{claim}
    Claim \ref{claim-two noncross 1-2} would imply $r_{P_1}\geq\frac{1}{10}r_P$. Combining this with Claim \ref{claim-two noncross 1-1}, we deduce that $(P_1,\gamma_1,\gamma_2)$ is reducible, proving the lemma.
    
    It remains to prove Claim \ref{claim-two noncross 1-2}. Let $\beta_i\subset U_i\ (i=1,2)$ be the equidistant curve passing through $\sigma_i(l)$, and let $V_i\subset U_i$ be the $1$-petal flower bounded by $\alpha_i,\beta_i$. Then $\sigma_i\subset V_i$ is a simple geodesic that joins the two boundaries of $V_i$. Let $V_i'$ denote the curved quadrilateral obtained by cutting $V_i$ along $\sigma_i$. Let $V'=V_1'\cup V_2'$ be the annulus obtained by gluing $\sigma_1^+$ with $\sigma_2^-$, $\sigma_1^-$ with $\sigma_2^+$. By Lemma \ref{lem-two noncross 11} below, $V'$ contains a collar of length $\frac{1}{2}l-100\epsilon\geq\frac{1}{3}l$.
    
    Write $S=N\cup V_1\cup V_2\cup W$ (disjoint union), where $W=S\backslash(N\cup V_1\cup V_2)$. Then $S_1=N\cup V'\cup W$. Since $V'$ contains a collar of length $\frac{1}{3}l$, we have $d_{S_1}(P_1,W)\geq\frac{1}{3}l$. $V'$ contains no branch point, while for any branch or boundary point $Q_1$ in $N$, we have $d_{S_1}(P_1,Q_1)\geq d_{N}(\alpha_1\cup\alpha_2,Q_1)\geq r_P(S)-5\epsilon>\frac{1}{3}l$. Claim \ref{claim-two noncross 1-2} follows from these observations.
\end{proof}

\begin{lemma}\label{lem-two noncross 11}
    Let $\alpha_1,\alpha_2\subset\H$ be geodesic segments from $P_1$ to $P_2$, $P_2$ to $P_3$ respectively, such that $\theta=\measuredangle(\overleftarrow{\alpha_1},\overrightarrow{\alpha_2})\in[\frac{1}{2}\pi,\frac{3}{2}\pi]$. Let $\sigma:[0,l]\rightarrow\H$ be a geodesic of length $l$ such that $\sigma(0)=P_2,\measuredangle(\overrightarrow{\alpha_2},\overrightarrow{\sigma})=\pi-\frac{1}{2}\theta$. Suppose that $l_{\alpha_1},l_{\alpha_2},l\leq\frac{1}{100}$. Let $A\in\mathrm{PSL}_2(\R)$ be a hyperbolic or parabolic element such that $A(P_1)=P_3$. Then the distance between the $A$-equidistant lines passing through $\sigma(0)$ and $\sigma(l)$ is at least $\frac{1}{2}l$.
\end{lemma}

\begin{proof}
    By Lemma \ref{cor-h1} one deduces that the angle between $\sigma$ and every $A$-equidistant line that intersects $\sigma$ is in the range $[\frac{1}{6}\pi,\frac{5}{6}\pi]$.
\end{proof}

\begin{figure}
    \centering

\tikzset{every picture/.style={line width=0.75pt}} %set default line width to 0.75pt        

\begin{tikzpicture}[x=0.75pt,y=0.75pt,yscale=-1,xscale=1]
%uncomment if require: \path (0,339); %set diagram left start at 0, and has height of 339

%Curve Lines [id:da6590092039583532] 
\draw    (73,109) .. controls (107,115) and (156,117) .. (192,108) ;
%Curve Lines [id:da31289710270663573] 
\draw    (74,180) .. controls (103,174) and (163,171) .. (194,178) ;
%Curve Lines [id:da5301422127029164] 
\draw    (110.97,114.7) .. controls (119.09,121.64) and (119.3,176.93) .. (110,176) ;
\draw [shift={(109,114)}, rotate = 58.46] [color={rgb, 255:red, 0; green, 0; blue, 0 }  ][line width=0.75]    (10.93,-4.9) .. controls (6.95,-2.3) and (3.31,-0.67) .. (0,0) .. controls (3.31,0.67) and (6.95,2.3) .. (10.93,4.9)   ;
%Curve Lines [id:da386480459146441] 
\draw  [dash pattern={on 4.5pt off 4.5pt}]  (109,114) .. controls (99,114) and (100,176) .. (110,176) ;
%Straight Lines [id:da7733347662339111] 
\draw    (118,148) -- (183,147) ;
\draw [shift={(156.5,147.41)}, rotate = 179.12] [color={rgb, 255:red, 0; green, 0; blue, 0 }  ][line width=0.75]    (10.93,-4.9) .. controls (6.95,-2.3) and (3.31,-0.67) .. (0,0) .. controls (3.31,0.67) and (6.95,2.3) .. (10.93,4.9)   ;
%Straight Lines [id:da8917738048038555] 
\draw    (233,146) -- (279,146) ;
\draw [shift={(281,146)}, rotate = 180] [color={rgb, 255:red, 0; green, 0; blue, 0 }  ][line width=0.75]    (10.93,-3.29) .. controls (6.95,-1.4) and (3.31,-0.3) .. (0,0) .. controls (3.31,0.3) and (6.95,1.4) .. (10.93,3.29)   ;
%Curve Lines [id:da8837018347038759] 
\draw    (315,103) .. controls (328,110) and (342,115) .. (358,116) ;
%Curve Lines [id:da18704232209529303] 
\draw    (315,178) .. controls (327,171) and (342,169) .. (359,167) ;
%Curve Lines [id:da11048484088304122] 
\draw    (358,116) .. controls (370,117) and (368,165) .. (359,167) ;
\draw [shift={(366.16,134.16)}, rotate = 87.81] [color={rgb, 255:red, 0; green, 0; blue, 0 }  ][line width=0.75]    (10.93,-3.29) .. controls (6.95,-1.4) and (3.31,-0.3) .. (0,0) .. controls (3.31,0.3) and (6.95,1.4) .. (10.93,3.29)   ;
%Curve Lines [id:da44347024460866125] 
\draw  [dash pattern={on 4.5pt off 4.5pt}]  (359,167) .. controls (351,166) and (348,120) .. (358,116) ;
%Curve Lines [id:da2634697547773208] 
\draw    (418,116) .. controls (458,112) and (474,112) .. (502,97) ;
%Straight Lines [id:da7390345132228898] 
\draw    (426,138) -- (482,144) ;
\draw [shift={(459.97,141.64)}, rotate = 186.12] [color={rgb, 255:red, 0; green, 0; blue, 0 }  ][line width=0.75]    (10.93,-3.29) .. controls (6.95,-1.4) and (3.31,-0.3) .. (0,0) .. controls (3.31,0.3) and (6.95,1.4) .. (10.93,3.29)   ;
%Straight Lines [id:da2116165551214757] 
\draw    (426,151) -- (482,144) ;
\draw [shift={(459.95,146.76)}, rotate = 172.87] [color={rgb, 255:red, 0; green, 0; blue, 0 }  ][line width=0.75]    (10.93,-3.29) .. controls (6.95,-1.4) and (3.31,-0.3) .. (0,0) .. controls (3.31,0.3) and (6.95,1.4) .. (10.93,3.29)   ;
%Curve Lines [id:da39990298976399286] 
\draw    (418,169) .. controls (452,168) and (484,179) .. (504,187) ;
%Curve Lines [id:da2788197300293531] 
\draw    (418,116) .. controls (423,115) and (428,131) .. (426,138) ;
%Curve Lines [id:da3630438865205742] 
\draw    (418,169) .. controls (423,168) and (428,158) .. (426,151) ;
%Curve Lines [id:da22504058862270016] 
\draw    (418,169) .. controls (411,166) and (409,118) .. (418,116) ;
\draw [shift={(412.35,149.26)}, rotate = 266.89] [color={rgb, 255:red, 0; green, 0; blue, 0 }  ][line width=0.75]    (10.93,-3.29) .. controls (6.95,-1.4) and (3.31,-0.3) .. (0,0) .. controls (3.31,0.3) and (6.95,1.4) .. (10.93,3.29)   ;

% Text Node
\draw (162,128.4) node [anchor=north west][inner sep=0.75pt]    {$\beta $};
% Text Node
\draw (84,136.4) node [anchor=north west][inner sep=0.75pt]    {$\alpha $};
% Text Node
\draw (120,133.4) node [anchor=north west][inner sep=0.75pt]    {$2\pi $};
% Text Node
\draw (245,122) node [anchor=north west][inner sep=0.75pt]   [align=left] {cut};
% Text Node
\draw (352,174.4) node [anchor=north west][inner sep=0.75pt]    {$\alpha ^{+}$};
% Text Node
\draw (401,173.4) node [anchor=north west][inner sep=0.75pt]    {$\alpha ^{-}$};
% Text Node
\draw (470,122.4) node [anchor=north west][inner sep=0.75pt]    {$\beta ^{+}$};
% Text Node
\draw (471,148.4) node [anchor=north west][inner sep=0.75pt]    {$\beta ^{-}$};
% Text Node
\draw (427,126.4) node [anchor=north west][inner sep=0.75pt]    {$2\pi $};

\end{tikzpicture}

    \caption{Cut along $\alpha,\beta$.}
    \label{fig: sur2}
\end{figure}
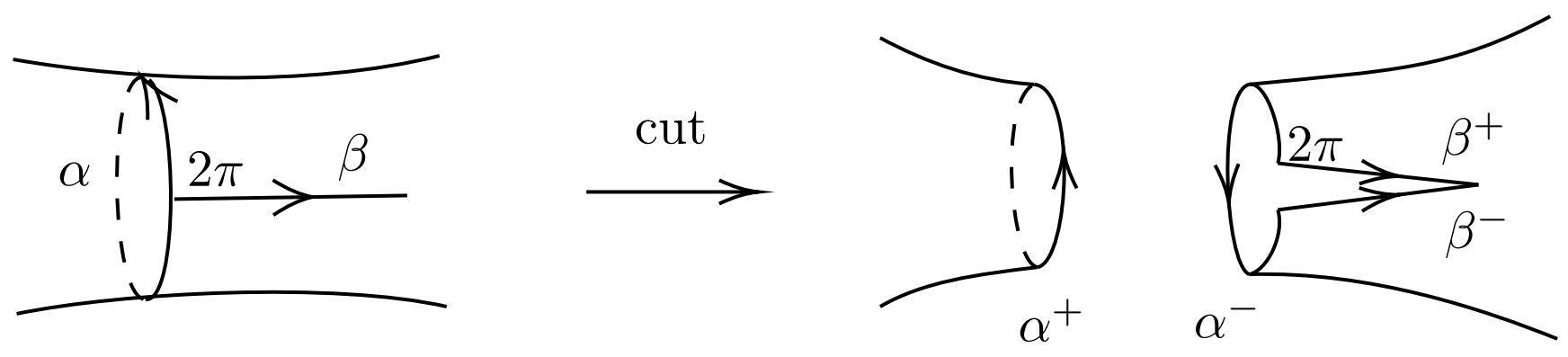

Before discussing the next case, we introduce a surgery which is a variant of moving along a pair of associated curves. The only difference is that in this case one curve in the pair is closed. Let $\alpha:[0,l]\rightarrow S$ be a simple geodesic loop based at $P$, and let $\beta:[0,l]\rightarrow S$ be a simple geodesic arc intersecting with $\alpha$ only at $P=\beta(0)$, such that $\measuredangle(\overrightarrow{\alpha},\overrightarrow{\beta})\equiv0\mod 2\pi$. Suppose that $\overrightarrow{\alpha},\overleftarrow{\alpha},\overrightarrow{\beta}$ is arranged in anti-clockwise order. Cut along $\alpha,\beta$. The resulting surface $S^{cut}$ has two new boundary curves, one of which is $\alpha^+$, while the other is composed of $\alpha^-,\beta^-,\beta^+$ (see Figure \ref{fig: sur2} for a topological picture). Glue back by identifying, for every $t\in[0,l]$, $\alpha^+(t)$ with $\beta^-(t)$, $\alpha^-(t)$ with $\beta^+(t)$. It is easy to check that the metric on $S^{cut}$ descends to a branched hyperbolic metric on the resulting surface. Let $S_1$ denote the resulting surface, equipped with the metric induced from $S^{cut}$.

\begin{lemma}
    $S_1\in[S].$
\end{lemma}

\begin{proof}
    $S_1$ is obtained from $S$ by applying two moves. $\alpha|_{[0,\frac{1}{2}l]}$ and $\beta|_{[0,\frac{1}{2}l]}$ is an associated pair. Let $S_0$ be the move of $S$ along $\alpha|_{[0,\frac{1}{2}l]},\beta|_{[0,\frac{1}{2}l]}$. The images of $\alpha|_{[\frac{1}{2}l,l]},\beta|_{[\frac{1}{2}l,l]}$ in $S_0$, denoted $\alpha_0,\beta_0$, is an associated pair on $S_0$. $S_1$ is the move of $S_0$ along $\alpha_0,\beta_0$.
\end{proof}

Suppose that $\measuredangle(\overrightarrow{\beta},\overrightarrow{\alpha})=2\pi$. Then one can check that the only branch point on $S_1$ lying on the image of $\alpha^+\cup\alpha^-\cup\beta^+\cup\beta^-$ is the point obtained by identifying $\beta^-(0),\beta^-(l),\alpha^+(0)=\alpha^+(l)$. Denote this point by $P_1$. By the same proof as that of Lemma \ref{lem-m12} and Lemma \ref{lem-m13}, one deduces that $r_{P_1}\geq r_P-l$ and $r_{\gamma}(S_1)\geq r_{\gamma}(S)-l,\ d(P_1,\gamma)\leq d(P,\gamma)+l$ for any boundary curve $\gamma$.

We now return to the discussion of two non-dividing loops.

\begin{lemma}\label{lem-two noncross 2}
    Suppose that $\measuredangle(\overleftarrow{\alpha_1},\overrightarrow{\alpha_2})\in(\frac{3}{2}\pi,2\pi]$. Then $(P,\gamma_1,\gamma_2)$ is reducible.
\end{lemma}

\begin{proof}
    By symmetry we may assume $l_{\alpha_1}\geq l_{\alpha_2}$. We first assume that $\measuredangle(\overleftarrow{\alpha_1},\overrightarrow{\alpha_2})\neq2\pi$. By Lemma \ref{lem-h5}, there is an embedded geodesic triangle $T\subset U_1$ with sides $\alpha_1,\beta,\sigma$ (clockwise ordered) such that $l_{\sigma}=l_{\alpha_2}$ and the angle between $\alpha_1$ and $\sigma$ is equal to $2\pi-\measuredangle(\overleftarrow{\alpha_1},\overrightarrow{\alpha_2})$. Orient $\sigma$ so that $\sigma(0)=P$. Then $\measuredangle(\overrightarrow{\sigma},\overrightarrow{\alpha_2})=2\pi$. So we are in the right position to apply the surgery introduced above. Cut along $\sigma,\alpha_2$ and identify $\sigma^+$ with $\alpha_2^-$, $\sigma^-$ with $\alpha_2^+$ (Figure \ref{fig: two noncross 2}). We get a surface $S_1\in[S]$ with a branched point $P_1$ satisfying $m(P_1)=m(P),\ r_{P_1}\geq r_P-5\epsilon,\ r_{\gamma_i}-5\epsilon\leq r_{\gamma_i}(S_1)\leq d_{S_1}(P_1,\gamma_i)\leq r_{\gamma_i}+5\epsilon$. It follows from Lemma \ref{lem-useful criterion} that $(P_1,\gamma_1,\gamma_2)$ is nice. Denote by $\beta_1$ the image of $\beta$ in $S_1$, which is a simple geodesic loop based at $P$. $\beta_1$ is non-dividing. Indeed, the connected set $N$, which contains $\gamma_1,\gamma_2$, is disjoint from $\beta_1$. By a computation of angles one sees that $\beta_1$ is trapped. Therefore by Lemma \ref{lem-noncross trap}, $(P,\gamma_1,\gamma_2)$ is reducible.

    If $\measuredangle(\overleftarrow{\alpha_1},\overrightarrow{\alpha_2})=2\pi$, we have $l_{\alpha_1}\neq l_{\alpha_2}$, otherwise $\rho(\alpha_1*\alpha_2)$ would be an elliptic element or the identity. We may assume $l_{\alpha_1}>l_{\alpha_2}$. The above argument goes through without change, except that the triangle $T$ is degenerate.
\end{proof}

\begin{figure}
    \centering

\tikzset{every picture/.style={line width=0.75pt}} %set default line width to 0.75pt        

\begin{tikzpicture}[x=0.75pt,y=0.75pt,yscale=-1,xscale=1]
%uncomment if require: \path (0,339); %set diagram left start at 0, and has height of 339

%Curve Lines [id:da3096447875410252] 
\draw    (153,166) .. controls (184,168) and (213,168) .. (236,165) ;
\draw [shift={(200.84,167.33)}, rotate = 179.46] [color={rgb, 255:red, 0; green, 0; blue, 0 }  ][line width=0.75]    (10.93,-4.9) .. controls (6.95,-2.3) and (3.31,-0.67) .. (0,0) .. controls (3.31,0.67) and (6.95,2.3) .. (10.93,4.9)   ;
%Curve Lines [id:da650920405170245] 
\draw    (153,166) .. controls (164,136) and (173,101) .. (172,69) ;
%Curve Lines [id:da33100832548168824] 
\draw    (236,165) .. controls (243,140) and (249,88) .. (249,61) ;
%Curve Lines [id:da4235090319311957] 
\draw    (186,123) .. controls (204,136) and (222,150) .. (236,165) ;
\draw [shift={(206.56,138.46)}, rotate = 38.77] [color={rgb, 255:red, 0; green, 0; blue, 0 }  ][line width=0.75]    (10.93,-4.9) .. controls (6.95,-2.3) and (3.31,-0.67) .. (0,0) .. controls (3.31,0.67) and (6.95,2.3) .. (10.93,4.9)   ;
%Curve Lines [id:da4066537916123646] 
\draw    (97,185) .. controls (127,161) and (156,140) .. (186,123) ;
%Curve Lines [id:da5600466851119742] 
\draw    (155,221) .. controls (183,197) and (206,180) .. (236,165) ;
%Curve Lines [id:da7103378987054448] 
\draw  [dash pattern={on 4.5pt off 4.5pt}]  (186,129) .. controls (201,137) and (221,154) .. (236,165) ;
%Curve Lines [id:da6596221725258998] 
\draw  [dash pattern={on 4.5pt off 4.5pt}]  (153,166) .. controls (165,150) and (177,139) .. (186,129) ;

% Text Node
\draw (176,175.4) node [anchor=north west][inner sep=0.75pt]    {$\alpha _{1}$};
% Text Node
\draw (210,124.4) node [anchor=north west][inner sep=0.75pt]    {$\alpha _{2}$};
% Text Node
\draw (196,146.4) node [anchor=north west][inner sep=0.75pt]    {$\sigma $};
% Text Node
\draw (172,141.4) node [anchor=north west][inner sep=0.75pt]    {$\beta $};
% Text Node
\draw (201,68.4) node [anchor=north west][inner sep=0.75pt]    {$U_{1}$};
% Text Node
\draw (112,186.4) node [anchor=north west][inner sep=0.75pt]    {$U_{2}$};

\end{tikzpicture}

    \caption{Cut along $\sigma,\alpha_2$.}
    \label{fig: two noncross 2}
\end{figure}

Combining Lemma \ref{lem-two noncross 1} and Lemma \ref{lem-two noncross 2} we get

\begin{lemma}\label{cor-two noncross}
    Suppose that $\measuredangle(\overleftarrow{\alpha_1},\overrightarrow{\alpha_2})\in[\frac{1}{2}\pi,2\pi]$. Then $(P,\gamma_1,\gamma_2)$ is reducible.
\end{lemma}

\begin{lemma}\label{lem-two noncross 3}
    Suppose that $\measuredangle(\overleftarrow{\alpha_1},\overrightarrow{\alpha_2})\in(2\pi,3\pi)$, and that there is an embedded geodesic triangle $T\subset S$ with sides $\alpha_1,\sigma,\beta$ (clockwise ordered), where $\sigma,\beta$ are simple geodesics, such that $\measuredangle(\overrightarrow{\sigma},\overrightarrow{\alpha_2})=2\pi,\ l_{\sigma}=l_{\alpha_2}$ and $\overline{T}\cap\alpha_2=\{P\}$. Then $(P,\gamma_1,\gamma_2)$ is reducible.
\end{lemma}

\begin{proof}
    The proof is the same as Lemma \ref{lem-two noncross 2}. Cut along $\sigma,\alpha_2$ and identify $\sigma^+$ with $\alpha_2^-$, $\sigma^-$ with $\alpha_2^+$ (Figure \ref{fig: two noncross 3}). On the resulting surface, $\beta$ is non-dividing and trapped. Therefore $(P,\gamma_1,\gamma_2)$ is reducible by Lemma \ref{lem-noncross trap}.
\end{proof}

\begin{figure}
    \centering

\tikzset{every picture/.style={line width=0.75pt}} %set default line width to 0.75pt        

\begin{tikzpicture}[x=0.75pt,y=0.75pt,yscale=-1,xscale=1]
%uncomment if require: \path (0,339); %set diagram left start at 0, and has height of 339

%Curve Lines [id:da9680752590262601] 
\draw    (120,117) .. controls (138,115) and (156,113) .. (172,115) ;
\draw [shift={(151.97,114.26)}, rotate = 177.31] [color={rgb, 255:red, 0; green, 0; blue, 0 }  ][line width=0.75]    (10.93,-4.9) .. controls (6.95,-2.3) and (3.31,-0.67) .. (0,0) .. controls (3.31,0.67) and (6.95,2.3) .. (10.93,4.9)   ;
%Curve Lines [id:da31232577663135297] 
\draw    (120,117) .. controls (119,93) and (122,72) .. (126,46) ;
%Curve Lines [id:da8826059001704406] 
\draw    (172,115) .. controls (172,91) and (174,70) .. (179,50) ;
%Curve Lines [id:da8823993124774034] 
\draw    (140,157) .. controls (150,144) and (165,129) .. (172,115) ;
\draw [shift={(152.5,141.9)}, rotate = 310.03] [color={rgb, 255:red, 0; green, 0; blue, 0 }  ][line width=0.75]    (10.93,-4.9) .. controls (6.95,-2.3) and (3.31,-0.67) .. (0,0) .. controls (3.31,0.67) and (6.95,2.3) .. (10.93,4.9)   ;
%Curve Lines [id:da15943921241995818] 
\draw  [dash pattern={on 4.5pt off 4.5pt}]  (138,152) .. controls (152,139) and (163,126) .. (172,115) ;
%Curve Lines [id:da44582334693713943] 
\draw    (140,157) .. controls (167,163) and (193,182) .. (207,194) ;
%Curve Lines [id:da6114690088868742] 
\draw    (172,115) .. controls (196,119) and (222,136) .. (236,151) ;
%Curve Lines [id:da4988891227900091] 
\draw  [dash pattern={on 4.5pt off 4.5pt}]  (120,117) .. controls (129,129) and (135,141) .. (138,152) ;

% Text Node
\draw (139,95.4) node [anchor=north west][inner sep=0.75pt]    {$\alpha _{1}$};
% Text Node
\draw (161,138.4) node [anchor=north west][inner sep=0.75pt]    {$\alpha _{2}$};
% Text Node
\draw (142,127.4) node [anchor=north west][inner sep=0.75pt]    {$\sigma $};
% Text Node
\draw (115,131.4) node [anchor=north west][inner sep=0.75pt]    {$\beta $};
% Text Node
\draw (144,46.4) node [anchor=north west][inner sep=0.75pt]    {$U_{1}$};
% Text Node
\draw (203,156.4) node [anchor=north west][inner sep=0.75pt]    {$U_{2}$};

\end{tikzpicture}

    \caption{Cut along $\sigma,\alpha_2$.}
    \label{fig: two noncross 3}
\end{figure}

\section{Proof of Proposition \ref{prop-rm}}

\subsection{Flower neighborhoods in nice triples}

Let $(P,\gamma_1,\gamma_2)$ be a nice triple such that $r_P\geq \epsilon_1,\ r_{\gamma_i}\leq2\delta_1\ (i=1,2)$. By Lemma \ref{lem-ff1}, flower neighborhoods exist for both $\gamma_1$ and $\gamma_2$.

\begin{lemma}\label{lem-rc1}
    Let $C$ be a flower neighborhood of $\gamma_1$ or $\gamma_2$, then exactly one petal of $C$ is dividing.
\end{lemma}

\begin{proof}
    Without lost of generality assume that $C$ is a flower neighborhood of $\gamma_1$. Let $\beta$ be a path from $\gamma_2$ to $\gamma_1$ of length $\leq 10\delta_1$ that does not pass through $P$. Let $Q$ be the first intersection point of $\beta$ with $\overline{C}$. Then $Q$ lies on some petal $\alpha$ of $C$. This implies that $\alpha$ is dividing and any other petal is non-dividing.
\end{proof}

\begin{corollary}\label{cor-rc1}
    If a flower neighborhood of $\gamma_1$ or $\gamma_2$ has at least three petals, then $(P,\gamma_1,\gamma_2)$ is reducible.
\end{corollary}

\begin{proof}
    Since exactly one petal is dividing, if there are at least three petals, there are two consecutive ones that are non-dividing. By the flatness condition, the angle between them is in the range $[\frac{1}{2}\pi,\frac{3}{2}\pi]$. Now apply Lemma \ref{lem-two noncross 1}.
\end{proof}

Let $C_1,C_2$ be flower neighborhoods of $\gamma_1,\gamma_2$ respectively. By Lemma \ref{lem-fd1} we have $C_1\cap C_2=\varnothing$. We have either $\overline{C_1}\cap\overline{C_2}=P$, or $\overline{C_1}\cap\overline{C_2}$ is a common petal of $C_1$ and $C_2$ (they can share at most one common petal since every common petal is dividing).

\begin{lemma}\label{lem-rc2}
    Suppose that there are $2$-petal flower neighborhoods $C_1,C_2$ of $\gamma_1,\gamma_2$ respectively, such that they share a common petal. Then $(P,\gamma_1,\gamma_2)$ is reducible.
\end{lemma}

\begin{proof}
    Denote their common petal by $\alpha$, and denote the other petal of $C_i$ by $\beta_i$. Then $\alpha$ is dividing, and $\beta_1,\beta_2$ are non-dividing. Orient $\beta_i$ so that $C_i$ lies on the right. By the flatness condition of $C_1,C_2$, we deduce that $\measuredangle(\overleftarrow{\beta_1},\overrightarrow{\beta_2})\in[2\pi-\frac{2}{5},2\pi+\frac{2}{5}]$. If $\measuredangle(\overleftarrow{\beta_1},\overrightarrow{\beta_2})\leq2\pi$ then $(P,\gamma_1,\gamma_2)$ is reducible by Lemma \ref{lem-two noncross 2}. It remains to consider the case $\measuredangle(\overleftarrow{\beta_1},\overrightarrow{\beta_2})>2\pi$. By symmetry we may assume $l_{\beta_1}\geq l_{\beta_2}$. By Lemma \ref{lem-fs2}, there exists an embedded geodesic triangle $T\subset C_1$ with sides $\beta_1,\sigma_1,\sigma_2$ (clockwise ordered), where $\sigma_1,\sigma_2$ are simple geodesics, such that $l_{\sigma_1}=l_{\beta_2}$ and the angle between $\beta_1,\sigma_1$ is equal to $\measuredangle(\overleftarrow{\beta_1},\overrightarrow{\beta_2})-2\pi$. The conditions of Lemma \ref{lem-two noncross 3} are verified. Therefore $(P,\gamma_1,\gamma_2)$ is reducible.
\end{proof}

\begin{corollary}\label{cor-rc2}
    Suppose that there are flower neighborhoods $C_1,C_2$ of $\gamma_1,\gamma_2$ such that $\overline{C_1}\cap\overline{C_2}\neq\{P\}$, then either $(P,C_1,C_2)$ is reducible, or $S$ is a model.
\end{corollary}

\begin{proof}
    As already mentioned, if $\overline{C_1}\cap\overline{C_2}\neq\{P\}$, then $\overline{C_1}\cap\overline{C_2}$ is a common petal. By Corollary \ref{cor-rc1} we may assume that both $C_1,C_2$ have at most $2$ petals. They cannot both have $1$ petal, since then $C_1\cup\alpha\cup C_2$ is an annulus and $P$ is not a branch point. If both $C_1,C_2$ have $2$ petals then $(P,C_1,C_2)$ is reducible by Lemma \ref{lem-rc2}. If one of them has $1$ petal while the other has $2$ petals, then $S$ satisfies the conditions of Model 1.
\end{proof}

Now suppose that $(P,\gamma_1,\gamma_2)$ is irreducible. By Corollary \ref{cor-rc1}, every flower neighborhood of $\gamma_1$ or $\gamma_2$ has at most two petals. We are going to choose specific flower neighborhoods of $\gamma_1,\gamma_2$.

\begin{definition}
    Let $i\in\{1,2\}$. We say that $\gamma_i$ is of the second type if it has a $\frac{1}{100}$-flower neighborhood with two petals. Otherwise $\gamma_i$ is of the first type.
\end{definition}

For $i=1,2$, we choose a flower neighborhood $C_i$ of $\gamma_i$ as follows.
\begin{itemize}
    \item If $\gamma_i$ is of the second type, let $C_i$ be any $\frac{1}{100}$-flower neighborhood with two petals. Denote its petals by $\alpha_i$ and $\beta_i$, where $\alpha_i$ is dividing and $\beta_i$ is non-dividing.
    \item If $\gamma_i$ is of the first type, let $C_i$ be a $1$-petal flower neighborhood of $\gamma_i$. The petal is denoted $\alpha_i$.
\end{itemize}

(Therefore $C_i$ has $k$ petals if it is of the $k$-st(nd) type.) Orient $\alpha_i,\beta_i$ so that $C_i$ lies on the right. If $\gamma_i$ is of the second type we define $\theta_i=\frac{1}{20}$. Otherwise we define $\theta_i=\frac{1}{200}$.

\begin{lemma}\label{lem-rc3}
    For $i=1,2$, $C_i$ admits a $\theta_i$-extension at every petal.
\end{lemma}

\begin{proof}
    If $\gamma_i$ is of the first type, it has no $\frac{1}{100}$-flower neighborhood with 2 petals. If $\gamma_i$ is of the second type, it has no flower neighborhood with 3 petals. In either case the lemma follows from Lemma \ref{lem-ft1}.
\end{proof}

By Corollary \ref{cor-rc2}, we may assume that $\overline{C_1}\cap\overline{C_2}=P$. By Lemma \ref{lem-rc3}, $C_i$ admits a $\theta_i$-extension at every petal. In particular it admits a $\theta_i$-extension at $\alpha_i$. Denote by $E_i=C_i\cup\mathring{\alpha}_i\cup T_i$ the $\theta_i$-extension of $C_i$ at $\alpha_i$. There are two cases, either $E_1,E_2$ intersect or they do not intersect.

\subsection{The case $E_1\cap E_2\neq\varnothing$}

\begin{lemma}\label{lem-case1}
    There are geodesic triangles $T_1'\subset T_1,T_2'\subset T_2$ based on $\alpha_1,\alpha_2$ respectively, with base angle less than $\frac{1}{20}$, such that they share a common tip $Q$. Moreover, $E_i':=C_i\cup\mathring{\alpha}_i\cup T_i'$ satisfy $E_1'\cap E_2'=\varnothing$. (See Figure \ref{fig: Two triangles with a common tip}.)
\end{lemma}

We play the same old trick. Let $E_i(t)$ be the $t\theta_i$-extension of $C_i\ (t\in(0,1])$. Gradually increase $t$ until $E_1(t)$ and $E_2(t)$ first touch. The touching point is the point $Q$ we seek.

\begin{proof}
    For $t\in(0,1]$, let $T_i(t)\subset T_i$ be the $t\theta_i$-cap based on $\alpha_i$, and let $E_i(t)=C_i\cup\mathring{\alpha}_i\cup T_i(t)$. Since $\overline{C_1}\cap\overline{C_2}=P$, we have $E_1(t)\cap E_2(t)=\varnothing$ for $t$ small enough, while the assumption $E_1\cap E_2\neq\varnothing$ implies $E_1(1)\cap E_2(1)\neq\varnothing$. Let $t_0=\sup\{t:E_1(t)\cap E_2(t)=\varnothing\}$. Then $E_1(t_0)\cap E_2(t_0)=\varnothing$ and we have either $(\partial T_1(t_0)\backslash\alpha_1)\cap \partial E_2(t_0)\neq\varnothing$ or $(\partial T_2(t_0)\backslash\alpha_2)\cap \partial E_1(t_0)\neq\varnothing$. Without loss of generality assume that $(\partial T_1(t_0)\backslash\alpha_1)\cap \partial E_2(t_0)\neq\varnothing$. Let $Q\in(\partial T_1(t_0)\backslash\alpha_1)\cap \partial E_2(t_0)$.
    \begin{claim}
        $Q\in\partial T_2(t_0)\backslash\alpha_2$.
    \end{claim}
    If $\gamma_2$ is of the first type then $\partial E_2(t_0)=\partial T_2\backslash\mathring\alpha_2$. Since $Q\neq P$, the claim follows. If $\gamma_2$ is of the second type then $\partial E_2(t_0)=(\partial T_2\backslash\alpha_2)\cup\beta_2$. Since $\beta_2$ is non-dividing, we have $Q\notin\beta_2$ \footnote{If $Q\in\beta_2$, then there is a path of length $\leq10\delta_1$ that joins $\gamma_1,\gamma_2$ and intersects $\beta_2$ transversally at $Q$. This implies that $\beta_2$ is dividing, a contradiction.}. The claim also follows.

    Let $T_i'\subset T_i$ be the geodesic triangle based on $\alpha_i$ with tip $Q$. Then $E_1'\cap E_2'=\varnothing$ because $E_i'\subset E_i(t_0)$ and $E_1(t_0)\cap E_2(t_0)=\varnothing$. This proves Lemma \ref{lem-case1}.
\end{proof}

\subsubsection{Interlude: exponential map at a point}

We pause a little to discuss the exponential map at a point, and obtain results of the same kind as Lemma \ref{cor-fs1}. These will be applied to $\exp_Q$ ($Q$ being the point in Lemma \ref{lem-case1}) to get structural results in a neighborhood of $\alpha_1,\alpha_2$.

Let $S$ be a branched hyperbolic surface. Let $P,Q\in S$, where $P$ is a branch point while $Q$ is not a branch point. Suppose that $d(Q,\partial S)\geq\epsilon_1$ and $d(Q,P_1)\geq\epsilon_1$ for any branch point $P_1$ different from $P$. Put the hyperbolic metric on $T_QS$, given by the pullback via $T_QS\cong T_0\D\xrightarrow{\exp_0}\D$. Then the exponential map
\[\exp_Q:T_QS\supset\mathcal{U}_Q\rightarrow S\]
is a local isometry. $\mathcal{U}_Q\cap N_{\epsilon_1}(0)$ is of the form $N_{\epsilon_1}(0)\backslash\bigcup_{k=1}^mI_k$ for finitely many segments $I_k$ of the form $I_k=(r_k,\epsilon_1]\cdot e^{i\varphi_k}$. Let $B=\{r_ke^{i\varphi_k}\}_{k=1}^m$, a finite set of points in $N_{\epsilon_1}(0)$. Then $\exp_Q(x)=P$ for every $x\in B$.

\begin{figure}
    \centering

\tikzset{every picture/.style={line width=0.75pt}} %set default line width to 0.75pt        

\begin{tikzpicture}[x=0.75pt,y=0.75pt,yscale=-1,xscale=1]
%uncomment if require: \path (0,339); %set diagram left start at 0, and has height of 339

%Straight Lines [id:da5737977297899685] 
\draw    (100,127) -- (200,227) ;
%Straight Lines [id:da1939955600573826] 
\draw    (206,95) -- (200,227) ;
%Straight Lines [id:da31017160941156896] 
\draw    (326,161) -- (200,227) ;
%Straight Lines [id:da23807020819968838] 
\draw    (206,95) -- (242,155) ;
%Straight Lines [id:da5801953423403787] 
\draw    (242,155) -- (326,161) ;
%Straight Lines [id:da21161188107993922] 
\draw    (147,144) -- (188,126) ;
%Straight Lines [id:da06561199262735862] 
\draw    (100,127) -- (147,144) ;
%Straight Lines [id:da8352976195437775] 
\draw    (188,126) -- (206,95) ;
%Straight Lines [id:da8171829486508602] 
\draw    (100,127) -- (206,95) ;
%Straight Lines [id:da5857344224690998] 
\draw    (206,95) -- (326,161) ;

% Text Node
\draw (222,173.4) node [anchor=north west][inner sep=0.75pt]    {$V_{i}$};
% Text Node
\draw (163,163.4) node [anchor=north west][inner sep=0.75pt]    {$V_{i-1}$};
% Text Node
\draw (141,116.4) node [anchor=north west][inner sep=0.75pt]    {$K_{i-1}$};
% Text Node
\draw (246,129.4) node [anchor=north west][inner sep=0.75pt]    {$K_{i}$};
% Text Node
\draw (205,139.4) node [anchor=north west][inner sep=0.75pt]    {$\sigma _{i}$};
% Text Node
\draw (201,70.4) node [anchor=north west][inner sep=0.75pt]    {$x_{i}$};
% Text Node
\draw (186,226.4) node [anchor=north west][inner sep=0.75pt]    {$0$};
% Text Node
\draw (287,184.4) node [anchor=north west][inner sep=0.75pt]    {$\sigma _{i+1}$};
% Text Node
\draw (107,165.4) node [anchor=north west][inner sep=0.75pt]    {$\sigma _{i-1}$};
% Text Node
\draw (76,108.4) node [anchor=north west][inner sep=0.75pt]    {$x_{i-1}$};
% Text Node
\draw (330,144.4) node [anchor=north west][inner sep=0.75pt]    {$x_{i+1}$};
% Text Node
\draw (268,215.4) node [anchor=north west][inner sep=0.75pt]    {$\cdots $};
% Text Node
\draw (87,193.4) node [anchor=north west][inner sep=0.75pt]    {$\cdots $};

\end{tikzpicture}

    \caption{The vertices of $K_i$ are points in $B$. $V$ can be think of as a star-shaped neighborhood with piecewise geodesic boundary and branch points as vertices.}
    \label{fig: exp point1}
\end{figure}
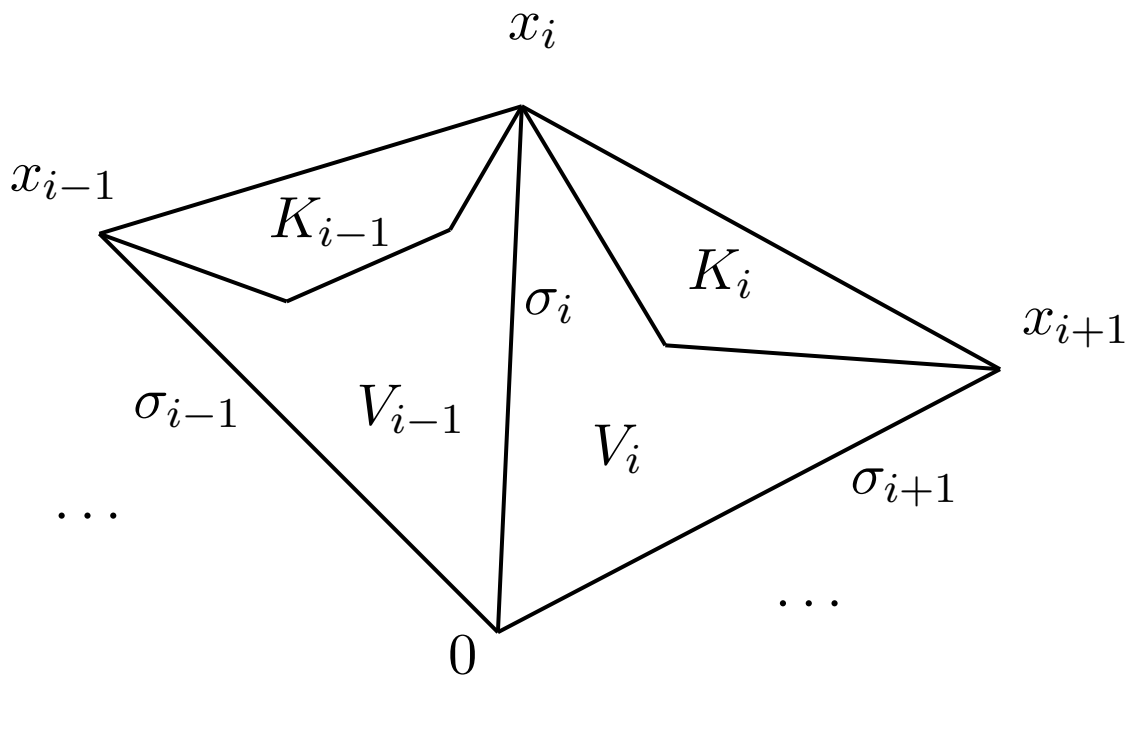

Let $\{x_1,\cdots,x_k\}\subset B\ (k\geq3)$ be a subset and let $\sigma_i$ be the geodesic segment from $0$ to $x_i$. Suppose that $\overrightarrow{\sigma_1},\cdots,\overrightarrow{\sigma_k}$ is arranged in clockwise order, and that $\measuredangle(\overrightarrow{\sigma_{i+1}},\overrightarrow{\sigma_{i}})<\pi\ (\forall i\in\Z/k\Z)$, so that $\sigma_i,\sigma_{i+1}$ span a geodesic triangle, which we denote by $T_i$. Let $K_i$ be the convex hull of $B\cap\overline{T_i}$, and let $V_i=\overline{T_i}\backslash K_i$. Let $V=\bigcup_iV_i$. Then $V\subset\mathcal{U}_Q$. Let $\sigma=\bigcup_i(\sigma_i\backslash\{x_i\})$. (See Figure \ref{fig: exp point1}.)

\begin{lemma}\label{lem-exp at a point}
    Notation as above. If $\exp_Q|_\sigma$ is injective, then $\exp_Q|_V$ is injective.
\end{lemma}

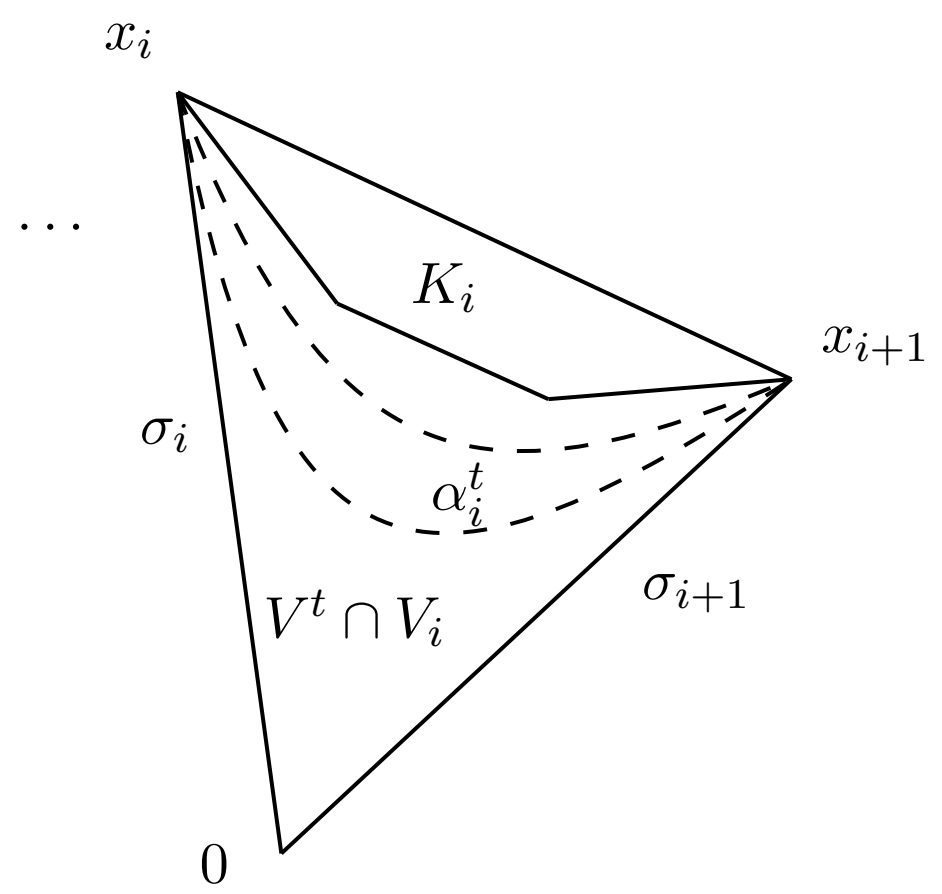
\begin{figure}
    \centering

\tikzset{every picture/.style={line width=0.75pt}} %set default line width to 0.75pt        

\begin{tikzpicture}[x=0.75pt,y=0.75pt,yscale=-1,xscale=1]
%uncomment if require: \path (0,339); %set diagram left start at 0, and has height of 339

%Straight Lines [id:da1939955600573826] 
\draw    (188,97) -- (214,288) ;
%Straight Lines [id:da31017160941156896] 
\draw    (342,169) -- (214,288) ;
%Straight Lines [id:da23807020819968838] 
\draw    (188,97) -- (228,150) ;
%Straight Lines [id:da5801953423403787] 
\draw    (281,174) -- (342,169) ;
%Straight Lines [id:da5857344224690998] 
\draw    (188,97) -- (342,169) ;
%Straight Lines [id:da9117562397549045] 
\draw    (228,150) -- (281,174) ;
%Curve Lines [id:da3185530519752483] 
\draw  [dash pattern={on 4.5pt off 4.5pt}]  (188,97) .. controls (227,198) and (266,201) .. (342,169) ;
%Curve Lines [id:da06609814587654839] 
\draw  [dash pattern={on 4.5pt off 4.5pt}]  (188,97) .. controls (214,232) and (253,227) .. (342,169) ;

% Text Node
\draw (245,138.4) node [anchor=north west][inner sep=0.75pt]    {$K_{i}$};
% Text Node
\draw (177,177.4) node [anchor=north west][inner sep=0.75pt]    {$\sigma _{i}$};
% Text Node
\draw (168,78.4) node [anchor=north west][inner sep=0.75pt]    {$x_{i}$};
% Text Node
\draw (192,284.4) node [anchor=north west][inner sep=0.75pt]    {$0$};
% Text Node
\draw (303,216.4) node [anchor=north west][inner sep=0.75pt]    {$\sigma _{i+1}$};
% Text Node
\draw (348,154.4) node [anchor=north west][inner sep=0.75pt]    {$x_{i+1}$};
% Text Node
\draw (146,126.4) node [anchor=north west][inner sep=0.75pt]    {$\cdots $};
% Text Node
\draw (250,188.4) node [anchor=north west][inner sep=0.75pt]    {$\alpha _{i}^{t}$};
% Text Node
\draw (208,220.4) node [anchor=north west][inner sep=0.75pt]    {$V^{t} \cap V_{i}$};

\end{tikzpicture}

    \caption{Take an increasing family of domains $V^t$, such that their boundaries are piecewise smooth, with vertices $x_i$, and the curvature at smooth points are outward pointing.}
    \label{fig: exp point2}
\end{figure}

\begin{proof}
    Take a continuous, increasing family of domains $\sigma\subset V^t\subset V\ (t\in(0,1])$ such that $V^1=V,\ \bigcap_{t\in(0,1]} V^t=\sigma$, and that for every $1\leq i\leq k,\ t\in(0,1)$, $\partial V^t\cap V_i$ is a smooth curve $\alpha_i^t$ connecting $x_i,x_{i+1}$ with outward pointing curvature (see Figure \ref{fig: exp point2}). By assumption $\exp_Q|_\sigma$ is injective, so that $\exp_Q|_{V^t}$ is injective for $t$ small enough. Let $t_0=\sup\{t\in(0,1]:\exp_Q|_{V^t}\text{ is injective}\}$. Then $\exp_Q|_{V^{t_0}}$ is injective. Suppose that $\exp_Q|_V$ is not injective, then $t_0<1$. Let $\alpha_i=\partial V^{t_0}\cap V_i$. Then $\bigcup_i\alpha_i$ is contained in the interior of $\mathcal{U}_Q$ (by convention $\alpha_i^t$ do not include their end points), on which $\exp_Q$ is a local diffeomorphism. Since the curvature of $\alpha_i$ is outward pointing, we deduce that $\exp_Q$ is injective on ${\bigcup}_i\alpha_i$, and that $d\exp_Q$ maps $\{\overleftarrow{\alpha_i},\overrightarrow{\alpha_i}\}_{i=1}^k$ to different tangent vectors at $P$. This implies that $\exp_Q|_{V^t}$ is injective for some $t>t_0$, a contradiction.
\end{proof}

\begin{figure}
    \centering

\tikzset{every picture/.style={line width=0.75pt}} %set default line width to 0.75pt        

\begin{tikzpicture}[x=0.75pt,y=0.75pt,yscale=-1,xscale=1]
%uncomment if require: \path (0,339); %set diagram left start at 0, and has height of 339

%Straight Lines [id:da8205412327211598] 
\draw    (106,49) -- (176,95) ;
%Straight Lines [id:da899342505638494] 
\draw    (176,95) -- (242,90) ;
%Straight Lines [id:da5294431821282404] 
\draw    (82,114) -- (176,95) ;
%Straight Lines [id:da4630766712318699] 
\draw    (176,95) -- (231,147) ;
%Straight Lines [id:da00010870578372013107] 
\draw    (82,114) -- (231,147) ;
%Straight Lines [id:da5823233658495346] 
\draw    (106,49) -- (242,90) ;
%Straight Lines [id:da12239765114916823] 
\draw    (322,141) -- (451,133) ;
%Straight Lines [id:da3443272085475454] 
\draw    (361,66) -- (451,133) ;
%Straight Lines [id:da06113844919540845] 
\draw    (361,66) -- (375,110) ;
%Straight Lines [id:da7751558577413996] 
\draw    (375,110) -- (451,133) ;
%Straight Lines [id:da3042271024181552] 
\draw    (322,141) -- (375,110) ;

% Text Node
\draw (154,109.4) node [anchor=north west][inner sep=0.75pt]    {$T'_{1}$};
% Text Node
\draw (167,74.4) node [anchor=north west][inner sep=0.75pt]    {$T'_{2}$};
% Text Node
\draw (212,114.4) node [anchor=north west][inner sep=0.75pt]    {$\sigma _{1}$};
% Text Node
\draw (104,92.4) node [anchor=north west][inner sep=0.75pt]    {$\sigma _{2}$};
% Text Node
\draw (117,67.4) node [anchor=north west][inner sep=0.75pt]    {$\sigma _{3}$};
% Text Node
\draw (202,95.4) node [anchor=north west][inner sep=0.75pt]    {$\sigma _{4}$};
% Text Node
\draw (141,139.4) node [anchor=north west][inner sep=0.75pt]    {$\alpha _{1}$};
% Text Node
\draw (163,49.4) node [anchor=north west][inner sep=0.75pt]    {$\alpha _{2}$};
% Text Node
\draw (150,85.4) node [anchor=north west][inner sep=0.75pt]    {$Q$};
% Text Node
\draw (373,141.4) node [anchor=north west][inner sep=0.75pt]    {$\alpha _{1}$};
% Text Node
\draw (370,112.4) node [anchor=north west][inner sep=0.75pt]    {$\sigma _{1} =\sigma _{4}$};
% Text Node
\draw (328,110.4) node [anchor=north west][inner sep=0.75pt]    {$\sigma _{2}$};
% Text Node
\draw (345,78.4) node [anchor=north west][inner sep=0.75pt]    {$\sigma _{3}$};
% Text Node
\draw (400,77.4) node [anchor=north west][inner sep=0.75pt]    {$\alpha _{2}$};

\end{tikzpicture}

    \caption{Two triangles with a common tip.}
    \label{fig: Two triangles with a common tip}
\end{figure}
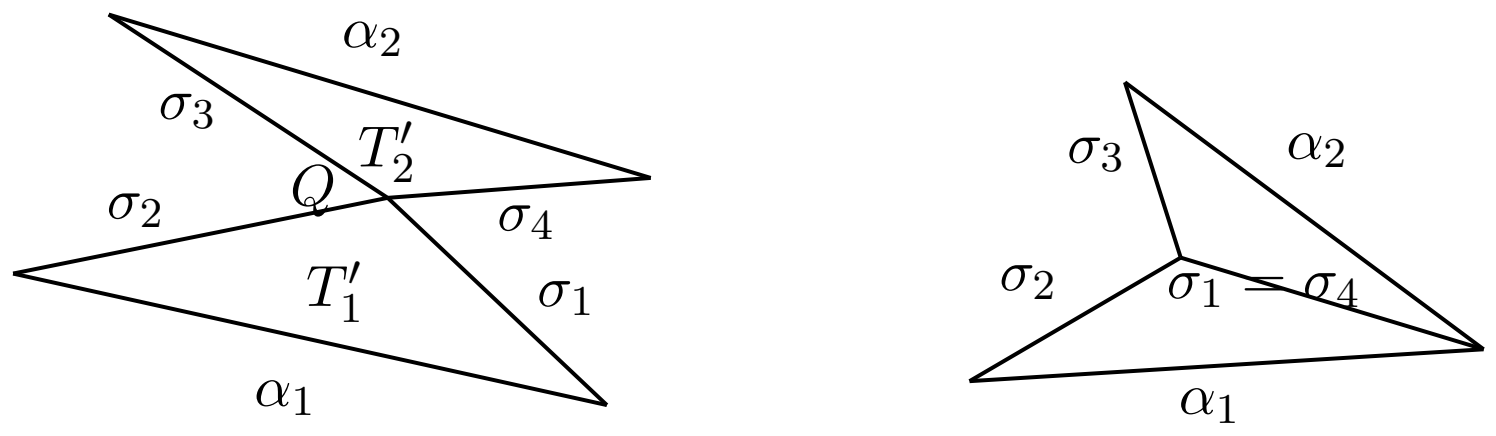

We now return to the discussion of the nice triple $(P,\gamma_1,\gamma_2)$. Recall that we have disjoint, embedded geodesic triangles $T_1',T_2'\subset S$ sharing a common tip $Q$, which is not a branch point since it lies in the interior of $T_1$. Let $\sigma_1,\sigma_2,\sigma_3,\sigma_4$ be the four geodesics from $Q$ to $P$, with $\sigma_1,\sigma_2$ sides of $T_1'$ and $\sigma_3,\sigma_4$ sides of $T_2'$, such that $\overrightarrow{\sigma_1},\overrightarrow{\sigma_2},\overrightarrow{\sigma_3},\overrightarrow{\sigma_4}$ are arranged in clockwise order (possibly with $\sigma_2=\sigma_3$ or $\sigma_1=\sigma_4$; see Figure \ref{fig: Two triangles with a common tip}). We have $\measuredangle(\overrightarrow{\sigma_{i+1}},\overrightarrow{\sigma_i})<\pi$ because the tip angles of $T_1',T_2'$ are  greater than $\frac{\pi}{2}$.

\subsubsection{Case 1: $\sigma_2=\sigma_3$ or $\sigma_1=\sigma_4$.}

Without loss of generality assume that $\sigma_1=\sigma_4$.

\begin{lemma}\label{lem-case11}
    $\alpha_1,\alpha_2$ are sides of an embedded geodesic triangle in $S$, whose third side is a simple geodesic loop based at $P$.
\end{lemma}

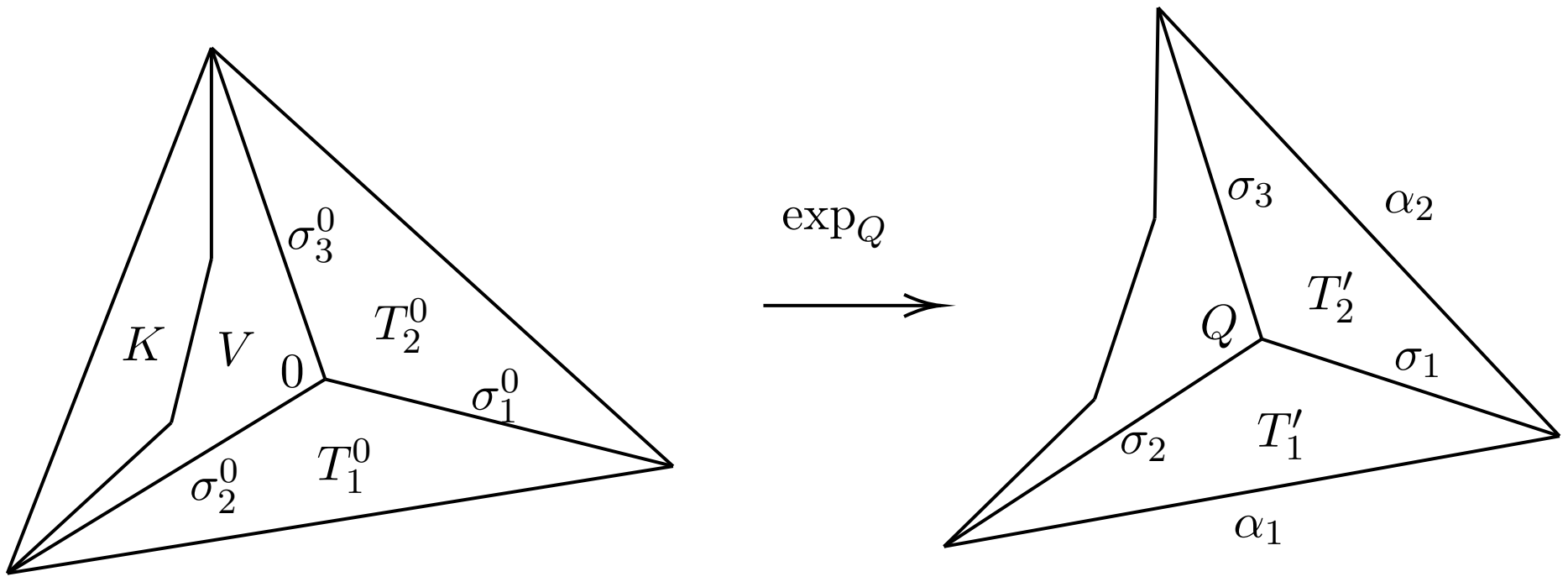
\begin{figure}
    \centering

\tikzset{every picture/.style={line width=0.75pt}} %set default line width to 0.75pt        

\begin{tikzpicture}[x=0.75pt,y=0.75pt,yscale=-1,xscale=1]
%uncomment if require: \path (0,339); %set diagram left start at 0, and has height of 339

%Straight Lines [id:da07164198880527495] 
\draw    (138,76) -- (276,201) ;
%Straight Lines [id:da869386428778581] 
\draw    (77,233) -- (276,201) ;
%Straight Lines [id:da8172892964403181] 
\draw    (138,76) -- (77,233) ;
%Straight Lines [id:da3226170184365308] 
\draw    (172,175) -- (276,201) ;
%Straight Lines [id:da12723653276292024] 
\draw    (138,76) -- (172,175) ;
%Straight Lines [id:da5330094505160677] 
\draw    (77,233) -- (172,175) ;
%Straight Lines [id:da008783313967596462] 
\draw    (138,76) -- (138,139) ;
%Straight Lines [id:da4315297870489311] 
\draw    (138,139) -- (126,188) ;
%Straight Lines [id:da5300076041897898] 
\draw    (77,233) -- (126,188) ;
%Straight Lines [id:da678043170697271] 
\draw    (303,153) -- (354,153) ;
\draw [shift={(356,153)}, rotate = 180] [color={rgb, 255:red, 0; green, 0; blue, 0 }  ][line width=0.75]    (10.93,-3.29) .. controls (6.95,-1.4) and (3.31,-0.3) .. (0,0) .. controls (3.31,0.3) and (6.95,1.4) .. (10.93,3.29)   ;
%Straight Lines [id:da41126222732734596] 
\draw    (421,64) -- (541,192) ;
%Straight Lines [id:da18627052350833562] 
\draw    (357,225) -- (541,192) ;
%Straight Lines [id:da35440047630337745] 
\draw    (452,163) -- (541,192) ;
%Straight Lines [id:da12040474482897712] 
\draw    (357,225) -- (452,163) ;
%Straight Lines [id:da6575701081455728] 
\draw    (421,64) -- (452,163) ;
%Straight Lines [id:da7994333624134127] 
\draw    (421,64) -- (420,127) ;
%Straight Lines [id:da9577606156484773] 
\draw    (420,127) -- (402,181) ;
%Straight Lines [id:da7677465431498848] 
\draw    (402,181) -- (357,225) ;

% Text Node
\draw (109.5,157.9) node [anchor=north west][inner sep=0.75pt]    {$K$};
% Text Node
\draw (138,159.4) node [anchor=north west][inner sep=0.75pt]    {$V$};
% Text Node
\draw (159,122.4) node [anchor=north west][inner sep=0.75pt]    {$\sigma _{3}^{0}$};
% Text Node
\draw (130,197.4) node [anchor=north west][inner sep=0.75pt]    {$\sigma _{2}^{0}$};
% Text Node
\draw (214,170.4) node [anchor=north west][inner sep=0.75pt]    {$\sigma _{1}^{0}$};
% Text Node
\draw (168,191.4) node [anchor=north west][inner sep=0.75pt]    {$T_{1}^{0}$};
% Text Node
\draw (185,149.4) node [anchor=north west][inner sep=0.75pt]    {$T_{2}^{0}$};
% Text Node
\draw (307,122.4) node [anchor=north west][inner sep=0.75pt]    {$\exp_{Q}$};
% Text Node
\draw (432,151.4) node [anchor=north west][inner sep=0.75pt]    {$Q$};
% Text Node
\draw (449,182.4) node [anchor=north west][inner sep=0.75pt]    {$T'_{1}$};
% Text Node
\draw (464,141.4) node [anchor=north west][inner sep=0.75pt]    {$T'_{2}$};
% Text Node
\draw (490,164.4) node [anchor=north west][inner sep=0.75pt]    {$\sigma _{1}$};
% Text Node
\draw (440,113.4) node [anchor=north west][inner sep=0.75pt]    {$\sigma _{3}$};
% Text Node
\draw (408,189.4) node [anchor=north west][inner sep=0.75pt]    {$\sigma _{2}$};
% Text Node
\draw (157,166.4) node [anchor=north west][inner sep=0.75pt]    {$0$};
% Text Node
\draw (442,214.4) node [anchor=north west][inner sep=0.75pt]    {$\alpha _{1}$};
% Text Node
\draw (487,117.4) node [anchor=north west][inner sep=0.75pt]    {$\alpha _{2}$};

\end{tikzpicture}

    \caption{All vertices of the polygon $V\cup T_1^0\cup T_2^0$ are mapped to $P$, and all sides (possibly further partitioned by points in $B_0$) are mapped to simple geodesic loops based at $P$. The lemma shows that in fact $K$ reduces to a line segment, and $V$ is the triangle spanned by $\sigma_2^0,\sigma_3^0$. Moreover, no interior point of the third side of $V$ is mapped to $P$.}
    \label{fig: lem-case11}
\end{figure}

\begin{proof}
    Let $\sigma_i^0\subset T_QS\ (i=1,2,3)$ be the geodesic segments such that $\exp_Q(\sigma_i^0)=\sigma_i$, and let $x_i$ be their endpoints other than $0$. Let $T\subset T_QS$ be the geodesic triangle with vertices $x_1,x_2,x_3$. We show that $\exp_Q$ is defined and injective on $T$. Then $\exp_Q(T)$ would be the desired triangle. Let $T_1^0$ (resp. $T_2^0$) be the geodesic triangle in $T_Q S$ spanned by $\sigma_1^0,\sigma_2^0$ (resp. $\sigma_1^0,\sigma_3^0$). Then $\exp_Q|_{T_i^0}$ is injective and $\exp_Q(T_i^0)=T_i'\ (i=1,2)$. Let $U\subset T_QS$ be the geodesic triangle spanned by $\sigma_2^0,\sigma_3^0$. Let $B=\{x\in\overline{U}\cap\mathcal{U}_Q:\exp_Q(x)=P\}$. Let $K$ be the convex hull of $B$ and let $V=\overline{U}\backslash K$. (See Figure \ref{fig: lem-case11}.)
    
    By Lemma \ref{lem-exp at a point}, $\exp_Q$ is injective on $V\cup T_1^0\cup T_2^0\cup\sigma_1^0$. $\partial K\cap\partial V$ is a piecewise geodesic with vertices in $B$. $B_0:=B\cap\partial V$ divides $\partial K\cap\partial V$ into geodesic arcs that are mapped to simple geodesic loops based at $P$ which are non-dividing (since $\alpha_1$ and $\alpha_2$ are connected by a path not passing through $\exp_Q(\partial V\cap\partial K)$).
    
    If $\#B\geq 3$, it is easy to see that $\#B_0\geq3$, so that $B_0$ divides $\partial K\cap\partial V$ into at least two geodesic arcs. Any two consecutive geodesic arcs would satisfy the conditions of Lemma \ref{cor-two noncross}, which implies that $(P,\gamma_1,\gamma_2)$ is reducible.
    
    Therefore $\#B=2$, consisting only of the two vertices of $U$ other than $0$. This implies that $T=V\cup T_1^0\cup T_2^0\cup\sigma_1^0\backslash\{x_1\}$, on which $\exp_Q$ is injective, and that $\exp_Q$ maps the third side of $U$ to a simple geodesic loop based at $P$. Therefore $\exp_Q(T)$ is the desired triangle.
\end{proof}

\begin{lemma}\label{lem-case12}
    Both $\gamma_1,\gamma_2$ are of the first type.
\end{lemma}

\begin{figure}
    \centering

\tikzset{every picture/.style={line width=0.75pt}} %set default line width to 0.75pt        

\begin{tikzpicture}[x=0.75pt,y=0.75pt,yscale=-1,xscale=1]
%uncomment if require: \path (0,339); %set diagram left start at 0, and has height of 339

%Straight Lines [id:da03672981700189315] 
\draw    (100,102) -- (167,115) ;
\draw [shift={(139.39,109.64)}, rotate = 190.98] [color={rgb, 255:red, 0; green, 0; blue, 0 }  ][line width=0.75]    (10.93,-3.29) .. controls (6.95,-1.4) and (3.31,-0.3) .. (0,0) .. controls (3.31,0.3) and (6.95,1.4) .. (10.93,3.29)   ;
%Straight Lines [id:da5495733478241059] 
\draw    (167,115) -- (271,102) ;
\draw [shift={(224.95,107.76)}, rotate = 172.87] [color={rgb, 255:red, 0; green, 0; blue, 0 }  ][line width=0.75]    (10.93,-3.29) .. controls (6.95,-1.4) and (3.31,-0.3) .. (0,0) .. controls (3.31,0.3) and (6.95,1.4) .. (10.93,3.29)   ;
%Straight Lines [id:da20691443994859682] 
\draw    (157,68) -- (167,115) ;
%Straight Lines [id:da2608843790726766] 
\draw    (157,68) -- (271,102) ;
%Curve Lines [id:da12440522671912202] 
\draw    (100,102) .. controls (107,127) and (108,147) .. (107,166) ;
%Curve Lines [id:da342220840774607] 
\draw    (271,102) .. controls (274,128) and (275,147) .. (273,168) ;

% Text Node
\draw (174,143.4) node [anchor=north west][inner sep=0.75pt]    {$C_{1}$};
% Text Node
\draw (213,113.4) node [anchor=north west][inner sep=0.75pt]    {$\alpha _{1}$};
% Text Node
\draw (183,87.4) node [anchor=north west][inner sep=0.75pt]    {$T$};
% Text Node
\draw (123,114.4) node [anchor=north west][inner sep=0.75pt]    {$\beta _{1}$};
% Text Node
\draw (148,84.4) node [anchor=north west][inner sep=0.75pt]    {$\eta $};
% Text Node
\draw (213,66.4) node [anchor=north west][inner sep=0.75pt]    {$\alpha _{2}$};
% Text Node
\draw (168,99.4) node [anchor=north west][inner sep=0.75pt]    {$\theta $};

\end{tikzpicture}

    \caption{Both $\beta_1$ and $\eta$ are non-dividing. If $\theta\leq\pi-\frac{1}{50}$, then $\beta_1,\eta$ imply reducibility. If $\theta>\pi-\frac{1}{50}$, then $C_2\cup\mathring\alpha_2\cup T$ is a $2$-petal flower neighborhood of $\gamma_2$ which shares a common petal with $C_1$, implying reducibility.}
    \label{fig: lem-case12}
\end{figure}
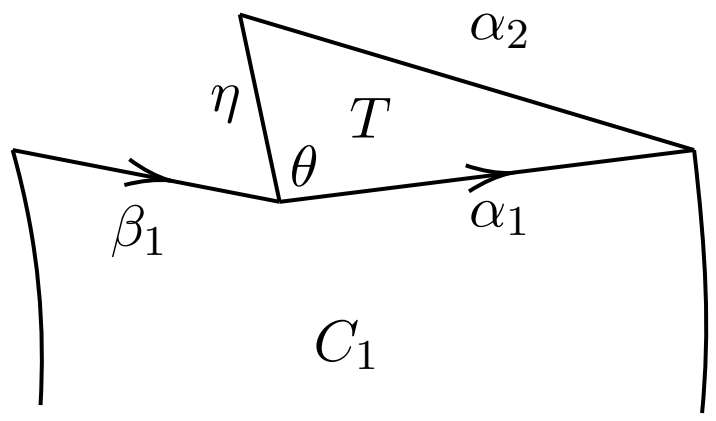

\begin{proof}
    Suppose that one of them, say $\gamma_1$, is of the second type. By the $\frac{1}{100}$-flatness of $C_1$, we have $\measuredangle(\overleftarrow{\beta_1},\overrightarrow{\alpha_1})\in[\pi-\frac{1}{50},\pi+\frac{1}{50}]$. Let $T$ be the geodesic triangle in Lemma \ref{lem-case11} and let $\eta$ be its third side. Let $\theta$ be the interior angle of $T$ opposite to $\alpha_2$. If $\theta\leq\pi-\frac{1}{50}$ then $\beta_1,\eta$ satisfies the conditions of Lemma \ref{cor-two noncross}, so that $(P,\gamma_1,\gamma_2)$ is reducible, a contradiction (see Figure \ref{fig: lem-case12}). Therefore $\theta>\pi-\frac{1}{50}$, so that the two interior angles of $T$ at the vertices of $\alpha_2$ are less than $\frac{1}{50}$. But $C_2$ admits a $\theta_2$ extension at $\alpha_2$. We deduce that $\theta_2<\frac{1}{50}$. Therefore $\gamma_2$ is of the first type, in particular $C_2$ is a $1$-petal flower. But then $C_2':=C_2\cup\mathring\alpha_2\cup T$ is a $2$-petal flower neighborhood of $\gamma_2$, which shares a common petal $\alpha_1$ with $C_1$. Therefore $(P,\gamma_1,\gamma_2)$ is reducible by Lemma \ref{lem-rc2}, a contradiction.
\end{proof}

By Lemma \ref{lem-case11} and Lemma \ref{lem-case12}, $S$ satisfies the conditions of Model 2.

\subsubsection{Case 2: $\sigma_2\neq\sigma_3$ and $\sigma_1\neq\sigma_4$.}

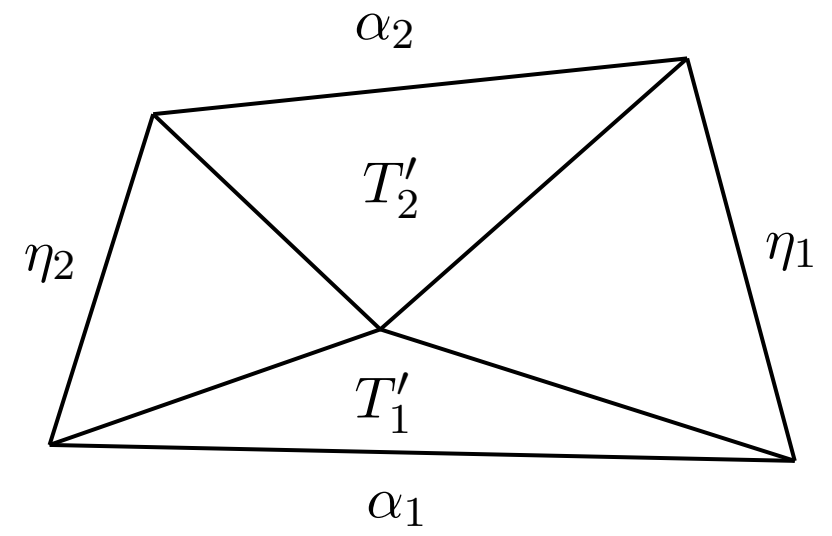
\begin{figure}
    \centering

\tikzset{every picture/.style={line width=0.75pt}} %set default line width to 0.75pt        

\begin{tikzpicture}[x=0.75pt,y=0.75pt,yscale=-1,xscale=1]
%uncomment if require: \path (0,339); %set diagram left start at 0, and has height of 339

%Straight Lines [id:da8251705101706773] 
\draw    (81,74) -- (138,128) ;
%Straight Lines [id:da8410389590605615] 
\draw    (138,128) -- (215,60) ;
%Straight Lines [id:da09799184949071171] 
\draw    (81,74) -- (215,60) ;
%Straight Lines [id:da24811092219783448] 
\draw    (138,128) -- (242,161) ;
%Straight Lines [id:da19874839073537032] 
\draw    (55,157) -- (138,128) ;
%Straight Lines [id:da3761894117673704] 
\draw    (81,74) -- (55,157) ;
%Straight Lines [id:da5552595252891289] 
\draw    (215,60) -- (242,161) ;
%Straight Lines [id:da8814255977682028] 
\draw    (55,157) -- (242,161) ;

% Text Node
\draw (130,137.4) node [anchor=north west][inner sep=0.75pt]    {$T'_{1}$};
% Text Node
\draw (132,83.4) node [anchor=north west][inner sep=0.75pt]    {$T'_{2}$};
% Text Node
\draw (133,167.4) node [anchor=north west][inner sep=0.75pt]    {$\alpha _{1}$};
% Text Node
\draw (130,47.4) node [anchor=north west][inner sep=0.75pt]    {$\alpha _{2}$};
% Text Node
\draw (47,105.4) node [anchor=north west][inner sep=0.75pt]    {$\eta _{2}$};
% Text Node
\draw (233,102.4) node [anchor=north west][inner sep=0.75pt]    {$\eta _{1}$};

\end{tikzpicture}

    \caption{If $\sigma_2\neq\sigma_3$ and $\sigma_1\neq\sigma_4$, then $\alpha_1,\alpha_2$ are opposite sides of an embedded geodesic quadrilateral $R$ (disjoint from $C_1,C_2$).}
    \label{fig: case21}
\end{figure}

\begin{lemma}
    There is an embedded geodesic quadrilateral $R\subset S$ with sides $\alpha_1,\eta_2,\alpha_2,\eta_1$ (clockwise ordered), where $\eta_1,\eta_2$ are simple geodesic loops based at $P$ (see Figure \ref{fig: case21}).
\end{lemma}

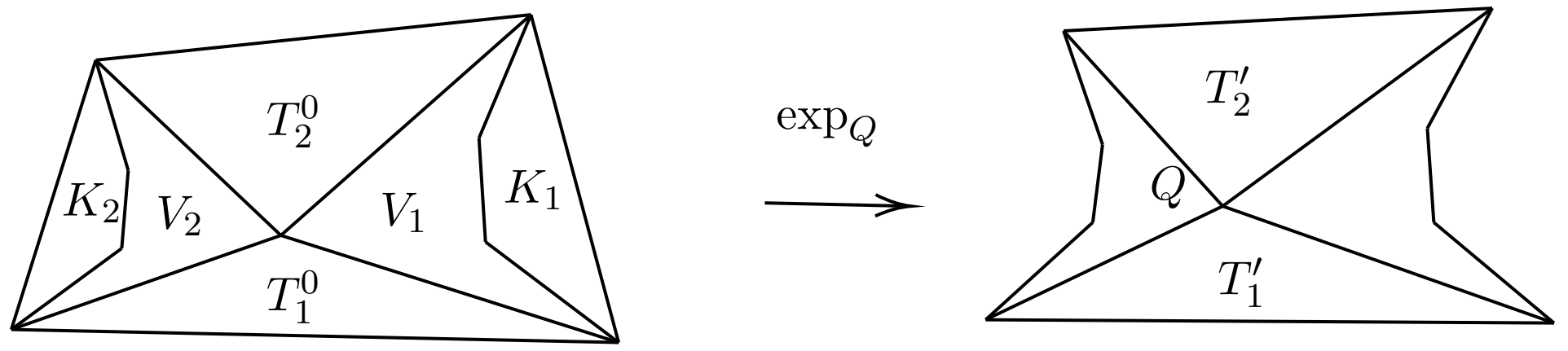
\begin{figure}
    \centering

\tikzset{every picture/.style={line width=0.75pt}} %set default line width to 0.75pt        

\begin{tikzpicture}[x=0.75pt,y=0.75pt,yscale=-1,xscale=1]
%uncomment if require: \path (0,339); %set diagram left start at 0, and has height of 339

%Straight Lines [id:da8251705101706773] 
\draw    (81,74) -- (138,128) ;
%Straight Lines [id:da8410389590605615] 
\draw    (138,128) -- (215,60) ;
%Straight Lines [id:da09799184949071171] 
\draw    (81,74) -- (215,60) ;
%Straight Lines [id:da24811092219783448] 
\draw    (138,128) -- (242,161) ;
%Straight Lines [id:da19874839073537032] 
\draw    (55,157) -- (138,128) ;
%Straight Lines [id:da3761894117673704] 
\draw    (81,74) -- (55,157) ;
%Straight Lines [id:da5552595252891289] 
\draw    (215,60) -- (242,161) ;
%Straight Lines [id:da8814255977682028] 
\draw    (55,157) -- (242,161) ;
%Straight Lines [id:da853327412157102] 
\draw    (81,74) -- (91,108) ;
%Straight Lines [id:da044055718296720614] 
\draw    (91,108) -- (89,132) ;
%Straight Lines [id:da5765221703656138] 
\draw    (89,132) -- (55,157) ;
%Straight Lines [id:da8011625591232401] 
\draw    (215,60) -- (199,98) ;
%Straight Lines [id:da4868218901930528] 
\draw    (199,98) -- (201,130) ;
%Straight Lines [id:da7626780579082479] 
\draw    (201,130) -- (242,161) ;
%Straight Lines [id:da22179857433392758] 
\draw    (287,118) -- (330,118.96) ;
\draw [shift={(332,119)}, rotate = 181.27] [color={rgb, 255:red, 0; green, 0; blue, 0 }  ][line width=0.75]    (10.93,-3.29) .. controls (6.95,-1.4) and (3.31,-0.3) .. (0,0) .. controls (3.31,0.3) and (6.95,1.4) .. (10.93,3.29)   ;
%Straight Lines [id:da6770872005163486] 
\draw    (379,65) -- (428,119) ;
%Straight Lines [id:da9717333550606836] 
\draw    (428,119) -- (511,58) ;
%Straight Lines [id:da8830144623134014] 
\draw    (379,65) -- (511,58) ;
%Straight Lines [id:da18227407886903302] 
\draw    (428,119) -- (530,155) ;
%Straight Lines [id:da21788532229770452] 
\draw    (355,154) -- (428,119) ;
%Straight Lines [id:da8704366215044338] 
\draw    (355,154) -- (530,155) ;
%Straight Lines [id:da5331057108723736] 
\draw    (379,65) -- (391,100) ;
%Straight Lines [id:da19321267150906074] 
\draw    (391,100) -- (388,124) ;
%Straight Lines [id:da24955479946655057] 
\draw    (388,124) -- (355,154) ;
%Straight Lines [id:da20756941676792295] 
\draw    (511,58) -- (491,95) ;
%Straight Lines [id:da7037045780306108] 
\draw    (491,95) -- (493,124) ;
%Straight Lines [id:da47378299337075547] 
\draw    (493,124) -- (530,155) ;

% Text Node
\draw (132,137.4) node [anchor=north west][inner sep=0.75pt]    {$T_{1}^{0}$};
% Text Node
\draw (132,83.4) node [anchor=north west][inner sep=0.75pt]    {$T_{2}^{0}$};
% Text Node
\draw (98,114.4) node [anchor=north west][inner sep=0.75pt]    {$V_{2}$};
% Text Node
\draw (69,110.4) node [anchor=north west][inner sep=0.75pt]    {$K_{2}$};
% Text Node
\draw (167,113.4) node [anchor=north west][inner sep=0.75pt]    {$V_{1}$};
% Text Node
\draw (205,106.4) node [anchor=north west][inner sep=0.75pt]    {$K_{1}$};
% Text Node
\draw (289,87.4) node [anchor=north west][inner sep=0.75pt]    {$\exp_{Q}$};
% Text Node
\draw (425,133.4) node [anchor=north west][inner sep=0.75pt]    {$T'_{1}$};
% Text Node
\draw (421,74.4) node [anchor=north west][inner sep=0.75pt]    {$T'_{2}$};
% Text Node
\draw (404,105.4) node [anchor=north west][inner sep=0.75pt]    {$Q$};

\end{tikzpicture}

    \caption{By the same argument as before, $K_1,K_2$ reduce to line segments.}
    \label{fig: case22}
\end{figure}

\begin{proof}
    The proof follows the same line as that of Lemma \ref{lem-case11}. Let $\sigma_i^0\subset T_QS\ (i=1,2,3,4)$ be the geodesic segments such that $\exp_Q(\sigma_i^0)=\sigma_i$. Let $T_1^0$ (resp. $T_2^0$) be the geodesic triangle in $T_QS$ spanned by $\sigma_1^0,\sigma_2^0$ (resp. $\sigma_3^0,\sigma_4^0$). Let $U_1$ (resp. $U_2$) be the geodesic triangle in $T_QS$ spanned by $\sigma_1^0,\sigma_4^0$ (resp. $\sigma_2^0,\sigma_3^0$). Let $B_i=\{x\in\overline{U_i}\cap\mathcal{U}_Q:\exp_Q(x)=P\}\ (i=1,2)$. Let $K_i$ be the convex hull of $B_i$ and let $V_i=\overline{U_i}\backslash K_i$ (see Figure \ref{fig: case22}). Then by Lemma \ref{lem-exp at a point}, $\exp_Q$ is injective on $V_1\cup V_2\cup T_1^0\cup T_2^0$. By the same argument as the proof of Lemma \ref{lem-case11}, we deduce that $\#B_1=\#B_2=2$, consisting only of the two vertices of $U_1,U_2$ other than $0$ respectively. Therefore $\exp_Q$ is injective on $U_1\cup U_2\cup T_1^0\cup T_2^0$, the image gives the desired quadrilateral.
\end{proof}

\begin{remark}
    $R$ need not be convex.
\end{remark}

Both $\eta_1,\eta_2$ are non-dividing, and $R$ lies in the same component of $N_{\frac{1}{2}}(P)\backslash\eta_i$ as $\gamma_1,\gamma_2$, in particular $\eta_1\neq\eta_2$. Orient $\eta_1,\eta_2$ so that $R$ lies on the right.

\begin{lemma}
    Both $\gamma_1$ and $\gamma_2$ are of the first type.
\end{lemma}

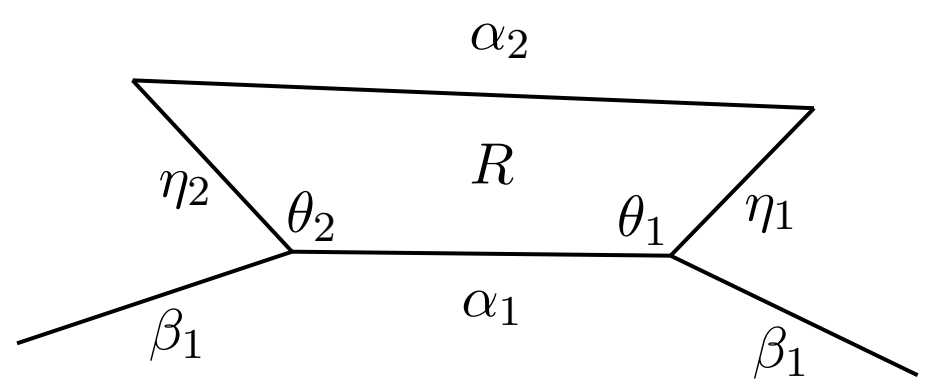
\begin{figure}
    \centering

\tikzset{every picture/.style={line width=0.75pt}} %set default line width to 0.75pt        

\begin{tikzpicture}[x=0.75pt,y=0.75pt,yscale=-1,xscale=1]
%uncomment if require: \path (0,339); %set diagram left start at 0, and has height of 339

%Straight Lines [id:da11540913237952921] 
\draw    (159,124) -- (254,125) ;
%Straight Lines [id:da047024082259646494] 
\draw    (90,147) -- (159,124) ;
%Straight Lines [id:da993633895265025] 
\draw    (254,125) -- (316,155) ;
%Straight Lines [id:da8211570622814811] 
\draw    (119,81) -- (159,124) ;
%Straight Lines [id:da12188509232727618] 
\draw    (290,88) -- (254,125) ;
%Straight Lines [id:da7291976221259564] 
\draw    (119,81) -- (290,88) ;

% Text Node
\draw (202,95.4) node [anchor=north west][inner sep=0.75pt]    {$R$};
% Text Node
\draw (156,107.4) node [anchor=north west][inner sep=0.75pt]    {$\theta _{2}$};
% Text Node
\draw (239,108.4) node [anchor=north west][inner sep=0.75pt]    {$\theta _{1}$};
% Text Node
\draw (124,102.4) node [anchor=north west][inner sep=0.75pt]    {$\eta _{2}$};
% Text Node
\draw (271,108.4) node [anchor=north west][inner sep=0.75pt]    {$\eta _{1}$};
% Text Node
\draw (200,132.4) node [anchor=north west][inner sep=0.75pt]    {$\alpha _{1}$};
% Text Node
\draw (202,65.4) node [anchor=north west][inner sep=0.75pt]    {$\alpha _{2}$};
% Text Node
\draw (121.5,136.9) node [anchor=north west][inner sep=0.75pt]    {$\beta _{1}$};
% Text Node
\draw (273,141.4) node [anchor=north west][inner sep=0.75pt]    {$\beta _{1}$};

\end{tikzpicture}

    \caption{If $\theta_i\leq\pi-\frac{1}{50}$ then $\beta_1,\eta_i$ imply reducibility.}
    \label{fig: case23}
\end{figure}

\begin{proof}
    The proof follows the same line as that of Lemma \ref{lem-case12}. Suppose that, say, $\gamma_1$ is of the second type. By the $\frac{1}{100}$-flatness of $C_1$, we have $\measuredangle(\overleftarrow{\beta_1},\overrightarrow{\alpha_1}),\measuredangle(\overrightarrow{\alpha_1},\overleftarrow{\beta_1})\in[\pi-\frac{1}{50},\pi+\frac{1}{50}]$. Let $\theta_1$ (resp. $\theta_2$) be the interior angle of $R$ between $\alpha_1$ and $\eta_1$ (resp. $\eta_2$)(see Figure \ref{fig: case23}). If $\theta_1\leq\pi-\frac{1}{50}$ then $\beta_1,\eta_1$ satisfies the conditions of Lemma \ref{cor-two noncross}, so that $(P,\gamma_1,\gamma_2)$ is reducible, a contradiction. Similarly, if $\theta_2\leq\pi-\frac{1}{50}$ then $\beta_1,\eta_2$ satisfies the conditions of Lemma \ref{cor-two noncross}, again a contradiction. We deduce that $\theta_1,\theta_2>\pi-\frac{1}{50}$. As a consequence the other two interior angles of $R$ are less than $\frac{1}{25}$. But $C_2$ admits a $\theta_2$ extension at $\alpha_2$. We deduce that $\theta_2<\frac{1}{25}$. Therefore $\gamma_2$ is of the first type, in particular $C_2$ is a $1$-petal flower. But then $C_2\cup\mathring\alpha_2\cup R$ is a $3$-petal flower neighborhood of $\gamma_2$, which is impossible by Corollary \ref{cor-rc1}.
\end{proof}

Let $\theta_i\ (1\leq i\leq 4)$ be the interior angle of $R$ at the vertex of $\sigma_i$ other than $Q$ (see Figure \ref{fig: case24}). We next show that $R$ is `almost a parallelogram', in the sense of (\ref{eq-c1})(\ref{eq-c2}) below.

\begin{figure}
    \centering

\tikzset{every picture/.style={line width=0.75pt}} %set default line width to 0.75pt        

\begin{tikzpicture}[x=0.75pt,y=0.75pt,yscale=-1,xscale=1]
%uncomment if require: \path (0,339); %set diagram left start at 0, and has height of 339

%Straight Lines [id:da5580952038203536] 
\draw    (52,81) -- (73,165) ;
\draw [shift={(60.8,116.21)}, rotate = 75.96] [color={rgb, 255:red, 0; green, 0; blue, 0 }  ][line width=0.75]    (10.93,-3.29) .. controls (6.95,-1.4) and (3.31,-0.3) .. (0,0) .. controls (3.31,0.3) and (6.95,1.4) .. (10.93,3.29)   ;
%Straight Lines [id:da7267319829182396] 
\draw    (52,81) -- (182,81) ;
\draw [shift={(110,81)}, rotate = 0] [color={rgb, 255:red, 0; green, 0; blue, 0 }  ][line width=0.75]    (10.93,-3.29) .. controls (6.95,-1.4) and (3.31,-0.3) .. (0,0) .. controls (3.31,0.3) and (6.95,1.4) .. (10.93,3.29)   ;
%Straight Lines [id:da45941323236037823] 
\draw    (73,165) -- (204,164) ;
\draw [shift={(144.5,164.45)}, rotate = 179.56] [color={rgb, 255:red, 0; green, 0; blue, 0 }  ][line width=0.75]    (10.93,-3.29) .. controls (6.95,-1.4) and (3.31,-0.3) .. (0,0) .. controls (3.31,0.3) and (6.95,1.4) .. (10.93,3.29)   ;
%Straight Lines [id:da5740298079190341] 
\draw    (182,81) -- (204,164) ;
\draw [shift={(194.54,128.3)}, rotate = 255.15] [color={rgb, 255:red, 0; green, 0; blue, 0 }  ][line width=0.75]    (10.93,-3.29) .. controls (6.95,-1.4) and (3.31,-0.3) .. (0,0) .. controls (3.31,0.3) and (6.95,1.4) .. (10.93,3.29)   ;
%Straight Lines [id:da5442646248259997] 
\draw    (297,78) -- (434,79) ;
%Straight Lines [id:da6495451526833292] 
\draw    (297,78) -- (319,169) ;
%Straight Lines [id:da5036770039635994] 
\draw    (434,79) -- (456,172) ;
%Straight Lines [id:da1874777882344556] 
\draw    (297,78) -- (371,129) ;
%Straight Lines [id:da7342290920698623] 
\draw    (319,169) -- (371,129) ;
%Straight Lines [id:da23398990008283316] 
\draw    (319,169) -- (456,172) ;

% Text Node
\draw (199,110.4) node [anchor=north west][inner sep=0.75pt]    {$\eta _{1}$};
% Text Node
\draw (36,118.4) node [anchor=north west][inner sep=0.75pt]    {$\eta _{2}$};
% Text Node
\draw (122,115.4) node [anchor=north west][inner sep=0.75pt]    {$R$};
% Text Node
\draw (129,174.4) node [anchor=north west][inner sep=0.75pt]    {$\alpha _{1}$};
% Text Node
\draw (108,62.4) node [anchor=north west][inner sep=0.75pt]    {$\alpha _{2}$};
% Text Node
\draw (179,146.4) node [anchor=north west][inner sep=0.75pt]    {$\theta _{1}$};
% Text Node
\draw (75,148.4) node [anchor=north west][inner sep=0.75pt]    {$\theta _{2}$};
% Text Node
\draw (61,89.4) node [anchor=north west][inner sep=0.75pt]    {$\theta _{3}$};
% Text Node
\draw (162,89.4) node [anchor=north west][inner sep=0.75pt]    {$\theta _{4}$};
% Text Node
\draw (357,139.4) node [anchor=north west][inner sep=0.75pt]    {$\sigma _{2}$};
% Text Node
\draw (348,96.4) node [anchor=north west][inner sep=0.75pt]    {$\sigma _{3}$};
% Text Node
\draw (319,137.4) node [anchor=north west][inner sep=0.75pt]    {$\phi _{1}$};
% Text Node
\draw (331,155.4) node [anchor=north west][inner sep=0.75pt]    {$\phi _{2}$};
% Text Node
\draw (305,93.4) node [anchor=north west][inner sep=0.75pt]    {$\phi _{3}$};
% Text Node
\draw (320,80.4) node [anchor=north west][inner sep=0.75pt]    {$\phi _{4}$};
% Text Node
\draw (283,117.4) node [anchor=north west][inner sep=0.75pt]    {$\eta _{2}$};
% Text Node
\draw (381,180.4) node [anchor=north west][inner sep=0.75pt]    {$\alpha _{1}$};
% Text Node
\draw (357,56.4) node [anchor=north west][inner sep=0.75pt]    {$\alpha _{2}$};
% Text Node
\draw (454,116.4) node [anchor=north west][inner sep=0.75pt]    {$\eta _{1}$};

\end{tikzpicture}

    \caption{Let $\theta_i$ be the interior angle of $R$ at the vertex of $\sigma_i$ other than $Q$. $\sigma_2$ divides $\theta_2$ into two angles $\phi_1,\phi_2$, while $\sigma_3$ divides $\theta_3$ into two angles $\phi_3,\phi_4$.}
    \label{fig: case24}
\end{figure}
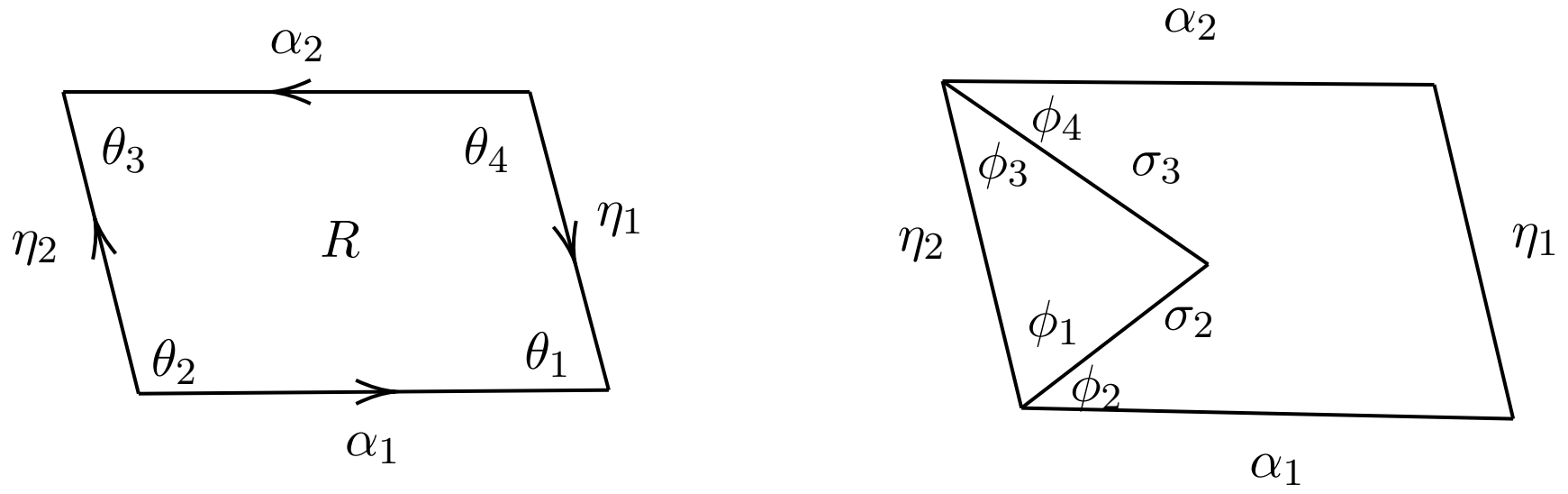

\begin{claim}
    $\theta_2+\theta_3\leq\pi+\frac{1}{10}.$
\end{claim}

\begin{proof}
    Recall that we have triangles $T_1',T_2'\subset R$ sharing a common tip $Q$ such that the base angles are less than $\frac{1}{20}$. $\sigma_2$ divides $\theta_2$ into two angles $\phi_1,\phi_2$, where $\phi_2$ is a base angle of $T_1'$, so that $\phi_2<\frac{1}{20}$ (see Figure \ref{fig: case24}). Similarly $\sigma_3$ divides $\theta_3$ into two angles $\phi_3,\phi_4$, where $\phi_4$ is a base angle of $T_2'$, so that $\phi_4<\frac{1}{20}$. Since $\phi_1,\phi_3$ are interior angles of the triangle with sides $\sigma_2,\sigma_3,\eta_2$, we deduce that $\phi_1+\phi_3<\pi$. Combining these we get $\theta_2+\theta_3\leq\pi+\frac{1}{10}$.
\end{proof}

By symmetry we have $\theta_1+\theta_4\leq\pi+\frac{1}{10}$. Since $\theta_1+\theta_2+\theta_3+\theta_4=2\pi-\mathrm{Area}(R)\geq2\pi-\epsilon$, we deduce that
\begin{align}\label{eq-c1}
    \theta_2+\theta_3,\ \theta_1+\theta_4\in[\pi-\frac{1}{9},\pi+\frac{1}{10}].
\end{align}

\begin{claim}
    $\theta_1+\theta_2\geq\pi-2\epsilon.$
\end{claim}

\begin{proof}
    By Lemma \ref{lem-h1} we have $\measuredangle(\overleftarrow{\alpha_1},\overrightarrow{\alpha_1})\in[\pi-2\epsilon,\pi+2\epsilon]$ (recall that $\gamma_1$ is of the first type, so $C_1$ is a $1$-petal flower). If $\theta_1+\theta_2<\pi-2\epsilon$, then $\pi-2\epsilon<\measuredangle(\overleftarrow{\eta_1},\overrightarrow{\eta_2})<2\pi$ (since $\measuredangle(\overleftarrow{\eta_1},\overrightarrow{\eta_2})=\measuredangle(\overleftarrow{\alpha_1},\overrightarrow{\alpha_1})+\theta_1+\theta_2$). Therefore $(P,\gamma_1,\gamma_2)$ is reducible by Lemma \ref{cor-two noncross}, a contradiction.
\end{proof}

By symmetry we have $\theta_3+\theta_4\geq\pi-2\epsilon$. Since $\theta_1+\theta_2+\theta_3+\theta_4<2\pi$, we deduce that
\begin{align}\label{eq-c2}
    \theta_1+\theta_2,\ \theta_3+\theta_4\in[\pi-2\epsilon,\pi+2\epsilon].
\end{align}

\begin{figure}
    \centering

\tikzset{every picture/.style={line width=0.75pt}} %set default line width to 0.75pt        

\begin{tikzpicture}[x=0.75pt,y=0.75pt,yscale=-1,xscale=1]
%uncomment if require: \path (0,339); %set diagram left start at 0, and has height of 339

%Straight Lines [id:da7116609909967058] 
\draw    (113,72) -- (224,72) ;
%Straight Lines [id:da7303993059753165] 
\draw    (142,159) -- (250,159) ;
%Straight Lines [id:da1183395847430262] 
\draw    (113,72) -- (142,159) ;
%Straight Lines [id:da8875873532651672] 
\draw    (224,72) -- (250,159) ;
%Straight Lines [id:da8696184722528741] 
\draw    (142,159) -- (224,72) ;
%Straight Lines [id:da8278011245668526] 
\draw    (338,164) -- (424,72) ;
%Straight Lines [id:da3799353332461295] 
\draw    (338,164) -- (471,166) ;
%Straight Lines [id:da11382776159487151] 
\draw    (424,72) -- (471,166) ;
%Straight Lines [id:da486825599199202] 
\draw    (424,72) -- (398,127) ;
%Straight Lines [id:da4741373886914926] 
\draw    (338,164) -- (398,127) ;

% Text Node
\draw (243,111.4) node [anchor=north west][inner sep=0.75pt]    {$\eta _{1}$};
% Text Node
\draw (103,114.4) node [anchor=north west][inner sep=0.75pt]    {$\eta _{2}$};
% Text Node
\draw (189,165.4) node [anchor=north west][inner sep=0.75pt]    {$\alpha _{1}$};
% Text Node
\draw (160,54.4) node [anchor=north west][inner sep=0.75pt]    {$\alpha _{2}$};
% Text Node
\draw (155,141.4) node [anchor=north west][inner sep=0.75pt]    {$\phi_1 $};
% Text Node
\draw (227,141.4) node [anchor=north west][inner sep=0.75pt]    {$\theta _{1}$};
% Text Node
\draw (197,122.4) node [anchor=north west][inner sep=0.75pt]    {$R_{1}$};
% Text Node
\draw (211,83.4) node [anchor=north west][inner sep=0.75pt]    {$\phi_{2}$};
% Text Node
\draw (167,105.4) node [anchor=north west][inner sep=0.75pt]    {$\beta $};
% Text Node
\draw (364,108.4) node [anchor=north west][inner sep=0.75pt]    {$\beta $};
% Text Node
\draw (399,170.4) node [anchor=north west][inner sep=0.75pt]    {$\alpha _{1}$};
% Text Node
\draw (453,116.4) node [anchor=north west][inner sep=0.75pt]    {$\eta _{1}$};
% Text Node
\draw (407,102.4) node [anchor=north west][inner sep=0.75pt]    {$\sigma _{4}$};
% Text Node
\draw (366,144.4) node [anchor=north west][inner sep=0.75pt]    {$\sigma _{2}$};

\end{tikzpicture}

    \caption{The angle between $\beta$ and $\sigma_2$ is less than $\frac{1}{9}$, while the angle between $\sigma_2$ and $\alpha_1$ is less than $\frac{1}{20}$. It follows that $\phi_1<\frac{1}{9}+\frac{1}{20}<\frac{1}{5}$.}
    \label{fig: case25}
\end{figure}
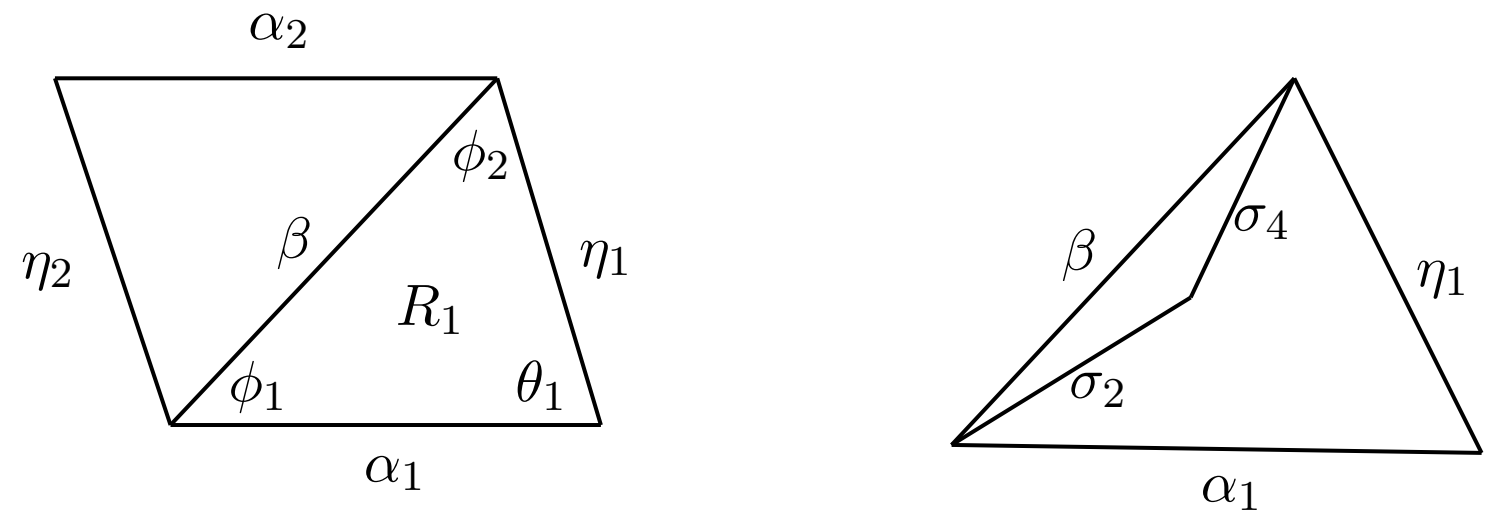

By symmetry we may assume that $l_{\eta_1}\geq l_{\eta_2},\ \theta_1\leq\theta_4$. By (\ref{eq-c1}) we get $\theta_1\leq\frac{\pi}{2}+\frac{1}{20}$. Since $R$ is `almost a parallelogram', this implies $\theta_3<\pi$. Therefore the common vertex of $\alpha_1,\eta_2$ can be joined with the common vertex of $\alpha_2,\eta_1$ by a geodesic $\beta$ in $R$. Let $R_1\subset R$ be the triangle with sides $\beta,\alpha_1,\eta_1$. Let $\phi_1,\phi_2$ be the interior angle of $R_1$ opposite to $\eta_1,\alpha_1$ respectively (see Figure \ref{fig: case25}).

\begin{claim}
\begin{align}\label{eq-cc}
    \phi_1\leq\frac{1}{5}.
\end{align}
\end{claim}

\begin{proof}
    Since the tip angles of $T_1',T_2'$ are greater than $\pi-\frac{1}{9}$, the angle formed by $\sigma_2,\sigma_4$ at $Q$ is greater than $\pi-\frac{1}{9}$. Considering the (possibly degenerate) triangle with sides $\sigma_2,\sigma_4,\beta$, we see that the angle between $\beta$ and $\sigma_2$ is less than $\frac{1}{9}$. Combining this with the fact that the base angle of $T_1'$ is less than $\frac{1}{20}$, we deduce that $\phi_1\leq\frac{1}{5}$. (See Figure \ref{fig: case25}.)
\end{proof}

Since $\theta_1+\phi_1+\phi_2=2\pi-\mathrm{Area}(R_1)\geq2\pi-\epsilon$, the estimates $\theta_1\leq\frac{\pi}{2}+\frac{1}{20},\ \phi_1\leq\frac{1}{5}$ imply
\begin{align}\label{eq-c3}
    \phi_2>1.
\end{align}
Combining $\measuredangle(\overleftarrow{\alpha_2},\overrightarrow{\alpha_2})\in[\pi-2\epsilon,\pi+2\epsilon]$ and (\ref{eq-c2}) we deduce that $\measuredangle(\overleftarrow{\eta_2},\overrightarrow{\eta_1})\in(2\pi-\frac{1}{10},2\pi+\frac{1}{10})$ (since $\measuredangle(\overleftarrow{\eta_2},\overrightarrow{\eta_1})=\measuredangle(\overleftarrow{\alpha_2},\overrightarrow{\alpha_2})+\theta_3+\theta_4$). If $\measuredangle(\overleftarrow{\eta_2},\overrightarrow{\eta_1})\leq2\pi$ then $(P,\gamma_1,\gamma_2)$ is reducible by Lemma \ref{lem-two noncross 2}. Now suppose $\measuredangle(\overleftarrow{\eta_2},\overrightarrow{\eta_1})>2\pi$.

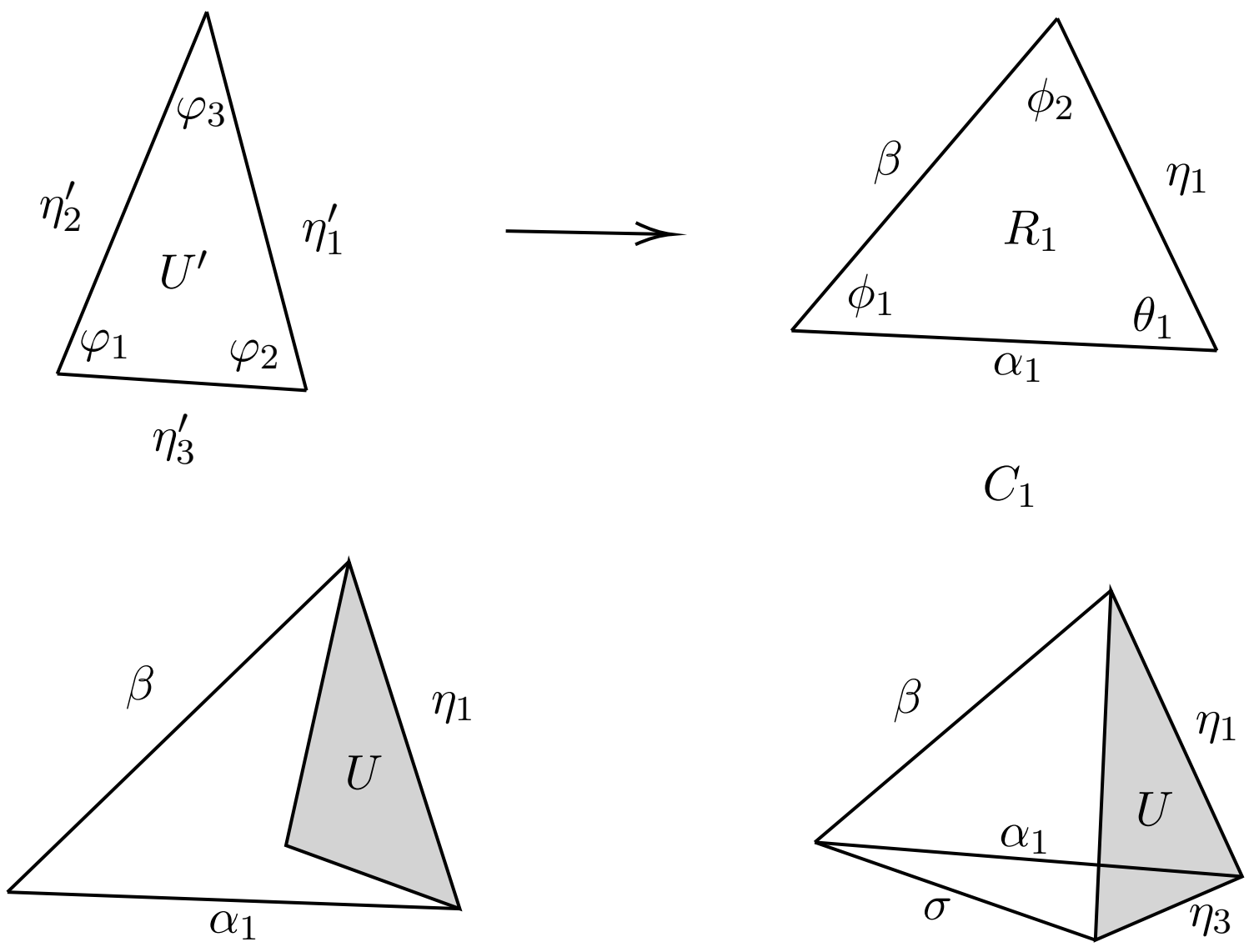
\begin{figure}
    \centering

\tikzset{every picture/.style={line width=0.75pt}} %set default line width to 0.75pt        

\begin{tikzpicture}[x=0.75pt,y=0.75pt,yscale=-1,xscale=1]
%uncomment if require: \path (0,339); %set diagram left start at 0, and has height of 339

%Straight Lines [id:da04801352407352921] 
\draw    (176,32) -- (206,146) ;
%Straight Lines [id:da18211489428391747] 
\draw    (176,32) -- (131,141) ;
%Straight Lines [id:da2562284899962325] 
\draw    (131,141) -- (206,146) ;

%Straight Lines [id:da6432870946024969] 
\draw    (266,98) -- (314,98.96) ;
\draw [shift={(316,99)}, rotate = 181.15] [color={rgb, 255:red, 0; green, 0; blue, 0 }  ][line width=0.75]    (10.93,-3.29) .. controls (6.95,-1.4) and (3.31,-0.3) .. (0,0) .. controls (3.31,0.3) and (6.95,1.4) .. (10.93,3.29)   ;
%Straight Lines [id:da20119581004926979] 
\draw    (432,34) -- (480,134) ;
%Straight Lines [id:da720527156339526] 
\draw    (432,34) -- (352,128) ;
%Straight Lines [id:da7975703967271837] 
\draw    (352,128) -- (480,134) ;

%Shape: Triangle [id:dp14614697322005632] 
\draw  [fill={rgb, 255:red, 155; green, 155; blue, 155 }  ,fill opacity=0.44 ] (218.79,197.52) -- (252.14,302.01) -- (199.84,282.95) -- cycle ;
%Straight Lines [id:da6093433599968483] 
\draw    (116,297) -- (218.79,197.52) ;
%Straight Lines [id:da00903512679047802] 
\draw    (116,297) -- (252.14,302.01) ;

%Shape: Triangle [id:dp5598553398094334] 
\draw  [fill={rgb, 255:red, 155; green, 155; blue, 155 }  ,fill opacity=0.42 ] (448.13,206.22) -- (487.62,292.32) -- (443.45,311.43) -- cycle ;
%Straight Lines [id:da15384699436632587] 
\draw    (448.13,206.22) -- (359,282) ;
%Straight Lines [id:da5644707876857565] 
\draw    (359,282) -- (443.45,311.43) ;
%Straight Lines [id:da5729718787852003] 
\draw    (359,282) -- (487.62,292.32) ;

% Text Node
\draw (216,254.4) node [anchor=north west][inner sep=0.75pt]    {$U$};
% Text Node
\draw (242,235.4) node [anchor=north west][inner sep=0.75pt]    {$\eta _{1}$};
% Text Node
\draw (150,227.4) node [anchor=north west][inner sep=0.75pt]    {$\beta $};
% Text Node
\draw (175,301.4) node [anchor=north west][inner sep=0.75pt]    {$\alpha _{1}$};
% Text Node
\draw (411,133.4) node [anchor=north west][inner sep=0.75pt]    {$\alpha _{1}$};
% Text Node
\draw (375,69.4) node [anchor=north west][inner sep=0.75pt]    {$\beta $};
% Text Node
\draw (414,90.4) node [anchor=north west][inner sep=0.75pt]    {$R_{1}$};
% Text Node
\draw (463,76.4) node [anchor=north west][inner sep=0.75pt]    {$\eta _{1}$};
% Text Node
\draw (453,116.4) node [anchor=north west][inner sep=0.75pt]    {$\theta _{1}$};
% Text Node
\draw (367,109.4) node [anchor=north west][inner sep=0.75pt]    {$\phi _{1}$};
% Text Node
\draw (421,50.4) node [anchor=north west][inner sep=0.75pt]    {$\phi _{2}$};
% Text Node
\draw (408,167.4) node [anchor=north west][inner sep=0.75pt]    {$C_{1}$};
% Text Node
\draw (160,102.4) node [anchor=north west][inner sep=0.75pt]    {$U'$};
% Text Node
\draw (203,88.4) node [anchor=north west][inner sep=0.75pt]    {$\eta '_{1}$};
% Text Node
\draw (124,81.4) node [anchor=north west][inner sep=0.75pt]    {$\eta '_{2}$};
% Text Node
\draw (158,151.4) node [anchor=north west][inner sep=0.75pt]    {$\eta '_{3}$};
% Text Node
\draw (136,126.4) node [anchor=north west][inner sep=0.75pt]    {$\varphi _{1}$};
% Text Node
\draw (181,129.4) node [anchor=north west][inner sep=0.75pt]    {$\varphi _{2}$};
% Text Node
\draw (165,56.4) node [anchor=north west][inner sep=0.75pt]    {$\varphi _{3}$};
% Text Node
\draw (454,265.4) node [anchor=north west][inner sep=0.75pt]    {$U$};
% Text Node
\draw (413,275.4) node [anchor=north west][inner sep=0.75pt]    {$\alpha _{1}$};
% Text Node
\draw (472,241.4) node [anchor=north west][inner sep=0.75pt]    {$\eta _{1}$};
% Text Node
\draw (381,231.4) node [anchor=north west][inner sep=0.75pt]    {$\beta $};
% Text Node
\draw (470,299.4) node [anchor=north west][inner sep=0.75pt]    {$\eta _{3}$};
% Text Node
\draw (390,297.4) node [anchor=north west][inner sep=0.75pt]    {$\sigma $};

\end{tikzpicture}

    \caption{We always have $\varphi_3<\phi_2$. If $\varphi_2\leq\theta_1$, then there is an embedding $U'\hookrightarrow R_1$. If $\varphi_2>\theta_1$, we construct a quadrilateral $R'$ such that the angle between $\eta_1$ and $\eta_3$ is equal to $\varphi_2$. Then the triangle formed by the sides $\eta_1,\eta_3$ satisfies the desired conditions.}
    \label{fig: case26}
\end{figure}

\begin{claim}\label{claim-c1}
    There is an embedded geodesic triangle $U\subset C_1\cup\alpha_1\cup R_1$ with sides $\eta_1,\eta_2',\eta_3$ (anti-clockwise ordered), such that $l_{\eta_2'}=l_{\eta_2}$ and the angle between $\eta_1,\eta_2'$ is equal to $\measuredangle(\overleftarrow{\eta_2},\overrightarrow{\eta_1})-2\pi$. Moreover, $\eta_2',\eta_3$ are simple geodesics. (See Figure \ref{fig: case26}.)
\end{claim}

Assuming Claim \ref{claim-c1}, the conditions of Lemma \ref{lem-two noncross 3} are satisfied. We conclude that $(P,\gamma_1,\gamma_2)$ is reducible.

It remains to prove Claim \ref{claim-c1}. Let $U'\subset\H$ be a geodesic triangle with sides $\eta_1',\eta_2',\eta_3'$ (anti-clockwise ordered) such that $l_{\eta_1'}=l_{\eta_1},\ l_{\eta_2'}=l_{\eta_2}$ and the angle between $\eta_1'$ and $\eta_2'$ is equal to $\measuredangle(\overleftarrow{\eta_2},\overrightarrow{\eta_1})-2\pi$. We show that there exists an orientation preserving isometric embedding $U'\hookrightarrow C_1\cup\alpha_1\cup R_1$ which maps the side $\eta_1'$ to $\eta_1$, and the sides $\eta_2',\eta_3'$ to simple geodesics. The image of $U'$ would be the desired triangle $U$.

Let $\varphi_i$ be the angle of $U'$ opposite to $\eta_i'$. By $\measuredangle(\overleftarrow{\eta_2},\overrightarrow{\eta_1})\in(2\pi,2\pi+\frac{1}{10})$, we get
\begin{align}\label{eq-c4}
    \varphi_3<\frac{1}{10}.
\end{align}
Combining (\ref{eq-c3}) and (\ref{eq-c4}) we get
\begin{align}
    \varphi_3<\phi_2.
\end{align}
If $\varphi_2\leq\theta_1$, then the angles of $U'$ at the two vertices of $\eta_1'$ is less than or equal to the corresponding angles of $R_1$ at the two vertices of $\eta_1$ respectively. It follows that there is an embedding $U'\hookrightarrow R_1$, proving the claim. (See Figure \ref{fig: case26}.)

Now suppose $\varphi_2>\theta_1$. Since $l_{\eta_2'}\leq l_{\eta_1'}$, we get $\varphi_2<\frac{\pi}{2}$. This implies
\begin{align}\label{eq-c6}
    \varphi_2-\theta_1<\frac{\pi}{2}.
\end{align}
By (\ref{eq-c4}), $\eta_3'$ cannot be the longest side of $U'$. So $l_{\eta_3'}\leq l_{\eta_1'}$. By (\ref{eq-cc}) and (\ref{eq-c3}), we get $\phi_1<\phi_2$, which implies $l_{\eta_1}\leq l_{\alpha_1}$. We obtain
\begin{align}\label{eq-c7}
    l_{\eta_3'}\leq l_{\alpha_1}.
\end{align}
By (\ref{eq-c6}),(\ref{eq-c7}) and Lemma \ref{lem-h5}, there is an embedded geodesic triangle $R_2\subset C_1$ with sides $\alpha_1,\eta_3,\sigma$ (clockwise ordered) such that $l_{\eta_3}=l_{\eta_3'}$ and the angle between $\alpha_1,\eta_3$ is equal to $\varphi_2-\theta_1$. We get a geodesic quadrilateral $R'=R_1\cup\alpha_1\cup R_2$, with sides $\beta,\eta_1,\eta_3,\sigma$, such that the interior angle between $\eta_1$ and $\eta_3$ is equal to $\varphi_2$. Since $l_{\eta_3}\leq l_{\alpha_1}$, the angle of $R_2$ opposite to $\eta_3$ is less than $\frac{\pi}{2}$. Combining this with (\ref{eq-cc}) we deduce that the angle of $R'$ between $\beta$ and $\sigma$ is less than $\pi$. Therefore $R'$ is convex. The geodesic triangle spanned by the sides $\eta_1$ and $\eta_3$ is the desired triangle $U$. This proves Claim \ref{claim-c1}.

We see that the case $\sigma_2\neq\sigma_3,\sigma_1\neq\sigma_4$ always end up with $(P,\gamma_1,\gamma_2)$ being reducible. So this case is in fact impossible by our assumption.

\subsection{The case $E_1\cap E_2=\varnothing$}

In this case we show that there is $S'\in[S]$ which satisfies the conditions of Model 3. We want both $\gamma_1$ and $\gamma_2$ to have $1$-petal flower neighborhoods which admit disjoint extensions. If $\gamma_i$ is of the first type then $C_i$ is the desired $1$-petal flower with a $\frac{1}{200}$-extension $E_i$. But if $\gamma_i$ is of the second type, we need to perform certain surgery so that $\gamma_i$ has a $1$-petal flower neighborhood, which admits a $\frac{1}{200}$-extension. Since $\beta_i$ is non-dividing, it is free by Lemma \ref{lem-noncross trap}, so that by Lemma \ref{lem-free} it is the petal of a $1$-petal flower $F_i$ of length $10\epsilon$. By Lemma \ref{lem-noncross trap2}, $F_i$ is in the component of $N_{\frac{1}{2}}(P)\backslash\beta_i$ not containing $\gamma_1,\gamma_2$. Therefore $F_i$ and $E_i$ lie in different components of $N_{\frac{1}{2}}(P)\backslash\beta_i$, in particular $E_i\cap F_i=\varnothing$. Our surgery will be supported in $E_i\cup\beta_i\cup F_i$. Let's focus on one index $i$, and simplify our notation.

Let $\gamma$ be a boundary curve of an $\epsilon$-good surface $S$. Let $C$ be a $\frac{1}{100}$-flower neighborhood of $\gamma$ with two petals $\alpha,\beta$, which admits a $\frac{1}{20}$-extension $E$ at $\alpha$. Suppose that $\beta$ is the petal of an embedded $1$-petal flower $F$ of length $10\epsilon$ disjoint from $E$. We are going to perform a surgery (supported in $E\cup\beta\cup F$), such that in the resulting surface, $\gamma$ has a $1$-petal flower neighborhood which admits a $\frac{1}{200}$-extension.

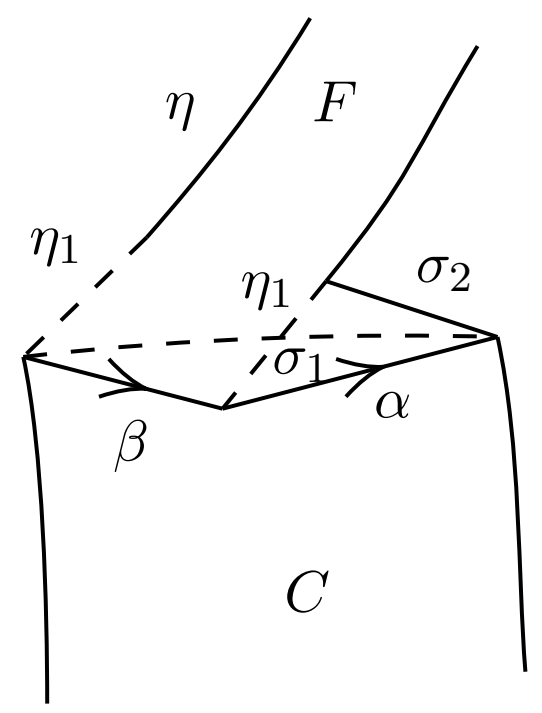
\begin{figure}
    \centering

\tikzset{every picture/.style={line width=0.75pt}} %set default line width to 0.75pt        

\begin{tikzpicture}[x=0.75pt,y=0.75pt,yscale=-1,xscale=1]
%uncomment if require: \path (0,339); %set diagram left start at 0, and has height of 339

%Curve Lines [id:da7897385295366213] 
\draw    (137,134) .. controls (142,157) and (143,197) .. (143,221) ;
%Curve Lines [id:da19194056975193763] 
\draw    (256,129) .. controls (262,157) and (261,188) .. (263,213) ;
%Curve Lines [id:da504256920289821] 
\draw    (213,115) .. controls (235,88) and (234,84) .. (251,56) ;
%Curve Lines [id:da1610921951086408] 
\draw    (168,104) .. controls (182,88) and (195,72) .. (209,49) ;
%Straight Lines [id:da7562609774917006] 
\draw    (137,134) -- (187,147) ;
\draw [shift={(167.81,142.01)}, rotate = 194.57] [color={rgb, 255:red, 0; green, 0; blue, 0 }  ][line width=0.75]    (10.93,-4.9) .. controls (6.95,-2.3) and (3.31,-0.67) .. (0,0) .. controls (3.31,0.67) and (6.95,2.3) .. (10.93,4.9)   ;
%Straight Lines [id:da9920448776861583] 
\draw    (187,147) -- (256,129) ;
\draw [shift={(227.31,136.49)}, rotate = 165.38] [color={rgb, 255:red, 0; green, 0; blue, 0 }  ][line width=0.75]    (10.93,-4.9) .. controls (6.95,-2.3) and (3.31,-0.67) .. (0,0) .. controls (3.31,0.67) and (6.95,2.3) .. (10.93,4.9)   ;
%Straight Lines [id:da45363103946755434] 
\draw    (213,115) -- (256,129) ;
%Curve Lines [id:da6902467425935204] 
\draw  [dash pattern={on 4.5pt off 4.5pt}]  (137,134) .. controls (179,129) and (223,128) .. (256,129) ;
%Straight Lines [id:da9576096442911047] 
\draw  [dash pattern={on 4.5pt off 4.5pt}]  (213,115) -- (187,147) ;
%Straight Lines [id:da6745144060230382] 
\draw  [dash pattern={on 4.5pt off 4.5pt}]  (168,104) -- (137,134) ;

% Text Node
\draw (198,130.4) node [anchor=north west][inner sep=0.75pt]    {$\sigma _{1}$};
% Text Node
\draw (234,107.4) node [anchor=north west][inner sep=0.75pt]    {$\sigma _{2}$};
% Text Node
\draw (223.5,141.4) node [anchor=north west][inner sep=0.75pt]    {$\alpha $};
% Text Node
\draw (201,186.4) node [anchor=north west][inner sep=0.75pt]    {$C$};
% Text Node
\draw (158,148.4) node [anchor=north west][inner sep=0.75pt]    {$\beta $};
% Text Node
\draw (171,66.4) node [anchor=north west][inner sep=0.75pt]    {$\eta $};
% Text Node
\draw (190,111.4) node [anchor=north west][inner sep=0.75pt]    {$\eta _{1}$};
% Text Node
\draw (137,100.4) node [anchor=north west][inner sep=0.75pt]    {$\eta _{1}$};
% Text Node
\draw (208,63.4) node [anchor=north west][inner sep=0.75pt]    {$F$};

\end{tikzpicture}

    \caption{$T$ is the triangle with sides $\alpha,\sigma_1,\sigma_2$. Perform a move along $\eta_1,\sigma_1$. $U$ is the triangle with sides $\alpha,\beta,\sigma$, where $\sigma$ is the horizontal dotted line in the figure. $\sigma$ and $\gamma$ bound a $1$-petal flower neighborhood $N$ of $\gamma$. Since $\eta$ and $\sigma_2$ both form angles greater than $\frac{1}{200}$ with $\sigma$, $N$ admits a $\frac{1}{200}$-extension at $\sigma$.}
    \label{fig: case 2}
\end{figure}

Let $P$ be the pistil of $C$, and let $T$ be the $\frac{1}{21}$-cap of $C$ at $\alpha$ (see Figure \ref{fig: case 2}). Denote by $\sigma_1,\sigma_2$ the other two sides of $T$ so that $\sigma_1,\sigma_2,\alpha$ is arranged in clockwise order on $\partial T$. Orient $\alpha,\beta$ so that $C$ lies on the right, and orient $\sigma_1$ so that $T$ lies on the right. Let $v$ be a tangent vector at $P$ such that $\measuredangle(v,\overrightarrow{\sigma_1})=2\pi$.
\begin{claim}\label{claim-c}
    There is a simple geodesic $\eta:[0,8\epsilon]\rightarrow\overline{F}$ of length $8\epsilon$ such that $\eta(0)=P,\overrightarrow{\eta}=v$.
\end{claim}

Since $C$ is $\frac{1}{100}$-flat, $\measuredangle(\overleftarrow\beta,\overrightarrow\alpha)\in[\pi-\frac{1}{50},\pi+\frac{1}{50}]$. By the definition of a cap, we have $\measuredangle(\overrightarrow\alpha,\overrightarrow{\sigma_1})=\frac{1}{21}$. Moreover, by Lemma \ref{lem-h2}, we get $\measuredangle(\overrightarrow{\beta},\overleftarrow{\beta})\in[\pi-2\epsilon,\pi+2\epsilon]$. Combining these estimates we get $\measuredangle(\overrightarrow\beta,v)\in[\frac{1}{40},\frac{1}{5}]$ (by the computation $\frac{1}{21}-\frac{1}{50}-2\epsilon>\frac{1}{40}$). Claim \ref{claim-c} then follows from Lemma \ref{cor-h1}.

Let $\eta_1=\eta|_{[0,l_{\sigma_1}]}$. Then $\eta_1,\sigma_1$ is an associated pair. Cut along $\eta_1,\sigma_1$ and identify $\eta_1^+$ with $\sigma_1^-$, $\eta_1^-$ with $\sigma_1^+$ (Figure \ref{fig: case 2}). Since $\measuredangle(\overleftarrow\beta,\overrightarrow\alpha)\in[\pi-\frac{1}{50},\pi+\frac{1}{50}]$, the other two angles of the hyperbolic triangle formed by $\alpha,\beta$ are less than $\frac{1}{50}$, which is less than $\measuredangle(\overrightarrow\beta,\overrightarrow{\eta})$ and the base angle of $T$. Therefore, there is an embedded geodesic triangle $U$ in the resulting surface with sides $\alpha,\beta$, whose third side $\sigma$ is a simple geodesic loop based at $P$. Then the annulus $N$ in the resulting surface bounded by $\sigma$ and $\gamma$ is a $1$-petal flower neighborhood of $\gamma$.
\begin{claim}
    $N$ admits a $\frac{1}{200}$-extension at $\sigma$.
\end{claim}
This is because the angles between $\sigma$ and $\eta$, $\sigma$ and $\sigma_2$ are greater than $\frac{1}{200}$ (by the computation $\frac{1}{40}-\frac{1}{50}=\frac{1}{200}$).

For $i=1,2$, if $\gamma_i$ is of the second type then we perform the above surgery in $E_i\cup\beta_i\cup F_i$. We obtain $S'\in[S]$ such that both $\gamma_1$ and $\gamma_2$ have $1$-petal flower neighborhoods which admit disjoint $\frac{1}{200}$-extensions. Moreover the two flower neighborhoods have the same pistil (which is the base point of the simple loop in $\partial E_i$ other than $\beta_i$, unchanged under the surgery). Therefore $S'$ satisfies the conditions of Model 3.

%\printbibliography
\bibliography{ref}

\end{document}